\newcommand{\bfb}{\mathbf b}
\newcommand{\bft}{\mathbf t}
\newcommand{\bfx}{\mathbf x}

\documentclass[preprint,12pt]{elsarticle}
\usepackage{amsmath,amssymb,amsthm}
\usepackage{graphicx}
\usepackage[margin=1in]{geometry} 
\usepackage{subfigure}
\usepackage{placeins}
\usepackage{multirow}
 \usepackage{url}
\numberwithin{equation}{section}
\usepackage{color}
\newcommand{\bfe}{\mathbf e}
\newtheorem{lem}{Lemma}

\begin{document}
\journal{J Computational and Applied Mathematics}
\begin{frontmatter}

\title{Non-negatively constrained least squares and parameter choice by the residual periodogram for the inversion of electrochemical impedance spectroscopy}

\author[label1]{Jakob Hansen}
\ead{jkhanse2@asu.edu}
\author[label1]{Jarom Hogue}
\ead{jdhogue@asu.edu}
\author[label1]{Grant Sander}
\ead{gksander@asu.edu}
\author[label1]{Rosemary A Renaut\corref{cor1}}
\ead{renaut@asu.edu}
\ead[url]{math.asu.edu/~rosie}
\author[label2]{Sudeep C Popat }
\ead{scp@asu.edu}
\cortext[cor1]{Corresponding Author: Rosemary Renaut, 480 965 3795 }
\address[label1]{School of Mathematical and Statistical Sciences, Arizona State University, Tempe, AZ 85287-1804, USA}
\address[label2]{Swette Center for Environmental Biotechnology, Biodesign Institute, Arizona State University, Tempe, AZ 85287, USA}
\begin{abstract}

The inverse problem associated with electrochemical  impedance spectroscopy  requiring the solution of a Fredholm integral equation of the first kind is considered. If the underlying physical model is not clearly determined, the inverse problem needs to be solved using a regularized linear least squares problem that is obtained from the discretization of the integral equation.     For this system,  it is shown that the model error can be made negligible by a  change of variables and by extending the effective range of quadrature. This change of variables serves as a right preconditioner that  significantly improves the condition of the system. Still, to obtain 
 feasible solutions the  additional constraint of non-negativity is  required. Simulations with artificial, but realistic, data demonstrate that the use of non-negatively constrained least squares  with a smoothing norm provides higher quality solutions than those obtained without the non-negative constraint.   Using higher-order smoothing norms also reduces the error in the solutions. The L-curve and residual periodogram parameter choice criteria, which are used for parameter choice with regularized  linear least squares, are successfully adapted to be used for the non-negatively constrained Tikhonov least squares problem. 
Although these results  have been verified within the context of the   analysis of electrochemical impedance spectroscopy, there is no reason to suppose that they would not be relevant within the broader framework of solving Fredholm integral equations  for other applications.

\end{abstract}


\end{frontmatter}

\section{Introduction}
We consider the numerical solution of ill-posed inverse problems that are motivated by measurements of electrochemical impedance spectra  from which   a model of the underlying physical reaction mechanisms is desired. There is extensive literature on a wide range of applications in which the same, or similar models can be applied. These include measurements for   solid oxide fuel cells   \cite{Enetal:10,Leetal:08,Leetal:10,Lietal:10,Lietal:11,Maetal:04,Scetal:02,Soetal:08}, microbial fuel cells \cite{CSUMS12}, as well as of physiological parameters,  and from a diverse range of dielectric models, \cite{Barsukov,Macdonald-quad, ln-drt, Ward:96}. In these applications the unknown distribution function of relaxation times (DRT) is related to a set of impedance measurements by   the Fredholm integral equation
\begin{equation}\label{modeleq}
Z(\omega) = R_0 + R_\text{pol} \int_0^\infty \frac{g(t)}{1+ i\omega t} \, dt,
\end{equation}
where $\omega$ is angular frequency, $t$ is   time, and $g(t)$ is the desired DRT with normalization  $\int_0^{\infty}g(t)dt=1$.

There are several models used to represent the individual processes of a DRT, many of which are  mostly used for the analysis of dielectric materials and are described in \cite{Barsukov}.   Several are directly applicable to the fuel cell modeling case, where they usually take the form of theoretical circuit components used in constructing equivalent circuit models. Equivalent circuit elements used for fuel cell modeling include the Cole-Cole (also known as RQ or ZARC) element, the Generalized Finite-Length Warburg element, and the Gerischer impedance \cite{Barsukov, Leonide-DRT,Maetal:04}. In analysis of specific fuel cell designs a log-normal form for the DRT has also been used \cite{ln-drt, CSUMS12}.   
Here we focus our investigations on the Cole-Cole DRT, which can be rendered temperature independent only in the limiting cases of $\beta\rightarrow 0, 1$, and the temperature independent lognormal DRT, denoted throughout by RQ and LN, respectively.  

The RQ impedance  is a generalization of  a simple parallel RC circuit and for a single process has an impedance given by
\begin{equation}\label{Coleimp}
Z_{\mathrm{RQ}}(\omega) = \frac{1}{1 + (i\omega t_0)^\beta}, 
\end{equation}
where $t_0$ is the point of maximum distribution, and $\beta$ is a shape parameter controlling the width of the distribution. The   corresponding DRT is 
\begin{equation}\label{DRTtCole}
g_{\mathrm{RQ}}(t) = \frac{1}{2\pi t}\frac{\sin \beta \pi}{\cosh\left(\beta \ln\left(\frac{t}{t_0}\right)\right) + \cos \beta \pi},
\end{equation}
which reduces  to the Dirac delta distribution  when $\beta = 1$, \cite{Barsukov}.
There is, however, no analytic  form for   the impedance   corresponding to the  log-normal DRT given by 
\begin{equation}\label{DRTtlog}
g_{\mathrm{LN}}(t) = \frac{1}{t\sigma\sqrt{2\pi}}\exp{\left(-\frac{(\ln(t) - \mu)^2}{2\sigma^2}\right)}.
\end{equation}

Although a number of options have been presented in the literature for geometrically assessing the parameterization of the DRT from  impedance data for a single physical process, e.g. as noted in \cite{Ward:96}, for given measured and noisy impedance data from multiple processes  there are effectively only two basic approaches that may be considered to estimate the DRT.  When a  specific analytic but parameter dependent form for the impedance is known, as in \eqref{Coleimp},  parametric nonlinear least squares (NLS) fitting may be used to determine the underlying parameters of the impedance and hence of the DRT, \cite{Macdonald-quad}. On the other hand, when no analytic representation of the impedance is available, as in \eqref{DRTtlog},  it is  still   possible, but more computationally expensive, to apply a parametric nonlinear fit by using direct numerical integration of \eqref{modeleq}. In either case, an alternative is to apply a linear least squares (LLS) fit directly to the DRT, but this is also challenging due to the general ill-posedness of the problem, e.g.  \cite{bjorck,hansen:1997-1,hansen:2007-1,hansen:2010-1,Vogel:02}. Both approaches, as well as the geometric analyses, have been extensively considered in the  literature, e.g. \cite{Barsukov}.  
When the model for the DRT is not known, perhaps when the physical process is not completely understood or the number of processes has not been determined, the only option is to fit directly to the DRT, without identifying its specific parameterization. 

Before further pursuing the LLS fit, we illustrate in Section~\ref{NLSfits} the use of direct NLS  fitting for a simple one-process example in order to emphasize the (self-evident) significance of the prior knowledge of the model. Assuming the wrong model leads to apparently robust data fitting, while at the same time potentially leading to incorrect conclusions about the DRT parameterization. With this conclusion we move in Section~\ref{sec:ls} to an analysis of the system describing the LLS fitting that arises when approximating \eqref{modeleq} discretely. The direct discretization of \eqref{modeleq} leads to two ill-conditioned systems of equations, for the real and imaginary parts separately.   Most  literature on the problem suggests the use of LLS for the systems obtained in this way, in conjunction with regularization to stabilize the  estimation of the  solution, \cite{Leetal:08,  Weese:92}.   In contrast, it was suggested in \cite{Maetal:04},  that rather than estimating the DRT in the given $t$-space, a transformation to $s-$space via $s=\log(t)$ would be preferable and that the resulting ill-posed system be solved using a non-negative least squares (NNLS)   algorithm, specifically imposing the constraint that the DRT is a positive distribution. In Section~\ref{model error} we investigate the modeling error that arises when using the $s-$space transformation, leading to new results that quantify the total modeling error due to discretization and truncation in \eqref{modeleq} for both real and imaginary terms. The results go 
   beyond those presented in \cite{CSUMS12} for the $t-$space formulation, by providing error estimates which are primarily determined by the kernel $h(\omega,t)=(1+i \omega t )^{-1}$,   only relying on standard smoothness and decay conditions for the DRT  functions.  
   
The numerical algorithms for the estimation of the DRT are discussed in Section~\ref{numerical}. First it is demonstrated that the $s-$transformation serves as a \textbf{right preconditioner}, leading to more stable estimation of the underlying basis for the solution when the time discretization is chosen appropriately in relation to the frequency measurements. Still the model remains ill-conditioned, and solution techniques using regularization are required, introducing the need for determination of a regularization parameter that weights the regularization term. Estimation of this regularization parameter for Tikhonov regularization is well-studied in the literature e.g. \cite{hansen:1997-1,Vogel:02}. We therefore give just a   brief and standard overview of regularized  LLS in Section~\ref{regsec}. On the other hand, the estimation of the regularization parameter in the context of non-negatively (NN) constrained Tikhonov regularization  is less well-studied. Thus to put the new work in context we focus on the estimation of the regularization parameter using the L-curve (LC) and residual periodogram (RP), or its normalized cumulative periodogram (NCP). This discussion leads to new regularization parameter estimation techniques for the NN-constrained Tikhonov LS  problem. In particular, the LC  and NCP  parameter estimation techniques  are extended to the  NNLS algorithm, with minimal additional algorithmic development, Section~\ref{nnlssec}.  

Finally, the theoretical developments are verified and evaluated for a set of impedance datasets,   consisting of two and three 
 physical processes, which are motivated by examples that are seen from practical configurations. The presented results justify both the use of inclusion of the non-negativity  constraint  for finding approximate DRTs from multiple physical processes, and the use of the  NCP and LC for estimating the regularization parameter within this context. 
 To verify that our approach is sufficiently general, numerical simulations use two different algorithms for the solution of the NNLS problem, namely the Matlab function \texttt{lsqnonneg} which is based on the algorithm in \cite{LawsonHansen:1995}, and the CVX software package, \cite{cvx,gb08}.
Our results verify that indeed the LC and NCP techniques can be used for optimal regularization parameter determination when solving the NN-constrained Tikhonov least squares problem.   The latter is of more general use for ill-posed systems of linear equations with  NN-constraints. Conclusions of the work are provided in Section~\ref{conc}. 
 
 \section{Parametric NLS Fitting: Distinguishing between models for one process}\label{NLSfits}
\begin{figure}[h]
\centering
\subfigure[DRTs:  $t-$space ]{\label{1a}\includegraphics[width=1.2in]{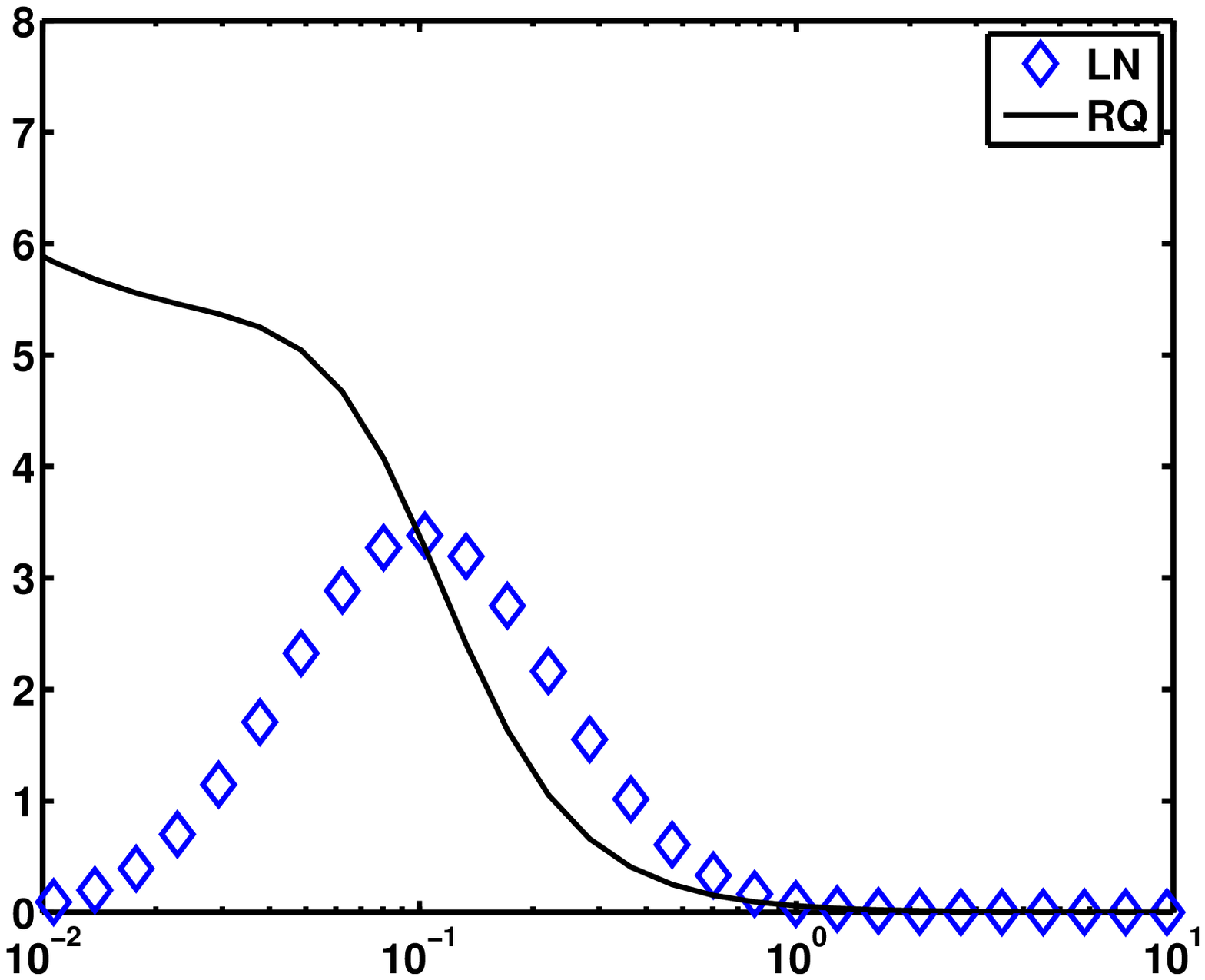}}
\subfigure[DRTs: $s-$space]{\label{1b}\includegraphics[width=1.2in]{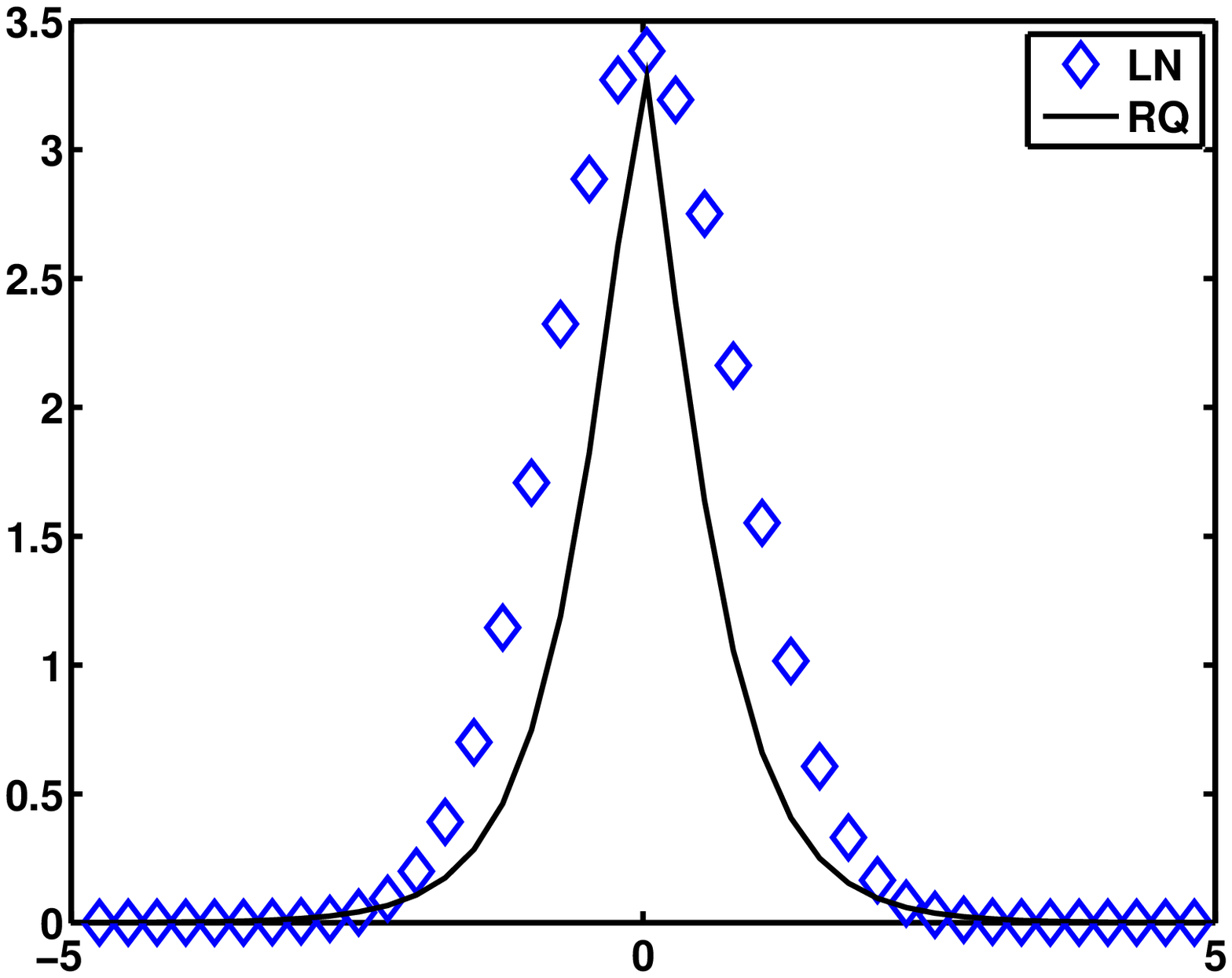}}
\subfigure[Nyquist plot]{ \label{1c}\includegraphics[width=1.2in]{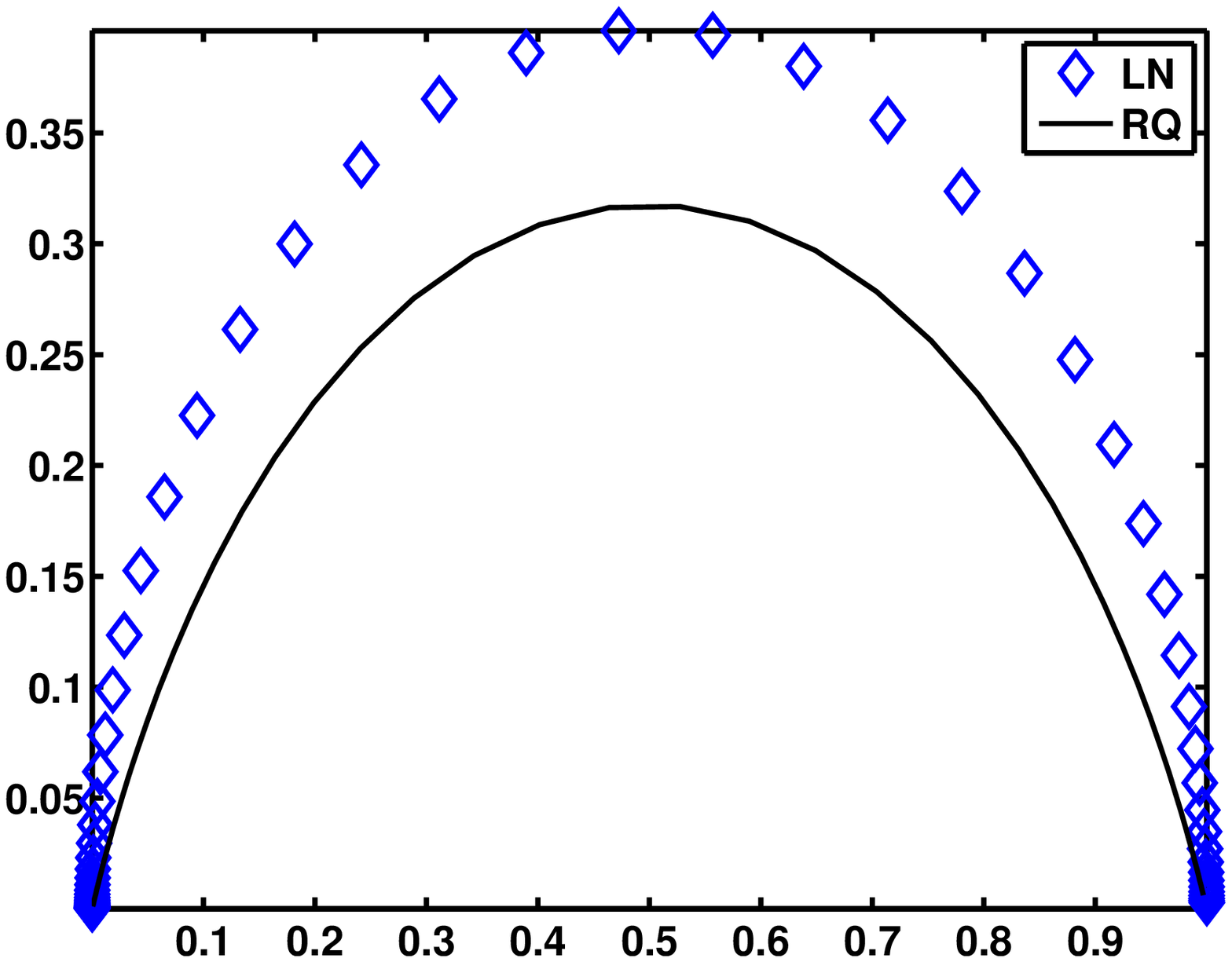}}
\subfigure[Components $Z_1$]{ \label{1d}\includegraphics[width=1.2in]{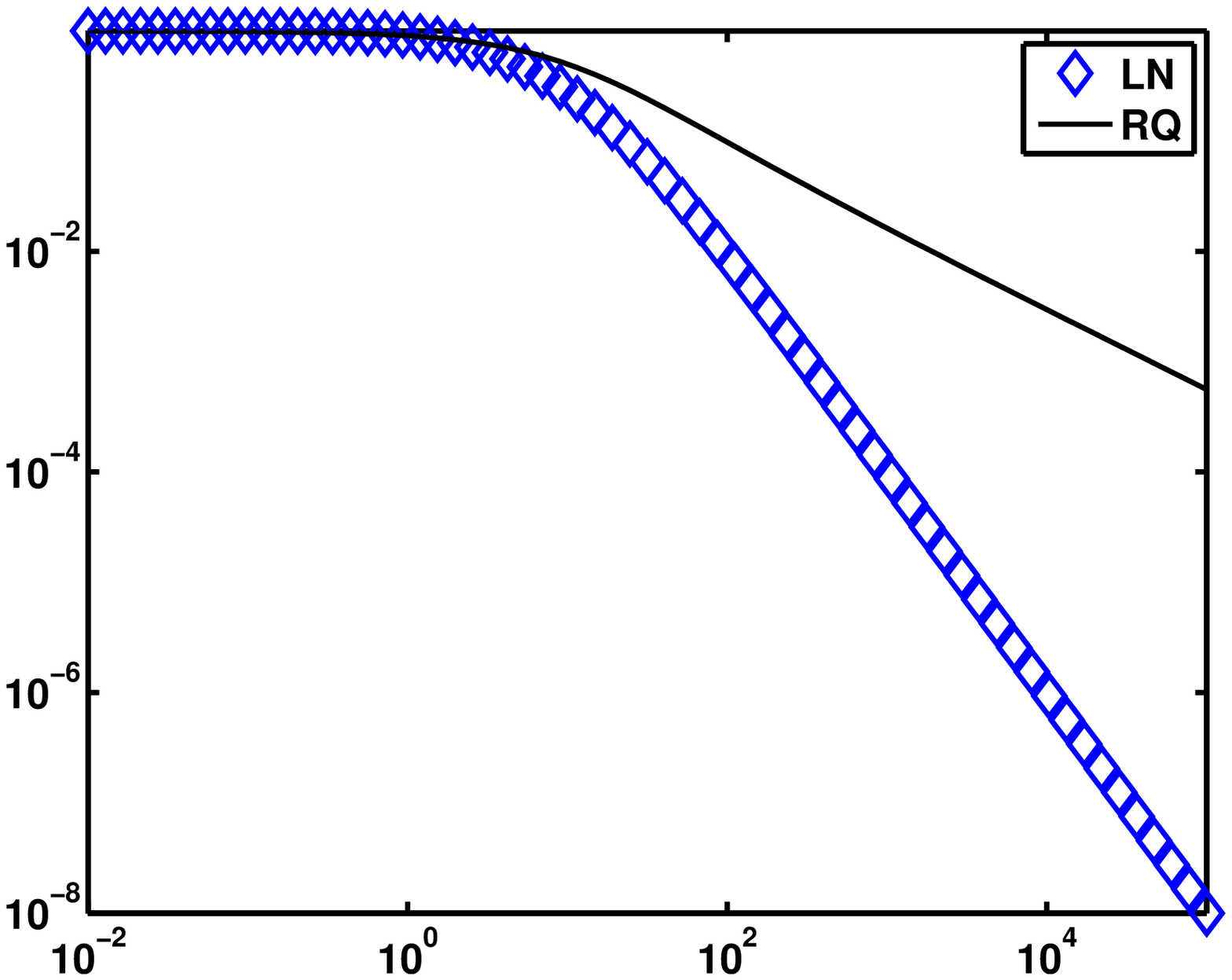}}
\subfigure[Components $Z_2$]{ \label{1e}\includegraphics[width=1.2in]{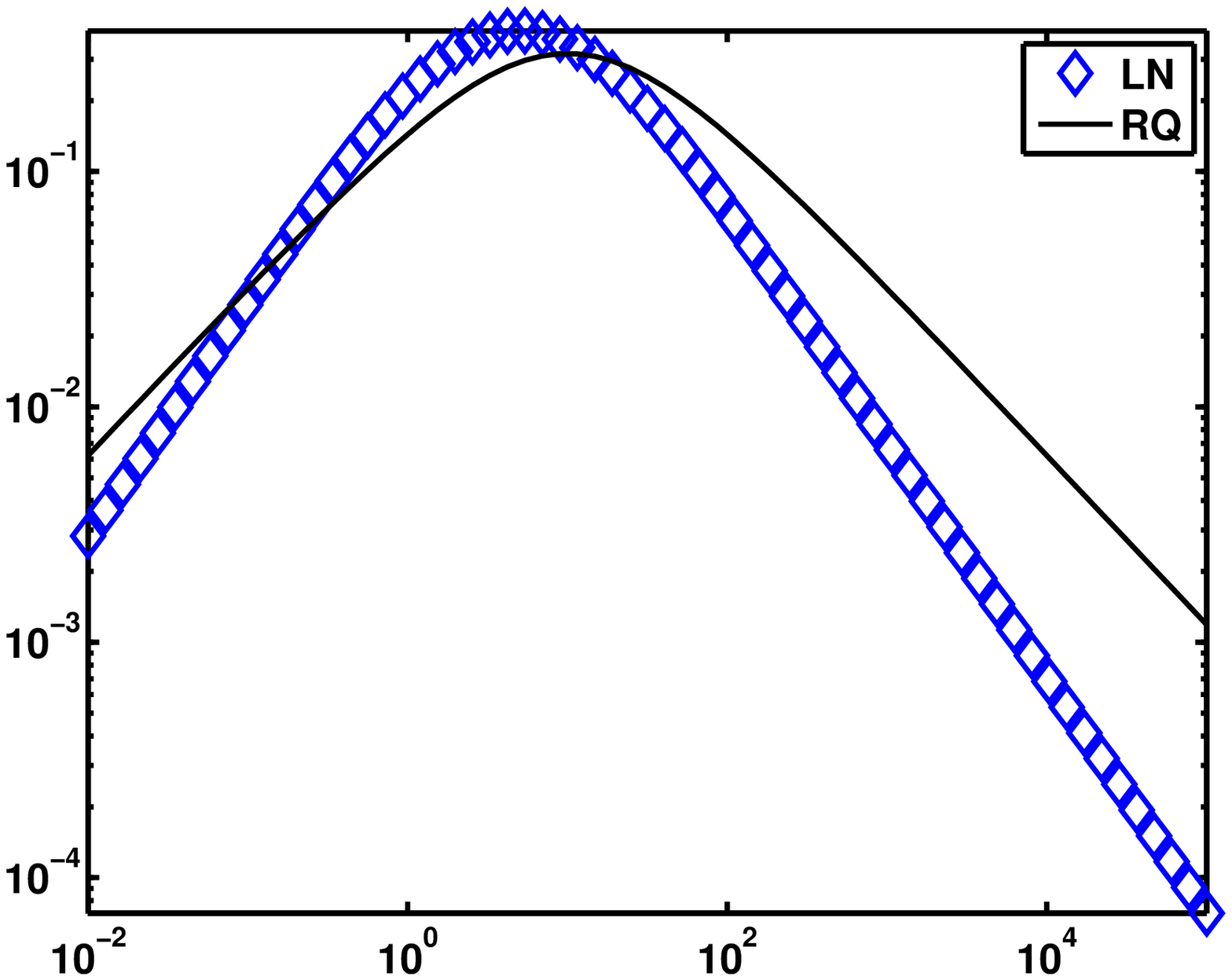}}
\caption{Simulated \textit{exact} data measured at $65$ logarithmically spaced points in $\omega$. In each case the solid line indicates the RQ functions and the $\diamond$ symbols the LN functions.}\label{fig1}
\end{figure}

 We investigate parametric nonlinear fitting to impedance data given by \eqref{modeleq} with $R_0=0$ and $R_{\mathrm{pol}}=1$  (without any loss of generality). For a given DRT, simulated data for $Z=Z_1 - i Z_2$ were generated for realistic frequencies $\omega \in [10^{-2}, 10^5]$ sampled logarithmically at $65$ points, providing the \textit{exact} fitting vector   $[Z_1;Z_2]$ of length $130$. The exact data were generated for the RQ and  LN  DRTs,  \eqref{DRTtCole} and \eqref{DRTtlog}, respectively, with  parameters $t_0=0.1$, $\beta=0.7203$ and $\sigma=\ln(2.3)$ chosen to  provide DRTs aligned with respect to the location and height of the peak in the $s=\ln(t/t_0)$ space,   see Figure~\ref{1b}, $\beta=(2/\pi)\arctan{((\sqrt{2\pi}/\sigma)\exp(-\sigma^2/2))}$, see Section~\ref{sec:nq} and \cite{Supp}. For the RQ impedance \eqref{Coleimp} was used to provide the vector $Z$,  whilst  for the LN impedance the data were generated by Matlab's \texttt{integrate()} function with the bounds $0$ and \texttt{Inf}. The resulting Nyquist plots  (complex plot of $Z$) and components $Z_1$ and $Z_2$ are quite similar, see Figures~\ref{1c}-\ref{1e}, and identification of the underlying model from such data, particularly when contaminated by noise, may not be possible. 
 
To simulate noisy measurements, white noise at a level $\eta_j$ was added according to   
\begin{equation}\label{noisy data}
\hat{Z}_{ij}  = Z + \eta_j \bfe_i,
\end{equation}
 where $\hat{Z}_{ij}$ denotes the noisy data, and the vector $\bfe_i$  is the $i^{\mathrm{th}}$ column of the array of size $130 \times 50$   generated using the Matlab function \texttt{randn}, corresponding to $50$ realizations of   white noise,  for noise levels
 $\eta_j$  $j=1:21$  logarithmically spaced between $10^{-6}$ and $10^{-2.5}$.  The choice of the highest noise level used here was determined by comparing the obtained values for $\hat{Z}$ as compared to data seen in practice.  The lowest noise level corresponds to effectively noise-free data, but avoids the {\it inverse crime} by assuring that the data used in the inversion were not exactly prescribed by the underlying forward model. Nonlinear fitting was performed using the Matlab function  \verb/lsqcurvefit()/ initialized with $\beta=0.8$, $\sigma=\ln(2)$, and $t_0=1/\omega_0$, where $\omega_0$ is obtained as the argument of the dominant peak value in $Z_2$, as is commonly used to find the value $t_0$  \cite{Barsukov}. Bounds   $0< t_0<100$, and $0.1<\beta, \sigma<1$  were imposed. Further, scaling of each DRT was introduced through a parameter $\alpha$ satisfying $0<\alpha<1.1$ initialized with $\alpha=1$. 
 
\begin{figure}[h]
\centering
\subfigure[Fitting to the RQ]{\label{2-a}\includegraphics[width=2.1in]{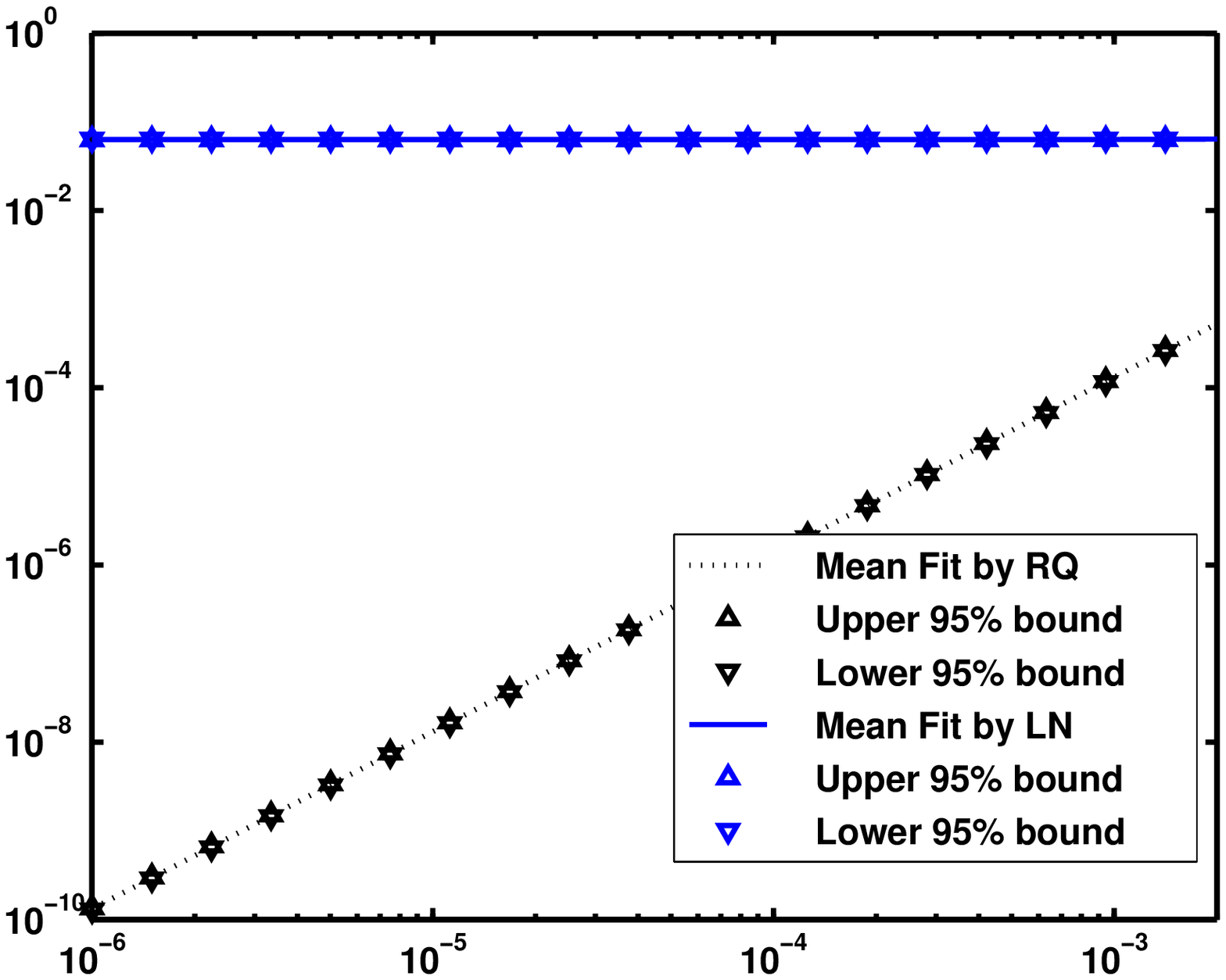}}
\subfigure[Fitting to the LN]{\label{2-b}\includegraphics[width=2.1in]{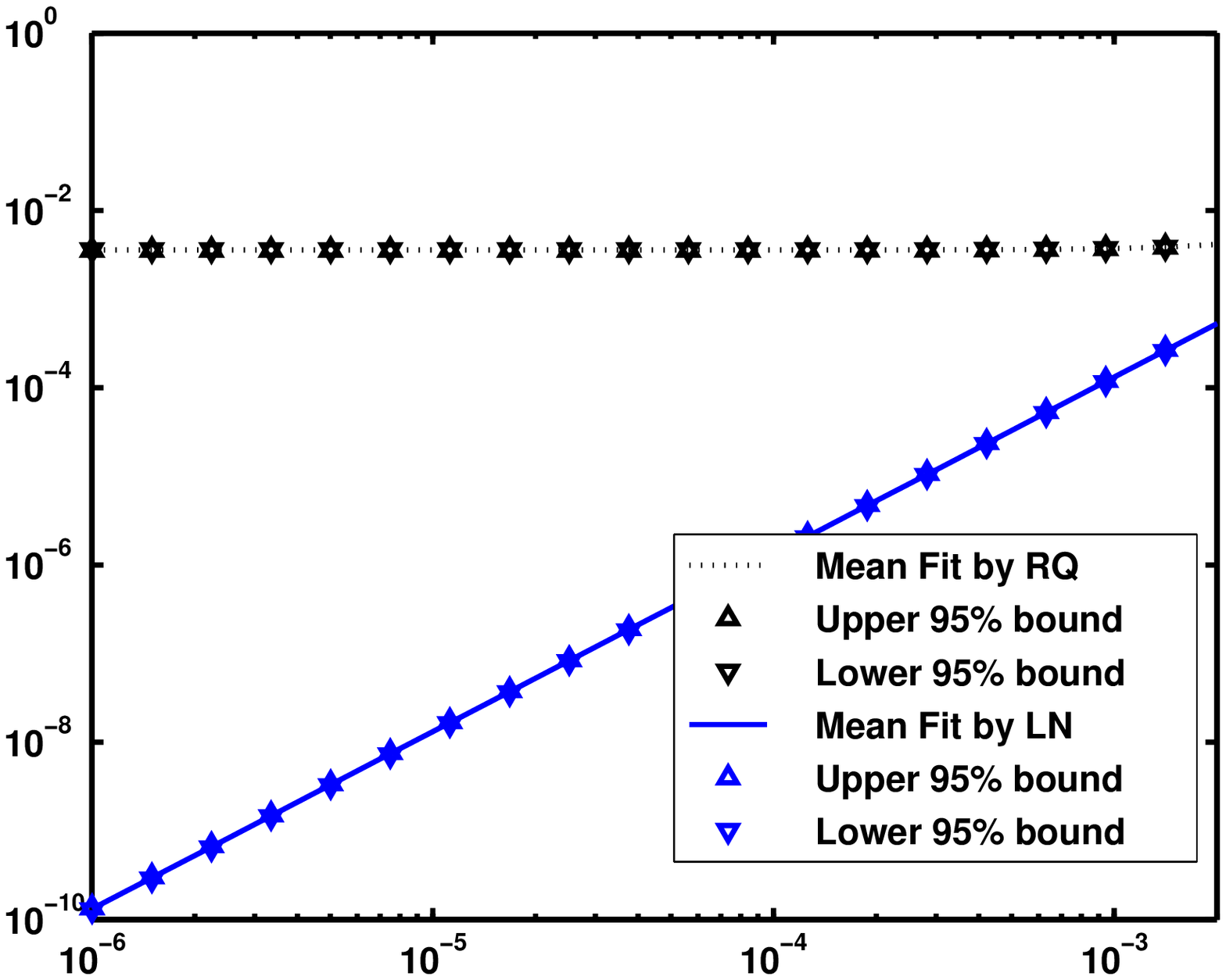}}
\caption{Residual norms fitting one process to the impedance spectrum of a DRT consisting of one   process with white noise in subfigure \ref{2-a} the RQ process and in subfigure \ref{2-b} the LN process. }\label{fig2}
\end{figure}

 Two fittings were performed for each $Z_{ij}$, one assuming the \textit{correct} DRT and one assuming the  \textit{incorrect} DRT, i.e. given $Z$ generated for the RQ DRT a fit was perfumed assuming the \textit{correct} RQ impedance  and the \textit{incorrect} LN impedance values. Similarly for the LN impedance values, \textit{correct} fitting was performed by fitting with a LN DRT, and \textit{incorrect} fitting by  an RQ DRT. When using the RQ DRT for fitting the analytic form of the impedance was used, whilst for the LN DRT all calculations used the Matlab  \texttt{integrate()}   function. Hence for each noise level and realization four fitting pairs were considered, RQ to RQ, LN to RQ, LN to LN, and RQ to LN. For fixed noise level $\eta_j$ and each fitting pair the mean and variance of the residual calculated over the $50$ realizations was calculated. Figure~\ref{fig2} demonstrates the imperfect residuals over all noise levels for the \textit{mismatched} fitting, and the increasing, but relatively smaller,  residuals for the \textit{matched} fittings. Further, fitting to the LN by  RQ yields a smaller residual than fitting   RQ by  LN. In Figure~\ref{fig2} the $95\%$ confidence bounds determined by the variance of the residual are very tight, indicating the robustness of the process. 

\begin{figure}[h!]
\centering
\subfigure[RQ \eqref{DRTtCole} by   LN \eqref{DRTtlog}]{\label{2a}\includegraphics[width=1.5in]{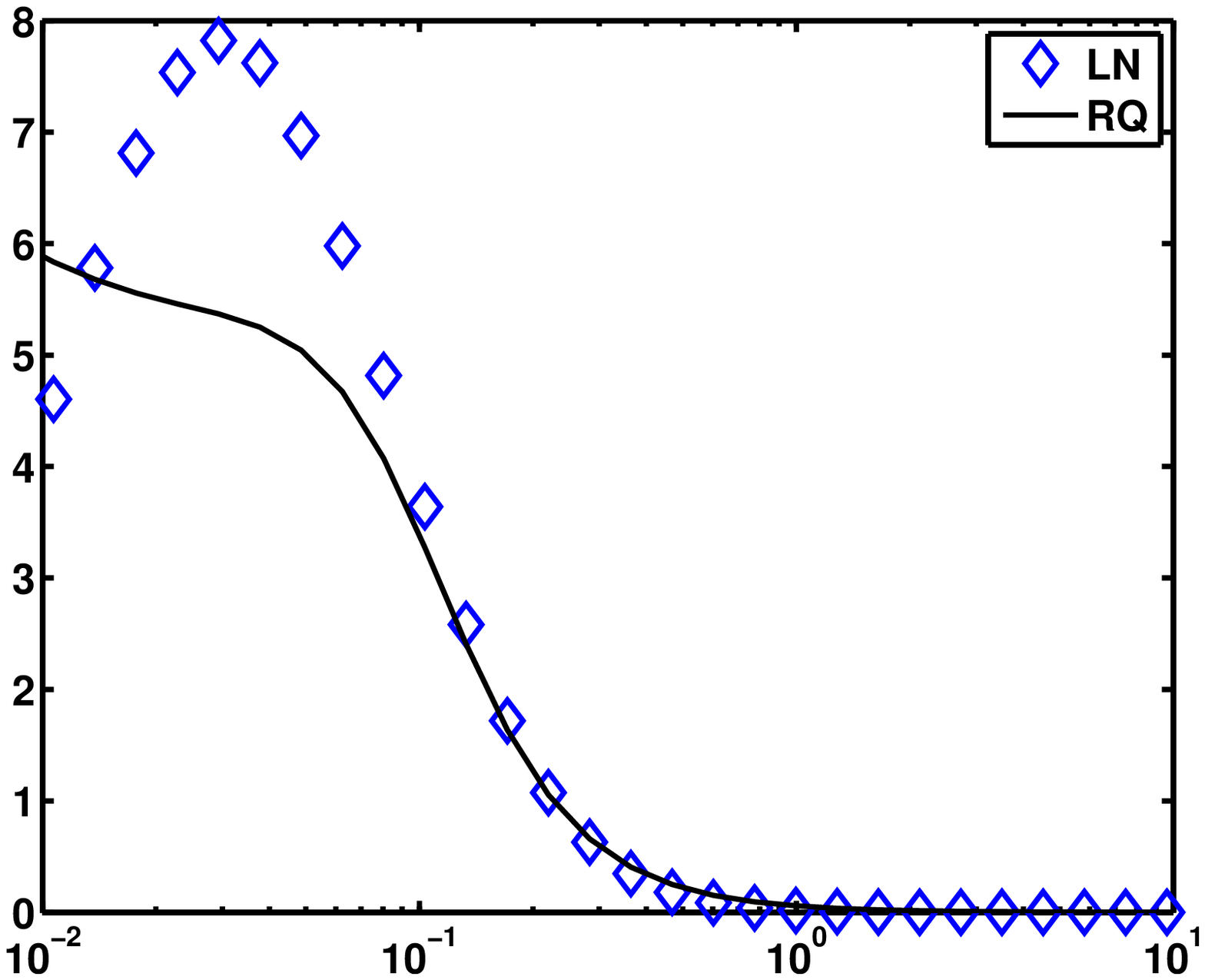}}
\subfigure[LN \eqref{DRTtlog} by   RQ \eqref{DRTtCole}]{\label{2b}\includegraphics[width=1.5in]{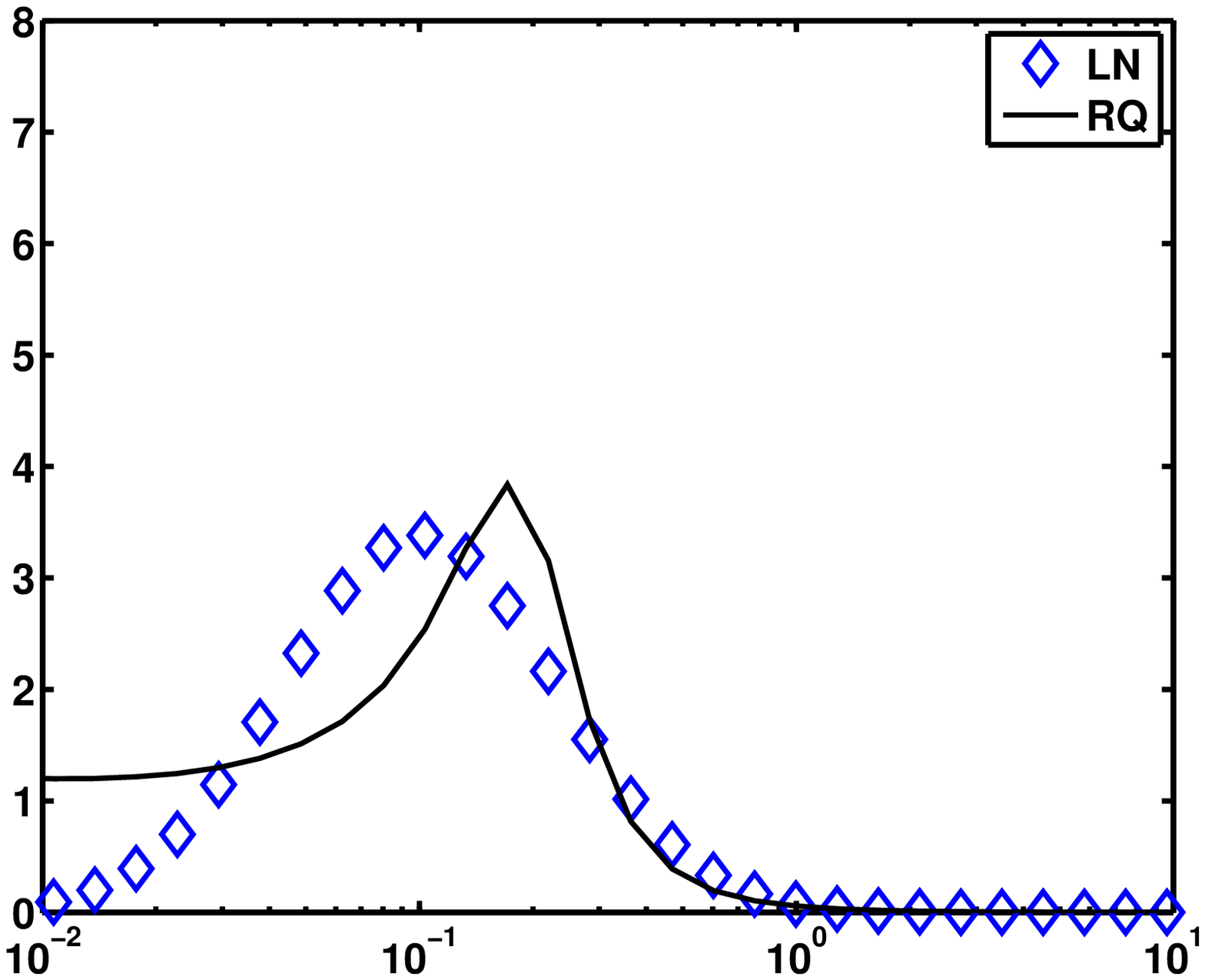}}
\subfigure[RQ  by  LN:  $s-$space]{\label{2d}\includegraphics[width=1.5in]{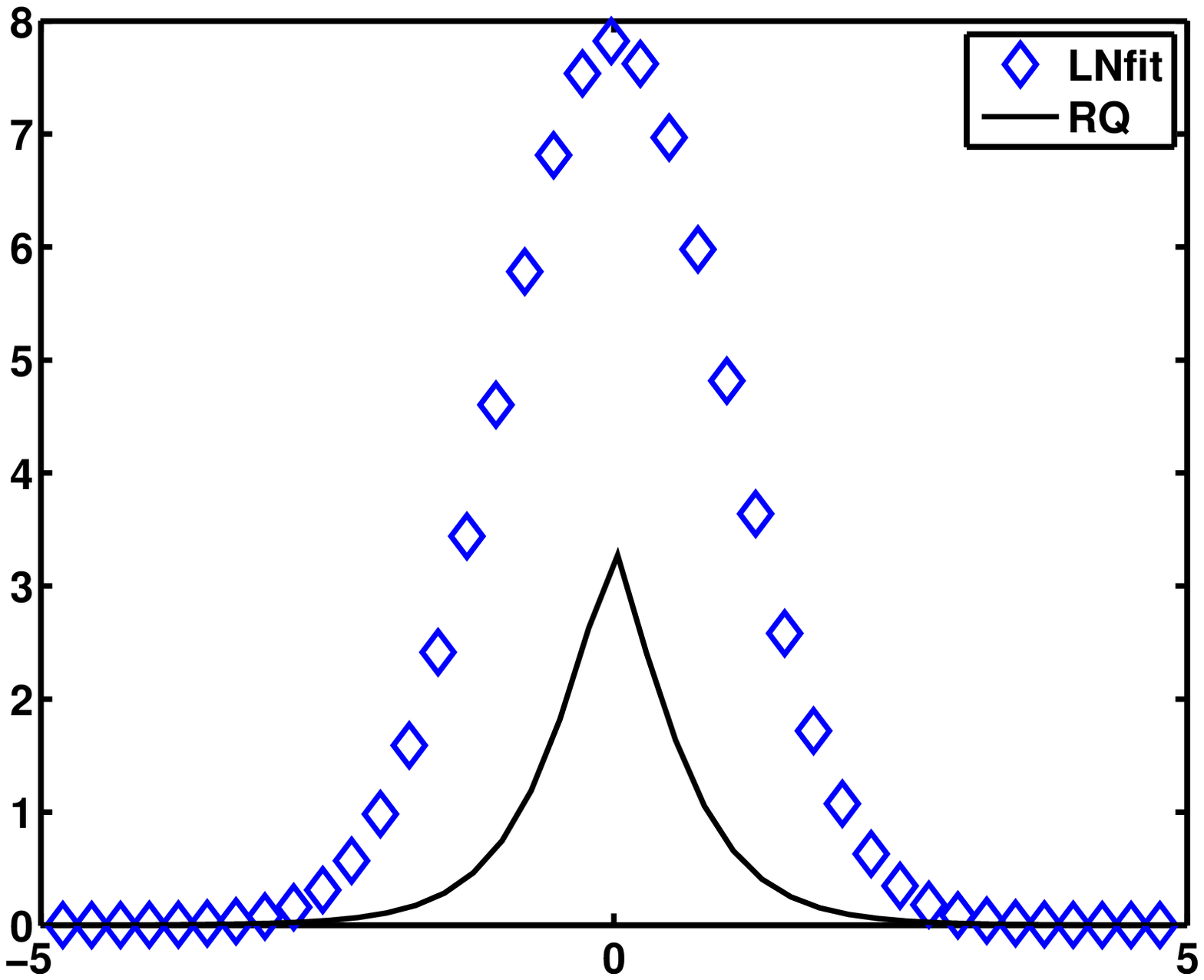}}
\subfigure[LN   by  RQ: $s-$space]{\label{2e}\includegraphics[width=1.5in]{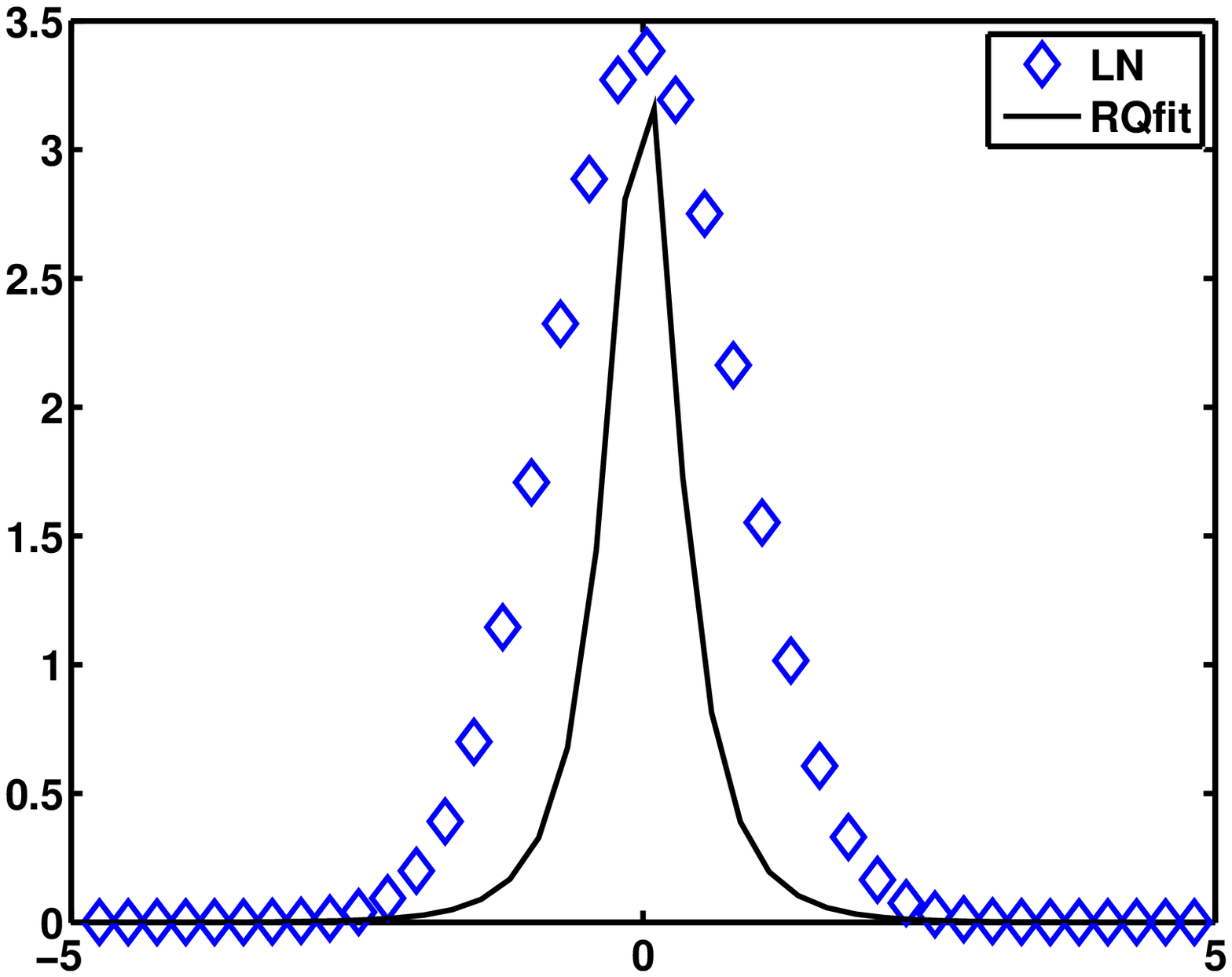}}
\caption{Fitting functions with the consistent mean values obtained and reported in Tables~\ref{RQ.tex}-\ref{LN.tex}. In Figures \ref{2a} and \ref{2d}  fitting the RQ with the LN, and in Figures  \ref{2b}-\ref{2e} fitting the LN with the RQ.
\label{Fits}}
\end{figure}

 To obtain a clearer picture of the quality of the fit, the mean and standard deviation in the estimates for the underlying parameters  for three different noise levels are given in  Tables~\ref{RQ.tex}-\ref{LN.tex}. In each case the \textit{matched} fitting does a good job of parameter estimation for all noise levels, while the  \textit{mismatched} fitting is consistently wrong.  Fitting the RQ impedance with the assumption of a LN DRT generates data that suggests the peak position has moved to the left, and the peak is relatively higher, Figure~\ref{2a}.  Fitting   the LN impedance by a RQ DRT generates a fit  moved to the right   with the height quite well-preserved, Figure~\ref{2b}. However, the processes are still aligned in the $s=\log(t/t_0)$ space, Figures~\ref{2d}-\ref{2e}, each plotted with respect to the identified $t_0$.
 
 We conclude that if the information on the underlying DRT is not known, fitting based on the parameters of the DRT itself will lead to misleading interpretation of the results. Consequently, LLS fitting which simply finds an estimate for the DRT, without finding its parameterization, is the only feasible option for understanding the physical processes when the precise model has not been determined.  We also observe that the fitting as seen in the $s-$space is much more informative in detecting differences and similarities of the DRTs. Of course these results are for the case of a single process. In practice multiple processes are generally exhibited and the impedance fitting is then also less robust in separating the linear combinations, even when the model is predetermined.

\begin{table}[h!] \begin{center}\begin{tabular}{|*{6}{c|}}\hline
&&&\multicolumn{3}{|c|}{Noise level $\eta$}\\ \hline 
 Fit &Parameter&True&$1e-6$&$5.6e-5$&$3.2e-3$\\ 
\hline RQ to RQ& $\beta$&$ 0.72$&$0.72$($5e-07$)&$0.72$($3e-05$)&$0.72$($1e-03$)\\
\hline RQ to RQ&  $t_0$&$0.10$&$0.10$($1e-07$)&$0.10$($6e-06$)&$0.10$($4e-04$)\\
\hline RQ to RQ&  $\alpha$&$1.00$&$1.00$($3e-07$)&$1.00$($1e-05$)&$1.00$($8e-04$)\\ \hline 
\hline LN to RQ&  $\sigma$&$ 0.83$&$1.00$($8e-16$)&$1.00$($3e-14$)&$1.00$($5e-09$)\\
\hline LN to RQ& $t_0$&$0.10$&$0.03$($4e-08$)&$0.03$($2e-06$)&$0.03$($1e-04$)\\
\hline LN to RQ&$\alpha$&$1.00$&$0.97$($2e-07$)&$0.97$($1e-05$)&$0.97$($7e-04$)\\
\hline
\end{tabular}\caption{Mean and standard deviation of absolute errors for obtained parameters for fitting to the RQ DRT.   \label{RQ.tex}}\end{center}  \end{table}

\begin{table}[h!] \begin{center}\begin{tabular}{|*{6}{c|}}\hline
&&&\multicolumn{3}{|c|}{Noise level $\eta$}\\ \hline 
 Fit &Parameter&True&$1e-6$&$5.6e-5$&$3.2e-3$\\ 
\hline RQ to LN&  $\beta$&$ 0.72$&$0.86$($4e-07$) &$0.86$($2e-05$)& $0.86$($1e-03$)\\
\hline RQ to LN& $t_0$&$0.10$&$0.20$($2e-07$) &$0.20$($1e-05$) &$0.20$($6e-04$)\\
\hline RQ to LN& $\alpha$&$1.00$&$1.01$($3e-07$) &$1.01$($1e-05$) &$1.01$($8e-04$)\\ \hline 
\hline LN to LN&$\sigma$&$ 0.83$&$0.83$($2e-06$) &$0.83$($1e-04$)& $0.83$($5e-03$)\\
\hline LN to LN&$t_0$&$0.10$&$0.10$($3e-07$) &$0.10$($2e-05$)& $0.10$($9e-04$)\\
 \hline LN to LN&$\alpha$&$1.00$&$1.00$($2e-07$) &$1.00$($1e-05$)& $1.00$($8e-04$)\\
\hline
\end{tabular}\caption{Mean and standard deviation of absolute errors for obtained parameters for fitting to the  LN DRT.  \label{LN.tex}}\end{center} \end{table} 
 
 \FloatBarrier
\section{Nonparametric Linear Least-Squares} 
\label{sec:ls}
\subsection{Numerical Quadrature}\label{sec:nq}

When a model for the physical system has not been established, a nonparametric method of estimating its DRT must be used. The most straightforward method of discretization, discussed at length in \cite{CSUMS12} with respect to the introduced model error, uses the trapezoidal rule for quadrature with logarithmically spaced points in time to generate  matrices  $A_1$ and $A_2$ approximating the real and imaginary integral operators    $h(\omega,t)=h_1(\omega,t)-i h_2(\omega,t)$ in \eqref{modeleq}. In \cite{Maetal:04} it was suggested to use a change of variables for the integration  before obtaining the quadrature formulae but no discussion or analysis   of the potential advantages or disadvantages was provided. Let $s = \ln t$, then  
 \begin{equation}\label{cov}
Z(\omega) = \int_0^{\infty} h(\omega,t) g(t) \,dt = \int_{-\infty}^\infty h(\omega,e^s) f(s) \,ds, \quad f(s):=t g(t).
\end{equation}
For the  DRTs \eqref{DRTtCole}-\eqref{DRTtlog} we obtain the functions\begin{align}f_{\mathrm{RQ}}(s) &=\frac{1}{  {2 \pi} } \frac{\sin(\beta \pi) }{(\cosh(\beta (s-\ln(t_0)))+\cos (\beta \pi))} \label{DRTsCole}
\\ \label{DRTslog}
f_{\mathrm{LN}}(s) &= \frac{1}{\sigma \sqrt{2 \pi}} \exp{\left(-\frac{(s-\mu)^2}{2\sigma^2}\right)}, \quad \ln(t_0)=\mu -\sigma^2, 
\end{align}
A motivation for this change of variables is to improve the interpretation of the graph of the function when plotted on the linear scale for $s$ as compared to the  logarithmic scale for $t$, see e.g. Figure~\ref{1b}.

In \cite[(8)-(10)]{CSUMS12}  formulae for the  trapezoidal quadrature weights $a_n$ in 
 \begin{align}\label{quadt}
 \int_0^{\infty} h(\omega,t) g(t)  \,dt &\approx  \int_{T_{\mathrm{min}}}^{T_{\mathrm{max}}} h(\omega,t) g(t) \,dt 
 \approx \sum_{n=1}^N a_n h(\omega,t_n) g(t_n),
 \end{align}
  show that $a_n$ are dependent on the logarithmic spacing for $t$, $t_{n+1}-t_{n}=(\Delta t)_{n+1}=t_{n} (10^{\Delta t} -1)$, where $\Delta t$ is constant.\footnote{We note the error in \cite{CSUMS12} which gives these  in terms of $\log$ rather than $\ln$.} 
The same rule   applied for the integral with respect to the $s$ variable, chosen so that $t_n=e^{s_n}$,  gives constant $\Delta s = s_{n+1}-s_n=\ln (t_{n+1})-\ln(t_n)= \ln(10) \Delta t$.
With the standard notation that the double prime on the summation indicates that first and last terms are halved, this yields,  with  $s_1=s_{\mathrm{min}} = \ln\left( {T_\mathrm{min}} \right)$ and $s_N=s_{\mathrm{max}} = \ln\left( {T_\mathrm{max}} \right)$ 
 \begin{align}\label{quads}
\int_{-\infty}^\infty h(\omega,e^s) f(s) \,ds \approx \int_{s_{\mathrm{min}}}^{s_{\mathrm{max}}} h(\omega,e^s) f(s) \,ds \approx \Delta s \sum_{n=1}^N{}^{''}  h(\omega,e^{s_n}) f(s_n).
 \end{align}
 
It is of interest to further investigate the impact of this change of variables on the condition of the resulting systems of equations and on the modeling error obtained from \eqref{quads} so as to justify the use of the $s-$space formulation rather than the $t-$space formulation. This  follows the similar investigation that was presented in \cite{CSUMS12} for \eqref{quadt}. The discretization requires the choice of values for $T_\mathrm{min}$ and $T_\mathrm{max}$, as well as the number of points $N$ used in the discretization of $f(s)$ or $g(t)$. Since in this problem $t$ and $\omega$ have a reciprocal relationship, as noted in, e.g., \cite{Leetal:08, Scetal:02}, we will assume the  range for $t$ is reciprocal to the given range for $\omega$, i.e., $T_\mathrm{max} = 1/\omega_\mathrm{min}$ and $T_\mathrm{min} = 1/\omega_\mathrm{max}$. Forthwith we will use $s_1$ and $s_N$ to denote   $s_{\mathrm{min}}$ and $s_{\mathrm{max}}$, and we  reiterate that $N$ depends on the number of samples for the impedance.

\subsection{Model error}\label{model error}
The model error involved in the discretization of the integral operators stems from two sources: the truncation of the improper integral and the approximation of the integral by a finite quadrature rule. It will be shown here that the model error can be reduced to a negligible level given reasonable assumptions on the DRT. 
\subsubsection{Quadrature error}\label{sec:quaderror}
Bounds for the quadrature error for the logarithmically-spaced trapezium rule \eqref{quadt}  applied for the lognormal $g(t)$ \eqref{DRTtlog} were shown in \cite[(23)]{CSUMS12}.   Because of the use of the variable spacing the error bound for each term of the quadrature varies with $t_n$, and thus   a rapidly decreasing integrand as $t\to\infty$ is necessary in order to control the error. It is well-known that the quadrature error for  the standard constant spacing composite trapezium rule is given by, where we use   $H(s)=h(\omega ,e^{s}) f(s)$ and $|H''(\zeta_n)| := \max_{s\in[s_n,s_{n+1}]} |H''(s)|$,  
\begin{align}  \label{quaderrors}
|E_\text{quad}| &\leq \sum_{n=1}^{N-1} \frac{|H''(\zeta_n)| (\Delta s)^3}{12}  = \frac{(\Delta s)^3}{12} \sum_{n=1}^{N-1} |H''(\zeta_n)|
=  \frac{(s_N-s_1)^3}{12N^2} |H''(\zeta)|
\end{align}
for some $\zeta \in [s_1,s_N] :=I_s$,  assuming appropriate continuity of $H(s)$ on $I_s$ \cite{Atkin}.

\begin{figure}[h!]
\centering
\subfigure[RQ: $t_0 = 0.01, \beta=0.5$.]{\includegraphics[width=2.1in]{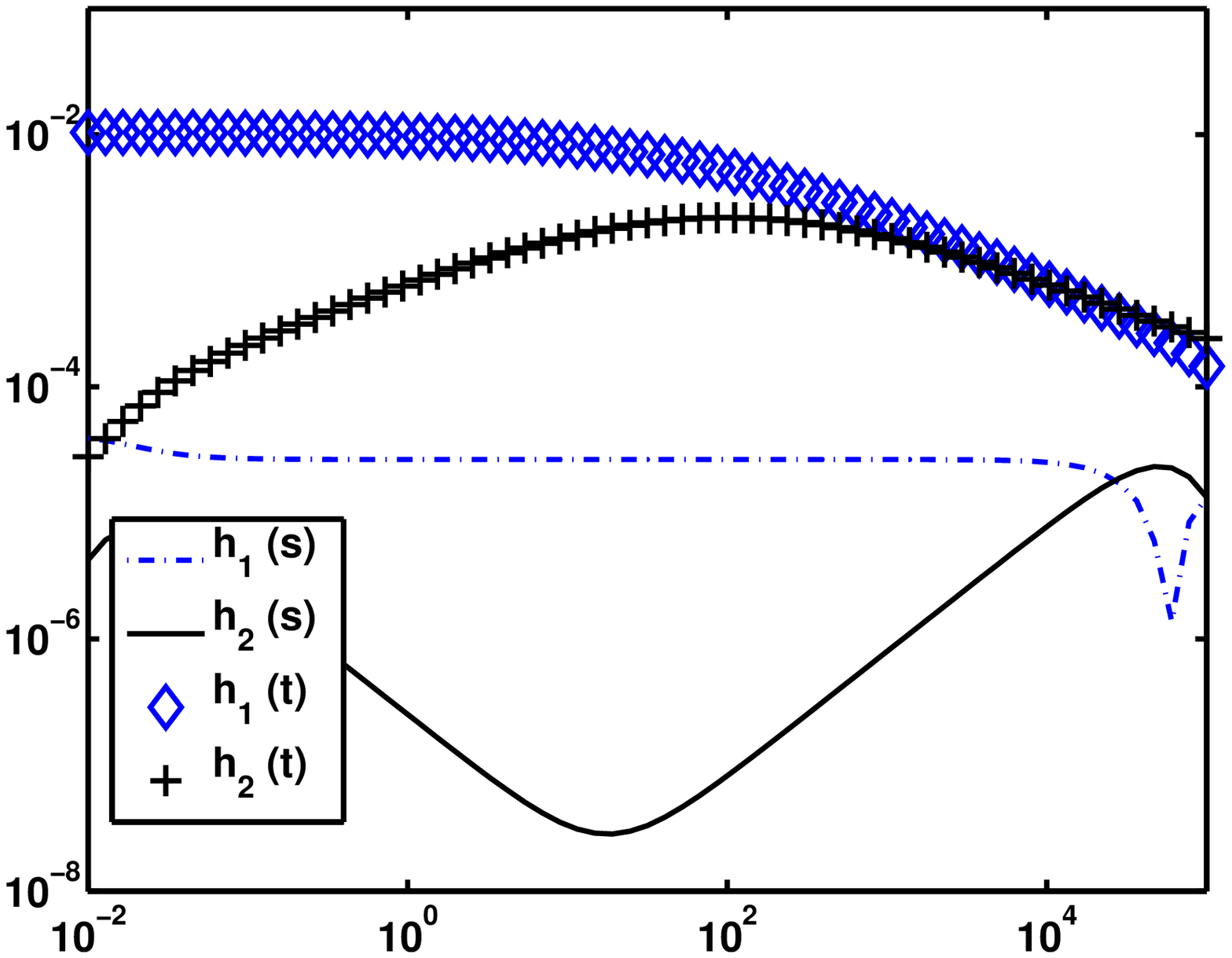}}
\subfigure[RQ: $t_0 = 0.1, \beta=0.5$.]{\includegraphics[width=2.1in]{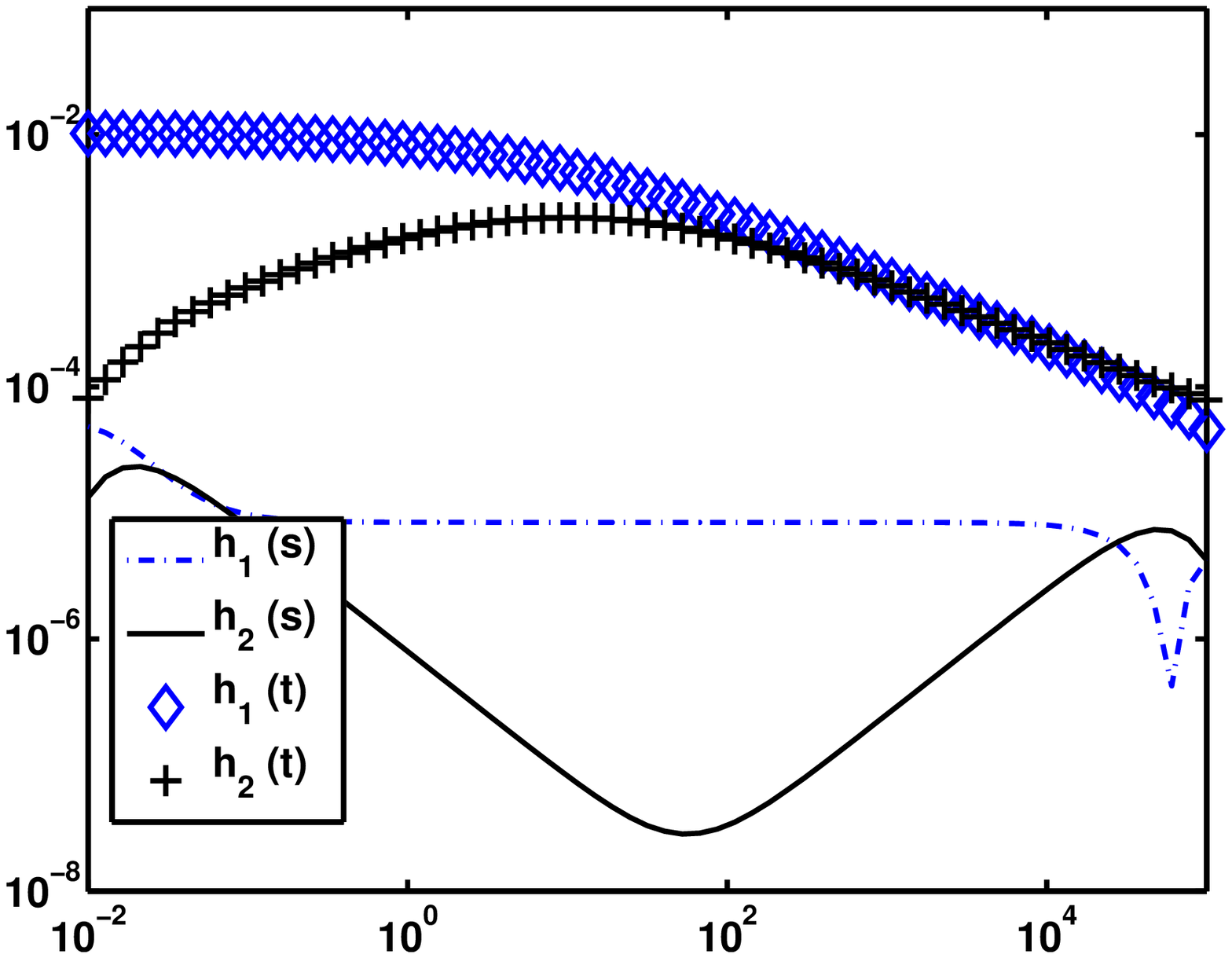}}
\subfigure[RQ: $t_0 = 1.0, \beta=0.5$.]{\includegraphics[width=2.1in]{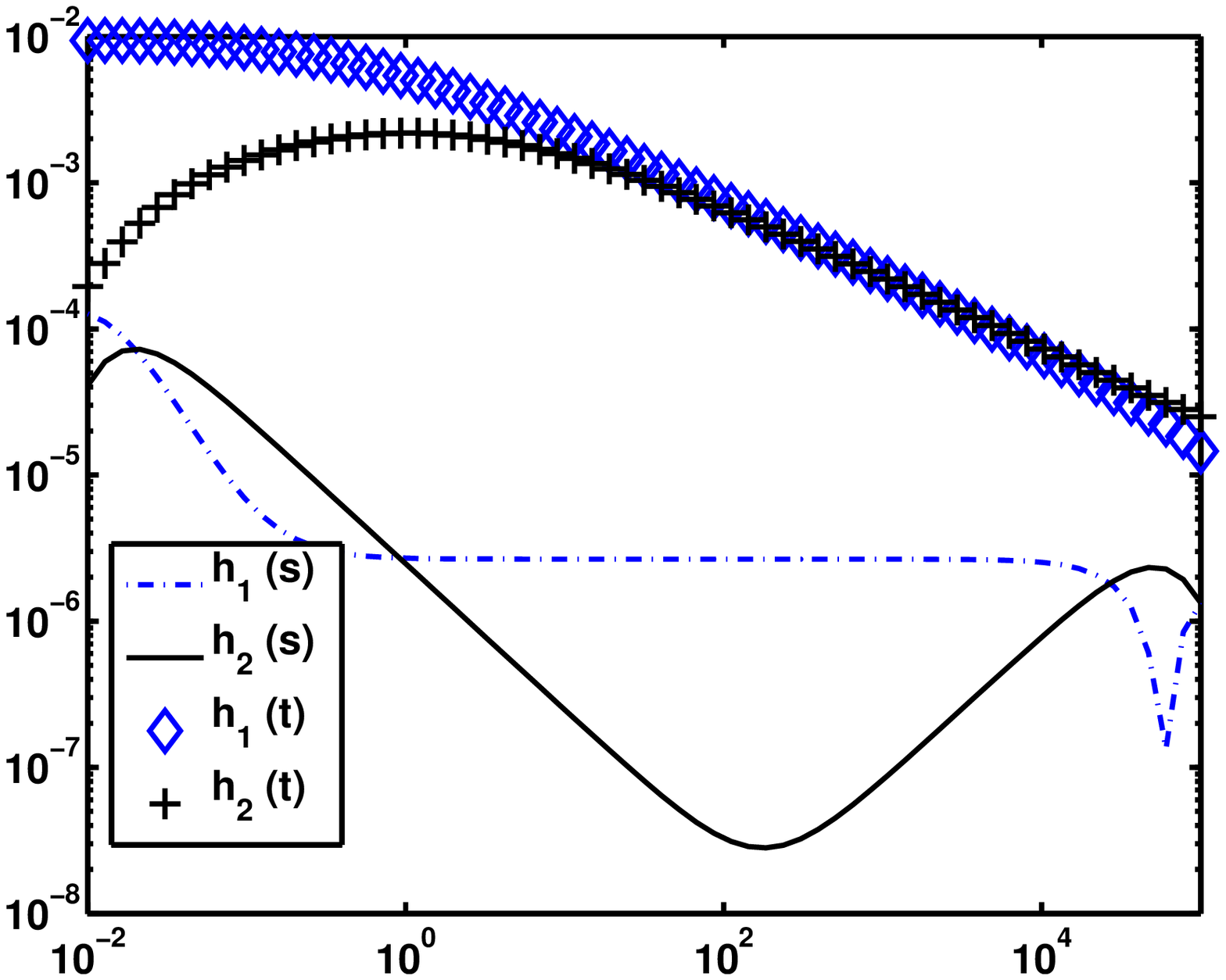}}
\caption{Quadrature error for a single RQ process with $N=65$ for quadrature in $s$ and $t$ as a function of $\omega$, plotted on a log-log scale, for the kernels corresponding to the real and imaginary components of the impedance. }\label{fig-quaderr65}
\end{figure}
\begin{figure}[h!]
\centering
\subfigure[LN: $t_0 = 0.01, \sigma=\ln(3)$.]{\includegraphics[width=2.1in]{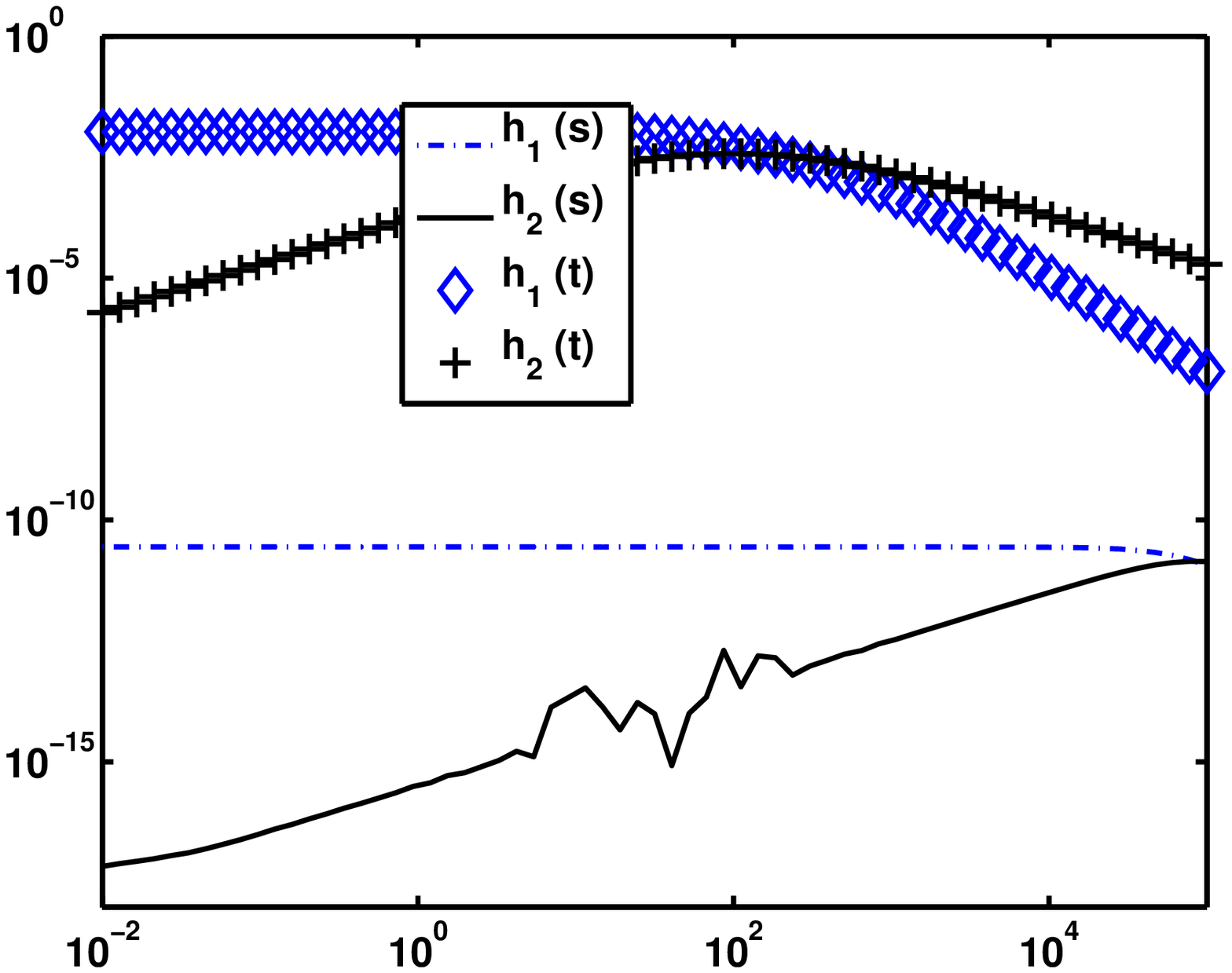}}
\subfigure[LN: $t_0 = 0.1,  \sigma=\ln(3)$.]{\includegraphics[width=2.1in]{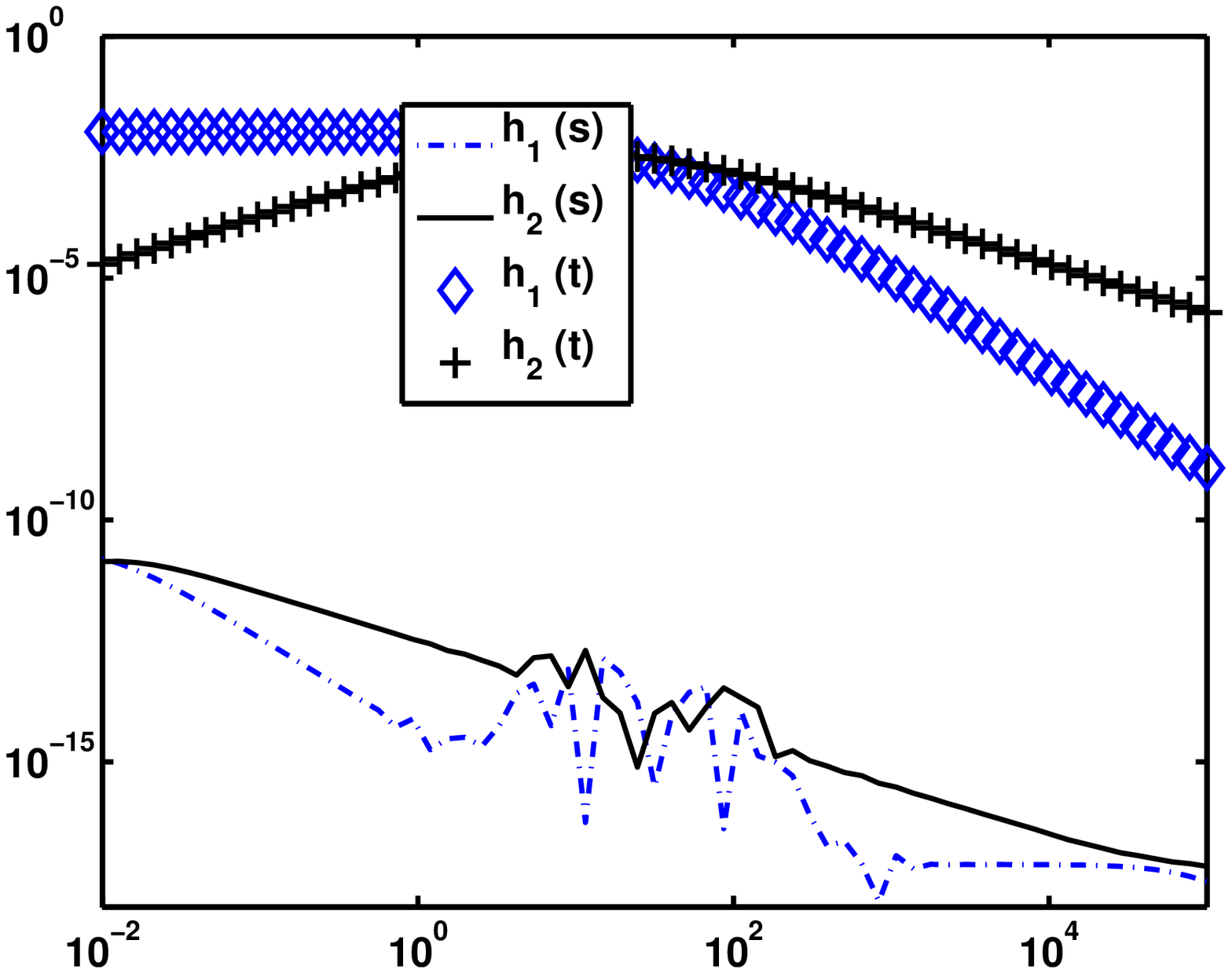}}
\subfigure[LN: $t_0 = 1.0,  \sigma=\ln(3)$.]{\includegraphics[width=2.1in]{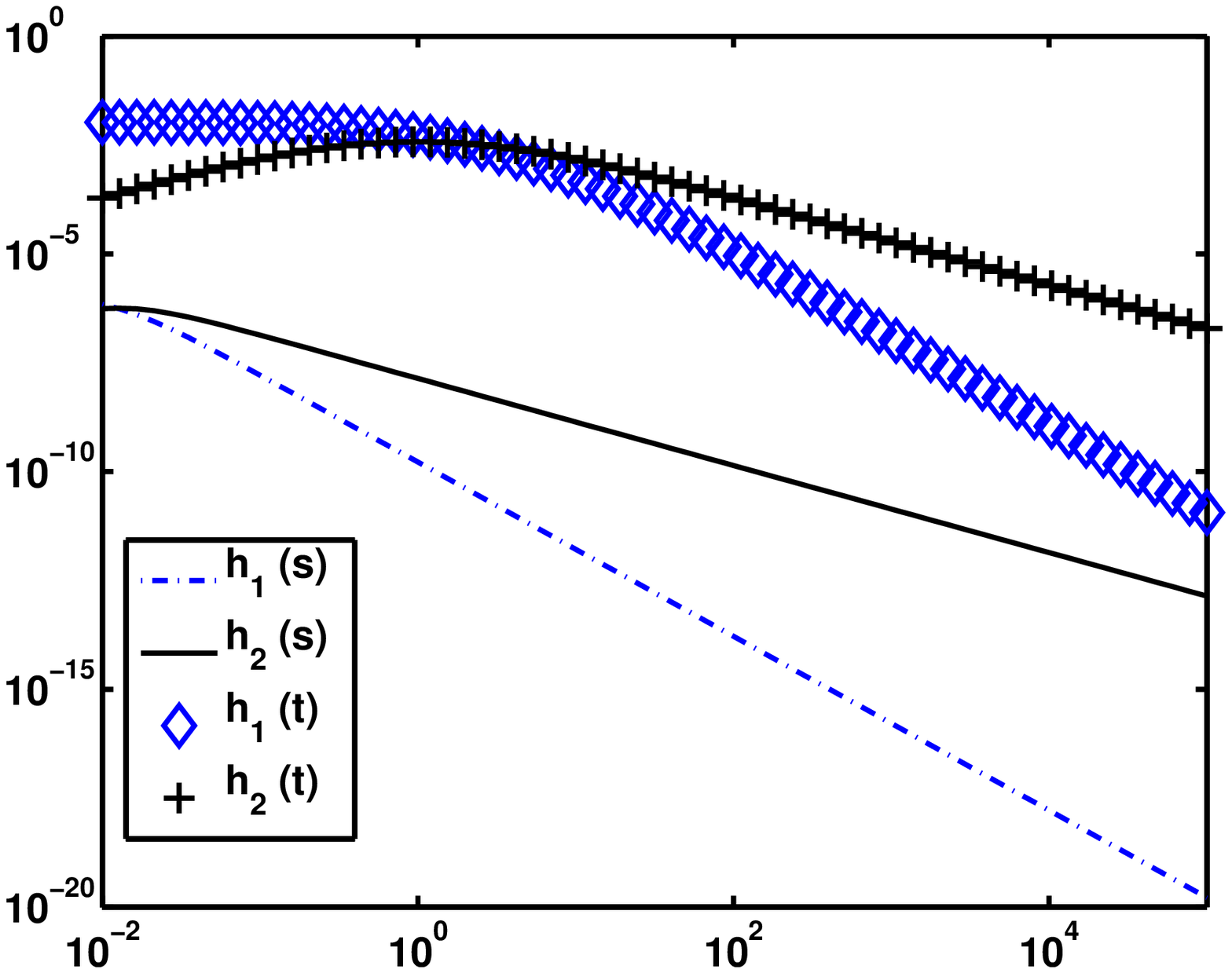}}
\caption{Quadrature error for a single LN process with $N=65$ for quadrature in $s$ and $t$ as a function of $\omega$, plotted on a log-log scale, for the kernels corresponding to the real and imaginary components of the impedance.}\label{fig-quaderrln65}
\end{figure}

  Figures~\ref{fig-quaderr65}-\ref{fig-quaderrln65} contrast the quadrature error as a function of $\omega$  using $N=65$ for \eqref{quadt} and \eqref{quads}, the $t$ and $s$ integrals, respectively, for the two DRTs  \eqref{DRTtCole}-\eqref{DRTtlog}. To obtain the error, an estimate of the {\it true} integral (to  machine epsilon) was calculated using the Matlab function \verb/integral()/. For each formulation of the DRT, the errors   were estimated for a single process for three choices of $t_0$, and a fixed spreading parameter, $\beta=0.5$ and $\sigma=\ln(3)$, respectively. The quadrature error for the integration  evaluated with respect to $s$   is so small as to be negligible for this problem; the error is always less than $10^{-4}$ and $10^{-5}$, for the RQ and LN cases, respectively, and always far less than the error  obtained for quadrature applied to the integration with respect to the variable $t$. We note further that for any integrand which has similar decay properties on the interval,  the same results will apply.

\subsubsection{Improved Quadrature}
Although we have shown that the error is small, the quadrature \eqref{quads} can be further improved as an approximation to the improper integral in \eqref{quads}. It was noted in \cite{Leetal:08} that the quadrature can be improved for the case of kernel $h_2$ by extrapolating the data outside the given interval by a straight line. Because the composite trapezium rule is based on linear interpolation on each interval, it is immediately possible to extend the range of the integration by an interval $\Delta s$, or indeed on any arbitrary interval,  on either side by an extrapolation that assumes  $H(s_N+\Delta s)=H(s_1-\Delta s)=0$, leading to the modification 
 \begin{align}\label{quads2}
\int_{-\infty}^\infty h(\omega,e^s) f(s) \,ds \approx \int_{s_1-\Delta s}^{s_N+\Delta s} h(\omega,e^s) f(s) \,ds \approx \Delta s \sum_{n=1}^N{}  h(\omega,e^{s_n}) f(s_n).
 \end{align}
Here the only change is that the weight on the first and last terms is no longer halved.  This extends the finite range for the quadrature, but does not handle the entire truncation. 
Suppose now that  $f(s)>0$ and $\lim_{s\rightarrow \pm \infty} f(s)=0$, then analytic integration 
gives
 \begin{align}\label{smaxRQ}
 \int_{s_N}^\infty f(s) h_k(\omega, e^{s }) \,ds &\le f(s_N) \int_{s_N}^\infty h_k(\omega, e^{s }) \,ds 
=f(s_N) r_{k,N}(\omega,e^{s_N}), \\
 r_{k,N}(\omega,e^{s_N}) &:= \begin{cases}  \frac{1}{2}\ln (1+(\omega e^{s_N})^{-2}) &k=1 \\
 \frac{\pi}{2} -\tan^{-1} (\omega e^{s_N})) & k=2, \end{cases}
 \end{align}
 while 
 \begin{align}
  \int_{-\infty}^{s_1} f(s) h_2(\omega, e^{s }) \,ds &\le f(s_1)   \int_{-\infty}^{s_1} h_2(\omega, e^{s }) \,ds  
 =f(s_1) r_{2,1}(\omega,e^{s_1})\\
r_{2,1}(\omega,e^{s_1})&:= \tan^{-1}(\omega e^{s_1}) .
 \end{align}
 For the real kernel  at the left hand end,  separating out the kernel integration in the same way is not useful. Instead we must use again the assumption that $f(s)$ decays fast enough that $f(s_1-\Delta s)=0$, (or equivalently that the entire integrand decays fast enough),  and set $r_{1,1}(\omega, e^{s_1})=\frac12 \Delta s\, h_1(\omega,e^{s_1}) $.
Putting these results together leads to the modified quadrature rule which more accurately accounts for the integration outside the  range determined by $I_s$, 
  \begin{align}\label{quads3}
\int_{-\infty}^\infty h_k(\omega,e^s) f(s) \,ds 
  \approx \Delta s \sum_{n=1}^{N}{''}  &h_k(\omega,e^{s_n}) f(s_n)+  f(s_N)r_{k,N}(\omega,e^{s_N})+f(s_1) r_{k,1}(\omega,e^{s_1}).
 \end{align}
 We emphasize the dependence here on the kernel of the integrand in obtaining this result, demonstrating the relative independence of the result on a sufficiently smooth   DRT. 
 \subsubsection{Truncation and Quadrature Error}\label{tqerror}
 The accuracy of \eqref{quads3}  as an approximation to the improper integral in \eqref{quads} depends on the quadrature error  \eqref{quaderrors} and  the truncation error 
 \begin{align*} 
 e_k(\omega)=\left|\int_{-\infty}^{s_1}h_k(\omega,s)f(s)\,ds +\int_{s_N}^{\infty}h_k(\omega,s)f(s)\,ds - f(s_1)r_{k,1}(\omega,s_1)-f(s_N)r_{k,N}(\omega,s_N)\right|.
\end{align*}
  Let $I_{\mathrm{L}}=[-\infty, s_1]$, and $I_{\mathrm{R}}=[s_N,\infty]$, and assume that $I_s$ has been chosen appropriately for the given data, so that $f(s)$ is monotonically decreasing on both $I_{\mathrm{L}}$ and $I_{\mathrm{R}}$, i.e. there exists $\epsilon>0$ such that $|f(s)|\le \epsilon$, $s\in I_{\mathrm{L}} \cup I_{\mathrm{R}}$. Then considering first $k=2$
  \begin{align*} 
e_2(\omega) 
&\le \left( \max_{s \in I_{\mathrm{L}}}\left |f(s) - f(s_1)\right | r_{2,1}(\omega,e^{s_1})\right)+ \left( \max_{s \in I_{\mathrm{R}}}\left |f(s) - f(s_N)\right | r_{2,N}(\omega,e^{s_N})\right) \\
&\le \epsilon  \left (r_{2,1}(\omega,e^{s_1}) +  r_{2,N}(\omega,e^{s_N}) \right ) = \epsilon(\frac{\pi}{2} + \tan^{-1}(\omega e^{s_1})-
\tan^{-1}(\omega e^{s_N}))\le \epsilon \pi
 \end{align*}
providing the total error for kernel $h_2$ 
\begin{align}\label{totalerrorkernel2}
E_2(\omega)\le \epsilon \pi +  \frac{(s_{N}-s_1)^3}{12N^2} |H_2''(\zeta)|.
\end{align}
For $k=1$  
  \begin{align*}\nonumber
e_1(\omega) &\le \left(\int_{-\infty}^{s_1}h_1(\omega,e^s)\,f(s)\,ds -\frac{\Delta s}{2} h_1(\omega,e^{s_1})f(s_1)\right) +  \left(f(s_N)  \frac{1}{2}\ln (1+(\omega e^{s_N})^{-2}) \right).
\end{align*}
For  the first term observe $h_1(\omega,e^s)<1$ for $s\in I_{\mathrm{L}}$ and that the   trapezoidal error bound can be applied. 
Then, assuming  $\int_{-\infty}^{s_1-\Delta s }f(s) \le \delta(f)$, and using   $\zeta \in [s_1-\Delta s, s_1]$, we obtain
\begin{align*}
\int_{-\infty}^{s_1}h_1(\omega,e^s)f(s)\,ds -\frac{\Delta s}{2} h_1(\omega,e^{s_1})f(s_1)&\le \frac{(\Delta s)^3}{12} |H_1^{\prime \prime}(\zeta)| +\frac{\Delta s}{2}f(s_1-\Delta s) + \delta(f),
\end{align*}
and
\begin{align*}
e_1(\omega) &\le \frac{(\Delta s)^3}{12} |H_1^{\prime \prime}(\zeta)| +\frac{\Delta s}{2}\epsilon + \delta(f)  + ( \frac{\epsilon}{2}\ln (1+(\omega e^{s_N})^{-2}) )\\
&\le \frac{\epsilon}{2} (\Delta s + \ln (2)) + \frac{(\Delta s)^3}{12}H_1^{\prime \prime}(\zeta) +\delta(f),
\end{align*}
where we use $\omega e^{s_N} = \omega/\omega_1$. Therefore, now with $\zeta \in I_{\mathrm{L}}$
\begin{align}\label{totalerrorkernel1}
E_1(\omega)\le  \frac{\epsilon}{2} (\Delta s + \ln (2)) +\delta(f)+ \frac{(s_{N}-s_1)^3}{12N^3}(N+1) |H_1''(\zeta)|.
\end{align}
In contrast to \eqref{totalerrorkernel2} this bound depends explicitly on the truncation error $\delta(f)$. On the other hand, for practical $N$ and $I_s$, 
$(\Delta s + \ln (2))<2\pi$.  The two bounds are thus very similar, both depending on the interval $I_s$, and the smoothness of the integrand. 

This analysis which examines the quadrature and truncation error together, contrasts the approach in \cite{CSUMS12} which found the bounds for each error separately, based explicitly on the use of $g_{\mathrm{LN}}(t)$. The same approach can be applied to examine the truncation error for $g_{\mathrm{RQ}}(t)$, and is presented in \ref{appA}. In particular the bound for $\delta(f)$,  which explicitly depends on $f$, is provided by \eqref{lowertruncbd} for $f_{\mathrm{RQ}}$ while the equivalent result for $f_{\mathrm{LN}}$ can be found in \cite{CSUMS12}. We deduce that the total model errors in \eqref{totalerrorkernel2}-\eqref{totalerrorkernel1} may be assumed small for the standard DRTs which decay quickly away from their centers at $t_0$ and have sufficiently small second derivatives. Moreover, these bounds improve on those in \cite{CSUMS12} through  the use of the improved quadrature in $s$ to take account of the extended range, thus reducing the size of  the error introduced by truncating the improper integral. Further, these results are largely independent of the specific DRT, for any DRT satisfying reasonable decay and smoothness assumptions in $s-$space. 

\section{Numerical Algorithms}\label{numerical}
We now turn to the numerical solution of the ill-posed system of equations defined by \eqref{quads3}. In the subsequent discussion, matrices created from the original formulation of the problem will be referred to without superscripts, while those created with the change of variable have a superscript $\mathrm{s}$.

\subsection{Right Preconditioning}\label{rightprec}
Suppose that the  measurements for the impedance are represented by the components of vector $\bfb$ and that the unknowns $f(s_n)$ are components of the vector $\bfx_1$. Then the matrix equation, $A^{\mathrm{s}} \mathbf{x}_1 \approx \bfb$,  for $\bfx_1$ is obtained from \eqref{quads} with $A^{\mathrm{s}}_{mn} = \Delta s{}^{''} \,h(\omega_m,e^{s_n})$, with the double prime indicating the halving for $n=1$ and $n=N$. With the improved quadrature indicated in \eqref{quads3} the components in $A^{\mathrm{s}}_{mn} $ are modified accordingly. In comparison, the discretization for $g(t_n) =f(s_n)/t_n$ is given by  $A \bfx \approx \bfb$ where $\bfx$ has components of $g(t)$ and $\mathrm{diag}(\bft) \, \bfx   =\bfx_1$. Contrasting these two formulations we see that the change of variables is effectively a  {\it right} preconditioning of the original system: $A \,\mathrm{diag}(1./\bft) \bfx_1 =\bfb$. However,  $A \,\mathrm{diag}(1./\bft) \ne  A^{\mathrm{s}}$. The entries in each column differ  in that the weights  $a_n$ in \eqref{quadt} after division by $t_n$ are proportional to $\sinh (\Delta s)$, excepting scale factors for $n=1$ and $n=N$, as compared to $\Delta s$ for matrix $A^{\mathrm{s}}$. When $\Delta s$ is small, for large enough $N$, however, $\Delta s \approx \sinh(\Delta s)$, so that the two matrices are nearly equal \cite{Supp}.  We have already shown that the change of variables impacts the modeling error, we now consider its impact as  a right preconditioner on the stability of the underlying system matrices. 
\subsection{Conditioning}\label{sec:conditioning}
\begin{figure}[h!]
\centering
\subfigure[Picard plot $A_1 \mathbf{x} \approx \mathbf{b}_1$]{\label{6a}\includegraphics[width=2.5in]{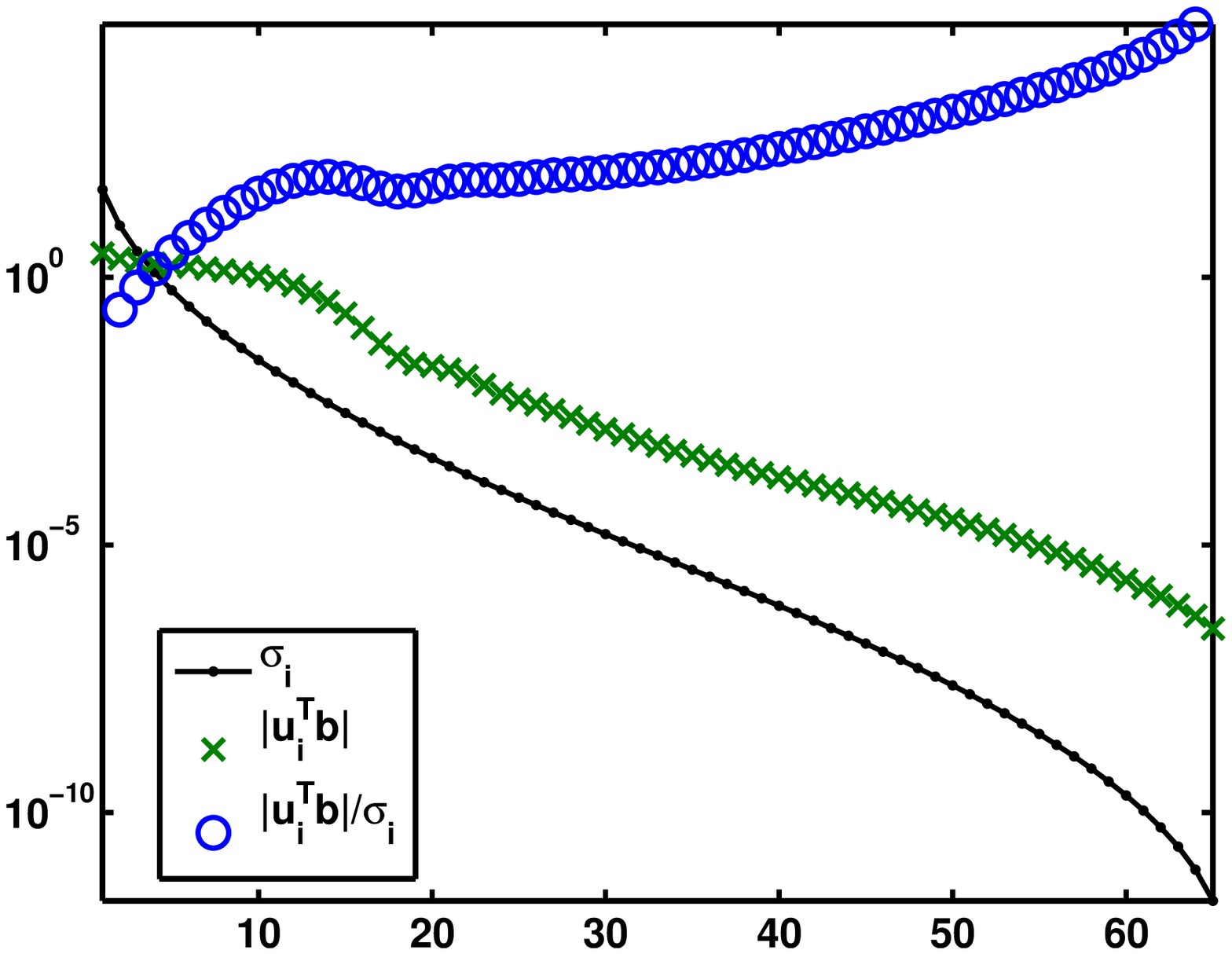}}
\subfigure[Picard plot $A_1^s \mathbf{x} \approx \mathbf{b}_1$]{\label{6b}\includegraphics[width=2.5in]{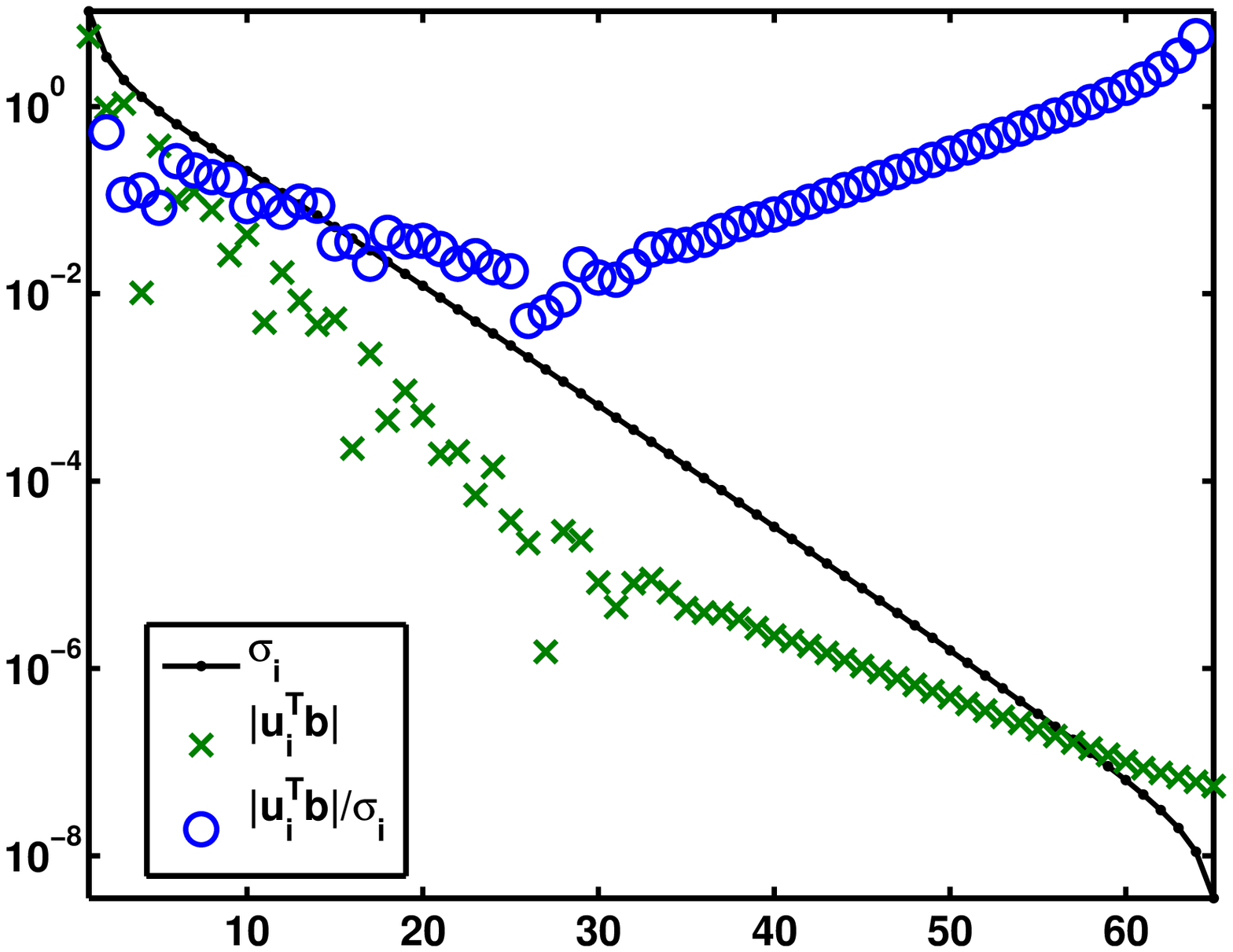}}
\caption{Rght hand side $\bfb$ generated with a single RQ process with $t_0 = 10^{-2}$ and $\beta = 0.8$, using $N=65$ and $\omega$ logarithmically spaced on the interval $[1e-2,1e+5]$. }\label{fig-picardplot}
\end{figure}

The stability of the solution of a system of linear equations $A\bfx \approx \bfb$ is well understood, see e.g. \cite{bjorck,Golub,hansen:1997-1}. Given the singular value decomposition (SVD), $A = U\Sigma V^T$, the na\"ive solution   is 
$
\mathbf{x} = V \Sigma^{-1} U^T \mathbf{b} = \sum_{i=1}^{N} {(\mathbf{u}_i^T \mathbf{b})}/{\sigma_i} \mathbf{v}_i, 
$
where $\mathbf{u}_i$ and $\mathbf{v}_i$ are the $i$th columns of $U$ and $V$, respectively, and $\sigma_i$ is the $i$th singular value of $A$. 
The sensitivity of this solution to errors (noise) in given measurements $\bfb$ can be examined via the Picard plot, \cite{hansen:2007-1}, which illustrates the values of $\sigma_i$, $|\mathbf{u}_i^T \mathbf{b}|$, and the ratio $|\mathbf{u}^T_i \mathbf{b}|/{\sigma_i}$. If   $\sigma_i$ decay faster than $|\mathbf{u}^T_i \mathbf{b}|$, a small perturbation in the true value of $\bfb$ will be amplified and the solution will be dominated by noise.  Figure~\ref{fig-picardplot} compares the Picard plots for a given simulation of $\bfb$ with $N=M=65$ measurements for kernel $h_1$, for both $t$ and $s$ quadrature, matrices  $A_1$ and $A_1^{\mathrm{s}}$, respectively.    For  $A_1$,  $|\mathbf{u}_i^T\mathbf{b}|/\sigma_i$ grow beginning at $i=1$, while for $A_1^{\mathrm{s}}$ they only begin to grow consistently around $i=28$.

\begin{figure}[h!]
\centering
\subfigure[$A_1$: $U$ above and $V$ below.]{\label{7a}\includegraphics[width=3.1in]{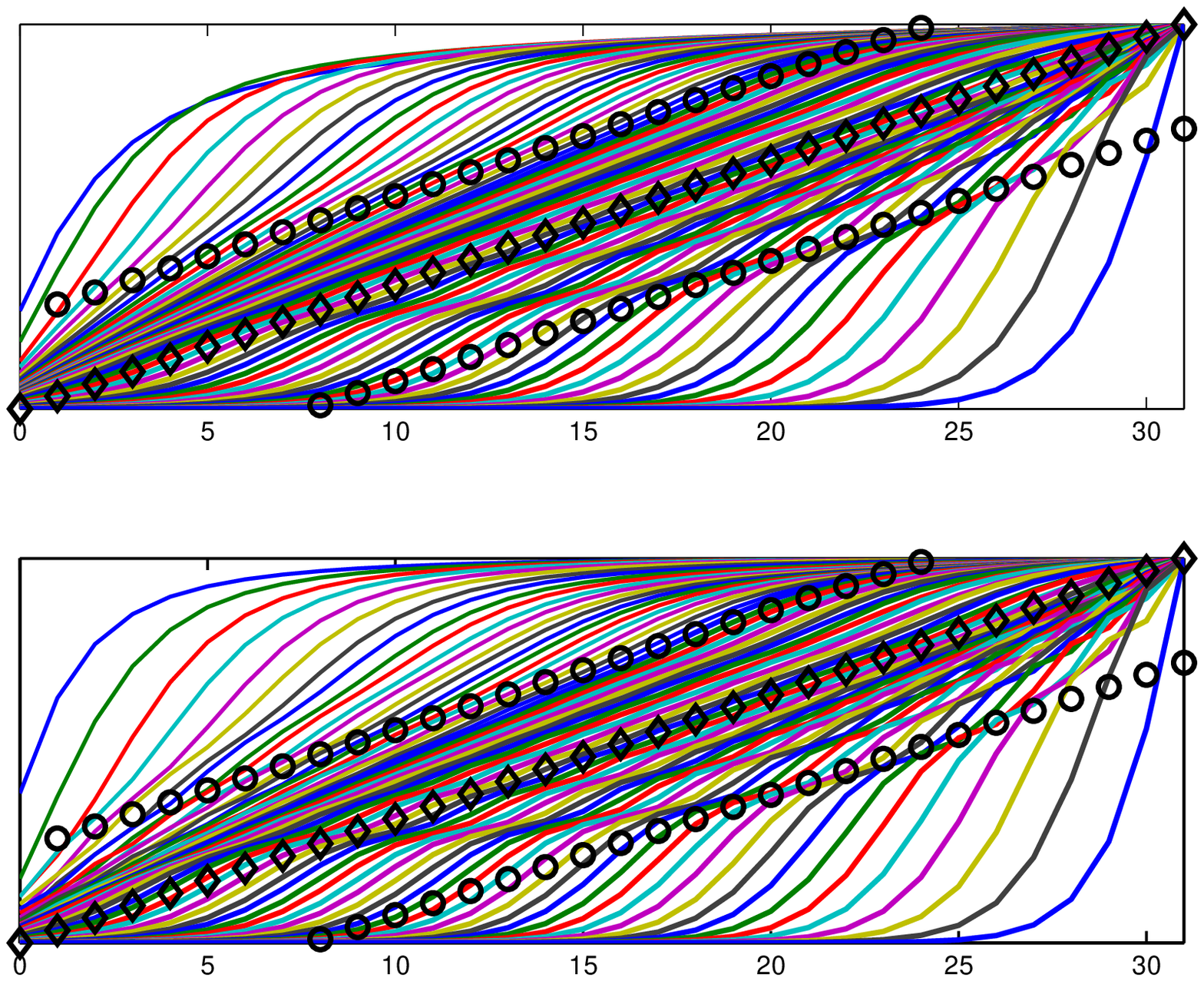}}
\subfigure[$A_1^{\mathrm{s}}$: $U$ above and $V$ below.]{\label{7b}\includegraphics[width=3.1in]{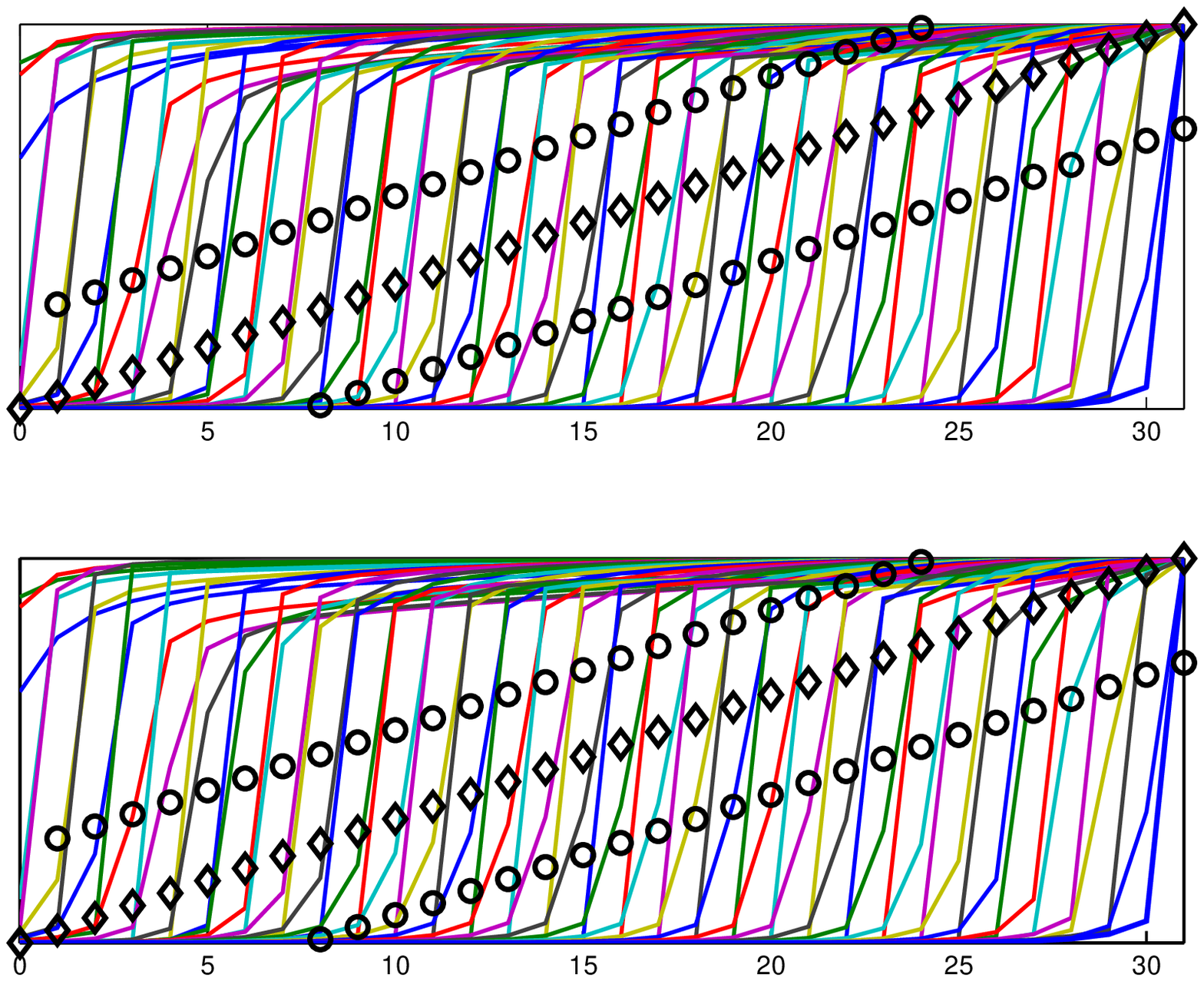}}
\caption{Normalized Cumulative Periodograms and Kolmogorov-Smirnov $95$\% confidence bounds for white noise for the matrices $A_1$ and $A_1^{\mathrm{s}}$. The NCPs for $A_2$ and $A_2^{\mathrm{s}}$ show a  similar separation of the frequency content of the respective basis vectors.}  \label{fig-NCPs}
\end{figure}

The sensitivity shown in Figure~\ref{fig-picardplot} can be exposed in part by examining the condition number, cond($A$), independent of the measurement $\bfb$.  A more informative analysis, however, considers not only the condition but also its impact on the calculation of the basis vectors  for the solution formed from the columns of the SVD matrices $U$ and $V$. It can be shown, e.g. \cite{hansen:2010-1},  that these columns predominantly approximate single frequency components. This frequency content can be visualized by forming the  normalized cumulative periodogram (NCP) for the vector regarded as a discrete time  sampling of a continuous function, as explained in the context of examining the residual vector in \cite{hansen:2006-1, rust:2008-1} and the basis vectors in \cite{CSUMS12}. Vectors which are primarily contaminated by noise have an NCP which falls within  Kolmogorov-Smirnov  bounds for a chosen confidence level, \cite{Fuller:96}. Figure~\ref{fig-NCPs} contrasts the NCPs of the matrices $A_1$ and $A_1^{\mathrm{s}}$, illustrating the better separation of the frequency content for the basis vectors formed from matrices $U$ and $V$ for matrix $A_1^{\mathrm{s}}$. This suggests  greater confidence in the use of a higher number of basis terms in the solution with the   right preconditioned matrix. A similar conclusion follows when comparing the matrices for the kernel $h_2$. 

Further demonstration of the impact of the right preconditioning is provided in Table~\ref{constable} which gives the condition of matrices for the $t$ and $s$ quadrature, for a selection of choices for $I_t=[T_{\mathrm{min}},T_{\mathrm{max}}]$ with $N=65$. It is clear that the choice for $I_\text{t}$, and hence $I_s$,   has a large impact on the conditioning of the problem. For instance, if the range for $t$ is too wide, the matrices will have several nearly linearly dependent columns, which greatly increases the condition number of the matrix, and increases the dimension of the numerical null space. Likewise, if the range for $t$ is too narrow, there will be nearly linearly dependent rows, again increasing the condition number of the matrix. When the reciprocal relationship $t=1/\omega$ is used to pick the  sampling in $t$  the optimal condition is obtained in both cases, as shown by the bold face row of the table.

In general there are relatively few measurements of the impedance for the particular biofuel application that can be used to find the solution. Typically there are on the order of $N=65$ usable values for each of $Z_1$ and $Z_2$, after estimates for the resistances $R_{\mathrm{pol}}$ and $R_0$ in \eqref{modeleq}  have been found.  It was shown in \cite{CSUMS12} that in order to make use of all available data, it is generally desirable to combine the matrices $A_1$ and $A_2$ into the stacked matrix 
$A_3 = [A_1 ; A_2 ]$ of size $2N \times N$, hence using systems associated with both kernels. Further, the overdetermined system can be replaced by a system using a square matrix $A_4$ of size $2N \times 2N$ through increasing the sampling used in the quadrature for the kernel equations, hence providing increased resolution in the discretization of the solution vector.   The same approach can be  used for the $s-$quadrature matrices $A_1^{\mathrm{s}}$ and $A_2^{\mathrm{s}}$ yielding matrices $A_3^{\mathrm{s}}$ and $A_4^{\mathrm{s}}$. In practice, therefore, it is the impact of the scaling on these augmented matrices which is also significant.

The results in Table~\ref{constable} are  given for the scaled matrices calculated using \eqref{quads}. To see that the improvement in the quadrature has little impact on the overall conditioning, the condition numbers of the matrices $A_1$ to $A_4$ for various quadrature rules, as noted by the respective equation numbers, are also   shown in Table~\ref{conditions}.  Overall, the right preconditioned matrices, regardless of selection of the system,  exhibit  better but still not ideal conditioning, so that the problem remains ill-posed and the solution   is subject to noise-contamination.

\begin{table}[ht]
\begin{center}
\begin{tabular}{|c|c|c|c|c|c|c|cccc}
\hline
$T_\text{min}$&$T_\text{max}$ & $A_1$ & $A_2$ & $A_1^{\mathrm{s}}$ & $A_2^{\mathrm{s}}$ \\
\hline
\multirow{3}{*}{$1e-6$}
&$1e+1$&$9.70e+17$&$6.28e+18$&$2.83e+17$&$1.90e+17$ \\\cline{2-6}
&$1e+2$&$6.47e+17$&$1.05e+20$&$2.99e+17$&$2.41e+17$ \\\cline{2-6}
&$1e+3$&$4.73e+18$&$7.45e+19$&$4.58e+20$&$6.22e+18$ \\\hline
\multirow{3}{*}{$\mathbf{1e-5}$}
&$1e+1$&$2.06e+19$&$1.74e+18$&$8.17e+17$&$6.40e+17$ \\\cline{2-6}
&$\mathbf{1e+2}$&$\mathbf{1.94e+13}$&$\mathbf{1.78e+13}$&$\mathbf{2.94e+09}$&$\mathbf{7.43e+07}$ \\\cline{2-6}
&$1e+3$&$2.20e+18$&$3.64e+19$&$2.12e+20 $&$1.87e+18$ \\\hline
\multirow{3}{*}{$1e-4$}
&$1e+1$&$5.05e+21$&$1.55e+20$&$4.38e+19$&$2.44e+18$ \\\cline{2-6}
&$1e+2$&$1.52e+21$&$1.99e+18$&$1.09e+20$&$3.90e+18$ \\\cline{2-6}
&$ 1e+3$&$6.49e+19$&$6.98e+18$&$1.30e+21$&$1.41e+20$ \\
\hline
\end{tabular}
\caption{Condition numbers for matrices $A_1$, $A_2$, $A_1^{\mathrm{s}}$, and $A_2^{\mathrm{s}}$ for $I_\text{t} = [T_{\mathrm{min}},T_{\mathrm{max}}]$ and $\omega$ logarithmically spaced on the interval $[1e-2,1e+5]$.}\label{constable}
\end{center}
\end{table}

\begin{table}[ht]\begin{center}\begin{tabular}{|*{7}{c|}}
\hline Quad&Equation&$A_1$&$A_2$&$A_3$&$A_4$\\
\hline t &\eqref{quadt}&$1.5e+13$&$1.4e+13$&$7.5e+12$&$4.1e+20$\\
\hline s&\eqref{quads}&$2.9e+09$&$7.4e+07$&$4.6e+08$&$2.3e+18$\\
\hline s&\eqref{quads2}&$3.1e+09$&$7.9e+07$&$5.0e+08$&$9.1e+17$\\
\hline s&\eqref{quads3}&$2.8e+09$&$7.4e+07$&$4.6e+08$&$9.0e+18$\\

\hline
\end{tabular}\caption{\label{conditions}Comparing condition number of matrices with different quadratures for the optimal selection of the nodes for $t_n$, using  $t=1/\omega$. In this case, as compared to Table~\ref{constable} the matrices are also reduced to size $N=64$, consistent with the loss of data when estimating $R_0$ and $R_{\text{pol}}$ in \eqref{modeleq}. 
}\end{center}\end{table}

 \FloatBarrier
\subsection{Regularization}\label{regsec}
Given the ill-conditioning of the system matrices, regularization is required in order to select an acceptable solution which provides  a reasonable fit to the measured data, but is at the same time controlled in its growth with respect to a chosen norm. 
There is an extensive literature discussing multiple formulations, e.g.  \cite{hansen:1997-1, hansen:2007-1, hansen:2010-1, LawsonHansen:1995,Vogel:02} and, for problems of this type, Tikhonov regularization is frequently applied  e.g. \cite{Leetal:08, CSUMS12, Weese:92}. The solution is then recast as the solution of the regularized problem   
\begin{equation}\label{regsoln}
\mathbf{x} = \arg\min \{ \|A\mathbf{x} - \mathbf{b}\|^2 + \lambda^2 \|L\mathbf{x}\|^2\} ,
\end{equation}
where $\lambda$ is a parameter affecting the amount of regularization applied and $L$ is a matrix chosen so the growth of $\mathbf x$ is controlled relative to the $L-$weighted norm. Typical choices for $L$ are approximations to derivative operators of order $0$, $1$, or $2$. For small problems the $\lambda$-dependent solution of \eqref{regsoln} can be expressed in terms of the SVD, ($L=I$), or  generalized singular value decomposition (GSVD), for other noninvertible choices for $L$, e.g.   \cite{hansen:1997-1}.

The choice of the regularization parameter $\lambda$ is a nontrivial but well-studied  problem. There are many algorithms available to choose this parameter. In general, these can be divided into two categories: those which require some knowledge of the characteristics and magnitude of the noise present in the right hand side, and those which do not. Techniques such as the discrepancy principle \cite{hansen:1997-1},  the unbiased predictive risk estimator (UPRE) \cite{Vogel:02} and the $\chi^2$ method \cite{MeRe:09} fall in the first category, while more heuristic methods such as the LC criterion, the quasi-optimality criterion, and generalized cross-validation (GCV) fall in the second, \cite{hansen:2007-1}. The NCP criterion applied to the residual for the data fit, which is  based on work in \cite{rust:1998-1} and then extended in  \cite{hansen:2006-1,rust:2008-1}, falls between these two categories. It requires that the noise be white, or  can be whitened,  but does not require an estimate of its magnitude.

The LC  criterion is commonly used when no information on the noise distribution is available. It relies on the fact that when plotted on a logarithmic scale, the weighted norm  $\|L\mathbf{x}\|$ and the data fidelity  norm $\|A\mathbf{x}-\mathbf{b}\|$ tend to form an L-shaped curve. That is, for a value of $\lambda$ near the corner, an increase in $\lambda$ would tend to increase the residual norm without reducing the solution norm much, while a decrease in $\lambda$ would increase the solution norm without much reducing the residual norm. The SVD and GSVD allow the simple construction of the L-curve, and the corner is generally determined by finding the point of maximum curvature, as implemented in \cite{hansen:2007-1}.  

The NCP approach is based on the assumption that for an appropriately chosen $\lambda$ the noise in the measurements is transferred to the residual vector so that it is completely dominated by white noise.  By calculating the NCP of the residual for suitably selected values of $\lambda$, the smallest value which produces a residual vector reasonably approximating white noise may be found. A common statistical test for white noise uses the Kolmogorov-Smirnov statistic to compare the residual NCP with the   cumulative distribution function of a uniform distribution. This test is described in more detail in, e.g., \cite{Fuller:96}.

\subsection{Non-negative Least-Squares}\label{nnlssec}
It was suggested  in \cite{Maetal:04} that  the solution of the DRT problem  should use the  additional information on the DRT, namely that 
$g(t)$ (and its cousin $f(s)$) are distribution functions satisfying $g(t) \geq 0$ for all $t$. Augmenting \eqref{regsoln} by this constraint, yields the NNLS formulation 
 \begin{equation}\label{nnls}
\bfx = \arg\min \{ ||A\mathbf{x} - \mathbf{b}||^2 + \lambda^2 ||L\mathbf{x}||^2,  \text{ s.t. } \mathbf{x} \geq 0\},
\end{equation}
for which efficient algorithms are available, given a specific choice for $\lambda$, including for instance the algorithm in \cite{LawsonHansen:1995}, which is implemented in Matlab. 
On the other hand, because we cannot immediately express the solution in terms of an expansion, many common parameter choice methods, such as the GCV and UPRE,  are less feasible for finding a suitable $\lambda$. For this small problem, $N\approx 65$, we explore the use of  the two parameter choice methods, the LC criterion and  the NCP analysis of the residual vector.

For both LC and NCP,  a standard approach for finding the optimal regularization parameter can still be applied. Specifically, solutions of \eqref{nnls} can be found for a range of values for $\lambda$. Then the   LC or NCP, respectively, is applied to assess the quality of each solution, as briefly described for \eqref{regsoln} in Section~\ref{regsec}. It is only necessary that an algorithm for the solution of \eqref{nnls} is available, e.g.    \texttt{lsqnonneg} which is the implementation of \cite{LawsonHansen:1995} in Matlab. To verify that our approach is indeed robust to the choice of the algorithm we also consider solutions using the CVX package, \cite{cvx,gb08}.
This  contrasts with the parameter choice method presented in \cite{Maetal:04} that relies on a  specific implementation of the  NNLS algorithm.  

\section{Simulation Results} 
The performance of the algorithms discussed and analyzed in Section~\ref{sec:ls} is investigated for  simulated data sets exhibiting properties seen in practical situations. The data sets are  described first in Section~\ref{datasets}. For each data set extensive comparisons have been performed using both matrices $A_3$ and $A_4$, in order to assess whether the better conditioning of $A_3$ or the better resolution offered by $A_4$ prevails as the optimum in each case. These examinations have demonstrated that the additional resolution of $A_4$  outweighs its worse conditioning as compared to $A_3$, (see Table~\ref{constable}). While in many cases the results using $A_3$ are quite comparable, as noise increases systems solved with matrix $A_4$ often provide better results. Thus in the presented results here we emphasize the comparisons between the LLS and NNLS methods and do not give results for systems with matrix $A_3$. Further results substantiating these comments are given in the supplementary materials \cite{Supp}.

\subsection{Data sets}\label{datasets}
\begin{figure}[!ht]
\centering
\subfigure[A: $g(t)$]{\includegraphics[width=1.3in]{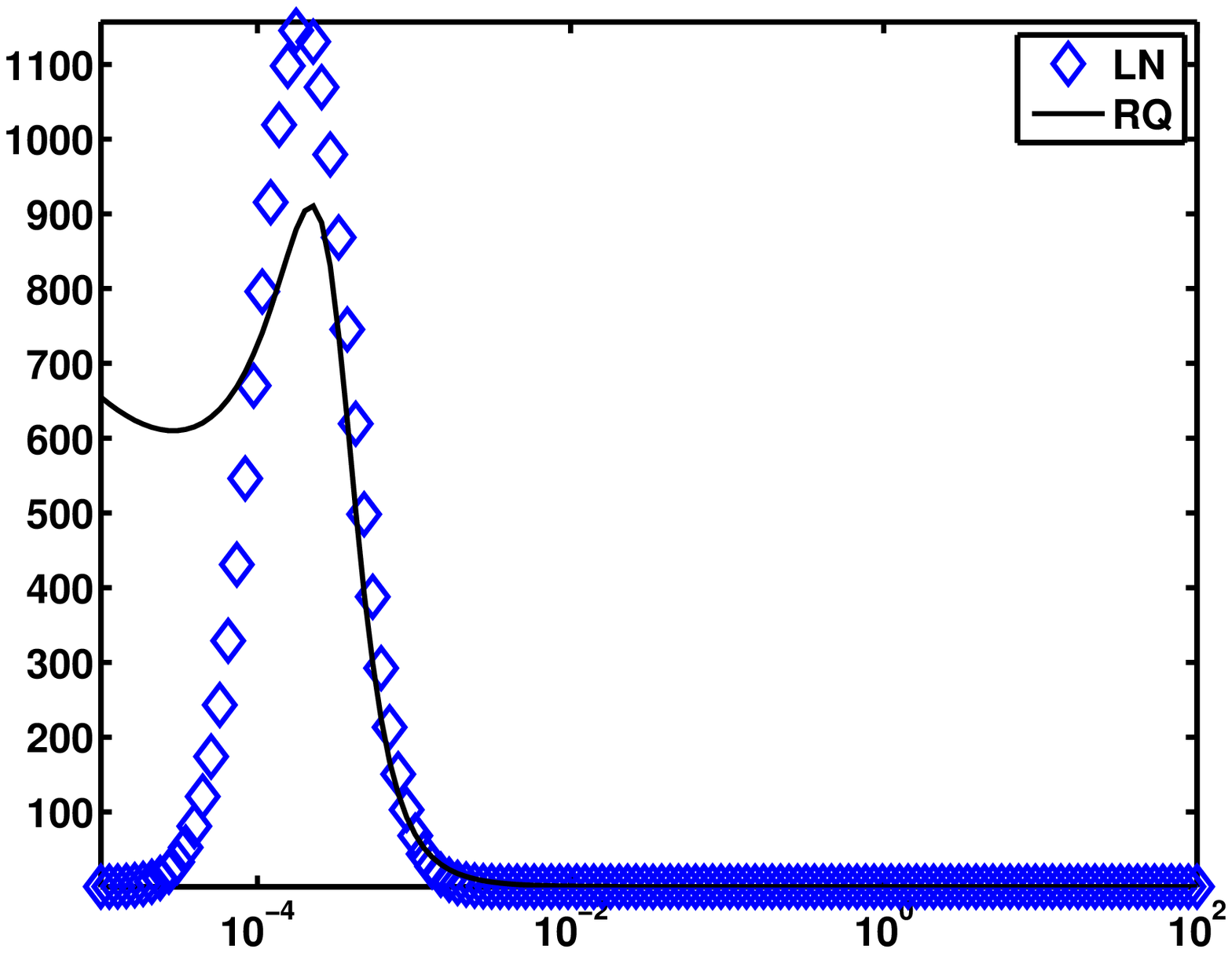}}\subfigure[$tg(t)$]{\includegraphics[width=1.3in]{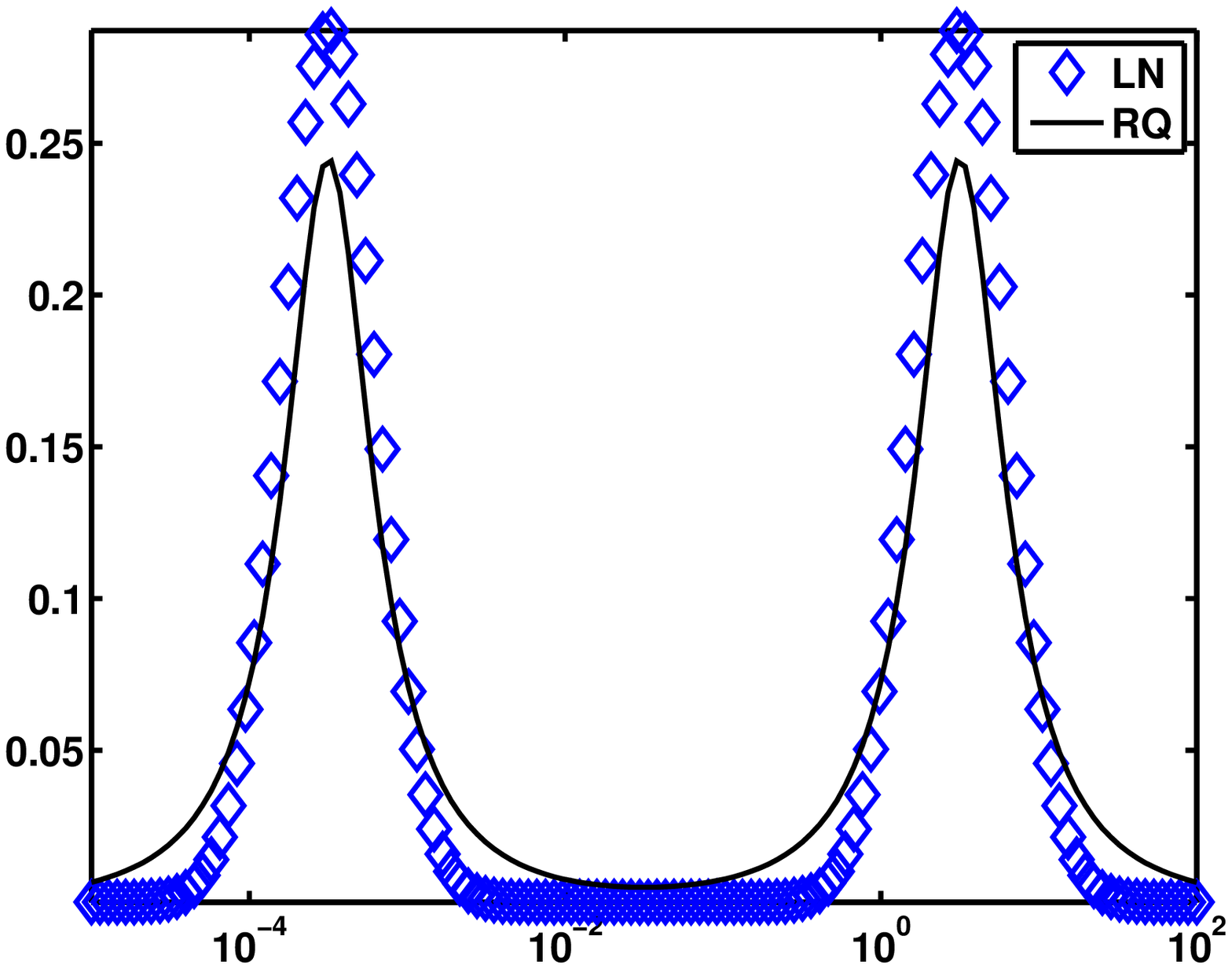}}\subfigure[Nyquist]{\includegraphics[width=1.3in]{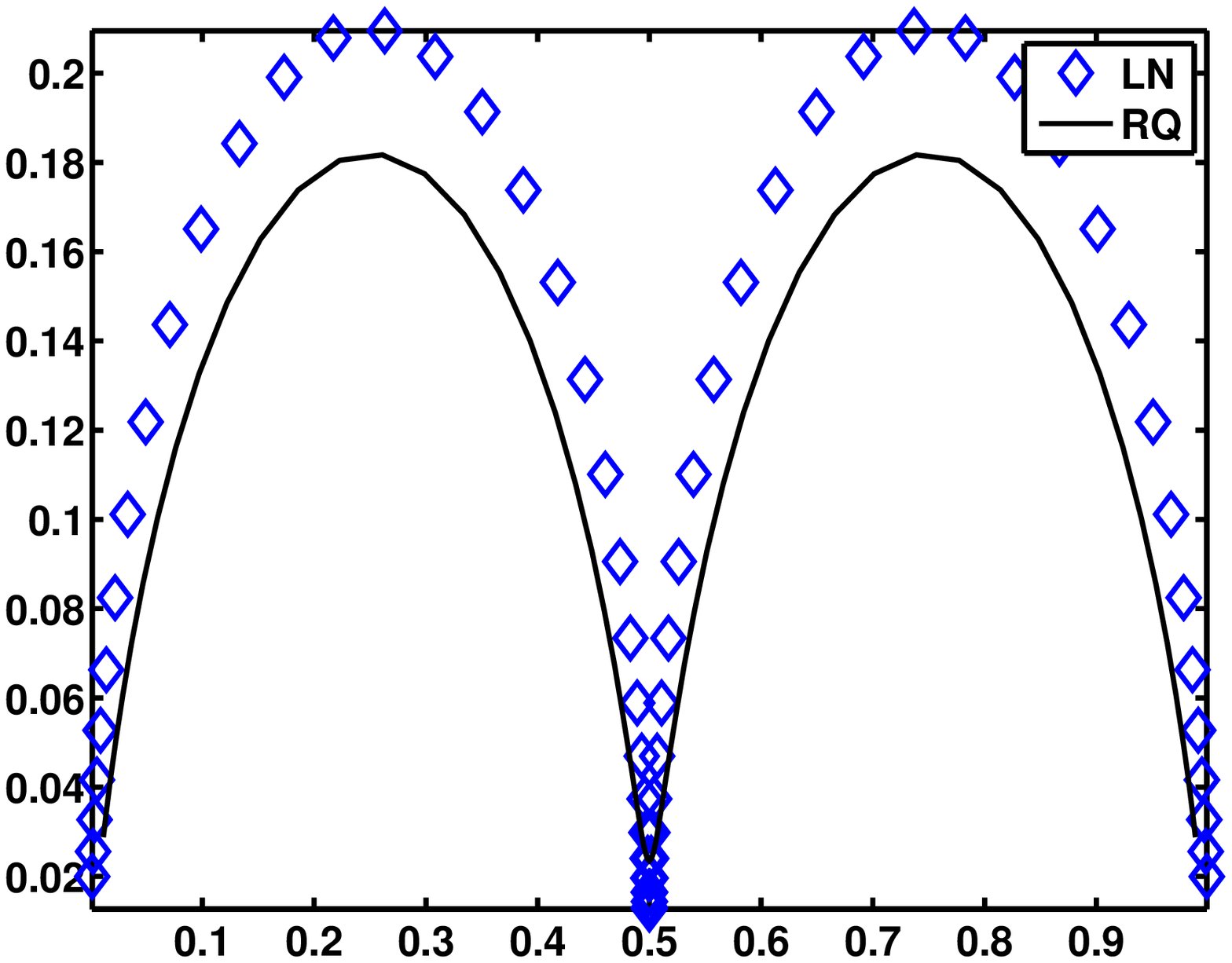}}\subfigure[$Z_1(\omega)$]{\includegraphics[width=1.3in]{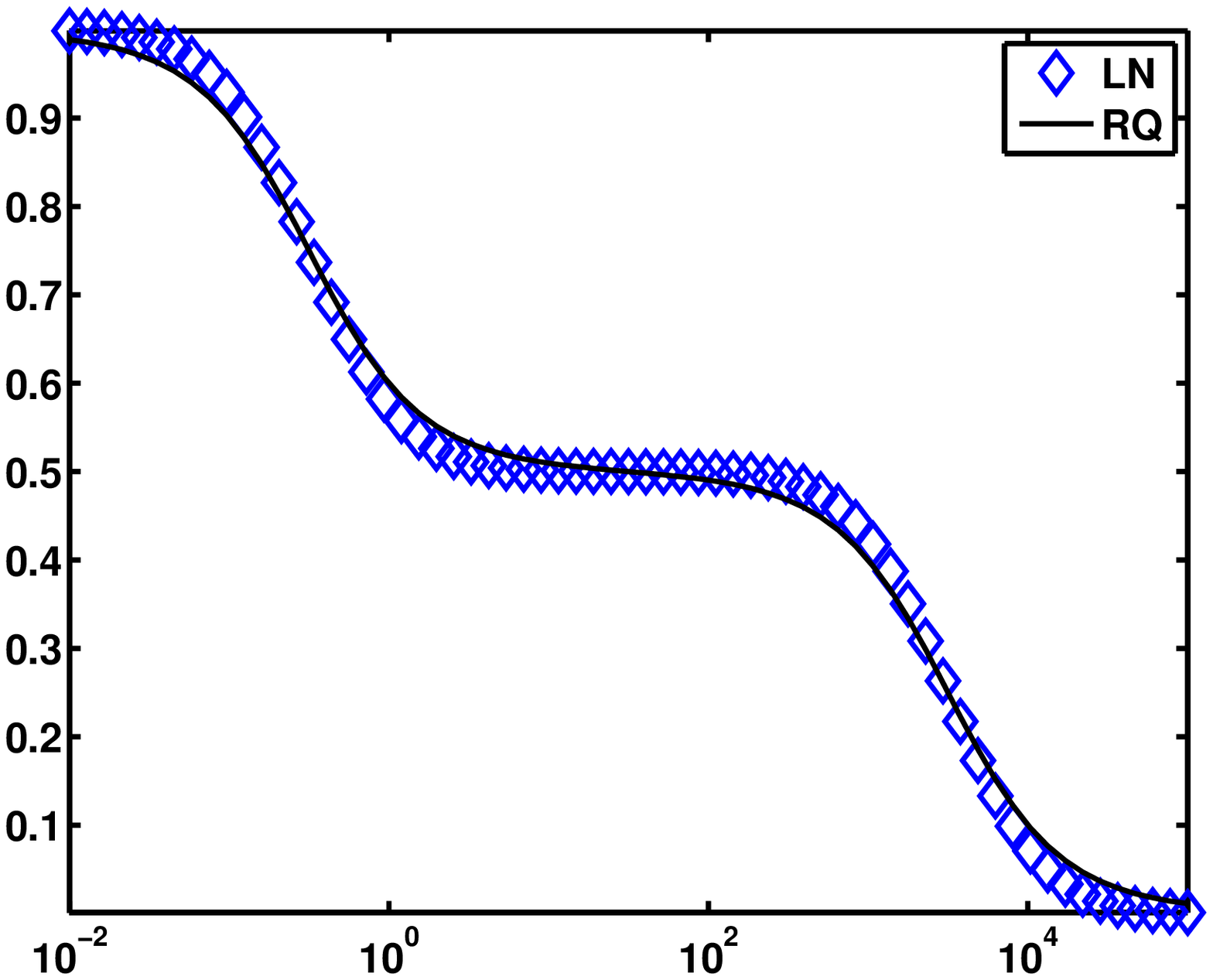}}\subfigure[$Z_2(\omega)$]{\includegraphics[width=1.3in]{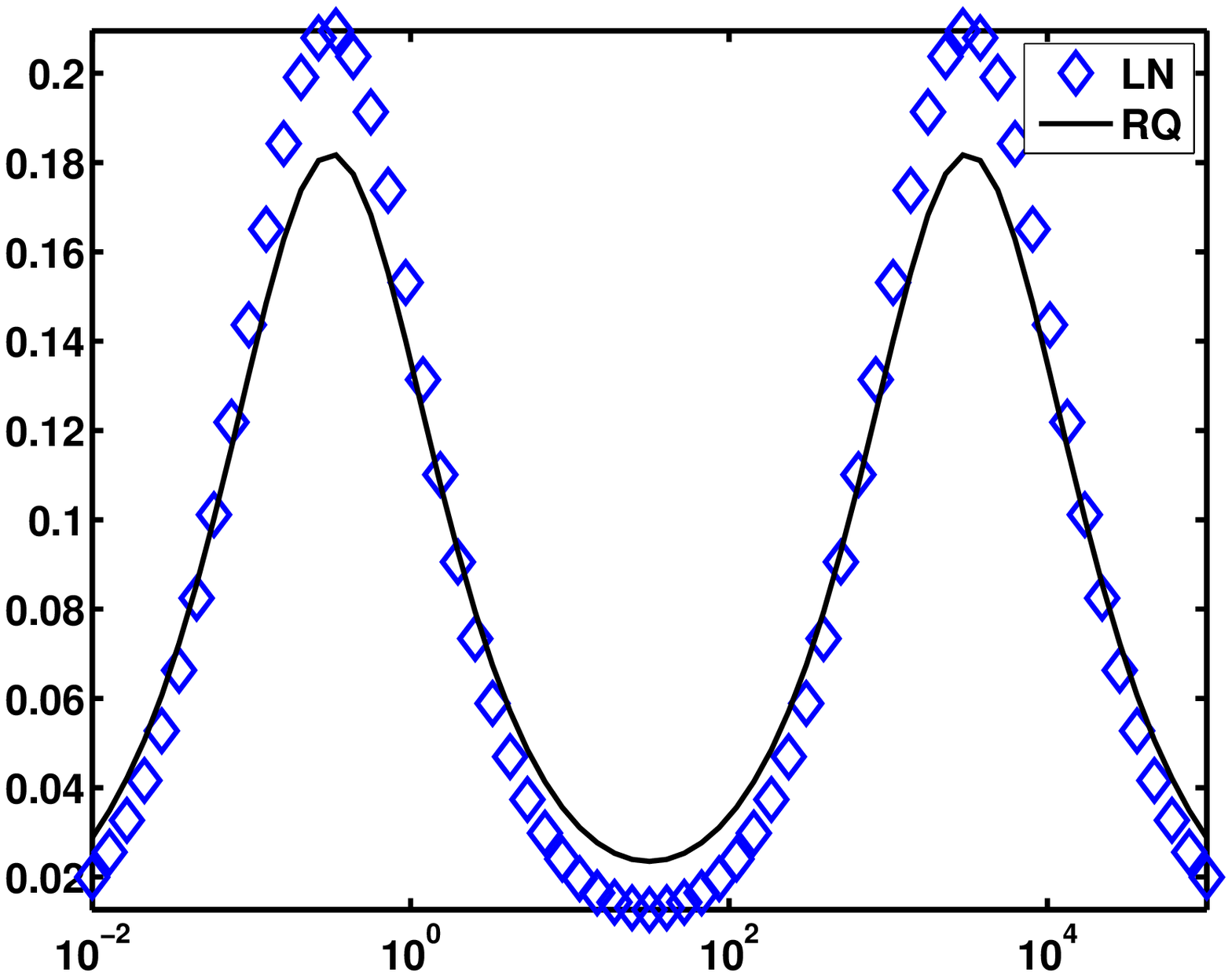}}

\subfigure[B: $g(t)$]{\includegraphics[width=1.3in]{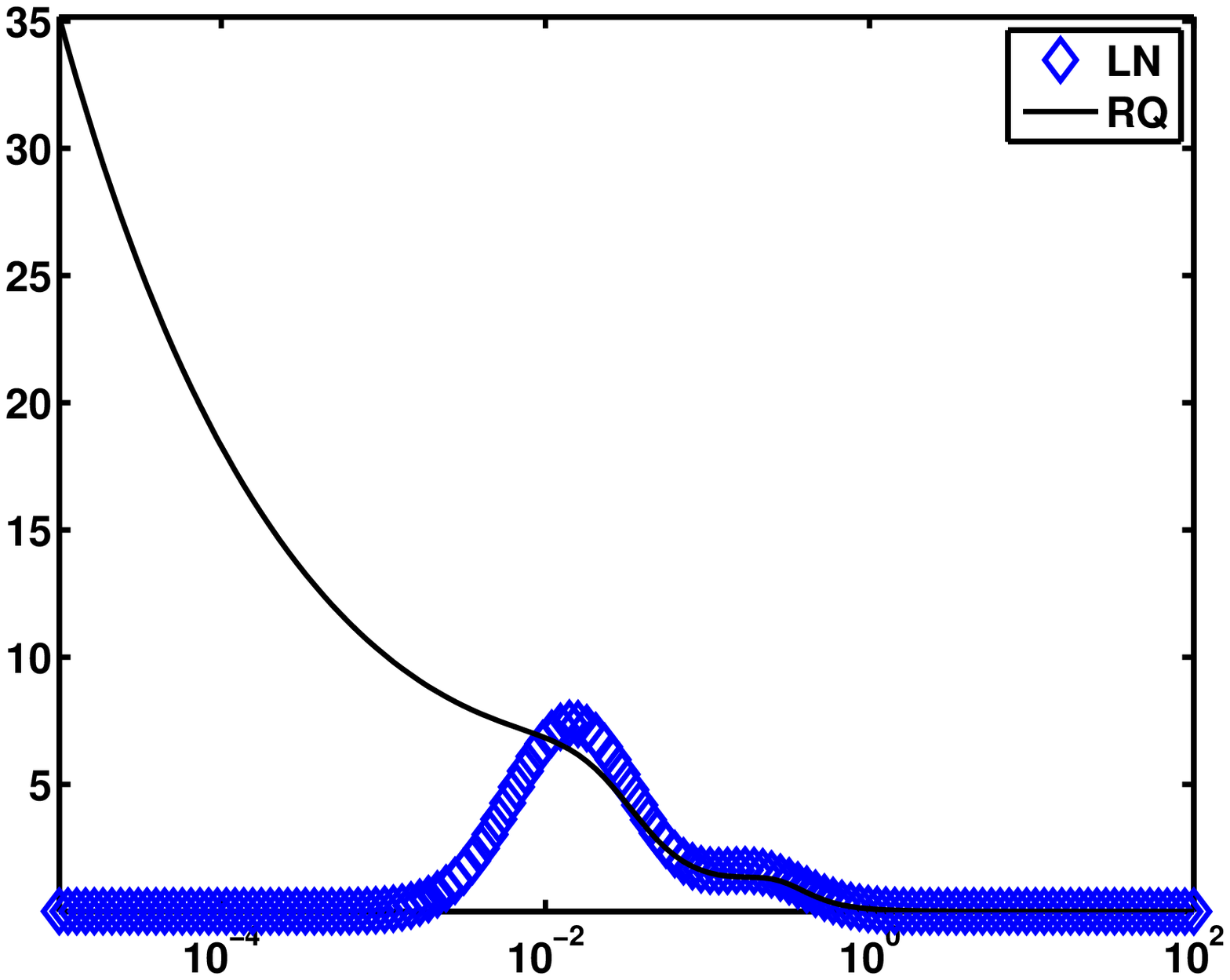}}\subfigure[$tg(t)$]{\includegraphics[width=1.3in]{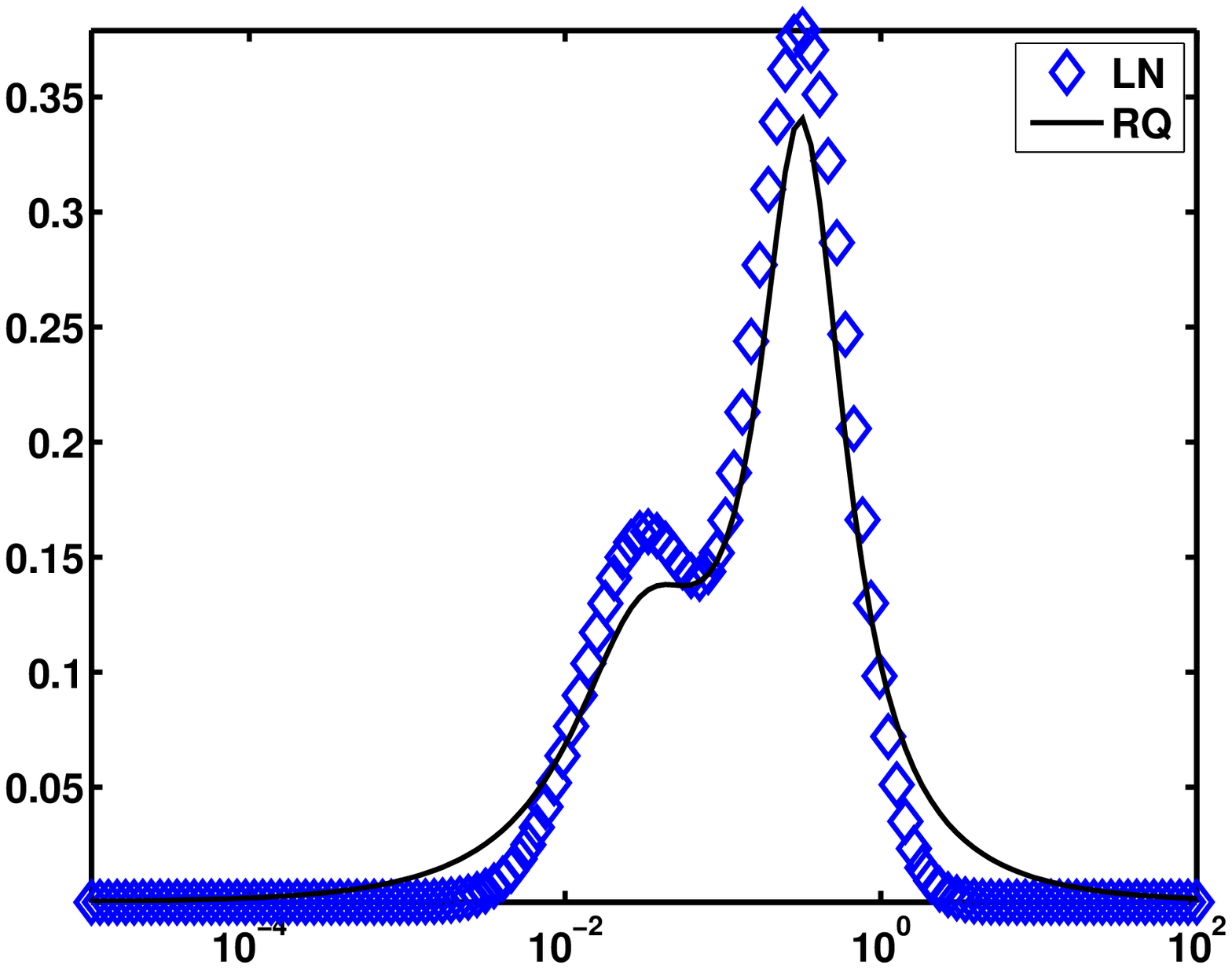}}\subfigure[Nyquist]{\includegraphics[width=1.3in]{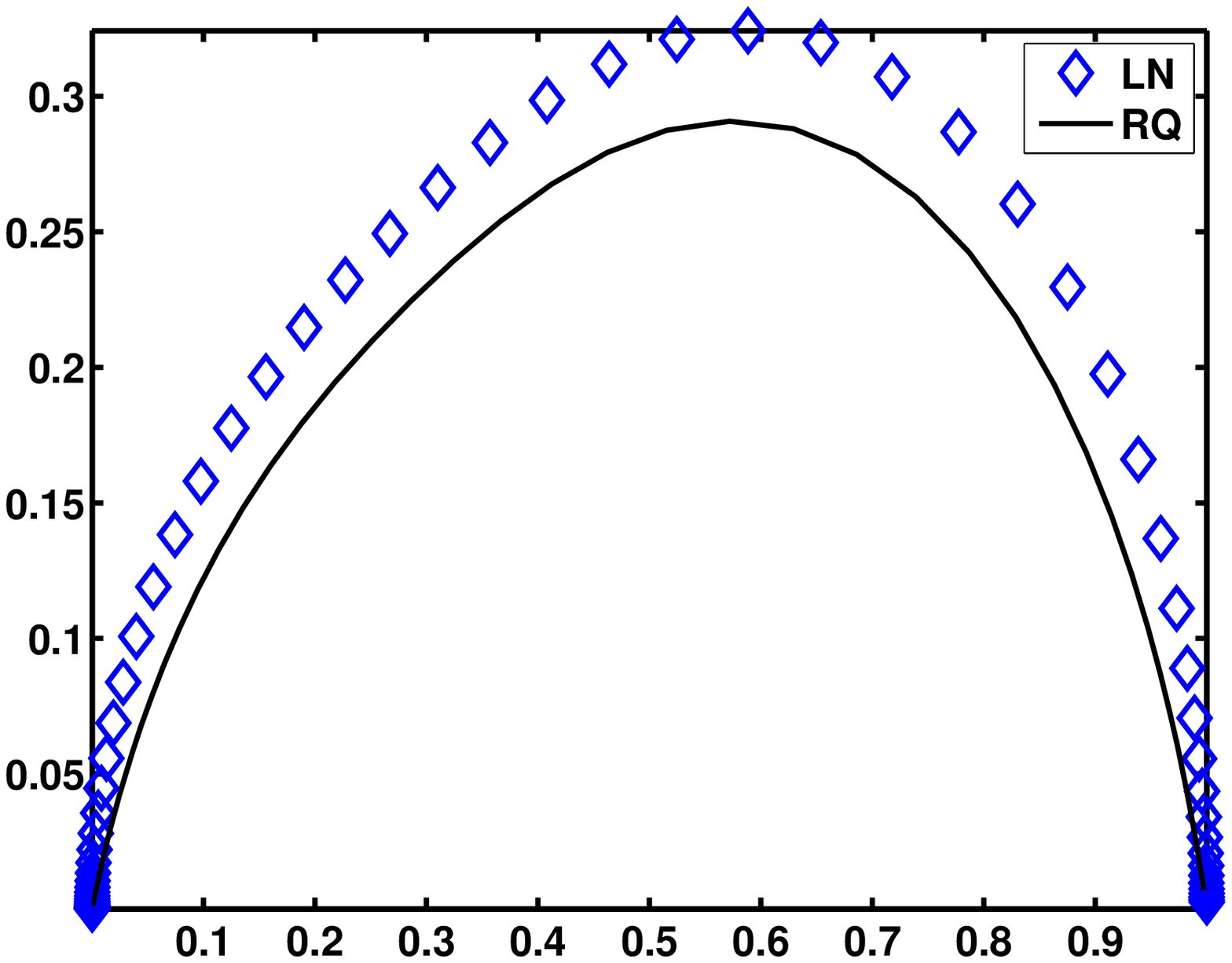}}\subfigure[$Z_1(\omega)$]{\includegraphics[width=1.3in]{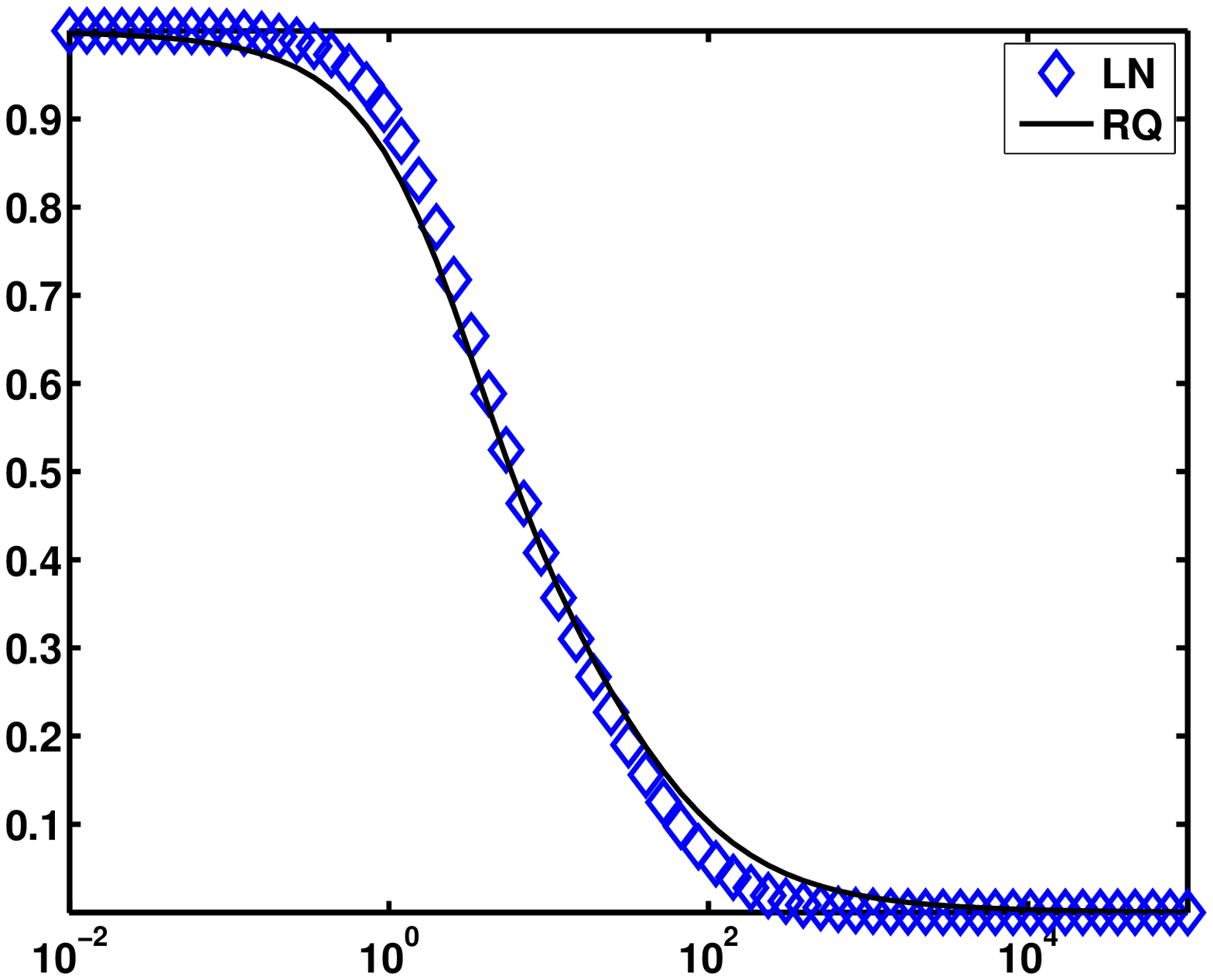}}\subfigure[$Z_2(\omega)$]{\includegraphics[width=1.3in]{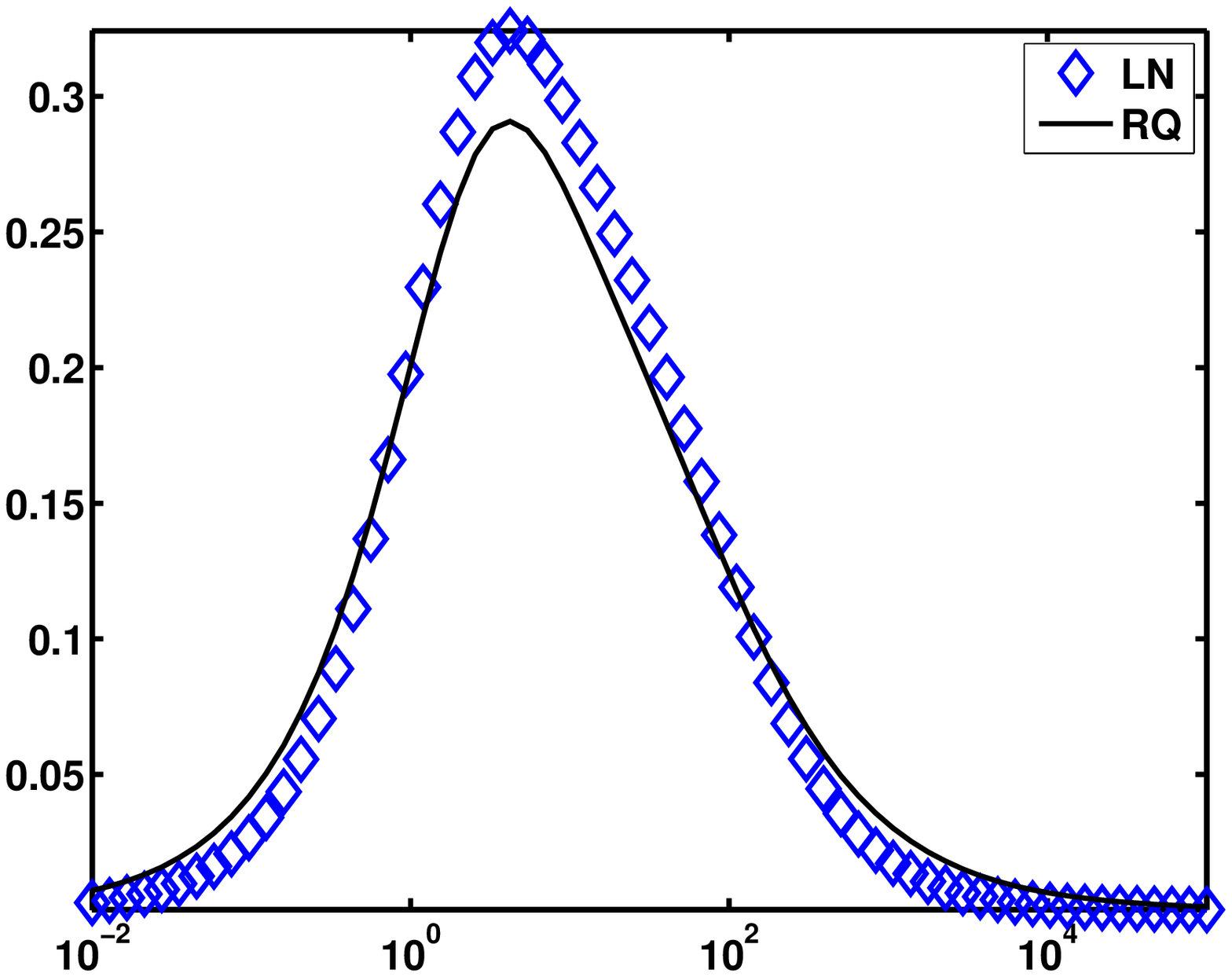}}

\subfigure[C: $g(t)$]{\includegraphics[width=1.3in]{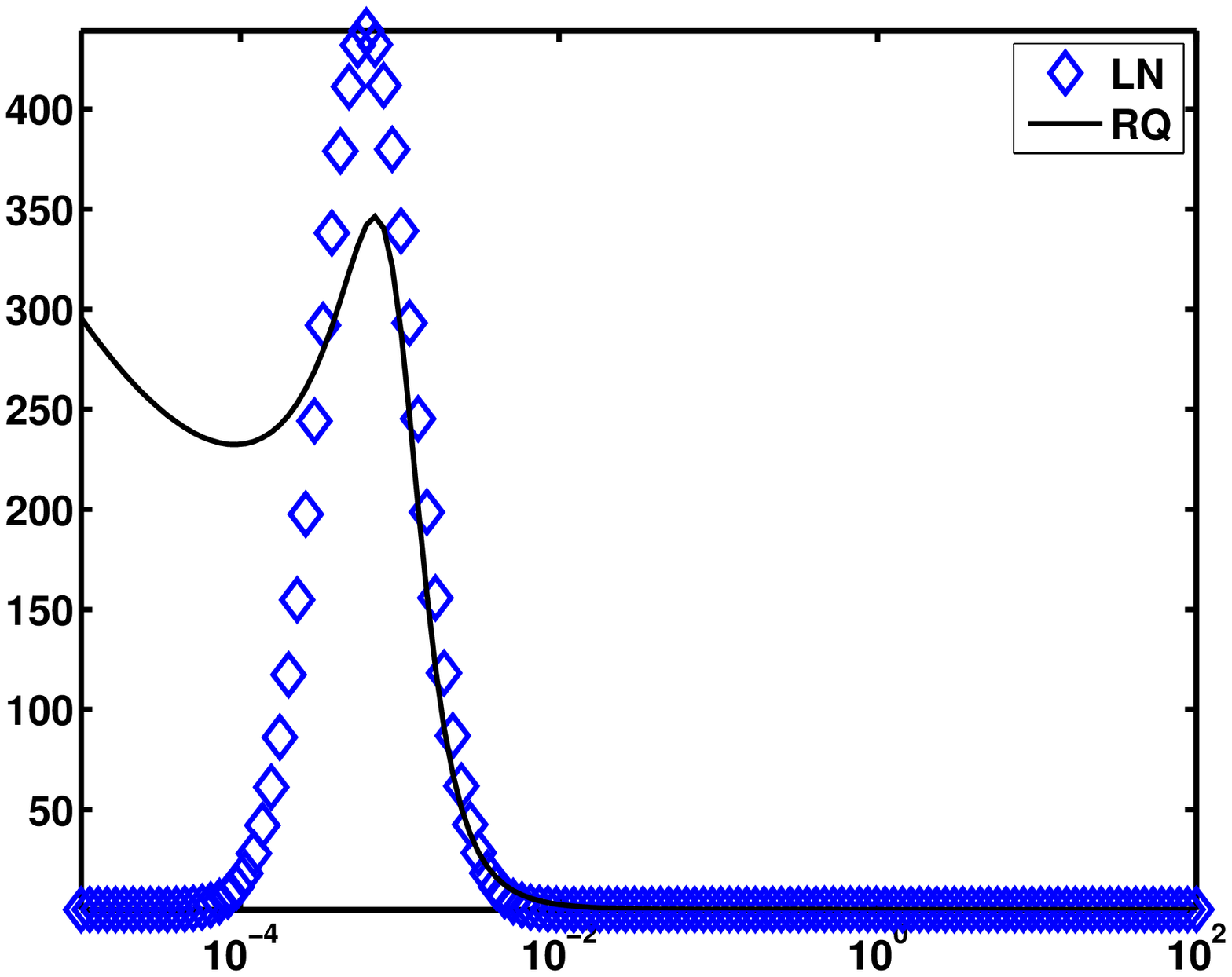}}\subfigure[$tg(t)$]{\includegraphics[width=1.3in]{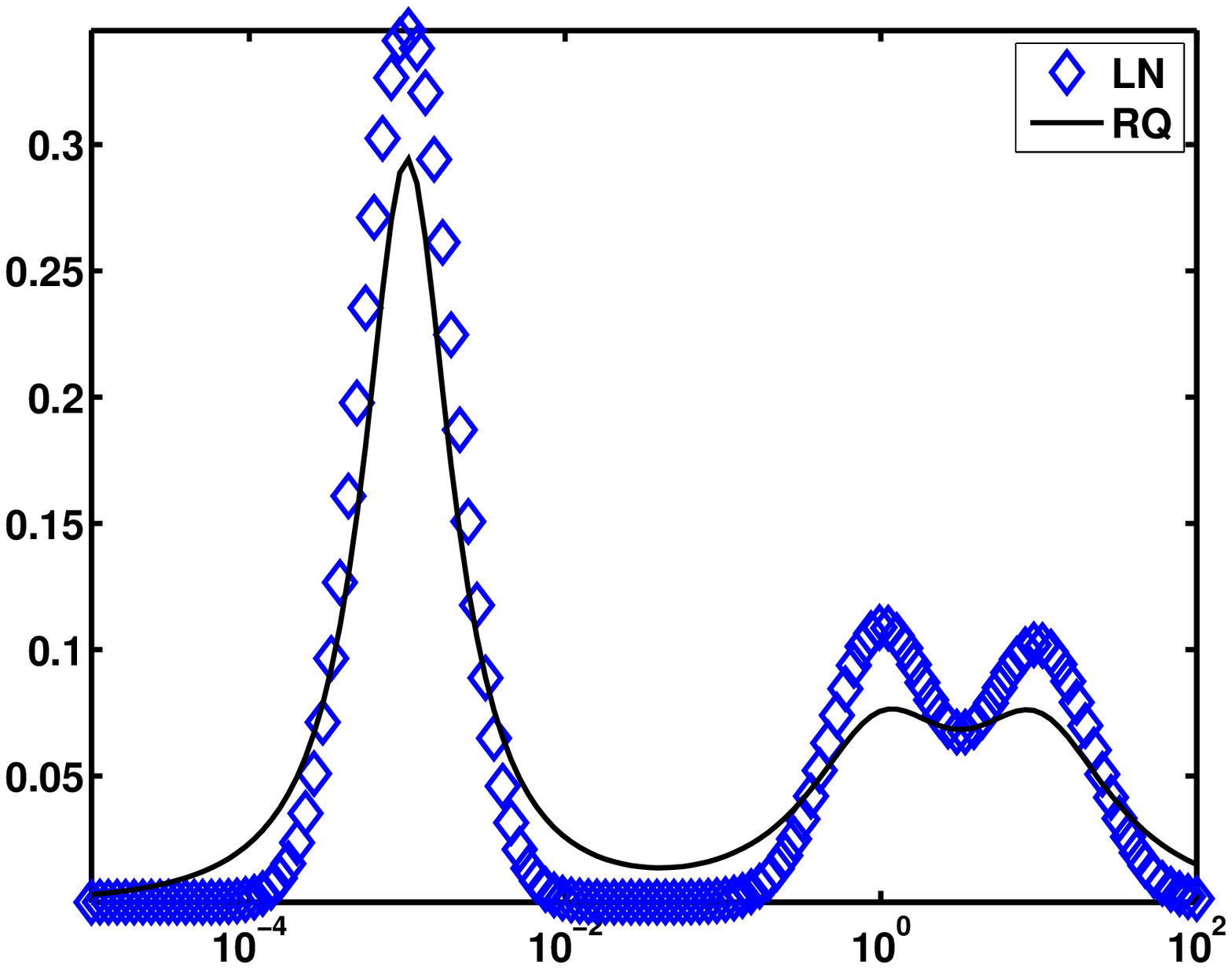}}\subfigure[Nyquist]{\includegraphics[width=1.3in]{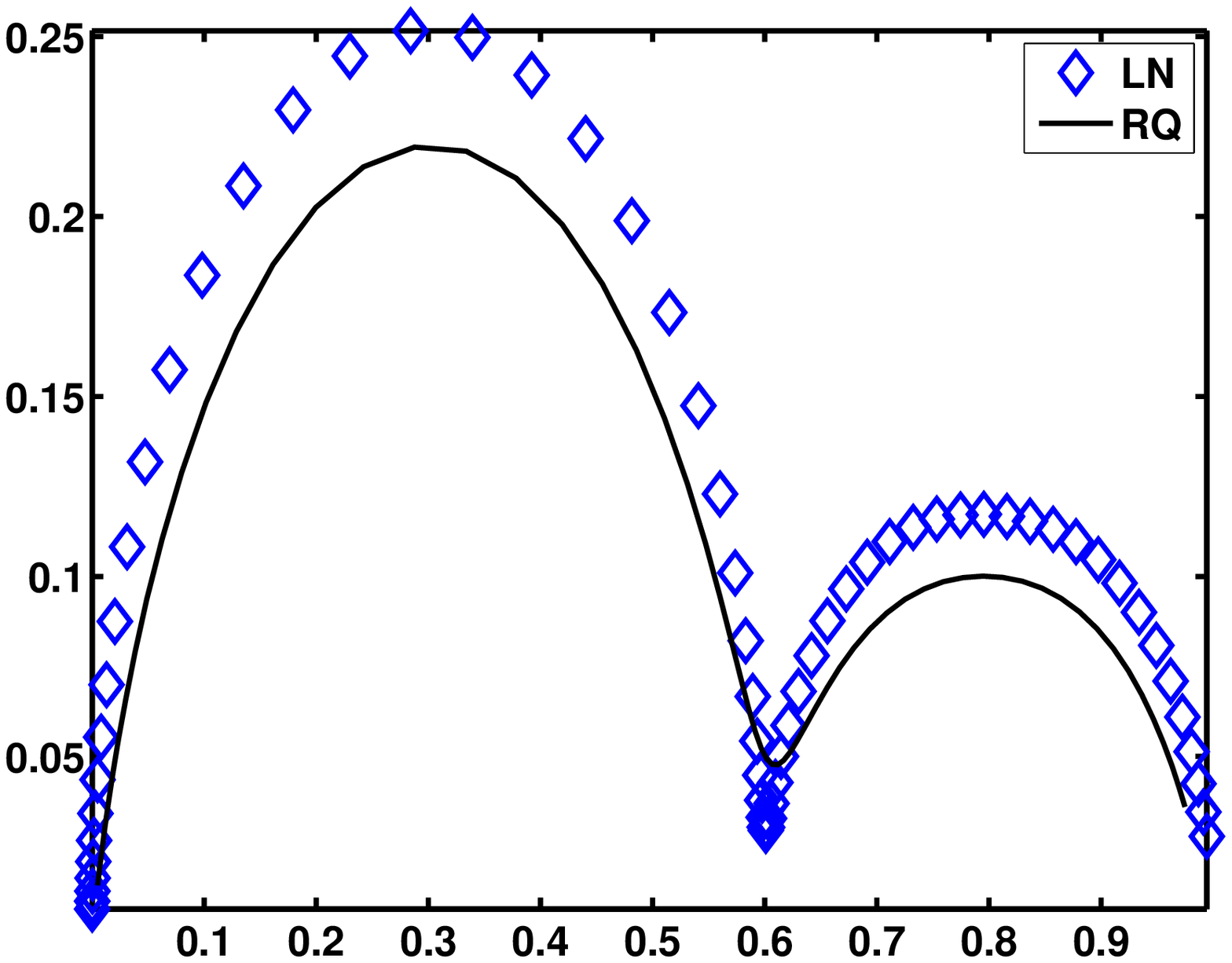}}\subfigure[$Z_1(\omega)$]{\includegraphics[width=1.3in]{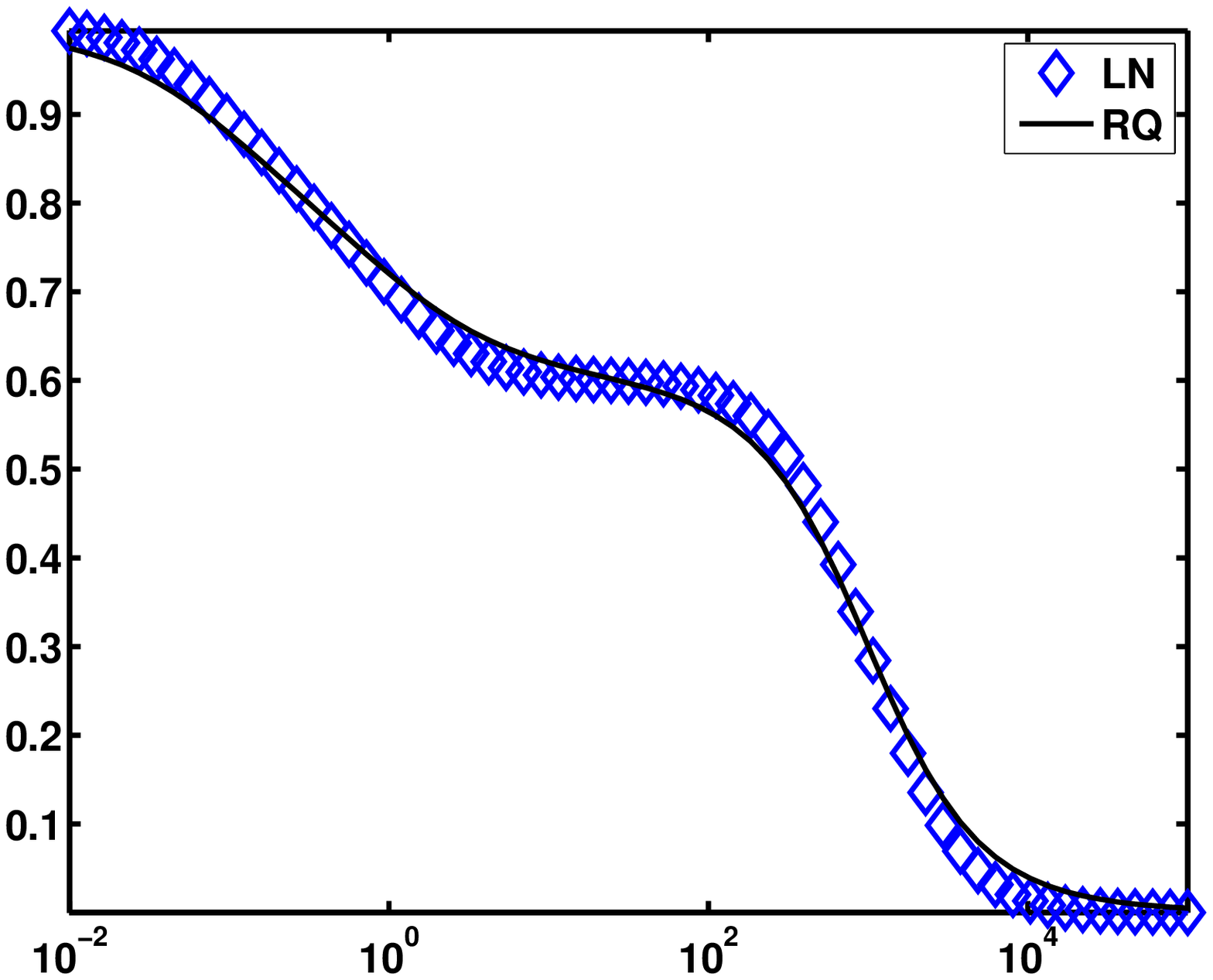}}\subfigure[$Z_2(\omega)$]{\includegraphics[width=1.3in]{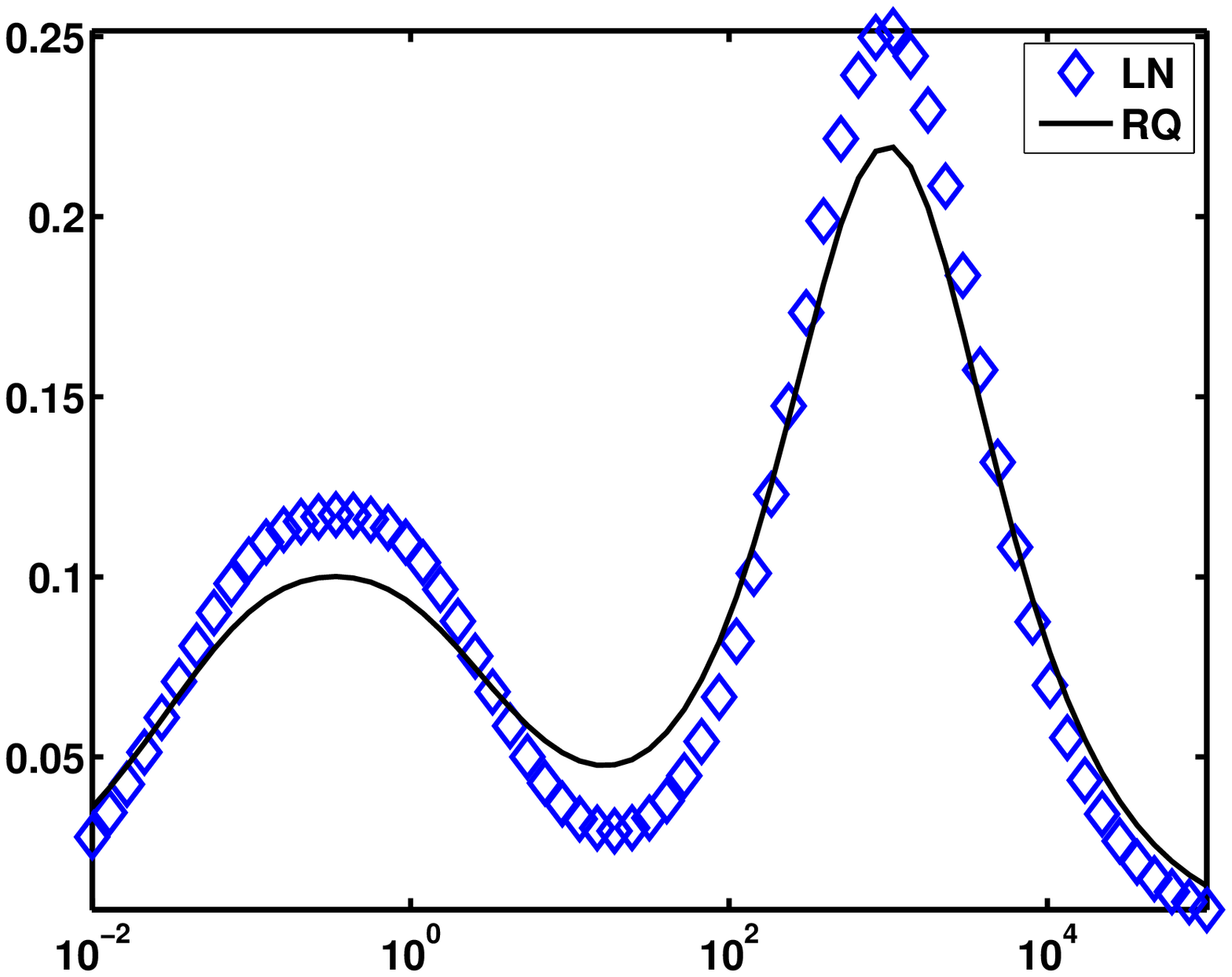}}
\caption{Visualization of the  simulation sets in Table~\ref{table:sim_params}, simulations $A$ to $C$ in  subfigures (a)-(e), (f)-(j) and (k)-(m), respectively. In each plot the RQ simulation is illustrated with the solid line and the LN simulation with the $\diamond$. Except for  the Nyquist plot, the scale of the $x-$ axis is logarithmic.  }
\label{fig:sim_figs}
\end{figure}

The parameters for the simulated data sets are detailed in  Table~\ref{table:sim_params}. Simulations A and B assume the existence of two underlying physical processes in the data; C assumes three processes. The RQ model uses the parameters $t_0$ and $\beta$, and the LN model uses the parameters $\mu=\ln(t_0)$ and $\sigma$. In each case the multiple processes are weighted by the weights $\alpha$. Note that the choice using $\mu=\ln(t_0)$ for the LN process provides data which are centered at $t_0$ in the $s-$space, and as can be seen in Figures~\ref{fig:sim_figs}(b), (g) and (l), are approximatelty aligned with the RQ model data. The set of figures for each simulation A to C demonstrates  the similarity of the chosen RQ and LN data, indicating the difficulty of distinguishing between these models from the impedance data alone,  Figures~\ref{fig:sim_figs}(c)-(e), (h)-(j) and (m)-(o).  In all plots, other than the Nyquist, the $y$-data are plotted against $x$ on a logarithmic scale.   

\begin{table}[h!]
\centering
\begin{tabular}{|c||c|c|c|c|} \hline
 Parameters  & $t_0$ & $\beta$ & $\sigma$ & $\alpha$ \\ \hline \hline
Simulation A & $[10^{-3.5},10^{0.5}]$ &$[0.8,0.8]$&$[\text{ln}(2),\text{ln}(2)]$&$[0.5,0.5]$ \\ \hline
Simulation B & $[10^{-1.5},10^{-0.5}]$ & $[0.7,0.8]$ & $[\text{ln}(2.4),\text{ln}(2)]$ & $[0.35,0.65]$ \\ \hline
Simulation C & $[10^{-3},1,10]$ & $[0.8,0.7,0.7]$ & $[\text{ln}(2),\text{ln}(2.1),\text{ln}(2.2)]$ & $[0.6,0.2,0.2]$ \\ \hline
\end{tabular}
\caption{Simulation parameters for simulated test data}
\label{table:sim_params}
\end{table}

The two processes in simulation Set A  have peaks that are far enough apart that the individuals processes are effectively separated. In such cases  distinct peaks are seen in the plot of $Z_2(\omega)$ at $\omega$ values which correspond to the reciprocals of the peak values in $t$.   
Such information can be used to verify the results of the fitting by examining the locations of the resulting peaks, \cite{Supp}. Moreover, this is 
  therefore a relatively well-behaved situation for which it should be possible to separate the underlying processes from measured data. In contrast, the two processes in set B are  close  but the existence of two peaks is not clear from $Z_2$. 
For the three processes in set C the first  is distinct in time from the latter two, which are close and overlapping for larger time. Again the impedance data does not clearly indicate the number of processes in the data.  

For testing the algorithms noise was added to the simulated values for $Z$ following \eqref{noisy data} as discussed in Section~\ref{NLSfits} for the examination of the NLS fitting. In all cases the data are sampled at $65$ points logarithmically spaced on $[\omega_{\mathrm{min}}, \omega_{\mathrm{max}}] = [10^{-2}, 10^5]$,  consistent with practical data, and with sampling of the DRT at   $t=1/\omega$ yielding the equal spacing in $s-$space at $130$ points for matrix $A_4$.  For the results presented here, the noise level $\eta=10^{-3}$ corresponding to $.1\%$ noise was chosen. Further results in \cite{Supp} use noise levels $1\%$ and $5\%$. While the higher noise levels, particularly around $1\%$ may be more consistent with practical data, the actual noise level for the measured data is unknown. On the other hand, even with the lower noise level the benefit of using the NNLS will become clear in Sections~\ref{parameterchoice}-\ref{lsnnls}.

\subsection{Parameter Choice Methods}\label{parameterchoice}
The  NCP or residual periodogram (RP) was presented in \cite{hansen:2006-1, rust:1998-1, rust:2008-1} as a parameter choice method for 
the Tikhonov LLS problem \eqref{regsoln}. The use of the NCP applied to the residual  for finding the optimal parameter for the NNLS problem  in \eqref{nnls} is, however,  to the best of our knowledge novel. We therefore first contrast the use of the LC and NCP for parameter choice in the context of the NNLS constrained problem, with the three standard  choices of zero ($L=I$), first ($L=L_1$) and second order ($L=L_2$)  derivative operators. These can  be obtained using the function \texttt{get\_l} in the Regularization toolbox \cite{hansen:2007-1}. The LC is also implemented in the Regularization toolbox, while for the NCP we use the  modification of the \texttt{ncp} in the toolbox, as used also for examining the frequency content of the basis vectors as shown in Section~\ref{regsec}. For both the NCP and L-Curve methods, solutions were found for $50$ choices of $\lambda$ logarithmically spaced between $10^{-3.5}$ and $10^{1.5}$. The optimal $\lambda$ for the LC was chosen using the corner of the LC, while for the 
NCP method, the Kolmogorov-Smirnov confidence level for white noise uses $p = 0.2$. Each simulation was tested over $50$ realizations of $.1\%$  white noise.  

For each noise realization the following information was recorded: the optimal solution obtained by the NCP and LC parameter choice methods, with the optimally found $\lambda_{\mathrm{NCP}}$ and $\lambda_{\mathrm{LC}}$, and the optimal solution over all $50$ choices for $\lambda$, with the respective $\lambda_{\mathrm{opt}}$, as measured with respect to the absolute error in the $s-$space.  The geometric means  of  $\lambda_{\mathrm{NCP}}$ and $\lambda_{\mathrm{LC}}$  were  calculated over all $50$ noise realizations. The absolute error for each choice of $\lambda$ was also recorded for each noise realization, and  the mean of these absolute errors taken to give an average error for a given $\lambda$ which can be visualized against $\lambda$. This follows the analysis presented in \cite{Maetal:04} for the examination of the optimal regularization parameter.  In the plots we thus show the average error against $\lambda$ indicated by the $\circ$ plot. On the same plot we indicate by the vertical lines the minimum $\lambda_{\mathrm{opt}}$, and the geometric means for  $\lambda_{\mathrm{NCP}}$ and $\lambda_{\mathrm{LC}}$, as the solid (red), dashed (green) and dot-dashed $\circ$ (blue) vertical lines, respectively.  For each simulation set the same procedure was performed for all smoothing norms $L$. To demonstrate the dependence of the obtained solution on the optimal parameter, an example  representative noise realization was chosen in each case and the solutions found using the chosen optimal parameters were compared with the exact solution. These are indicated by the solid line (black), $\diamond$ (red) , $\times$ (green)  and $\circ$ (blue), for the exact, $\lambda_{\mathrm{opt}}$, $\lambda_{\mathrm{NCP}}$, and $\lambda_{\mathrm{LC}}$ solutions, respectively. 

\begin{figure}[!ht]
 \centering
\subfigure[$L=I$]{\includegraphics[width=1.7in]{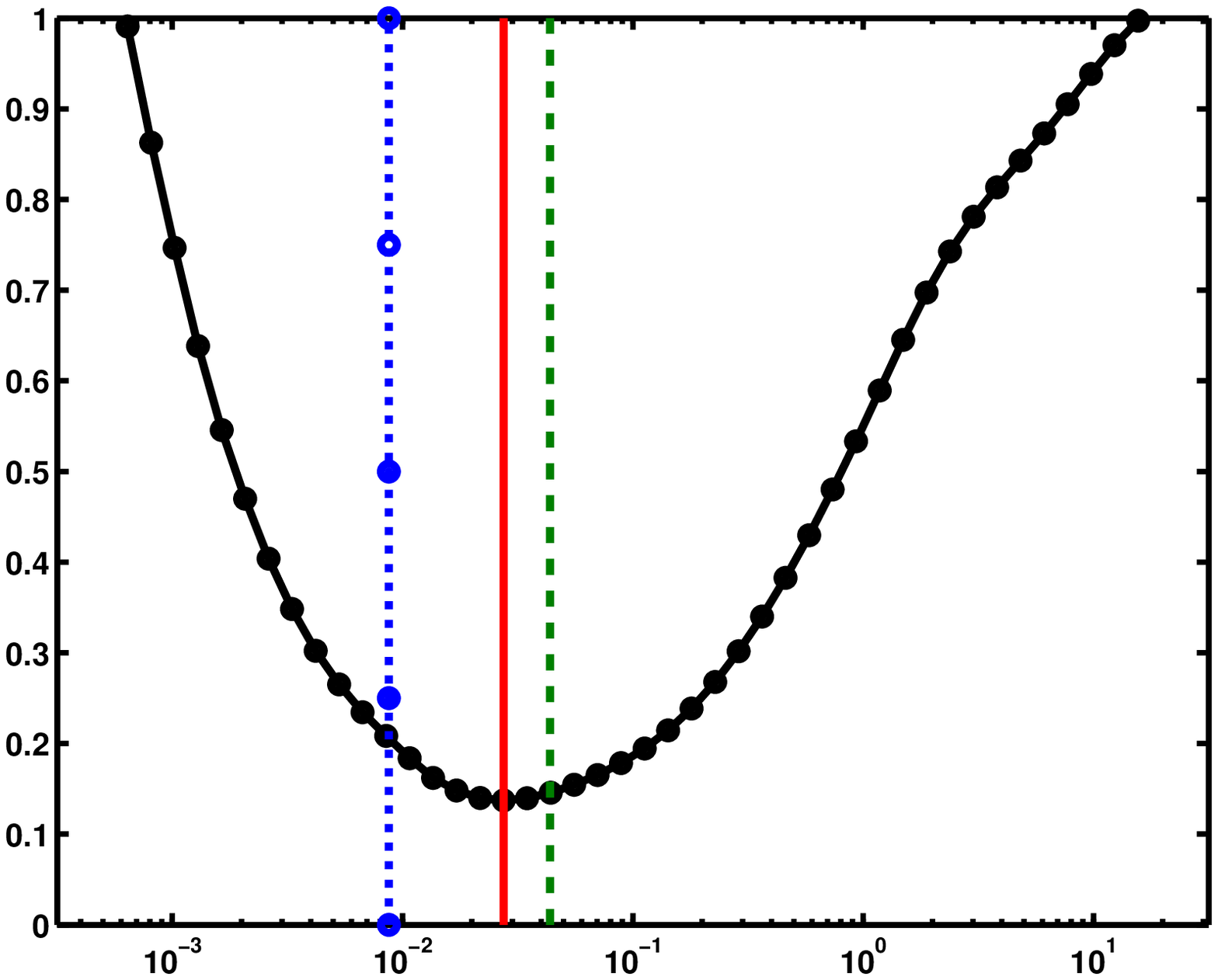}}
\subfigure[$L=L_1$]{\includegraphics[width=1.7in]{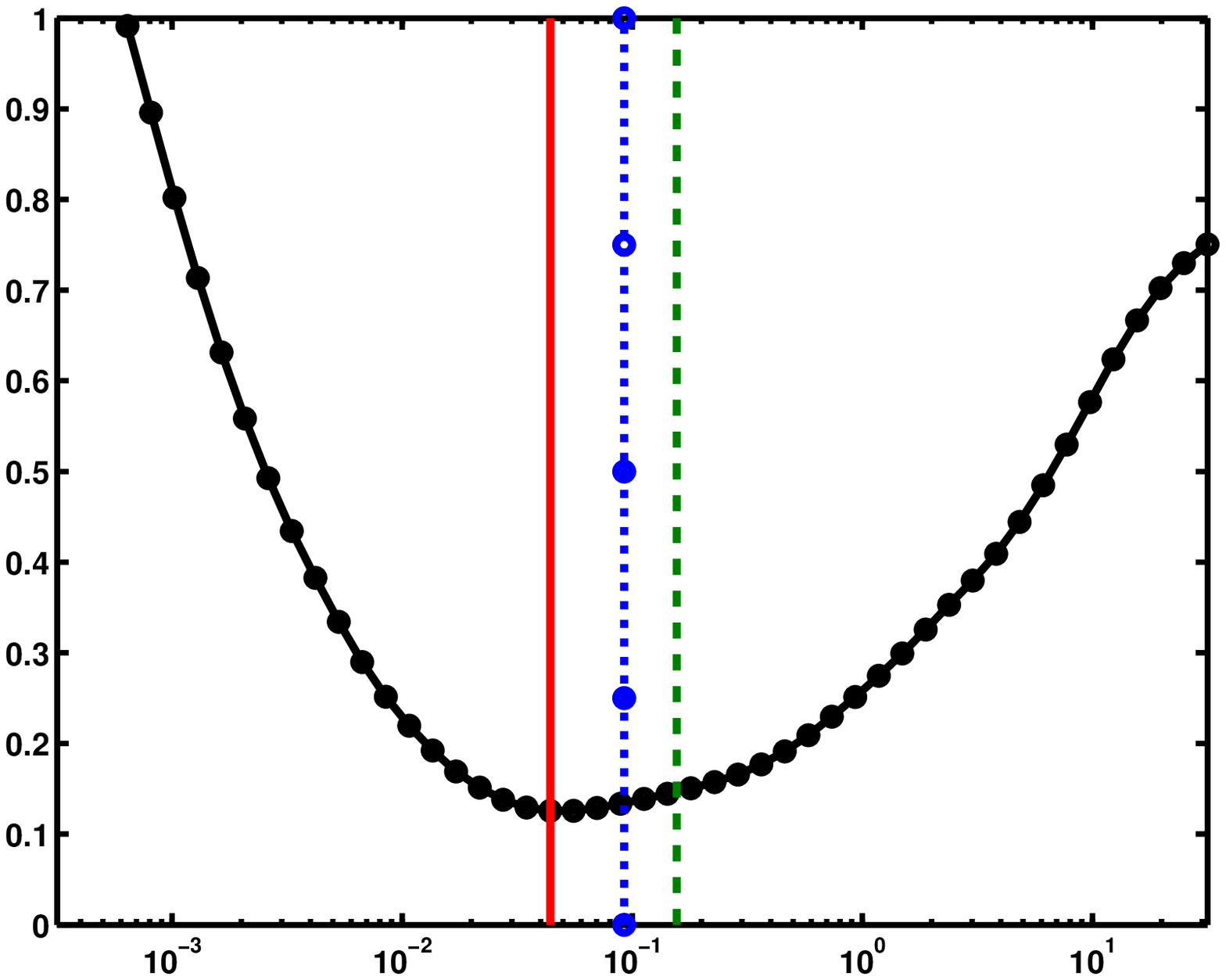}}
\subfigure[$L=L_2$]{\includegraphics[width=1.7in]{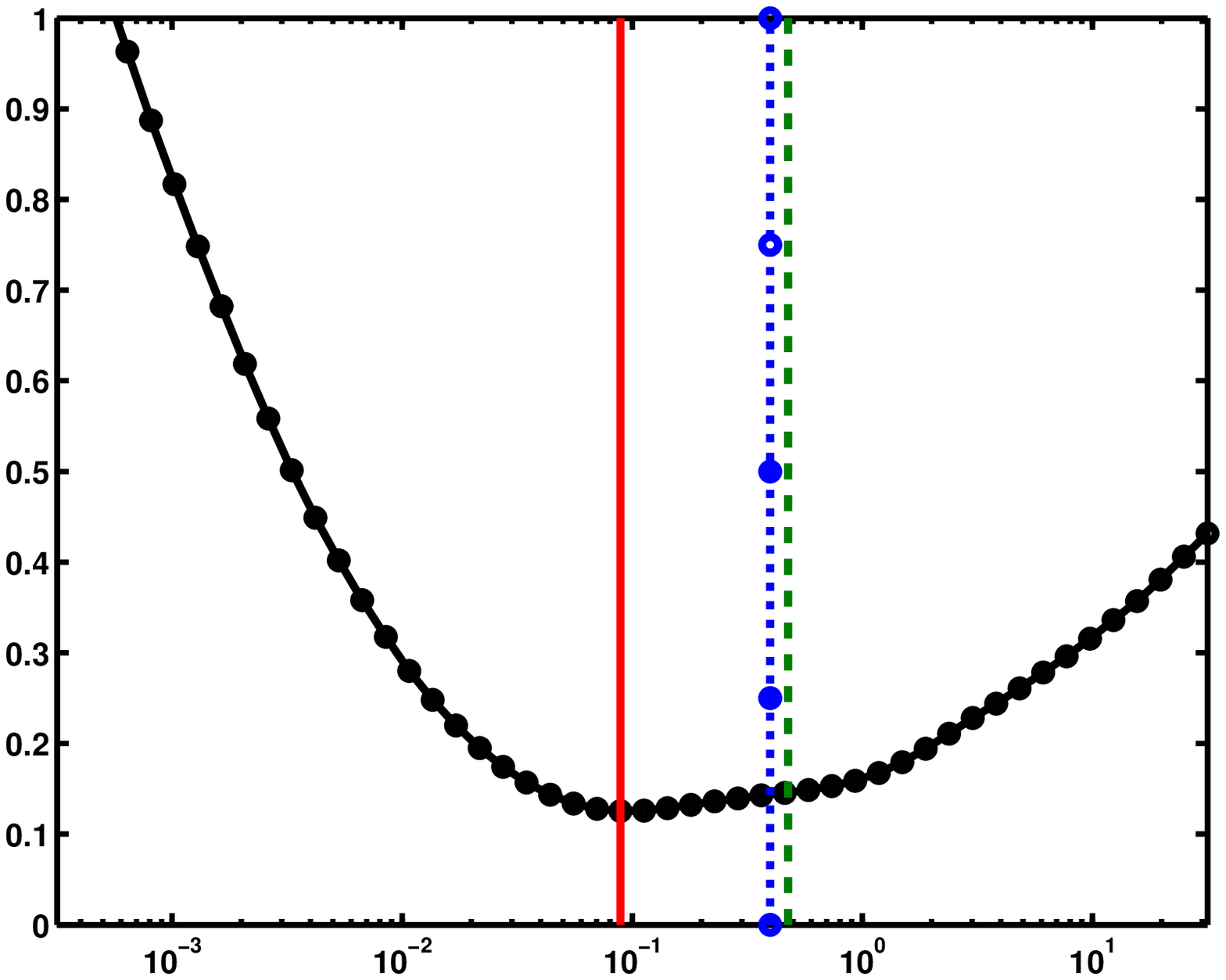}}
\subfigure[$L=I$]{\includegraphics[width=1.7in]{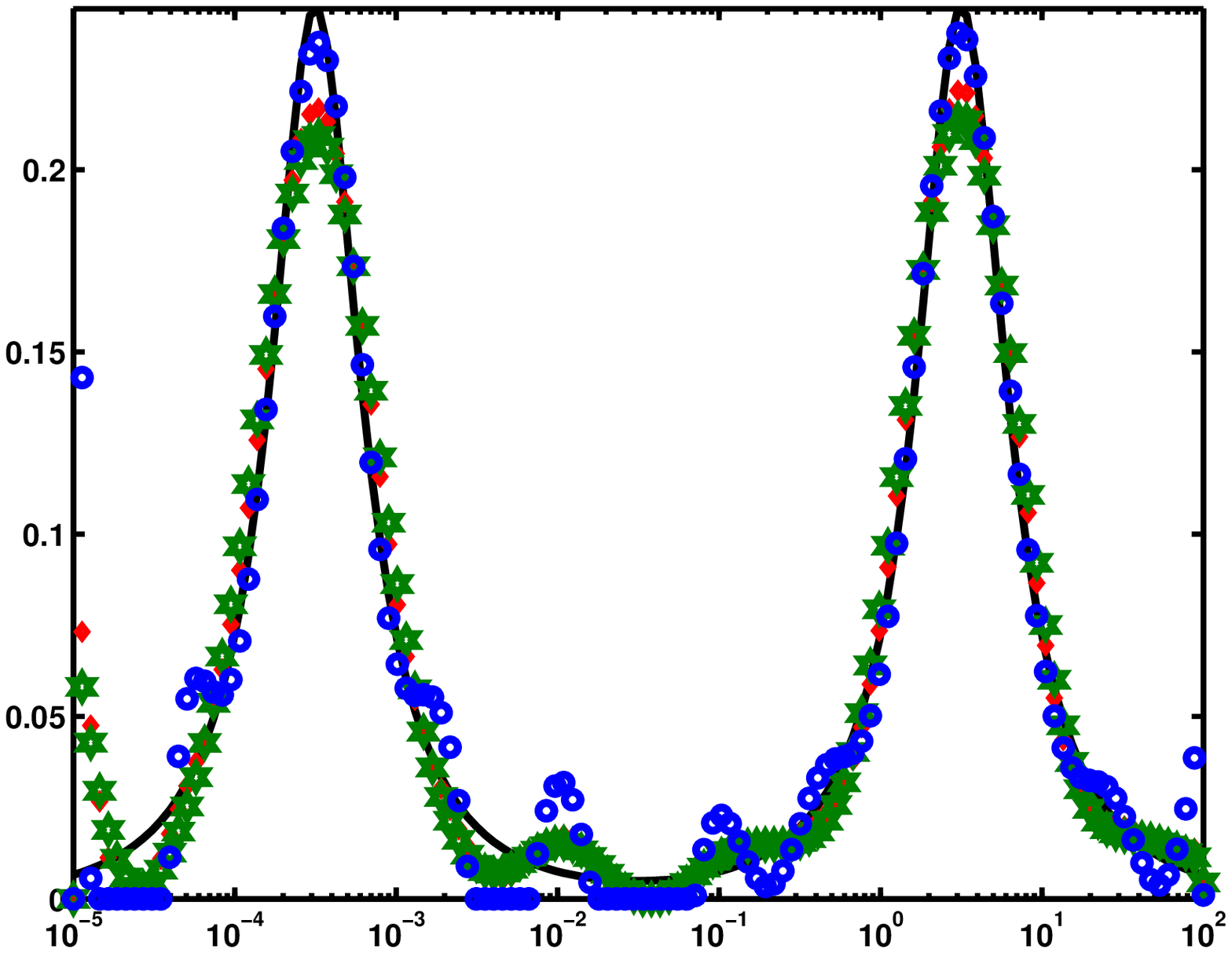}}
\subfigure[$L=L_1$]{\includegraphics[width=1.7in]{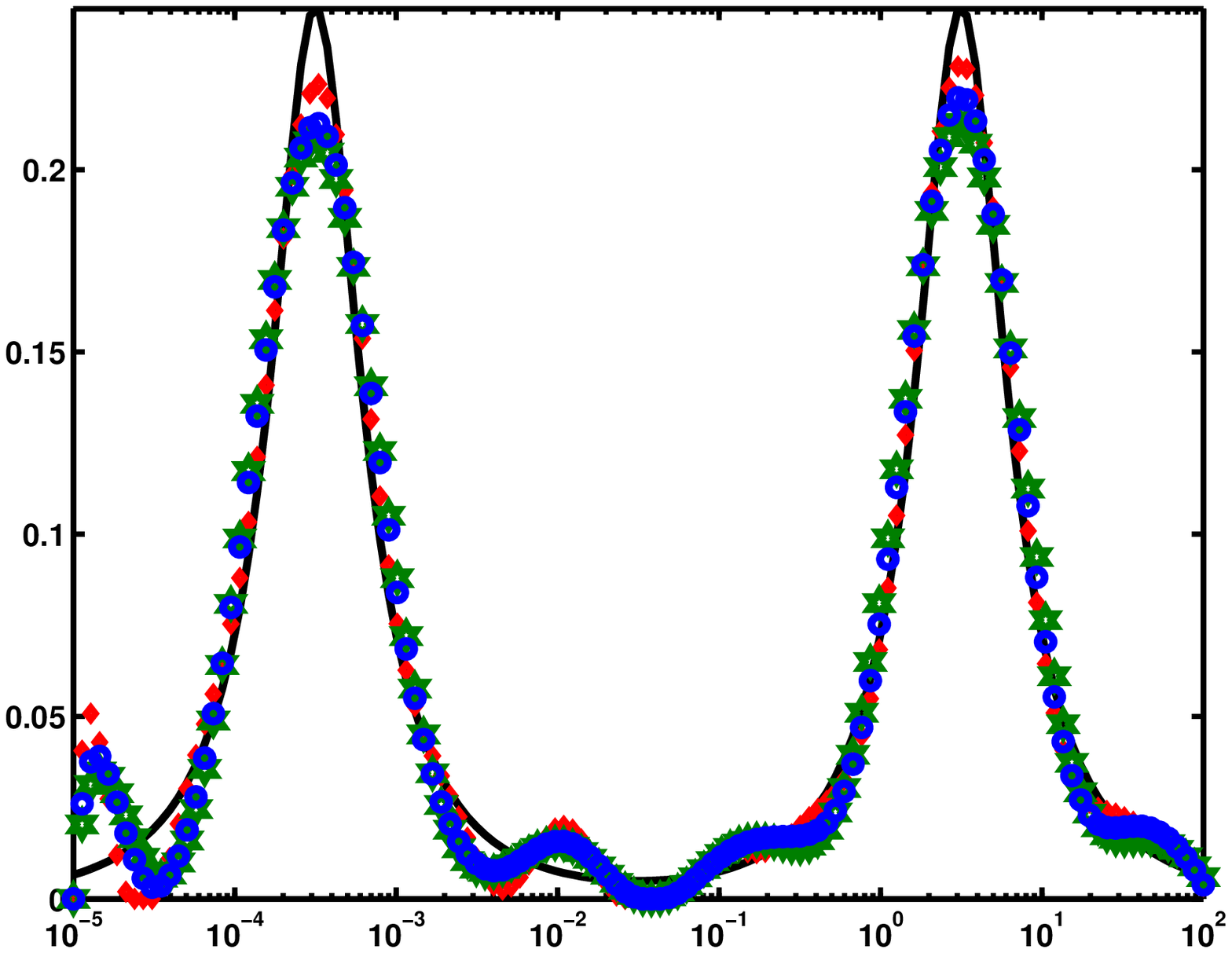}}
\subfigure[$L=L_2$]{\includegraphics[width=1.7in]{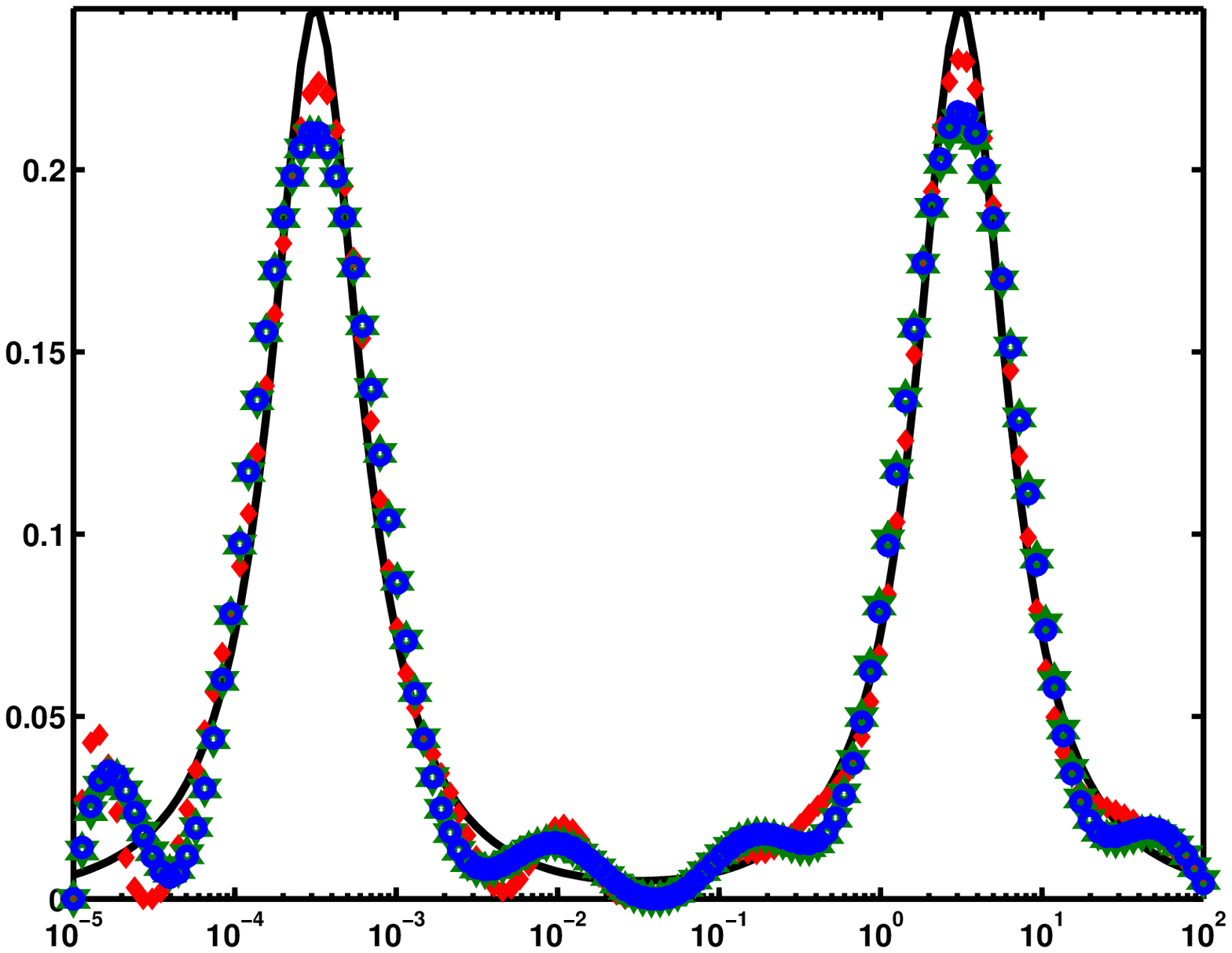}}
\caption{Mean error and example  \texttt{lsqnonneg}  NNLS solutions.  $.1\%$ noise, RQ-A data set, matrix $A_4$.}
\label{fig-lambdachoiceRQ1A4LN}
\end{figure}

\begin{figure}[!ht]
  \centering
\subfigure[$L=I$]{\includegraphics[width=1.7in]{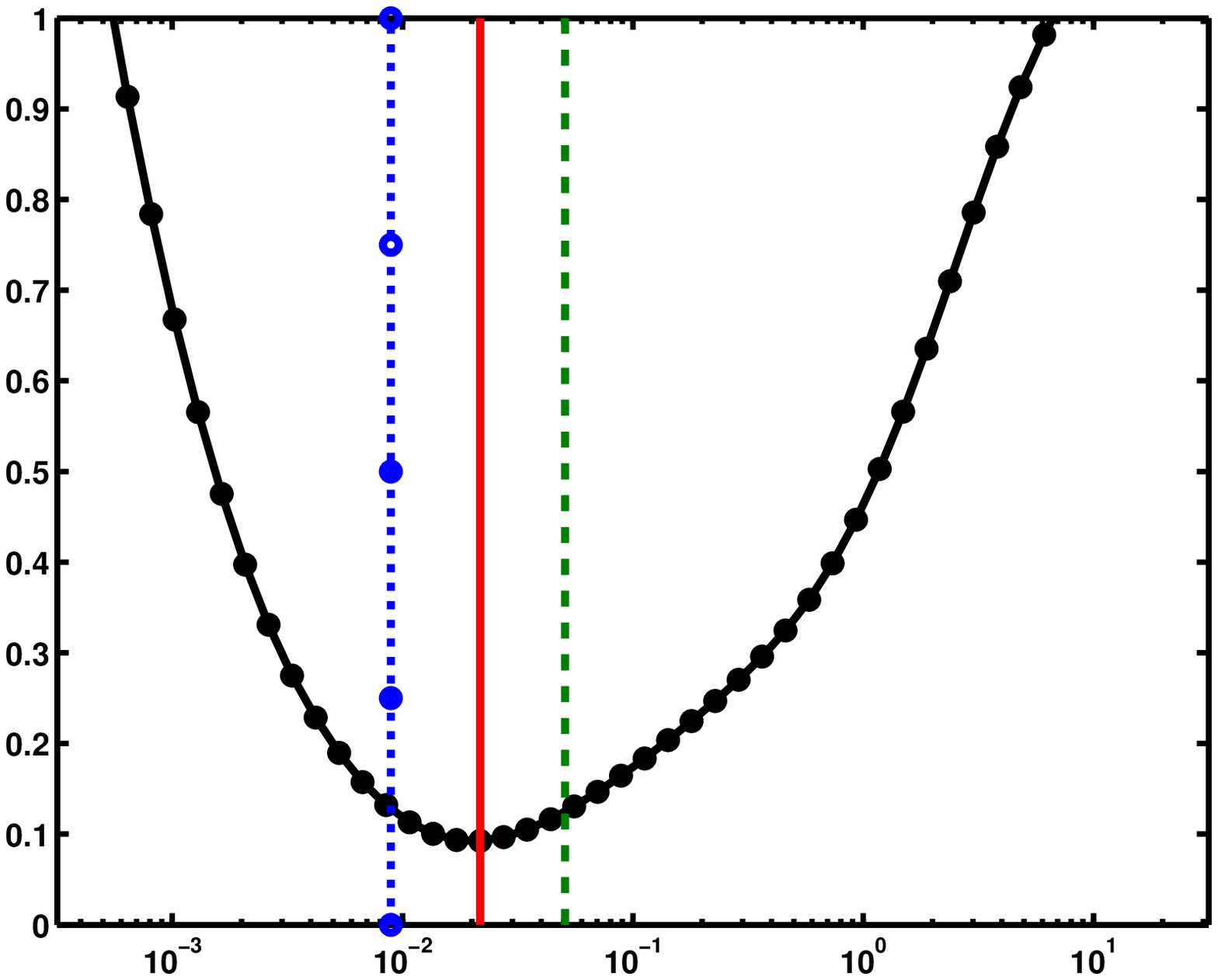}}
\subfigure[$L=L_1$]{\includegraphics[width=1.7in]{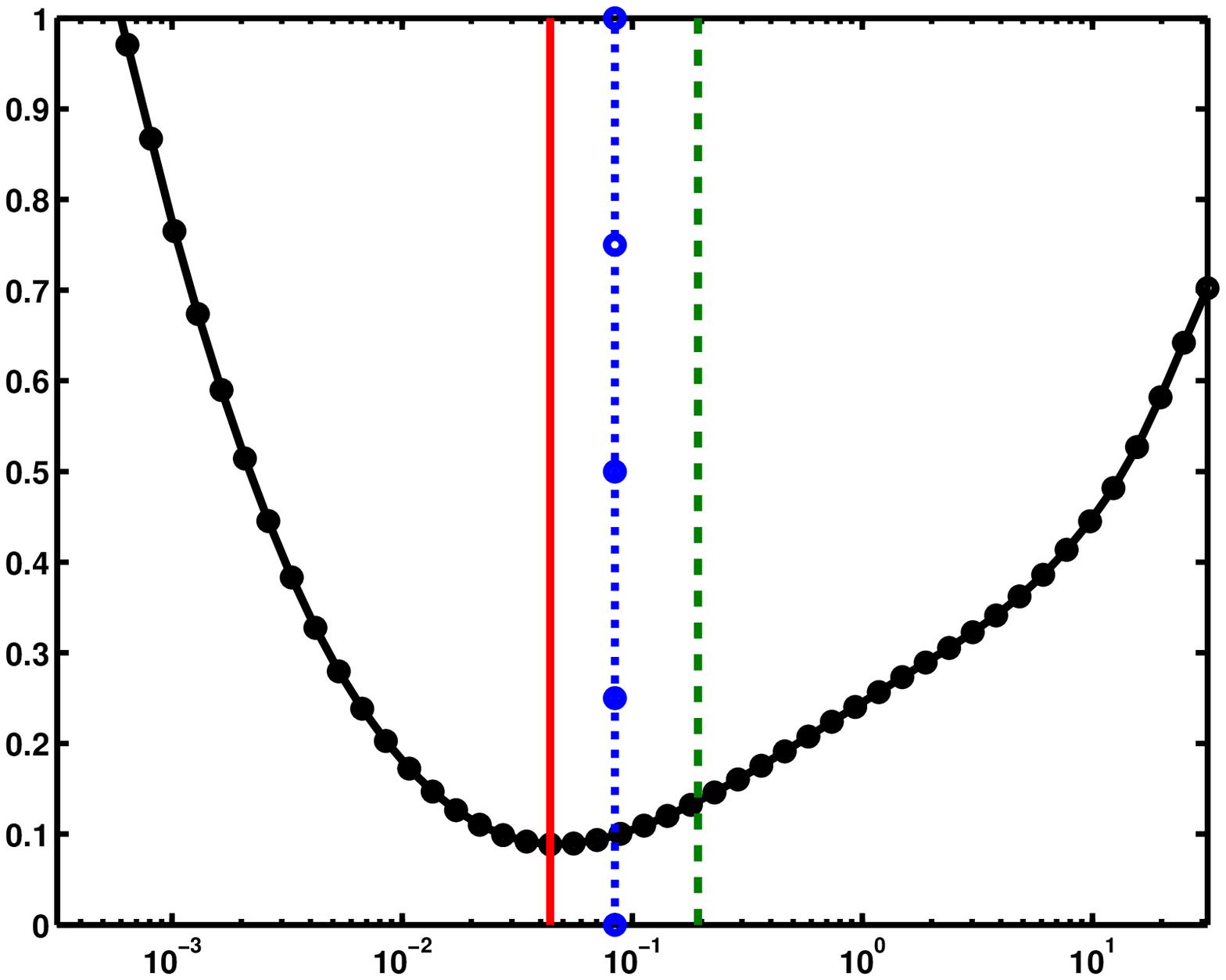}}
\subfigure[$L=L_2$]{\includegraphics[width=1.7in]{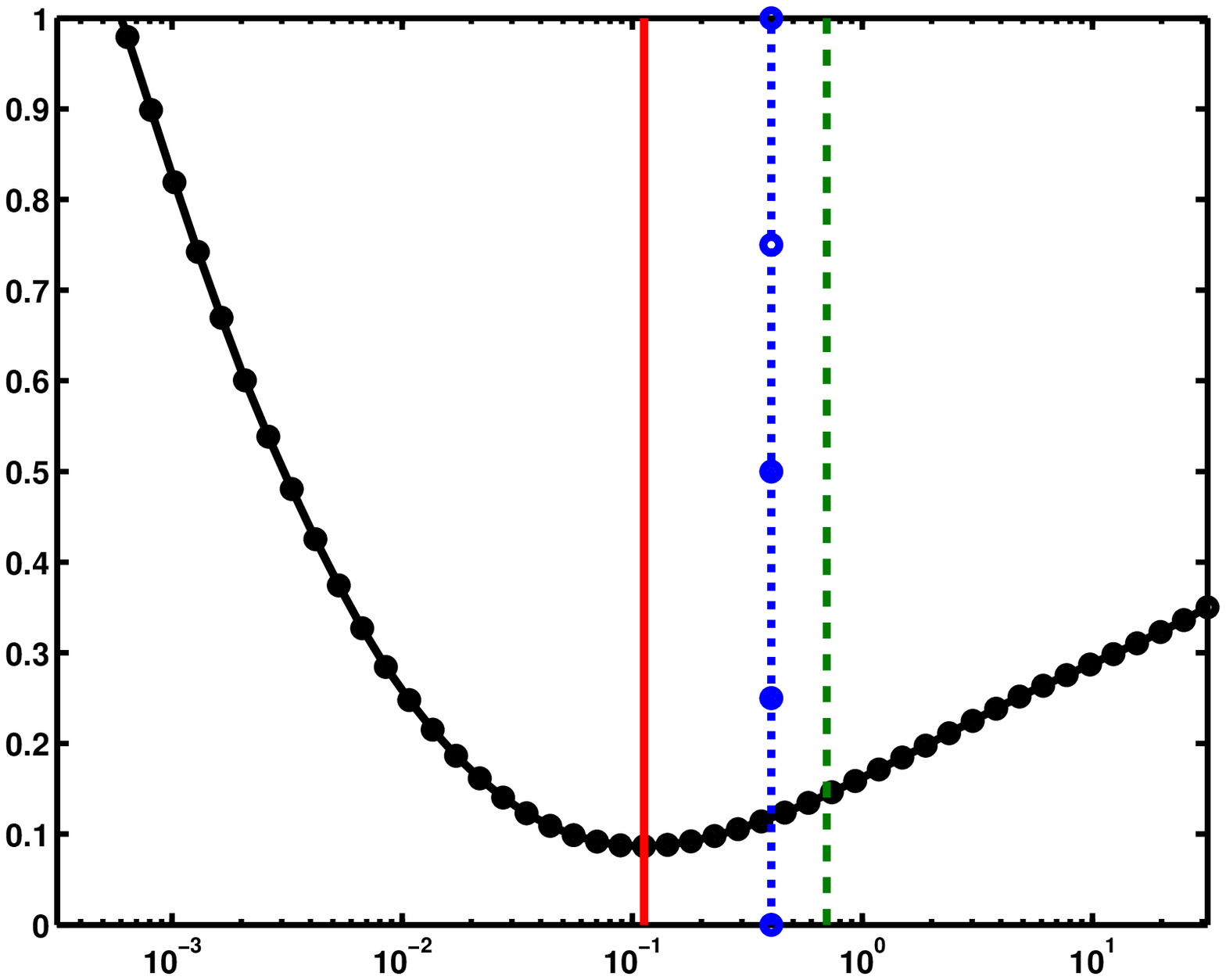}}
\subfigure[$L=I$]{\includegraphics[width=1.7in]{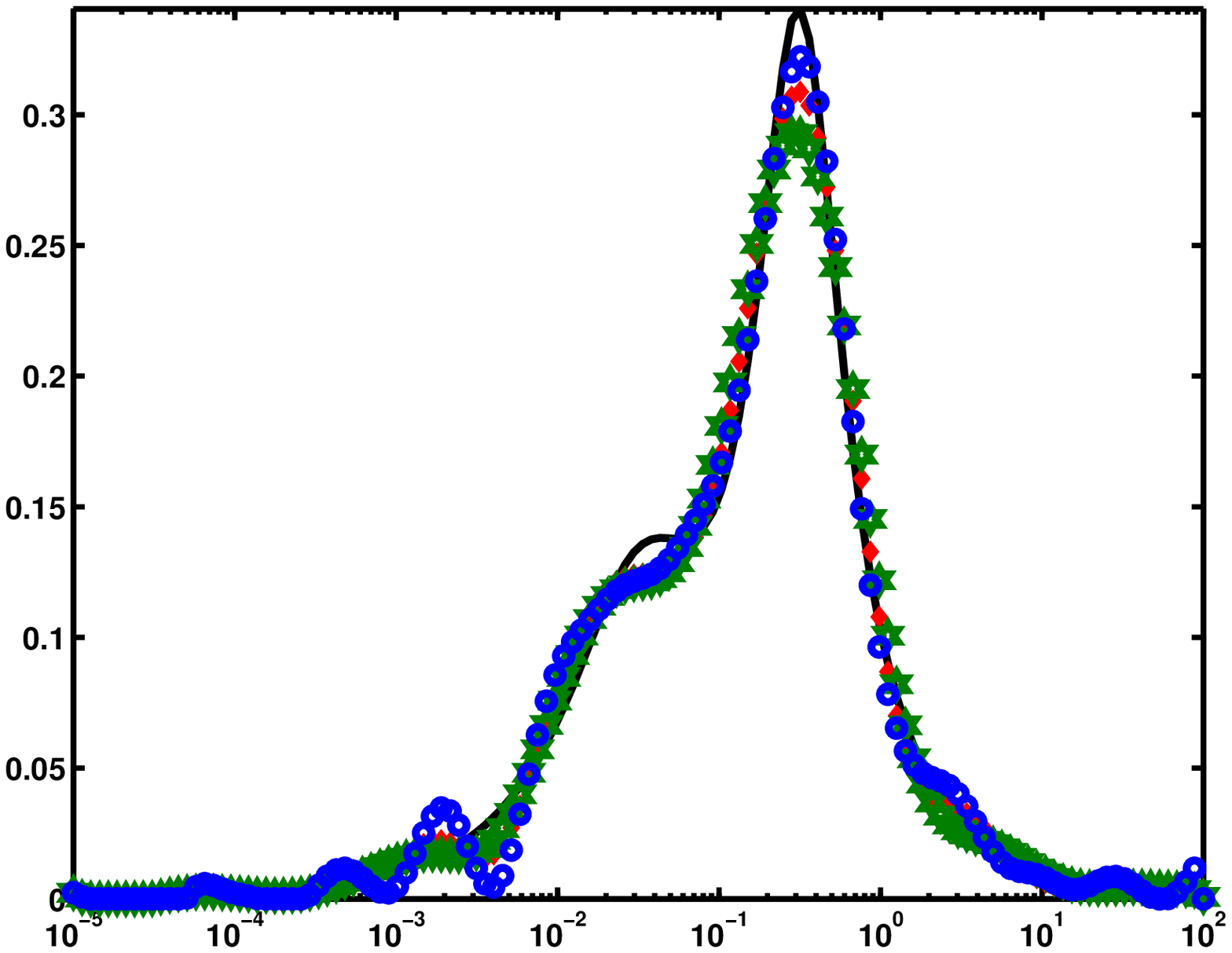}}
\subfigure[$L=L_1$]{\includegraphics[width=1.7in]{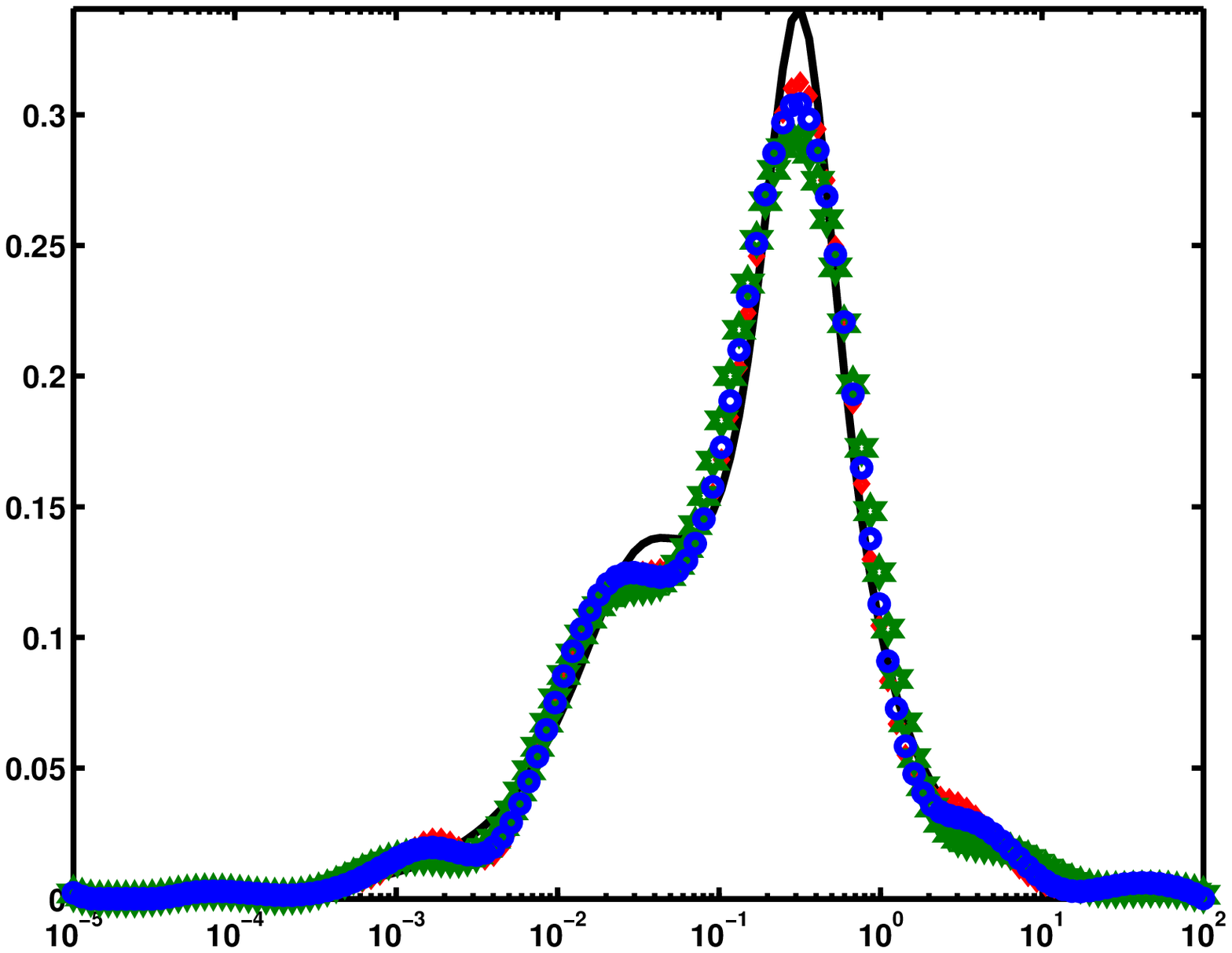}}
\subfigure[$L=L_2$]{\includegraphics[width=1.7in]{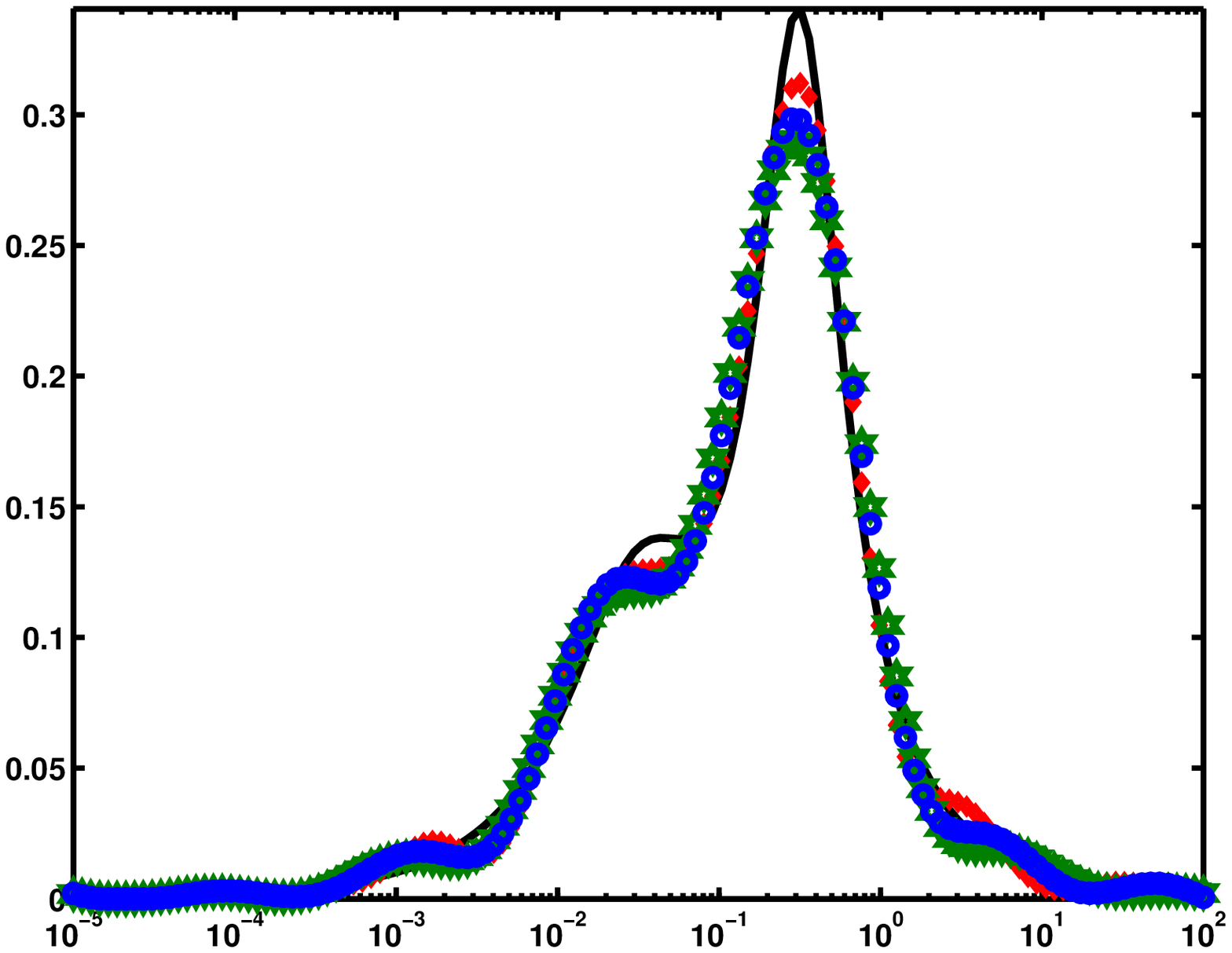}}
\caption{Mean error and example \texttt{lsqnonneg}  NNLS solutions. $.1\%$ noise, RQ-B data set, matrix $A_4$.}
\label{fig-lambdachoiceRQ5A4LN}
\end{figure}

\begin{figure}[!ht]
  \centering
\subfigure[$L=I$]{\includegraphics[width=1.7in]{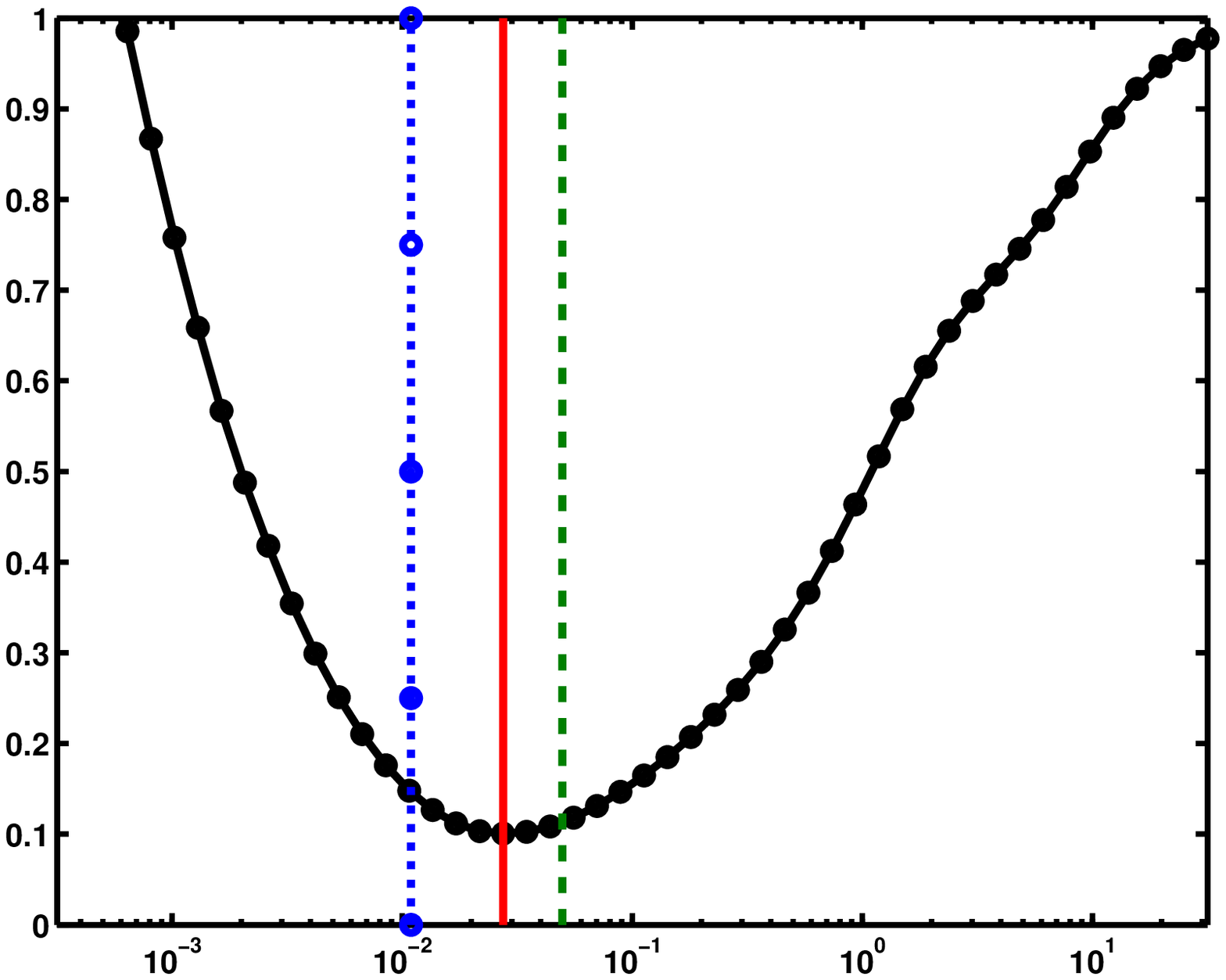}}
\subfigure[$L=L_1$]{\includegraphics[width=1.7in]{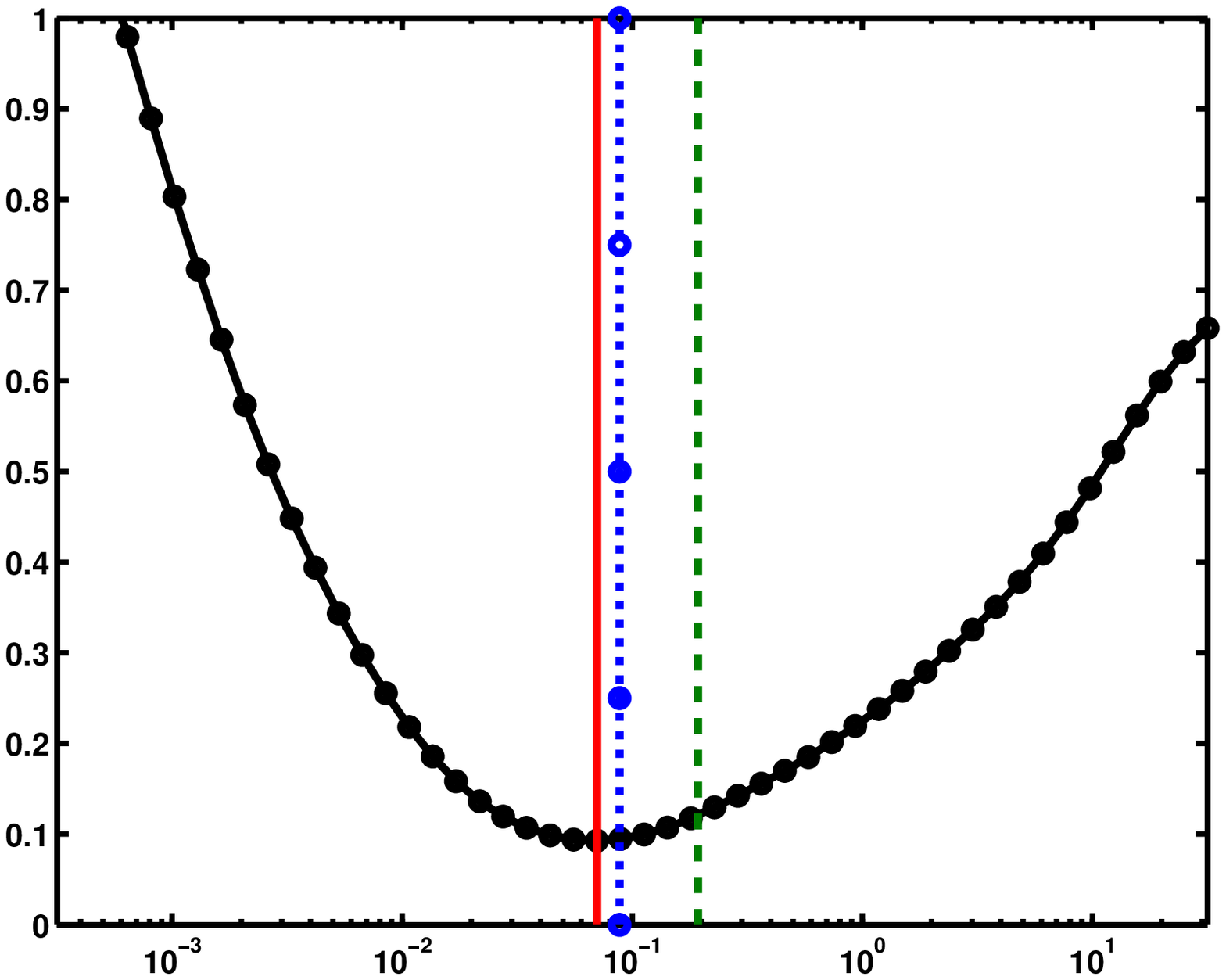}}
\subfigure[$L=L_2$]{\includegraphics[width=1.7in]{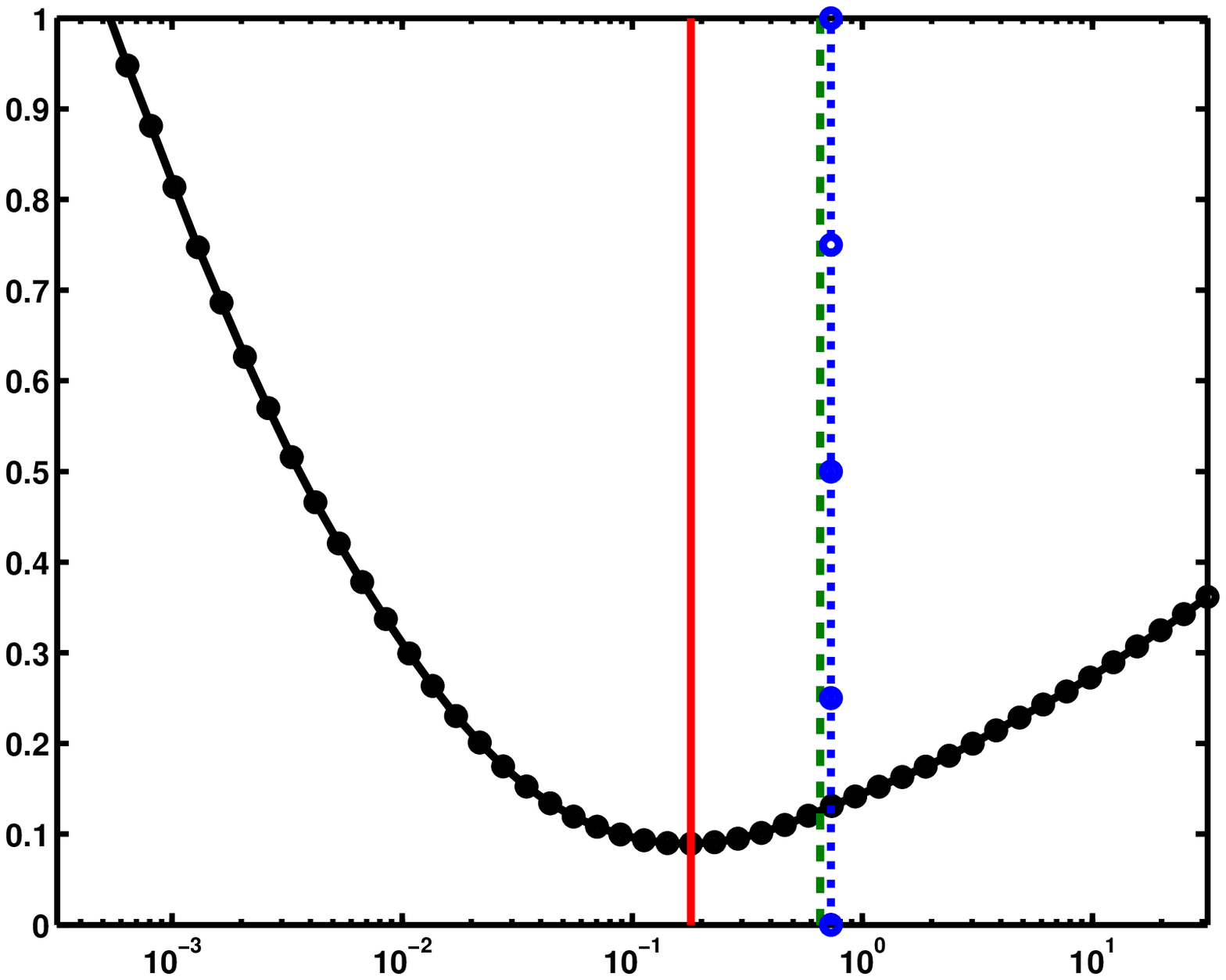}}
\subfigure[$L=I$]{\includegraphics[width=1.7in]{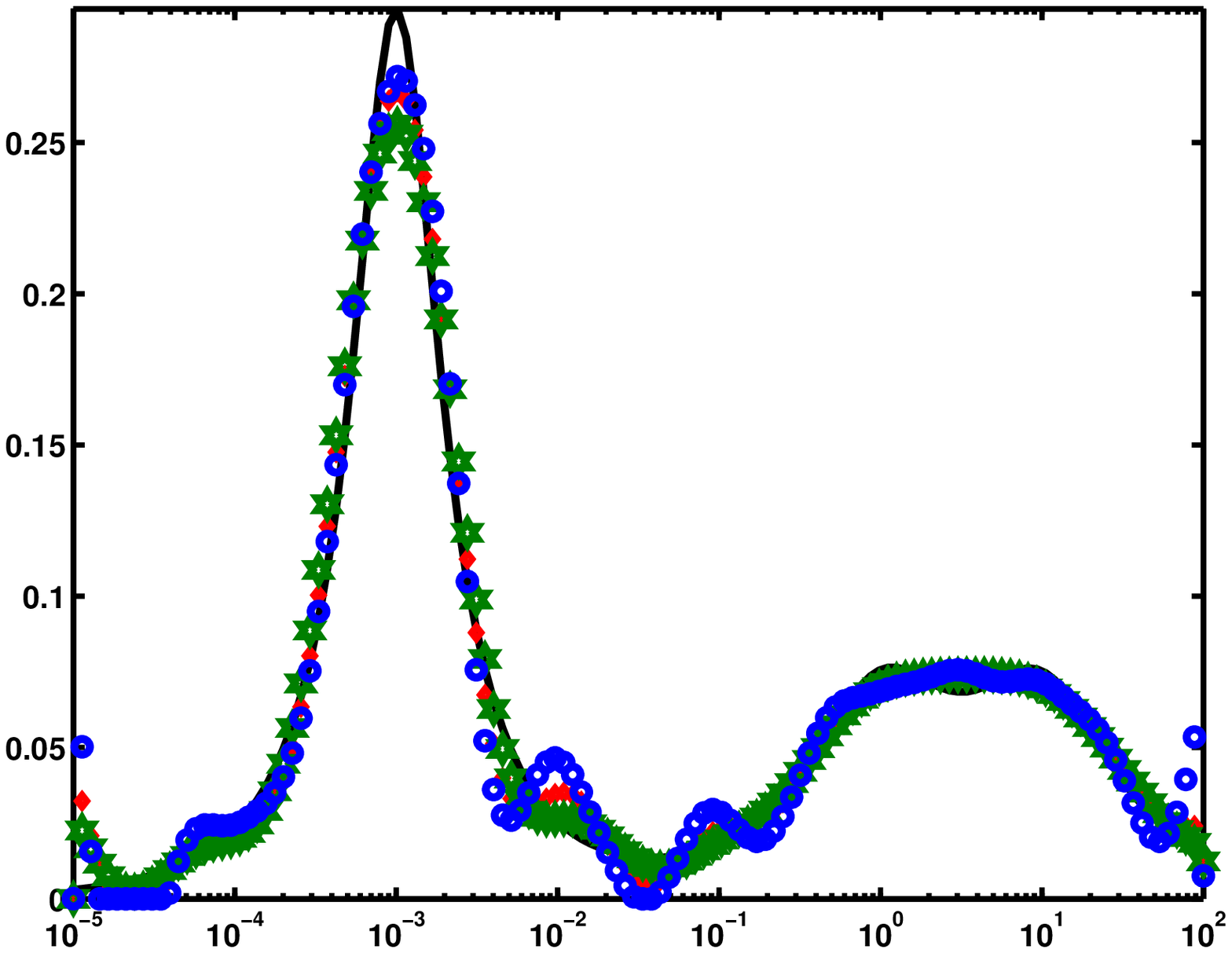}}
\subfigure[$L=L_1$]{\includegraphics[width=1.7in]{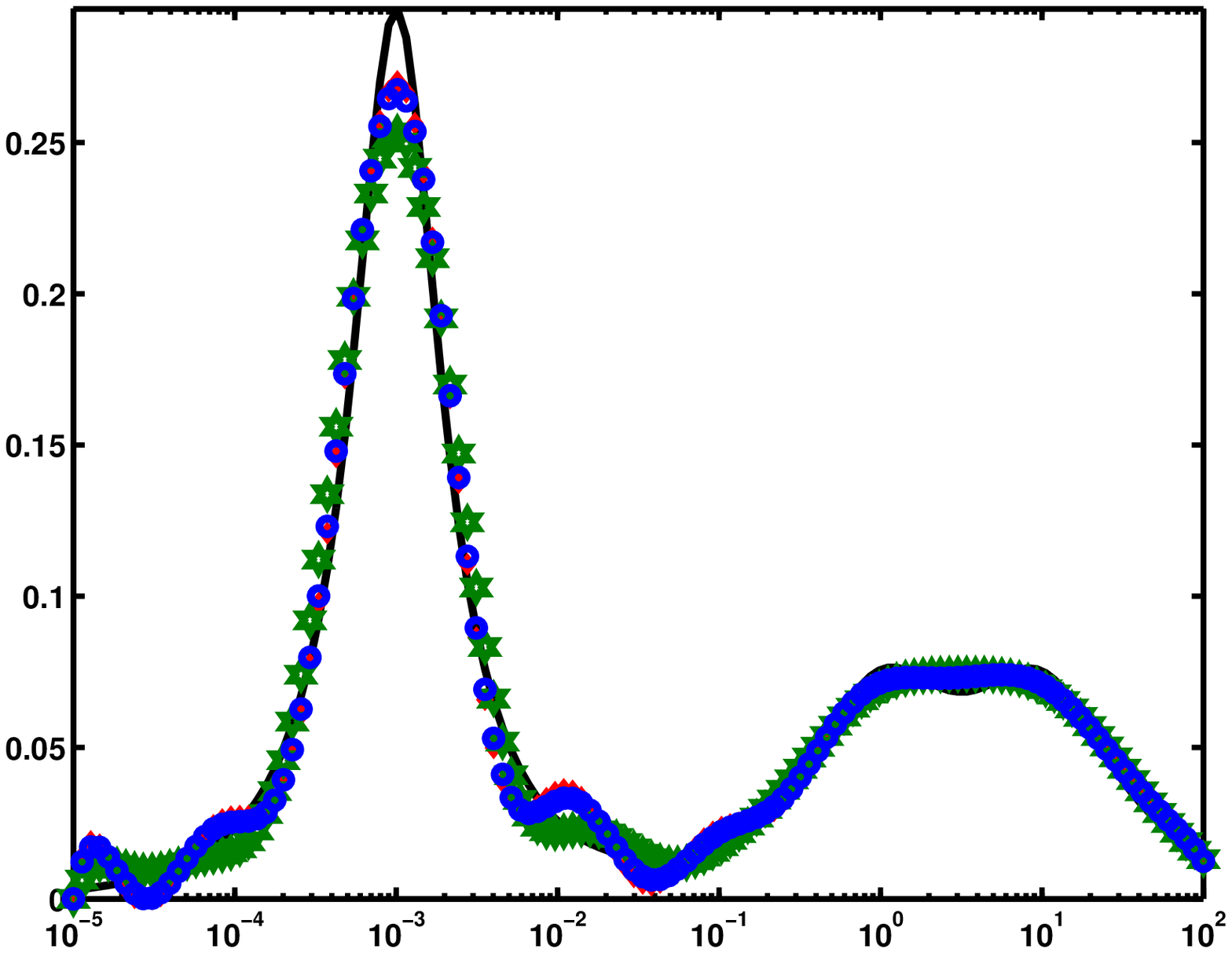}}
\subfigure[$L=L_2$]{\includegraphics[width=1.7in]{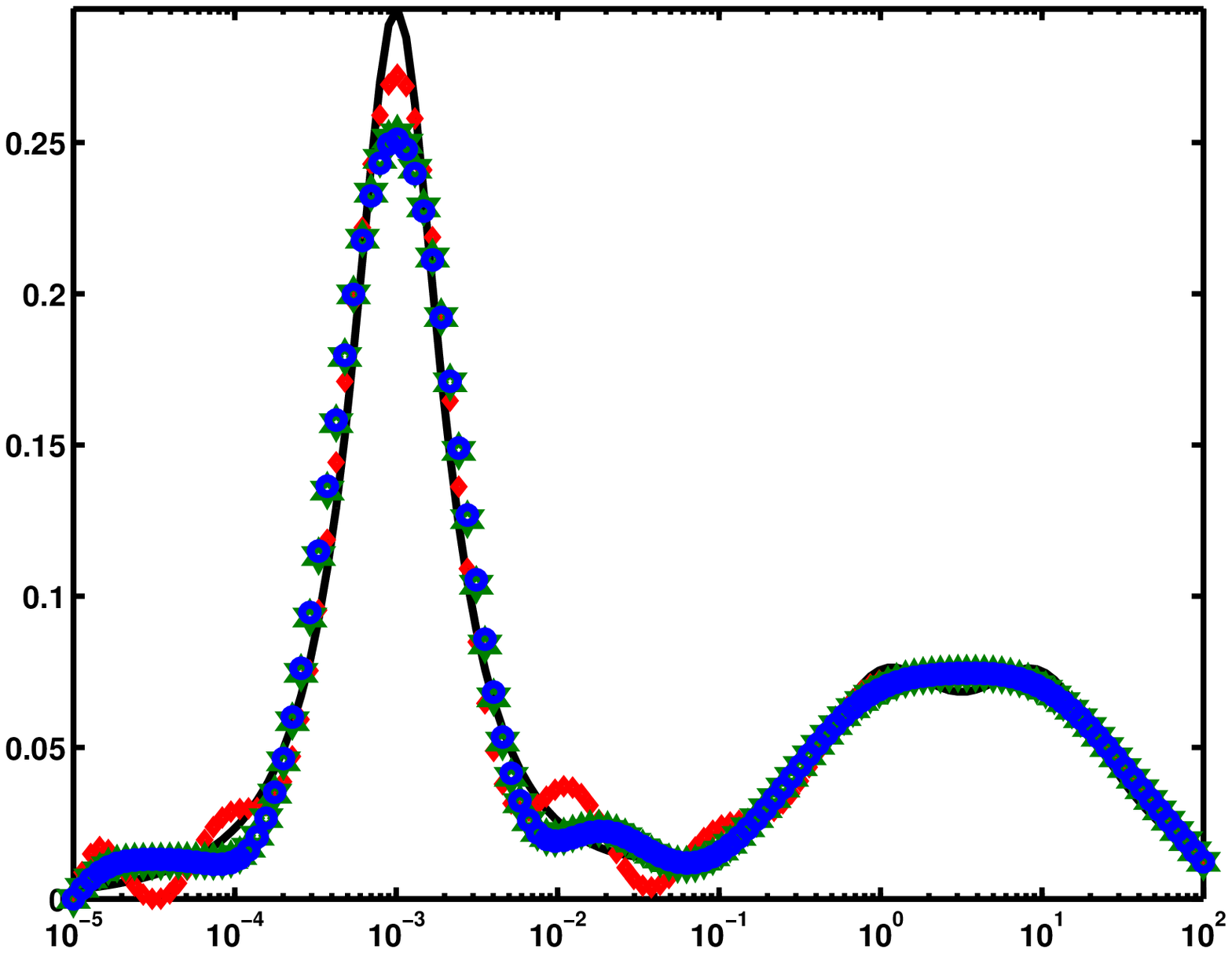}}
\caption{Mean error and example \texttt{lsqnonneg}  NNLS solutions.  $.1\%$ noise, RQ-C data set, matrix $A_4$.}
\label{fig-lambdachoiceRQ6A4LN}
\end{figure}

\begin{figure}[!ht]
  \centering
\subfigure[$L=I$]{\includegraphics[width=1.7in]{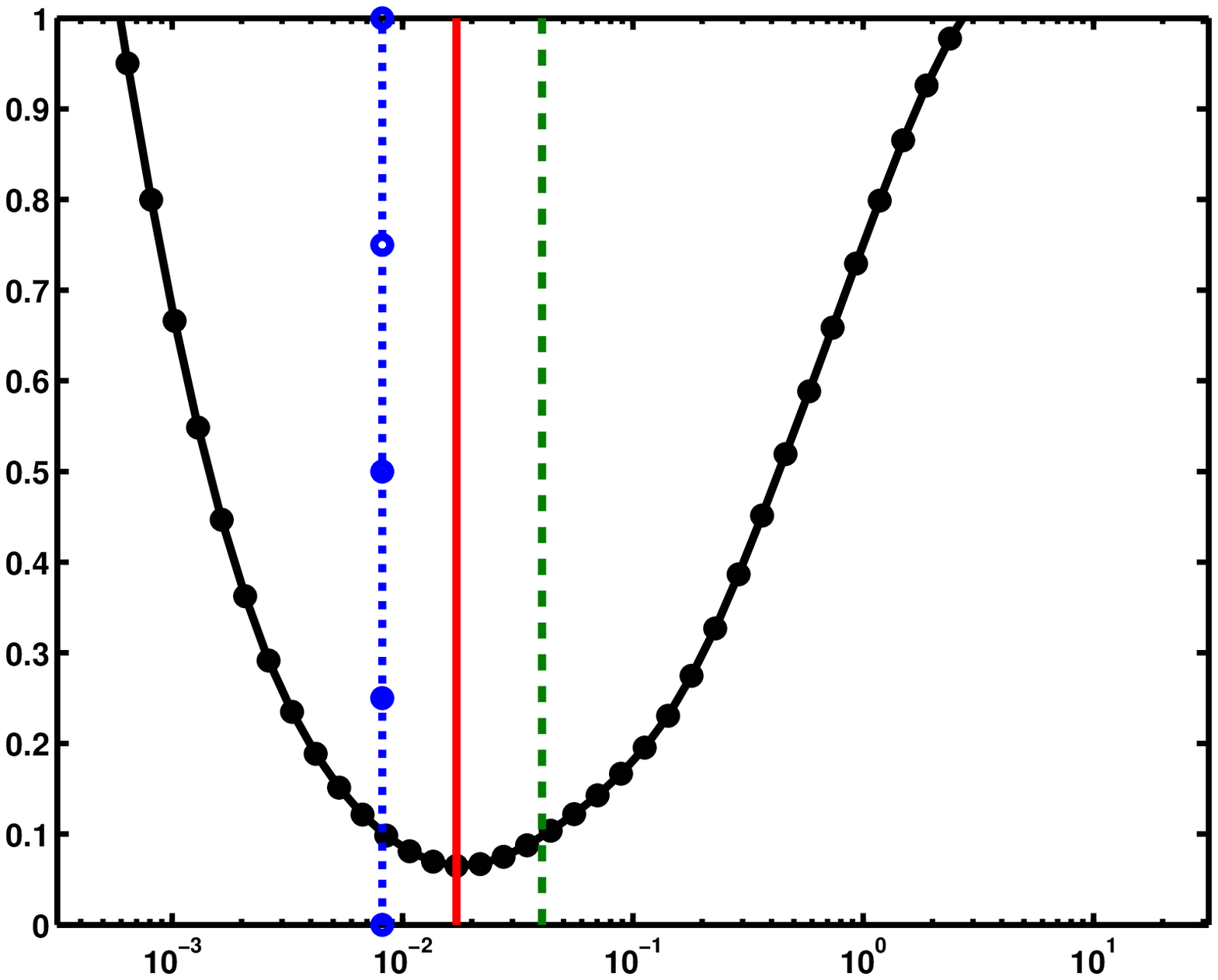}}
\subfigure[$L=L_1$]{\includegraphics[width=1.7in]{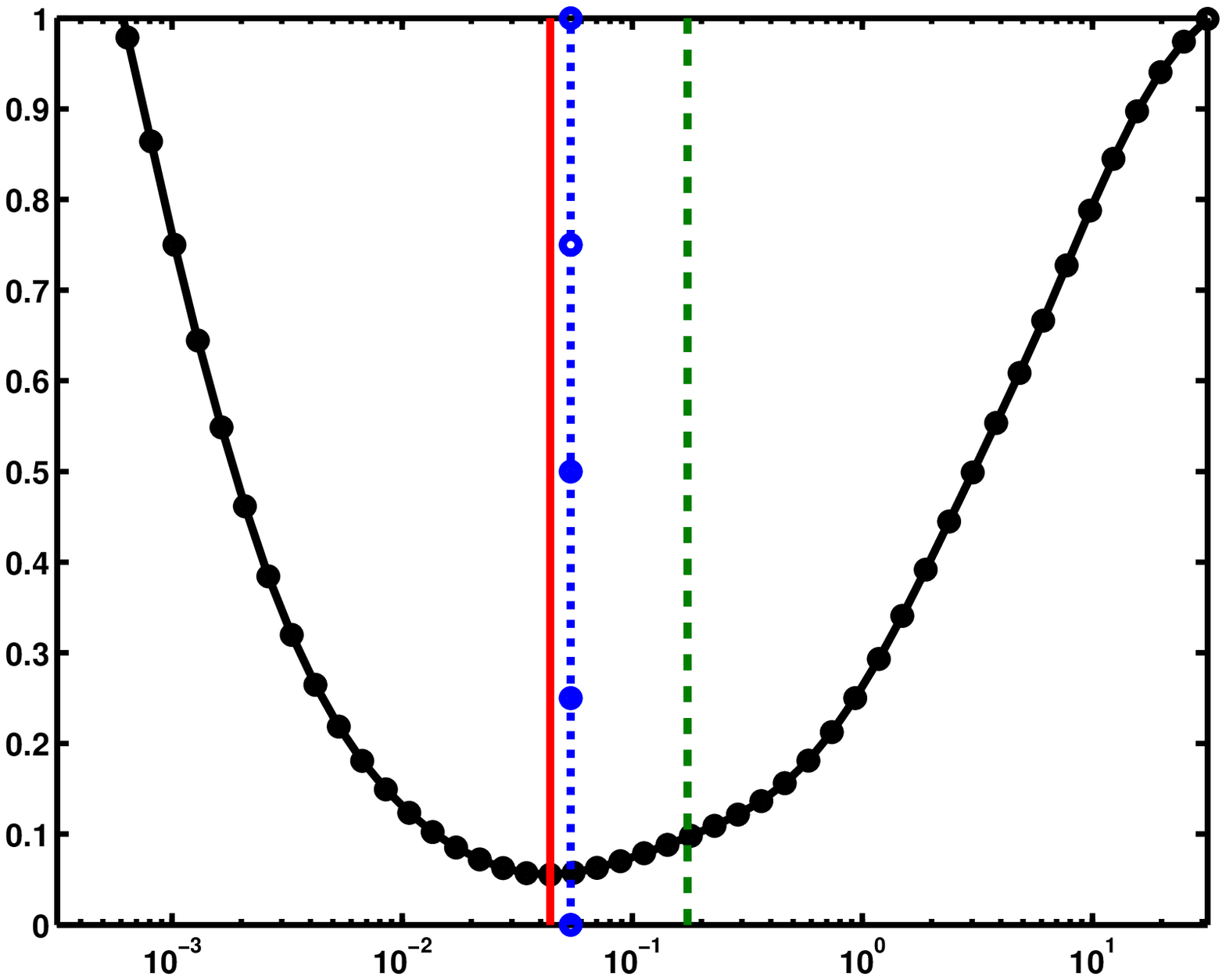}}
\subfigure[$L=L_2$]{\includegraphics[width=1.7in]{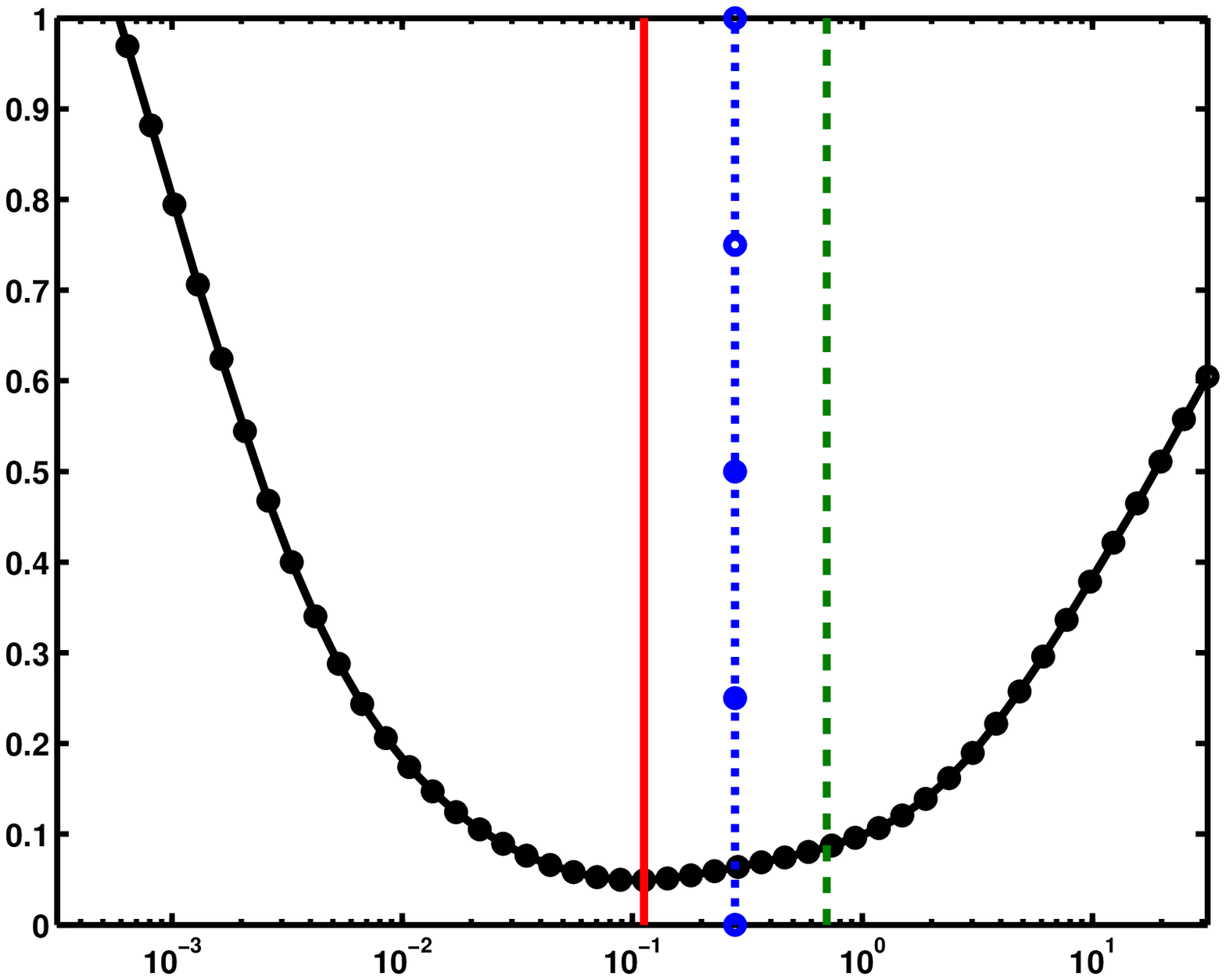}}
\subfigure[$L=I$]{\includegraphics[width=1.7in]{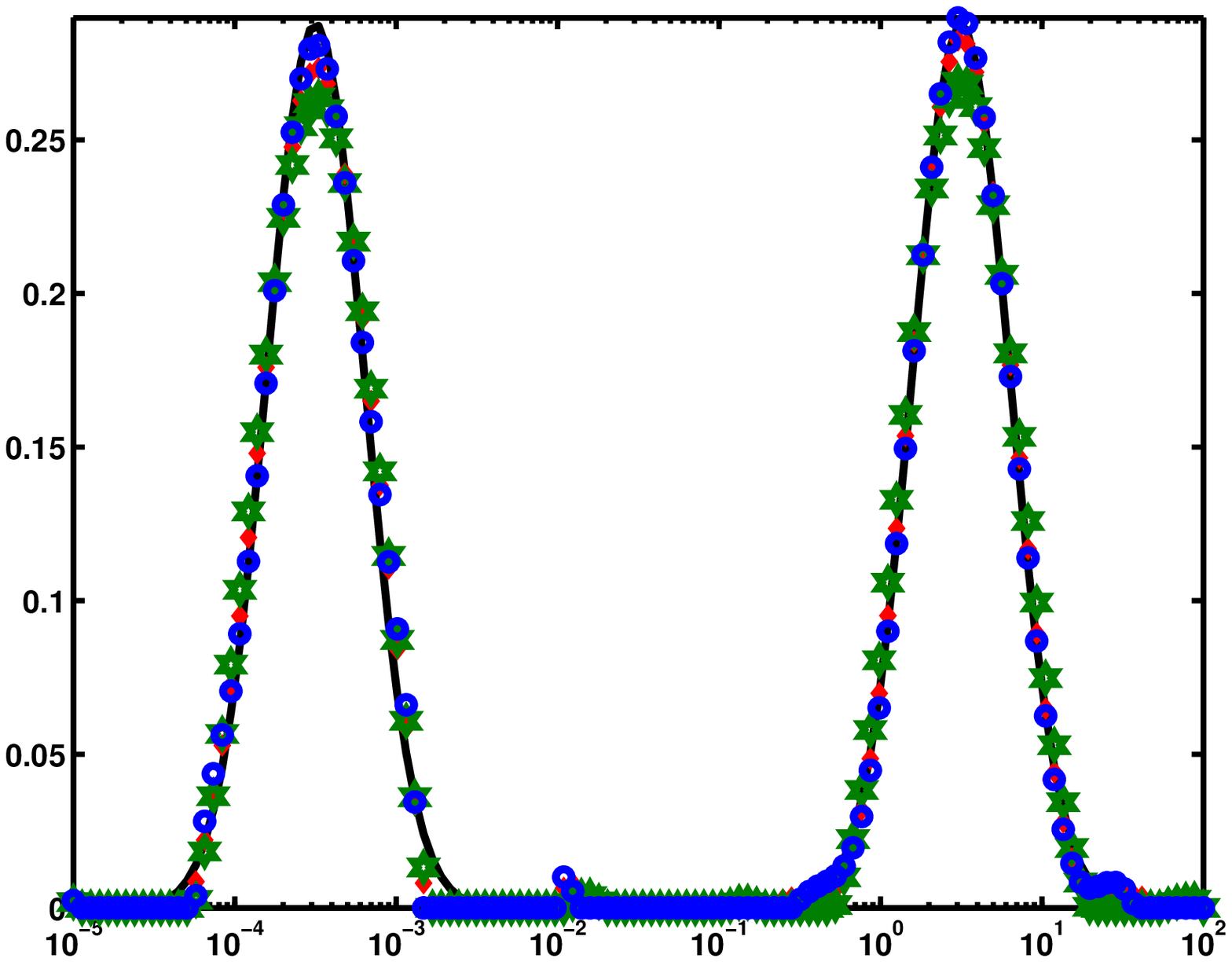}}
\subfigure[$L=L_1$]{\includegraphics[width=1.7in]{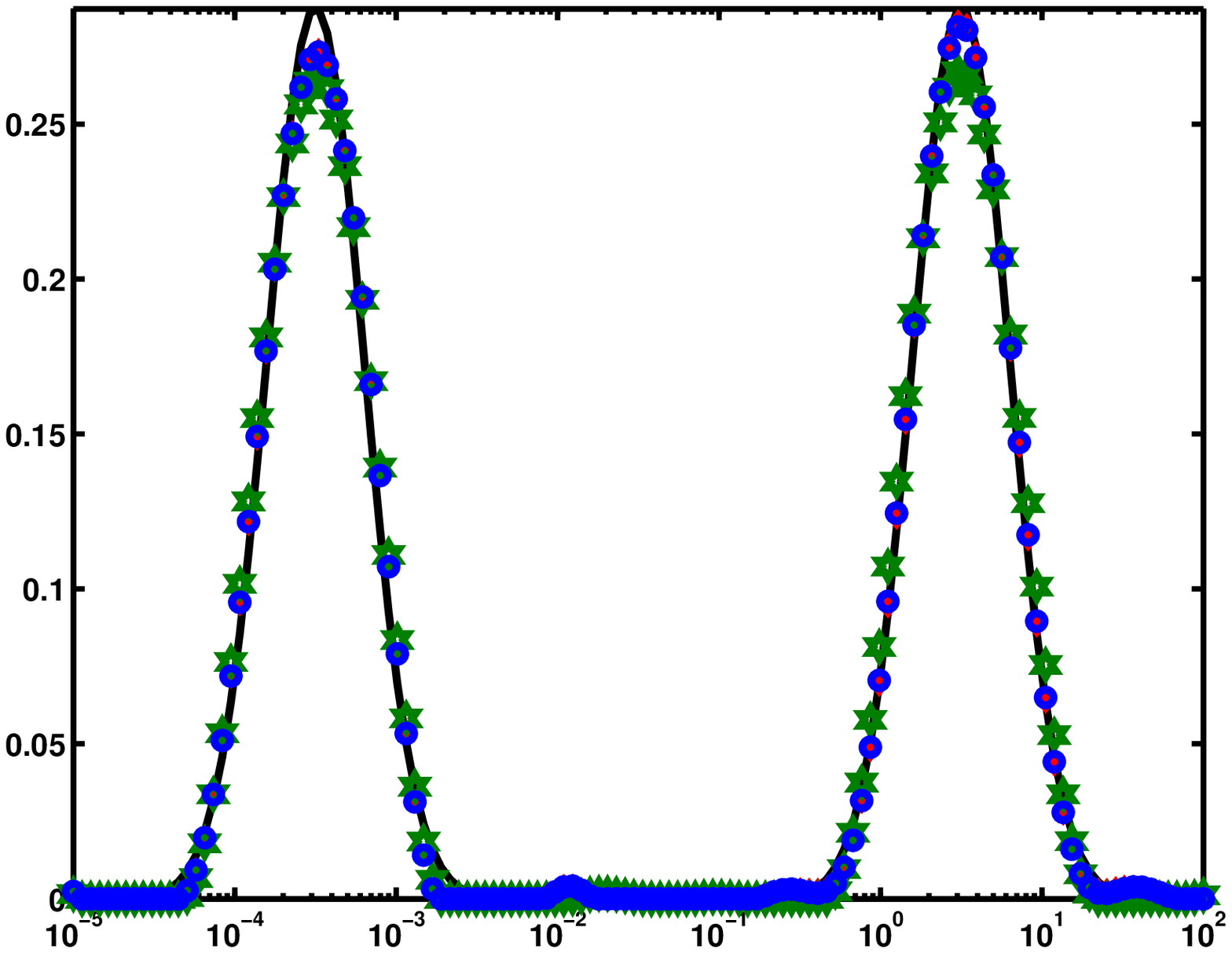}}
\subfigure[$L=L_2$]{\includegraphics[width=1.7in]{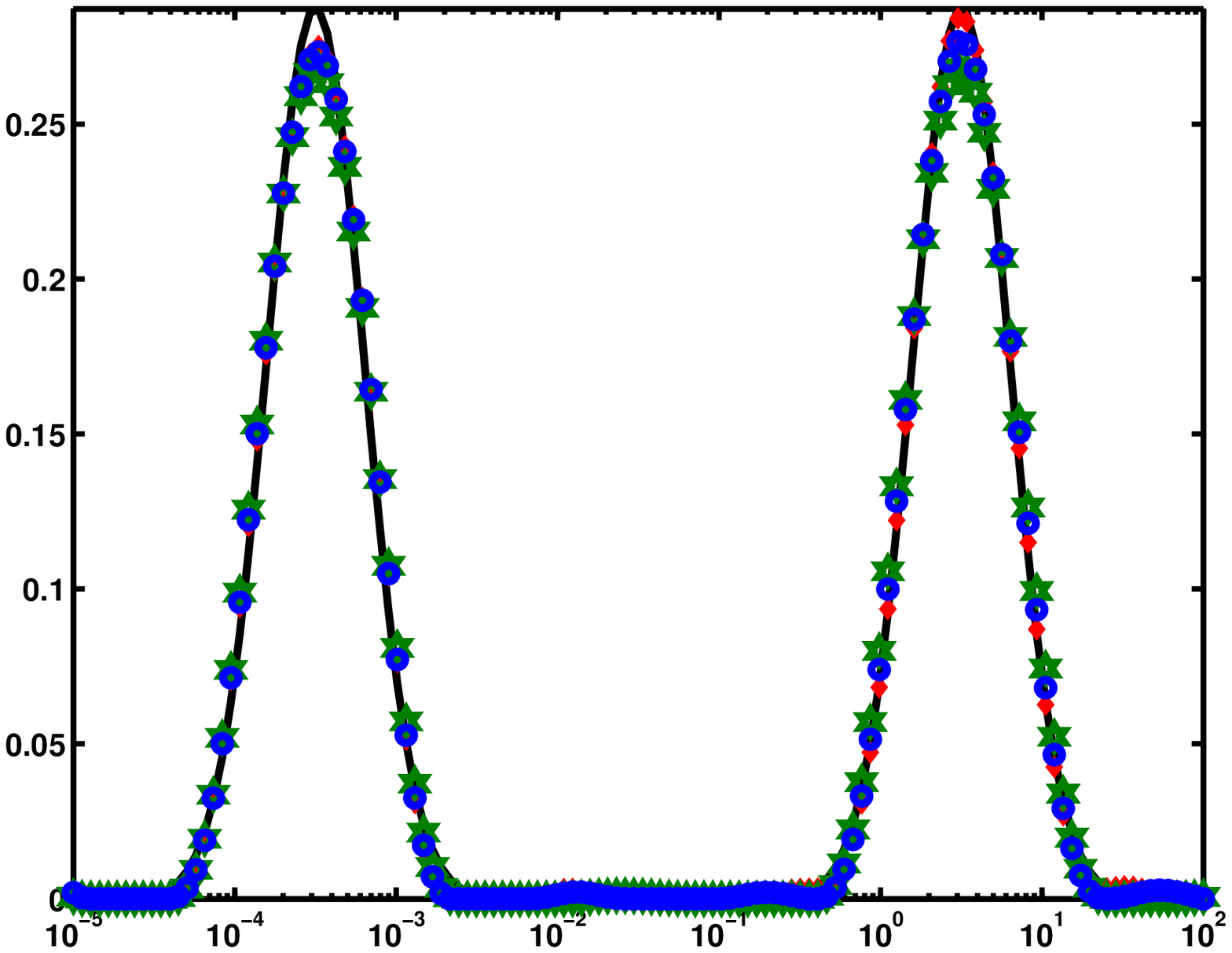}}
\caption{Mean error and example \texttt{lsqnonneg} NNLS solutions.  $.1\%$ noise, LN-A data set, matrix $A_4$.}
\label{fig-lambdachoiceLN2A4LN}
\end{figure}

\begin{figure}[!ht]
  \centering
\subfigure[$L=I$]{\includegraphics[width=1.7in]{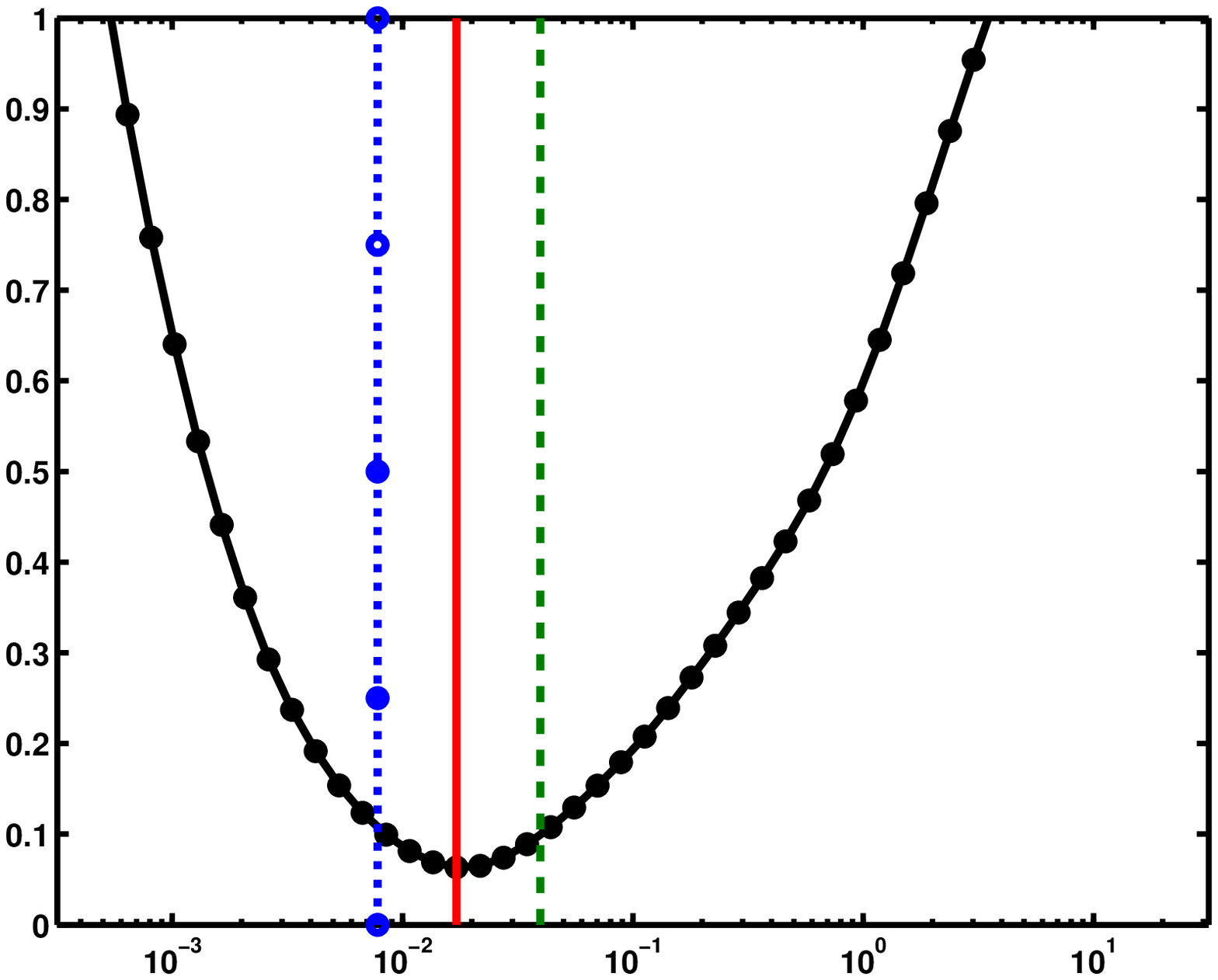}}
\subfigure[$L=L_1$]{\includegraphics[width=1.7in]{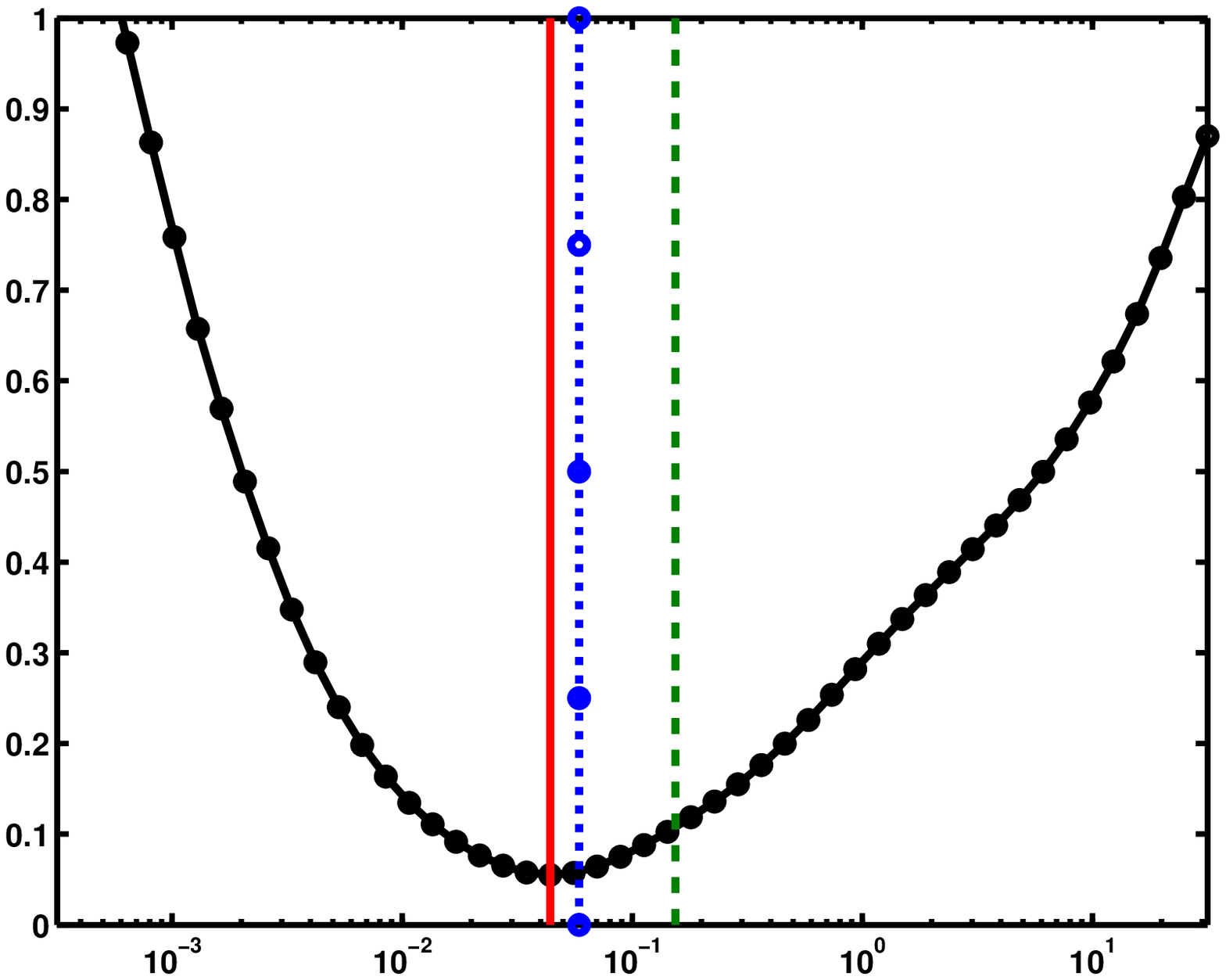}}
\subfigure[$L=L_2$]{\includegraphics[width=1.7in]{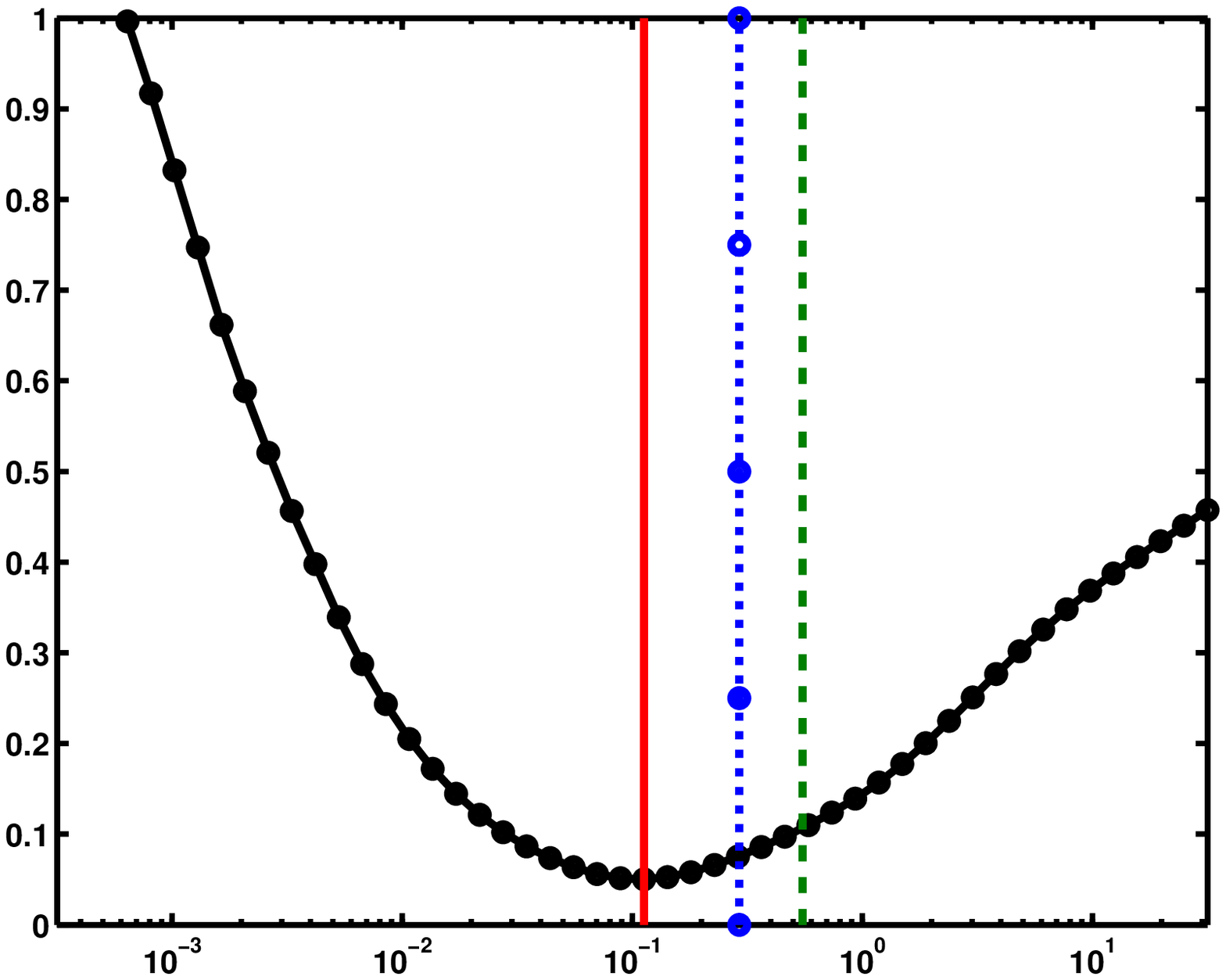}}
\subfigure[$L=I$]{\includegraphics[width=1.7in]{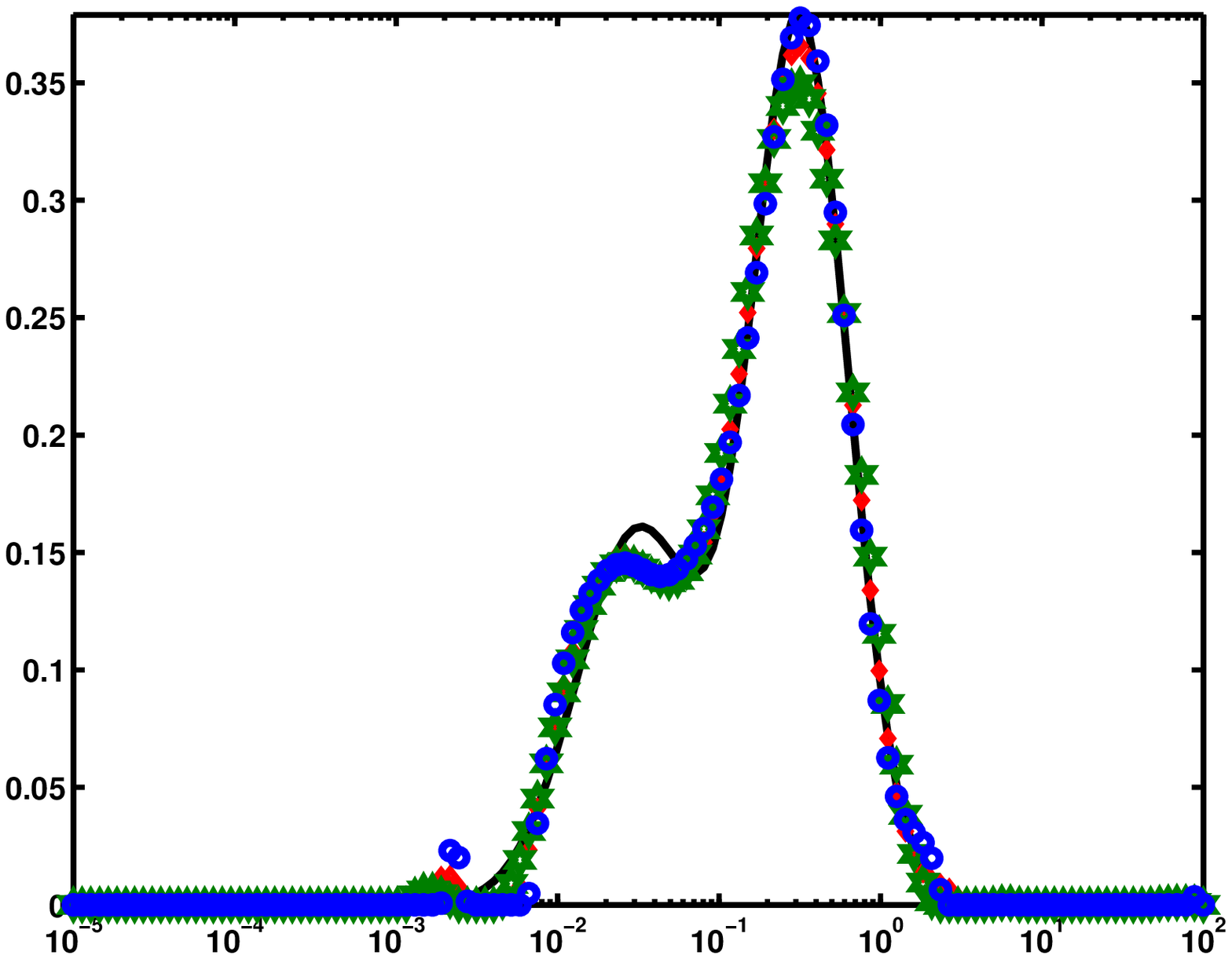}}
\subfigure[$L=L_1$]{\includegraphics[width=1.7in]{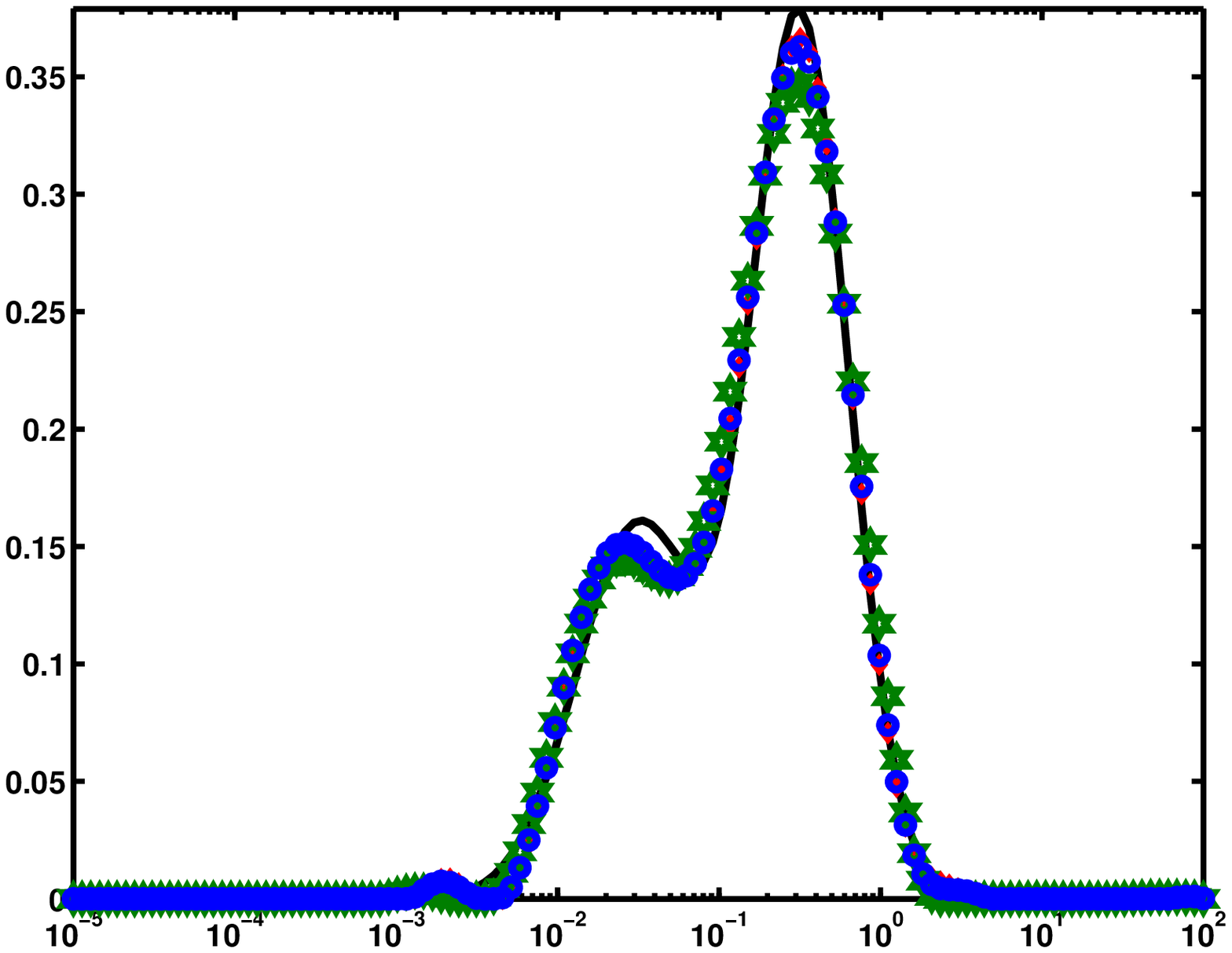}}
\subfigure[$L=L_2$]{\includegraphics[width=1.7in]{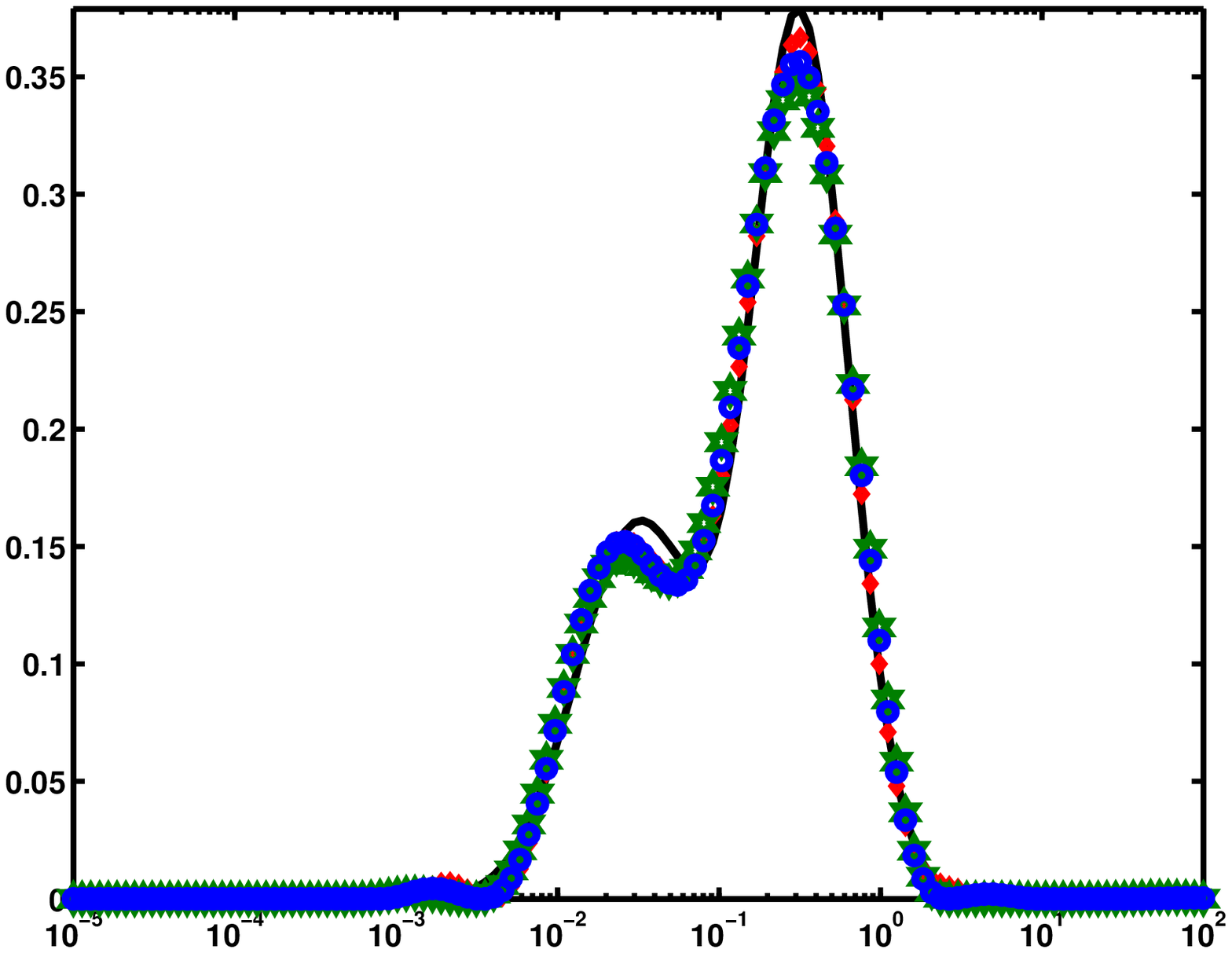}}
\caption{Mean error and example \texttt{lsqnonneg} NNLS solutions.  $.1\%$ noise, LN-B data set, matrix $A_4$.}
\label{fig-lambdachoiceLN5A4LN}
\end{figure}

\begin{figure}[!ht]
  \centering
\subfigure[$L=I$]{\includegraphics[width=1.7in]{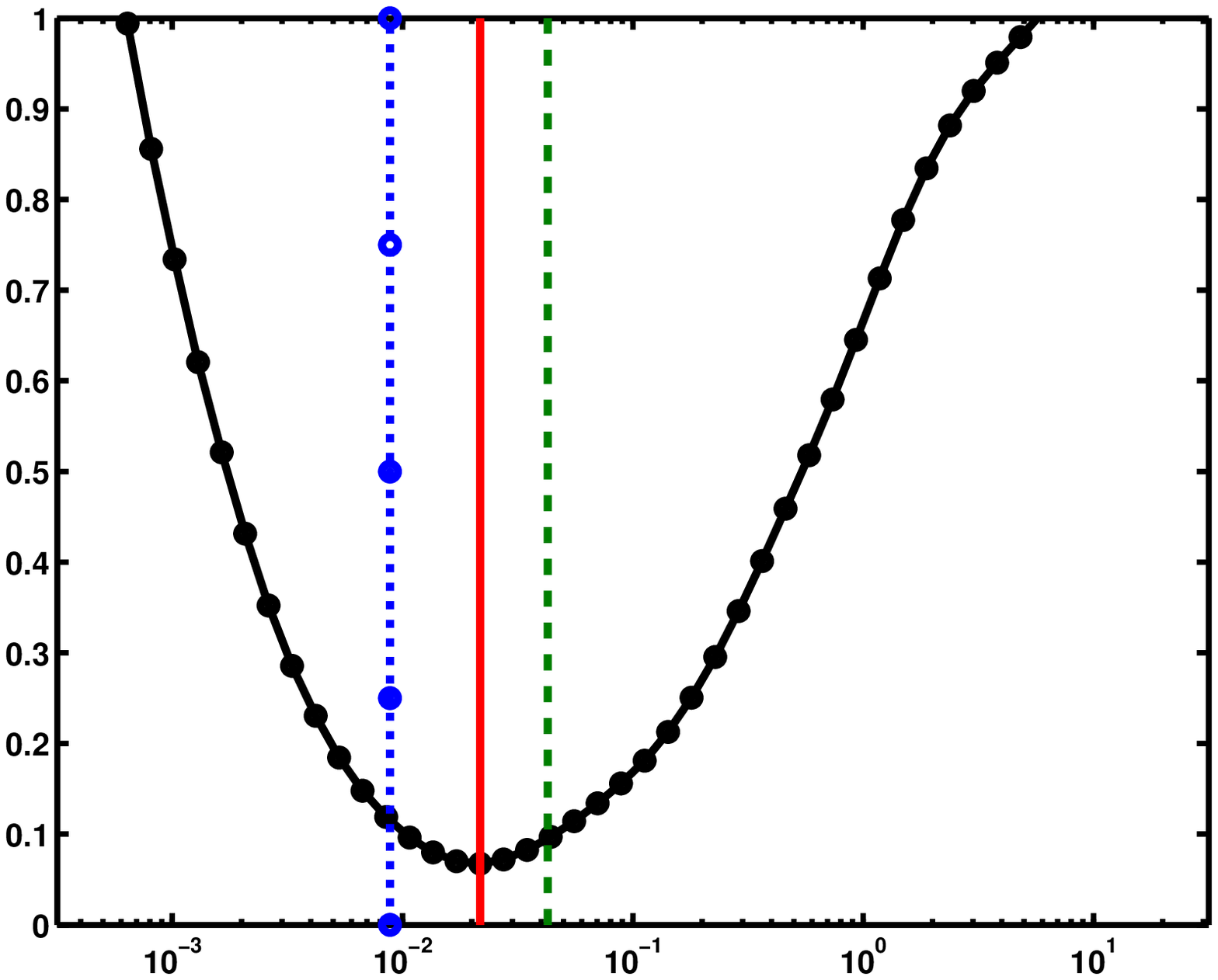}}
\subfigure[$L=L_1$]{\includegraphics[width=1.7in]{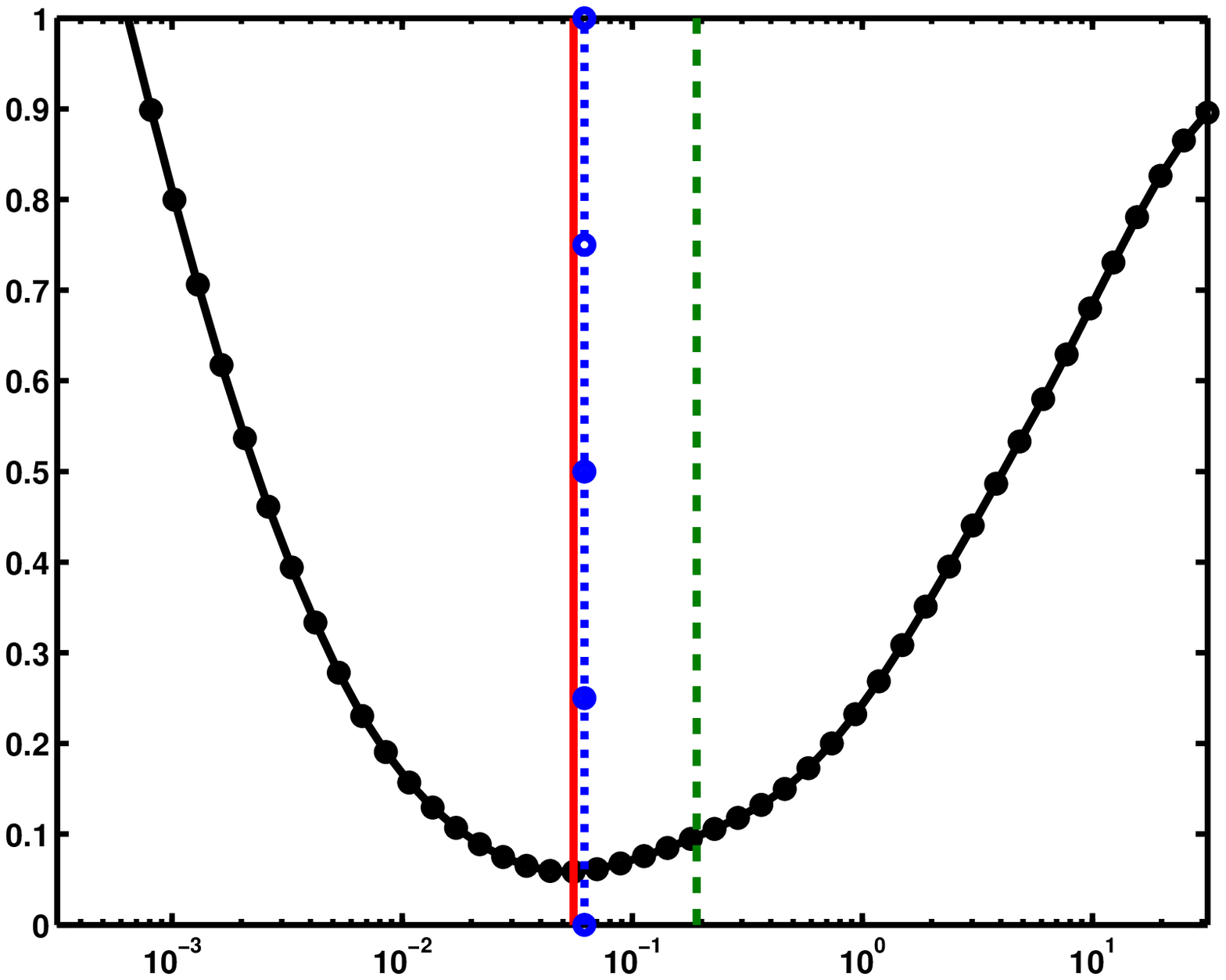}}
\subfigure[$L=L_2$]{\includegraphics[width=1.7in]{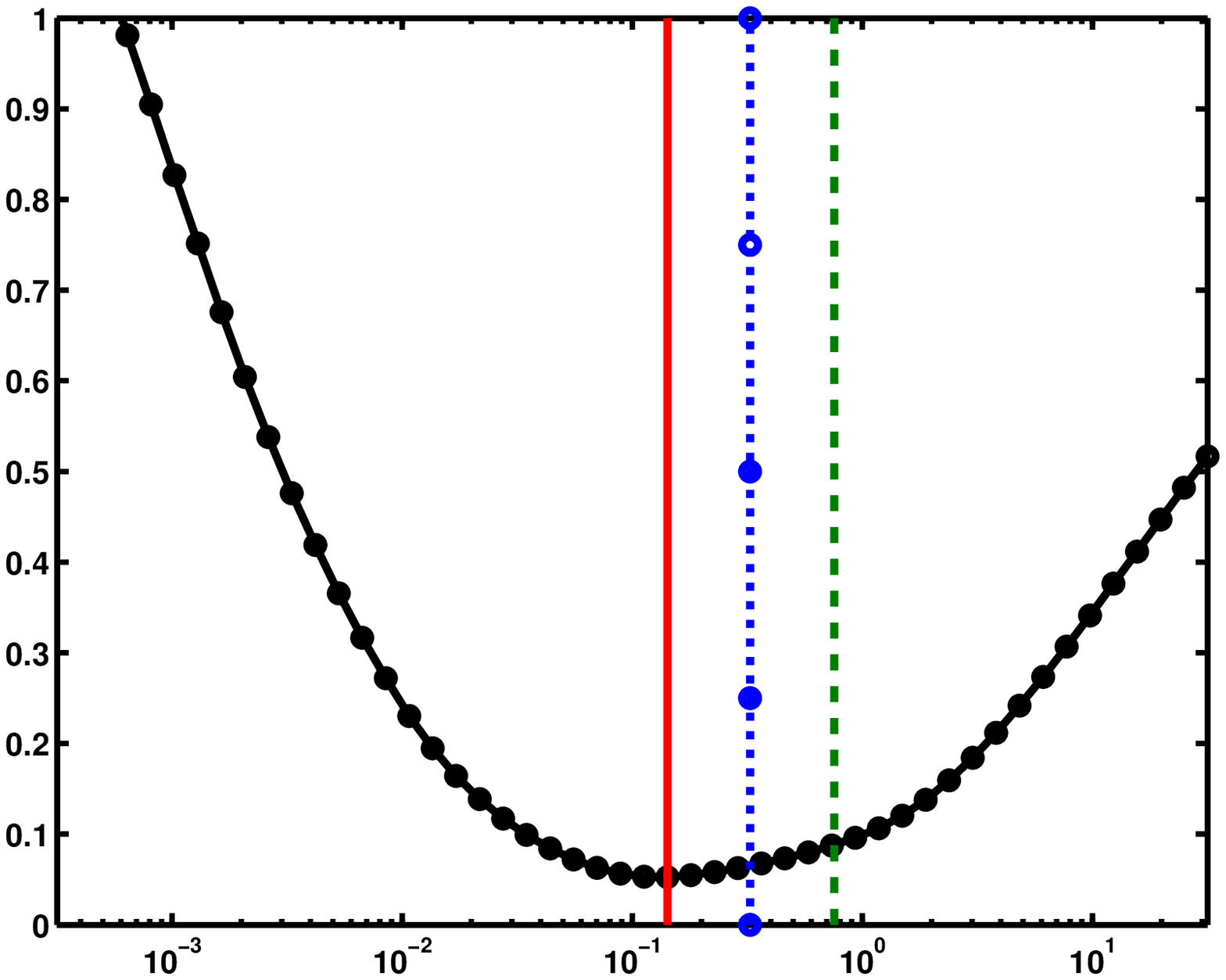}}
\subfigure[$L=I$]{\includegraphics[width=1.7in]{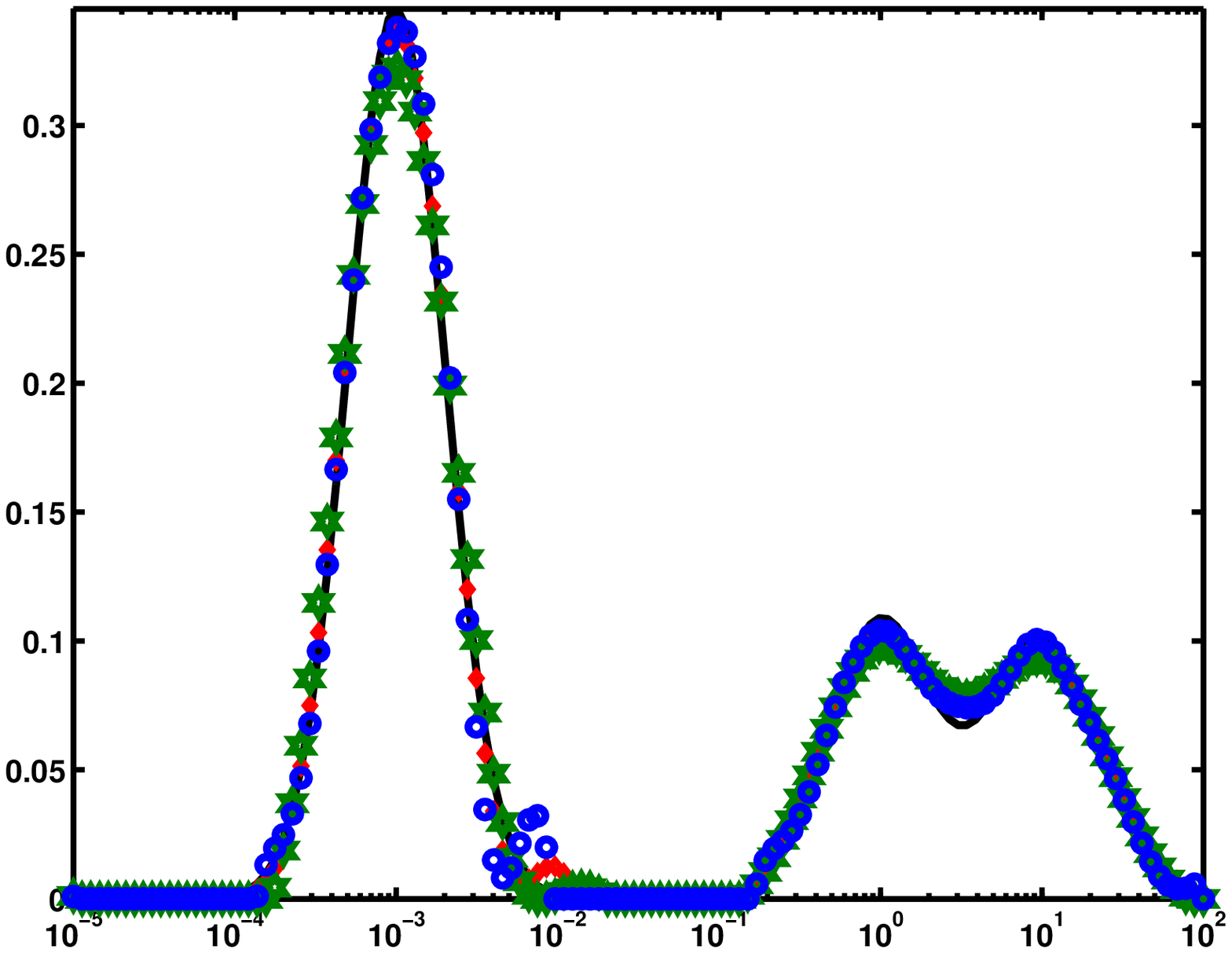}}
\subfigure[$L=L_1$]{\includegraphics[width=1.7in]{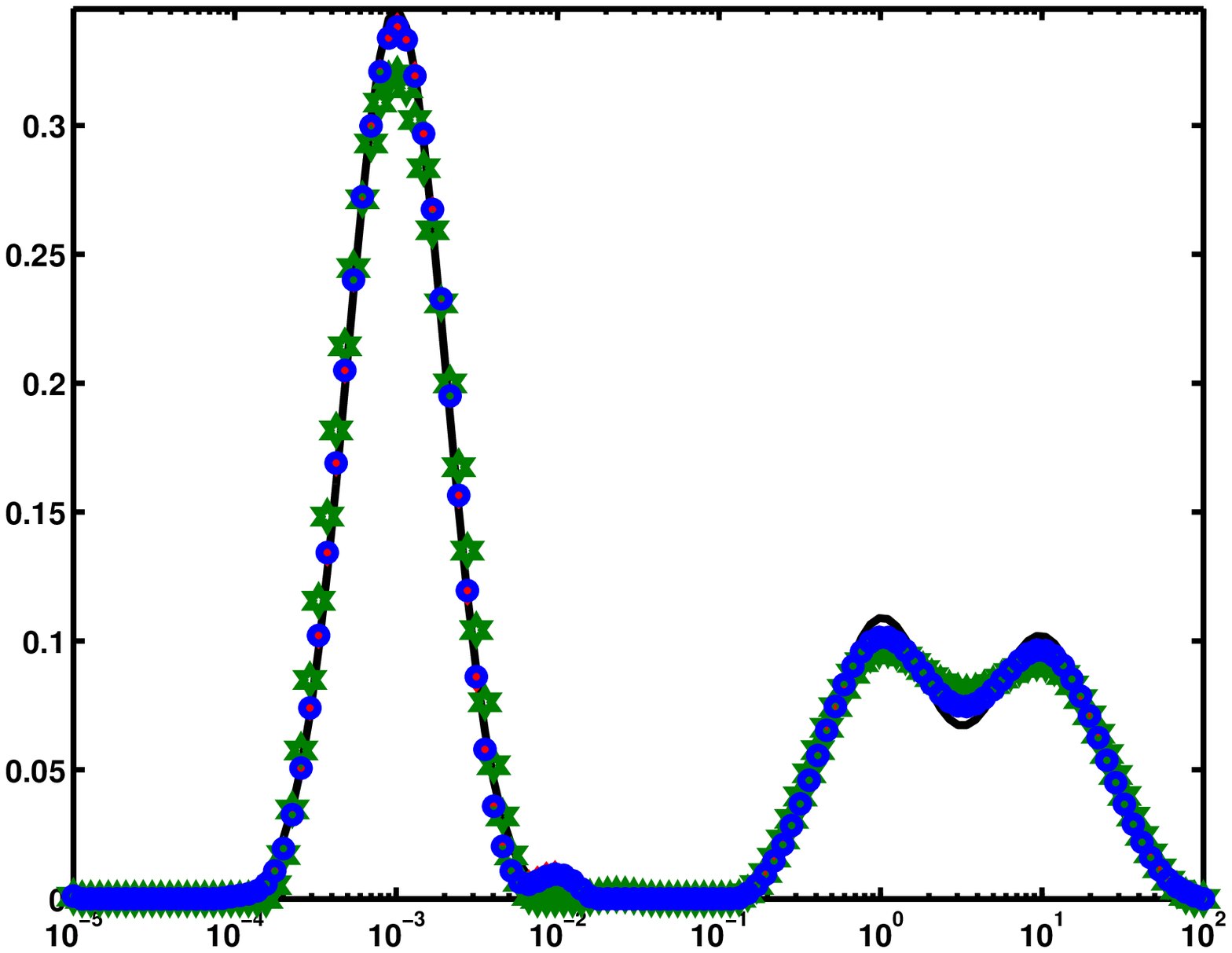}}
\subfigure[$L=L_2$]{\includegraphics[width=1.7in]{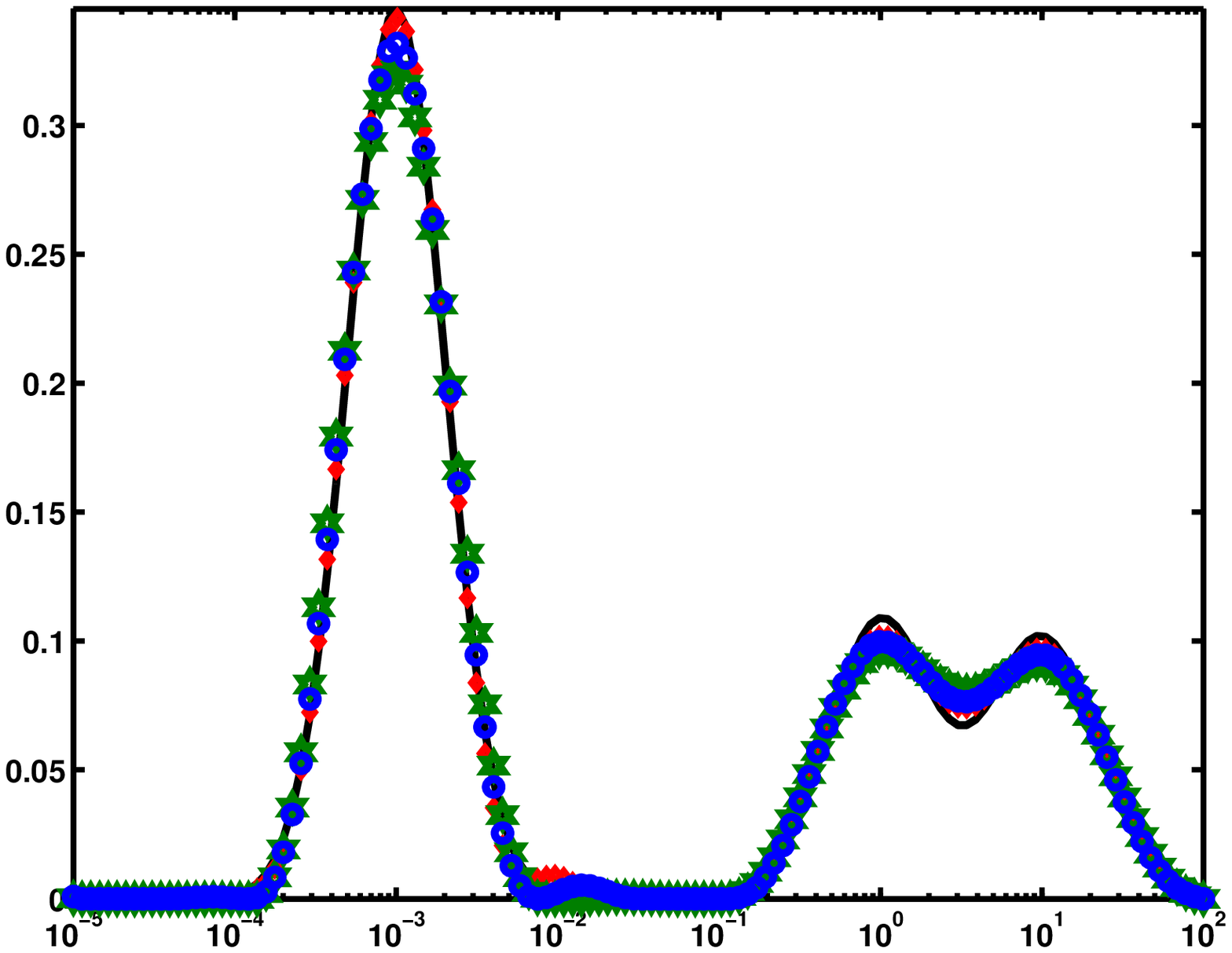}}
\caption{Mean error and example \texttt{lsqnonneg} NNLS solutions.  $.1\%$ noise, LN-C data set, matrix $A_4$.}
\label{fig-lambdachoiceLN6A4LN}
\end{figure}

The results are illustrated in Figures~\ref{fig-lambdachoiceRQ1A4LN}-\ref{fig-lambdachoiceLN6A4LN}. Figures (a)-(c) in each case indicate the mean error results for the different smoothing norms, and (d)-(f) demonstrate the sensitivity, or lack thereof, of the solution to the choice of $\lambda$ near the optimum. We see that the results are remarkably consistent; the results with the identity weighted norm are generally less robust, while overall the NCP parameter choice marginally outperforms the LC. On the other hand, we also conclude that the use of either parameter choice method is robust in terms of representing the \textit{optimal} but practically {\it unknown} solution. Thus the NCP is to be preferred for finding a suitable regularization parameter, and the $L_1$ operator provides a compromise between over smoothing (a reduced peak) by $L_2$ and under smoothing by $L=I$. 
Additional results for higher noise levels are provided in \cite{Supp}. There it is shown that for noise levels at $5\%$ the results deteriorate significantly in terms of the ability to accurately determine the number of processes, that the LC results are significantly under smoothed and that the extra resolution of matrix $A_4$ is most obviously worthwhile, a result that is not at all clear by examination of the relative error. Overall, it is also clear that the LN processes are better resolved than those given by the RQ. This is not surprising given the graph of the RQ DRT for small $t$, see e.g. Figure~\ref{1a}, \ref{2a}. 
 
\begin{figure}[!ht]
 \centering
\subfigure[$L=I$]{\includegraphics[width=1.7in]{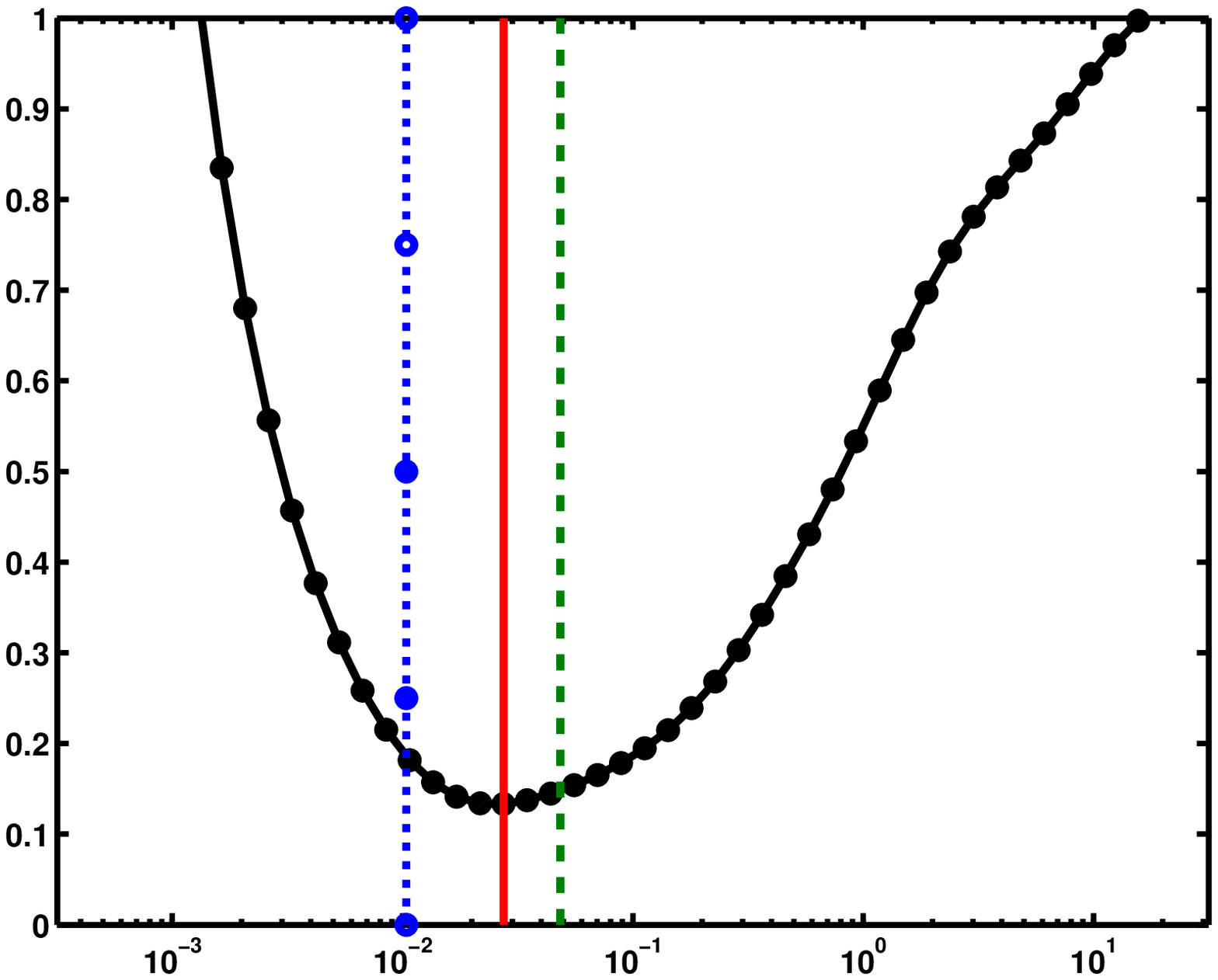}}
\subfigure[$L=L_1$]{\includegraphics[width=1.7in]{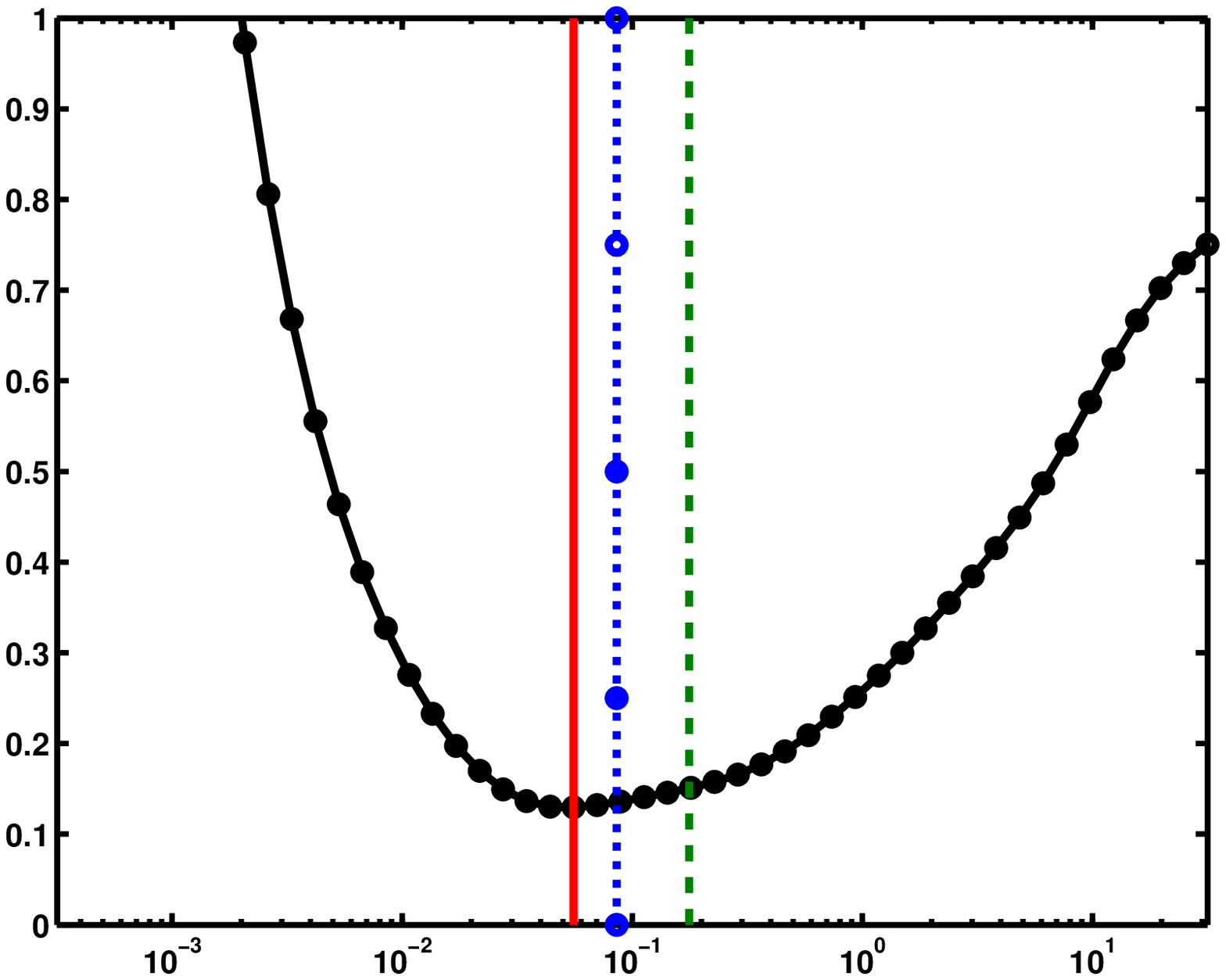}}
\subfigure[$L=L_2$]{\includegraphics[width=1.7in]{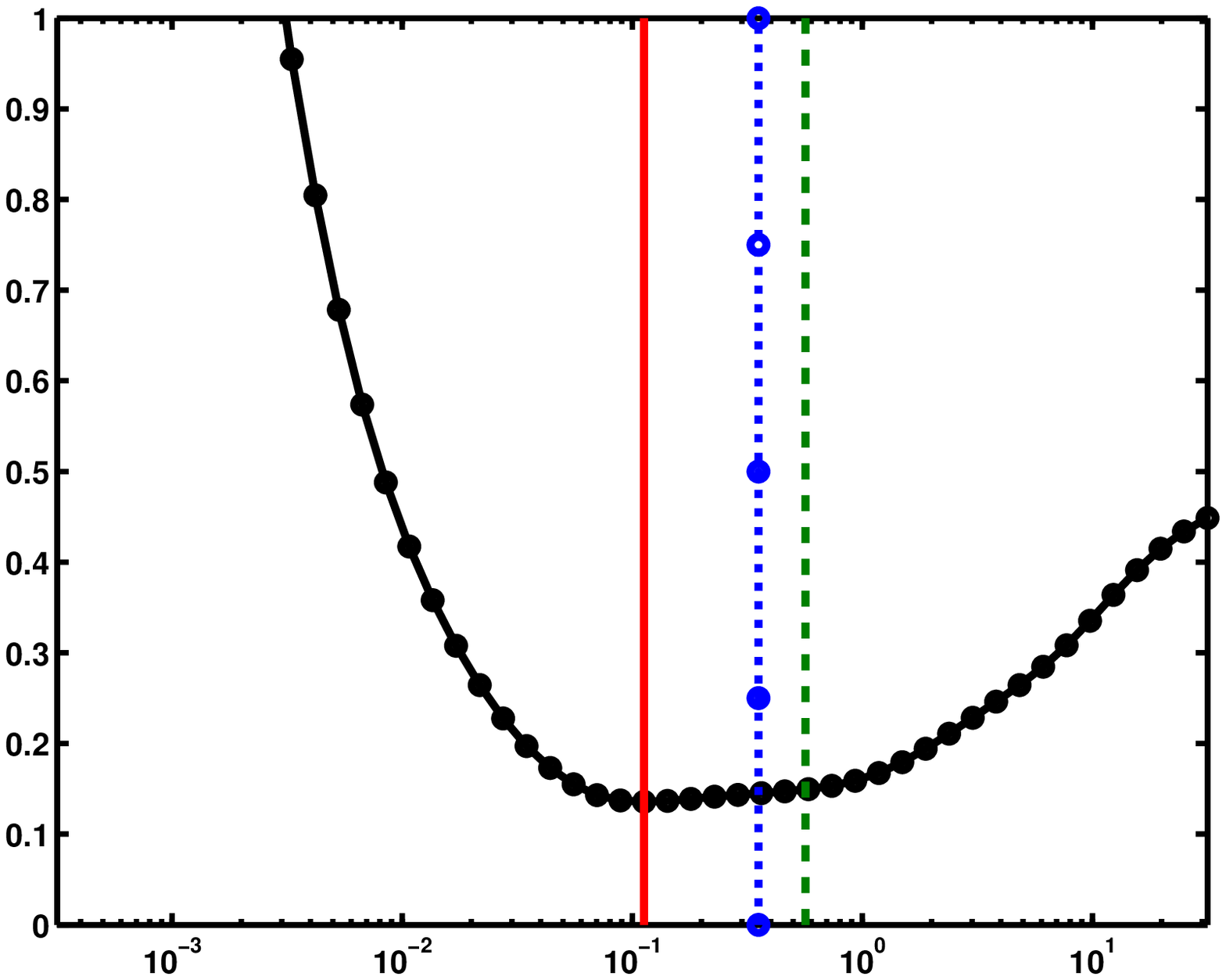}}
\subfigure[$L=I$]{\includegraphics[width=1.7in]{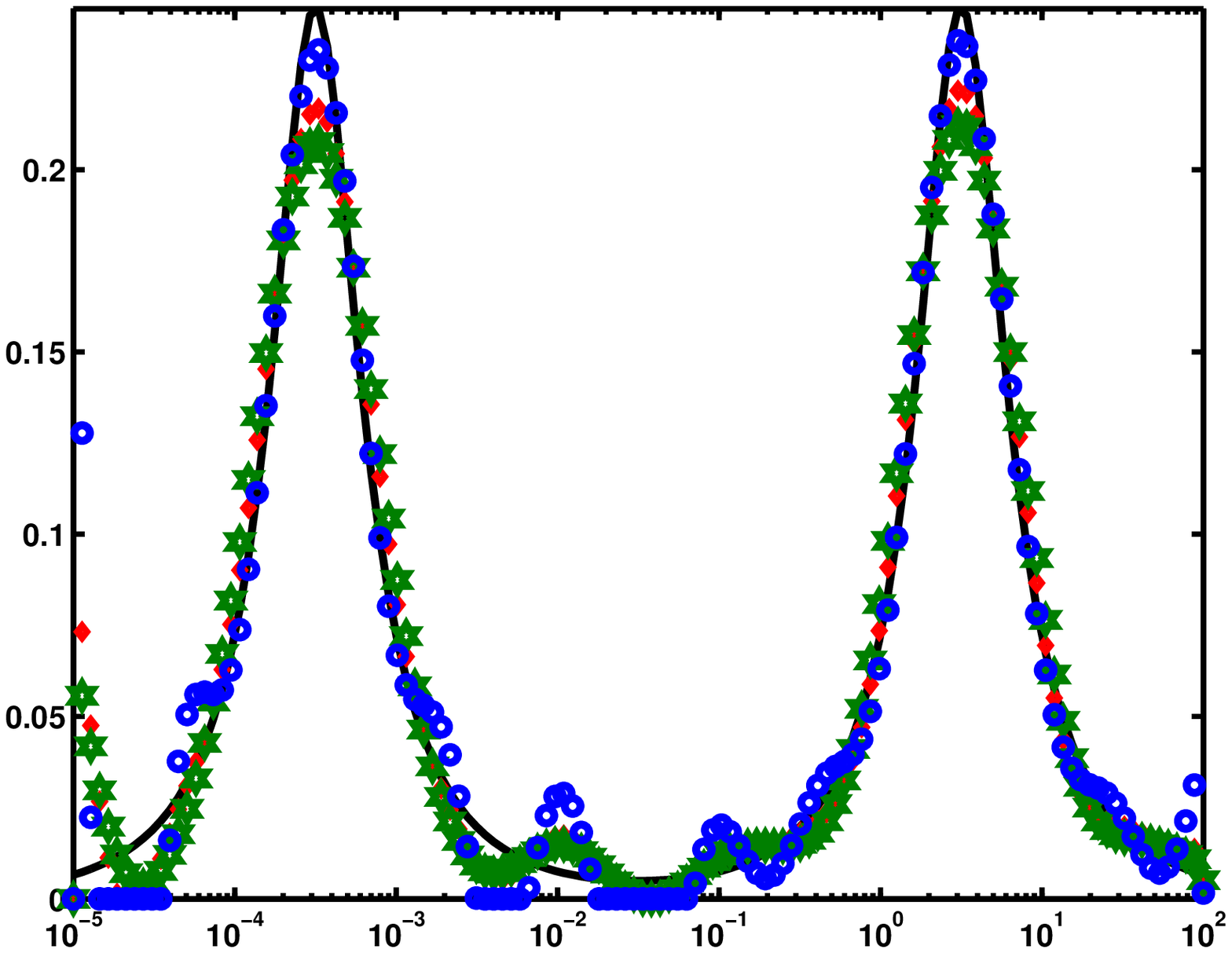}}
\subfigure[$L=L_1$]{\includegraphics[width=1.7in]{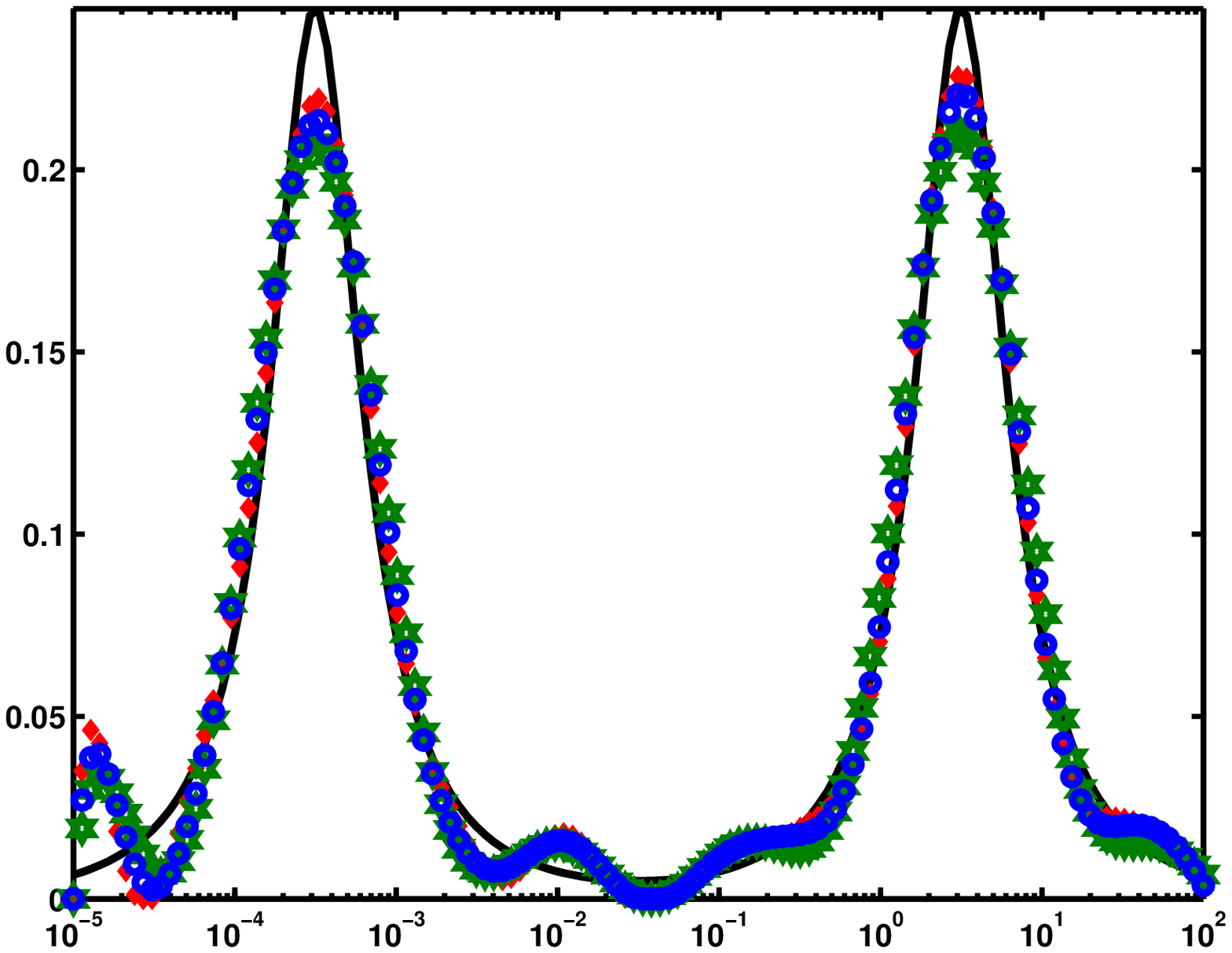}}
\subfigure[$L=L_2$]{\includegraphics[width=1.7in]{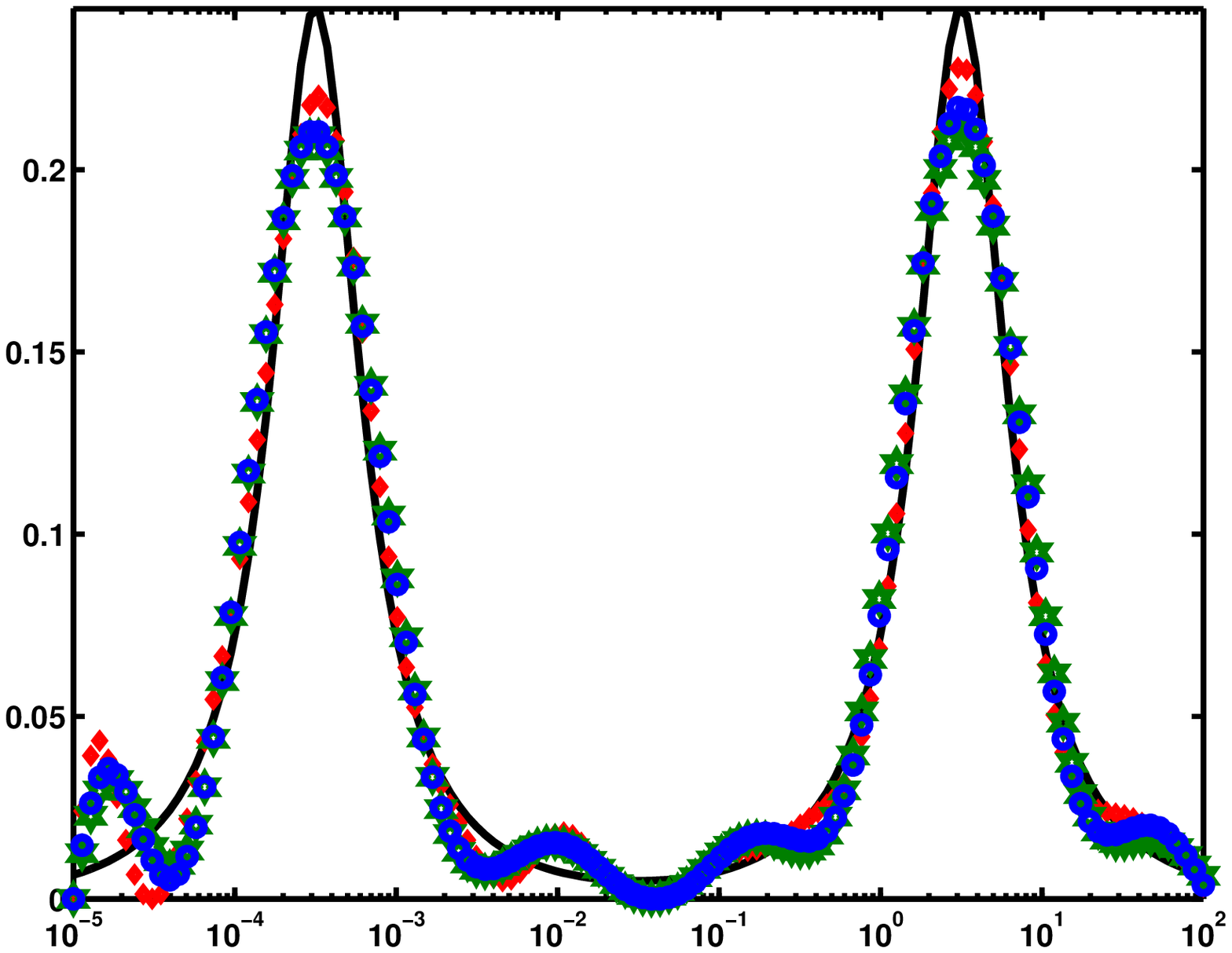}}
\caption{Mean error and example CVX NNLS solutions.  $.1\%$ noise, RQ-A data set, matrix $A_4$.}
\label{fig-lambdachoiceRQ1A4LNCVX}
\end{figure}

\begin{figure}[!ht]
  \centering
\subfigure[$L=I$]{\includegraphics[width=1.7in]{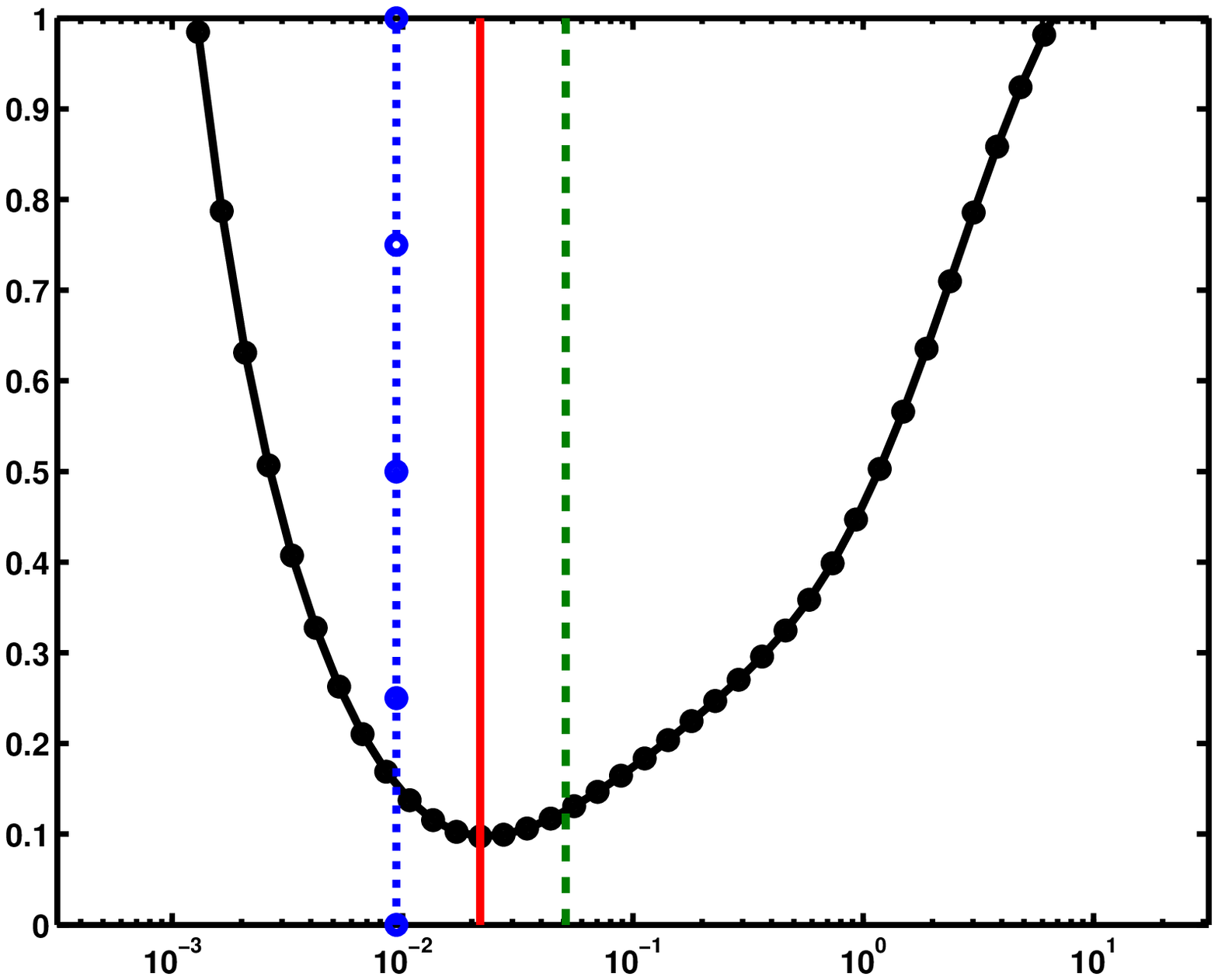}}
\subfigure[$L=L_1$]{\includegraphics[width=1.7in]{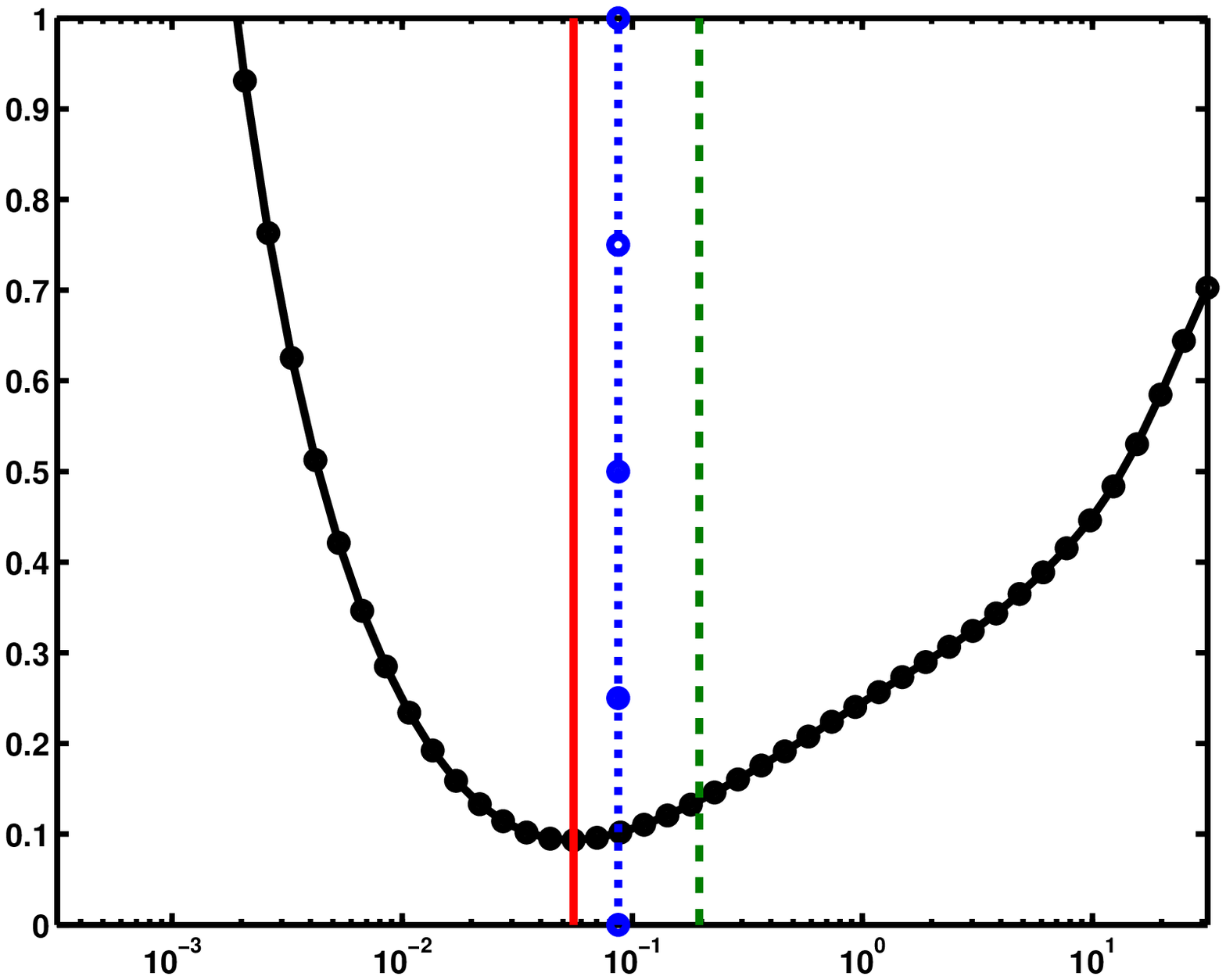}}
\subfigure[$L=L_2$]{\includegraphics[width=1.7in]{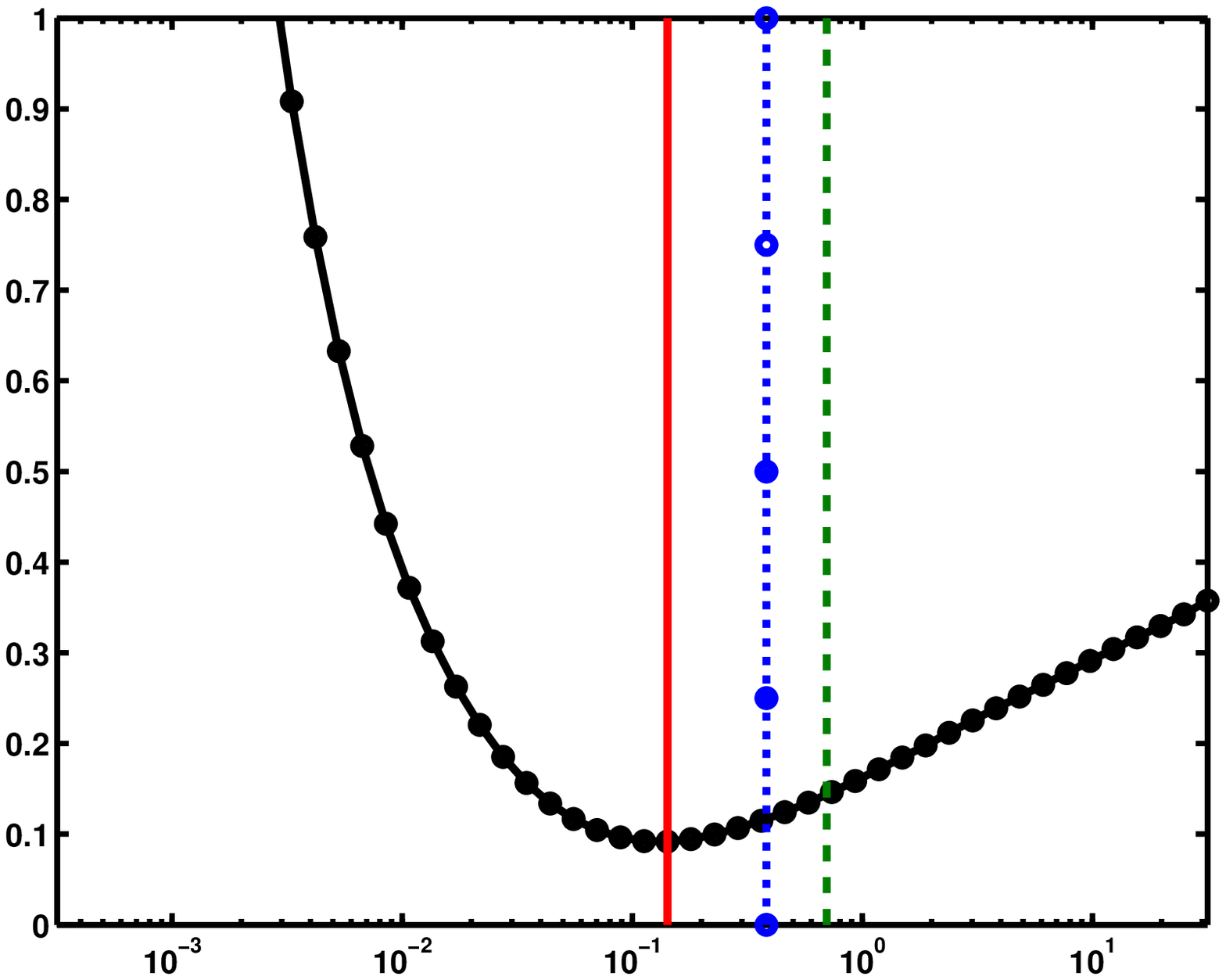}}
\subfigure[$L=I$]{\includegraphics[width=1.7in]{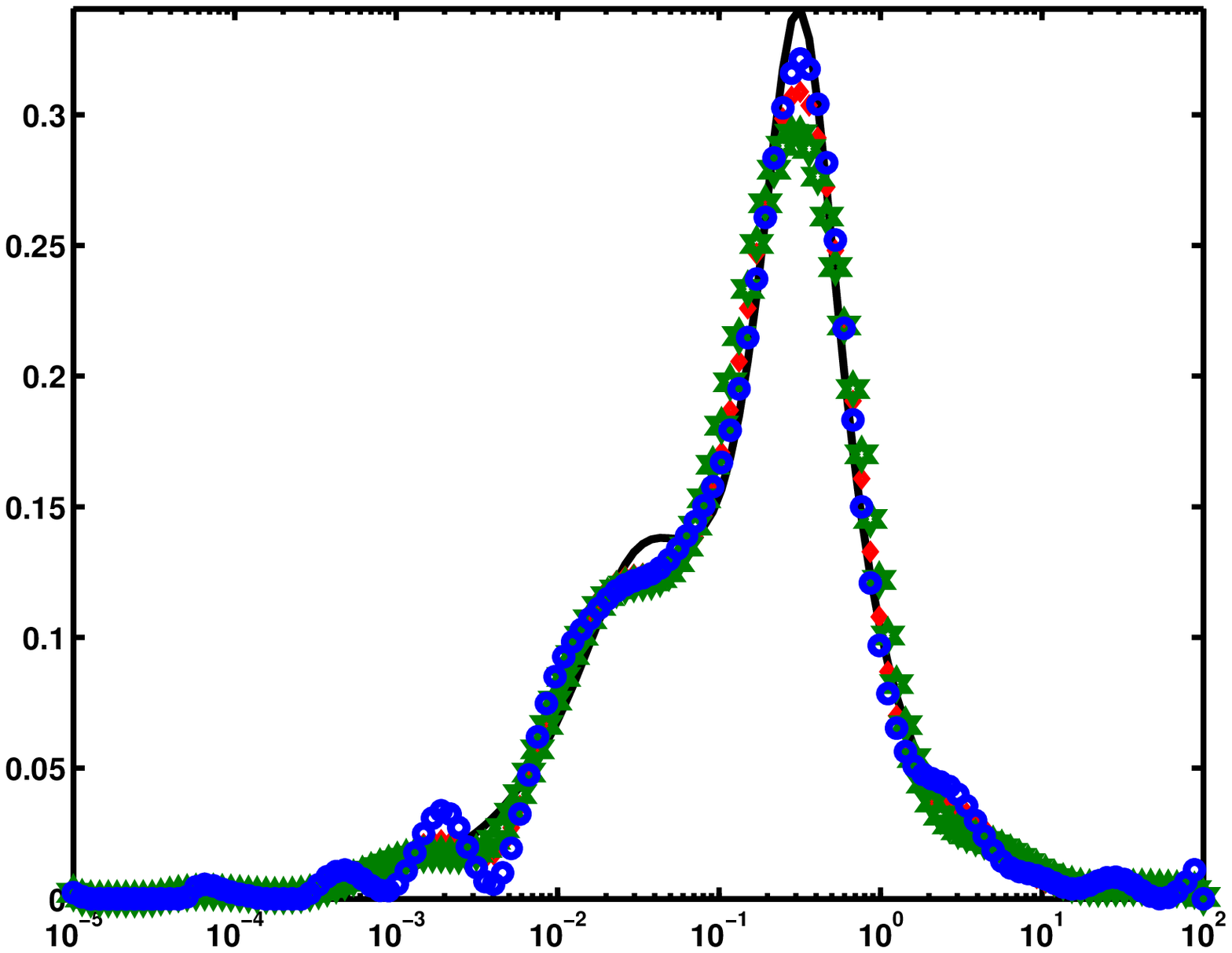}}
\subfigure[$L=L_1$]{\includegraphics[width=1.7in]{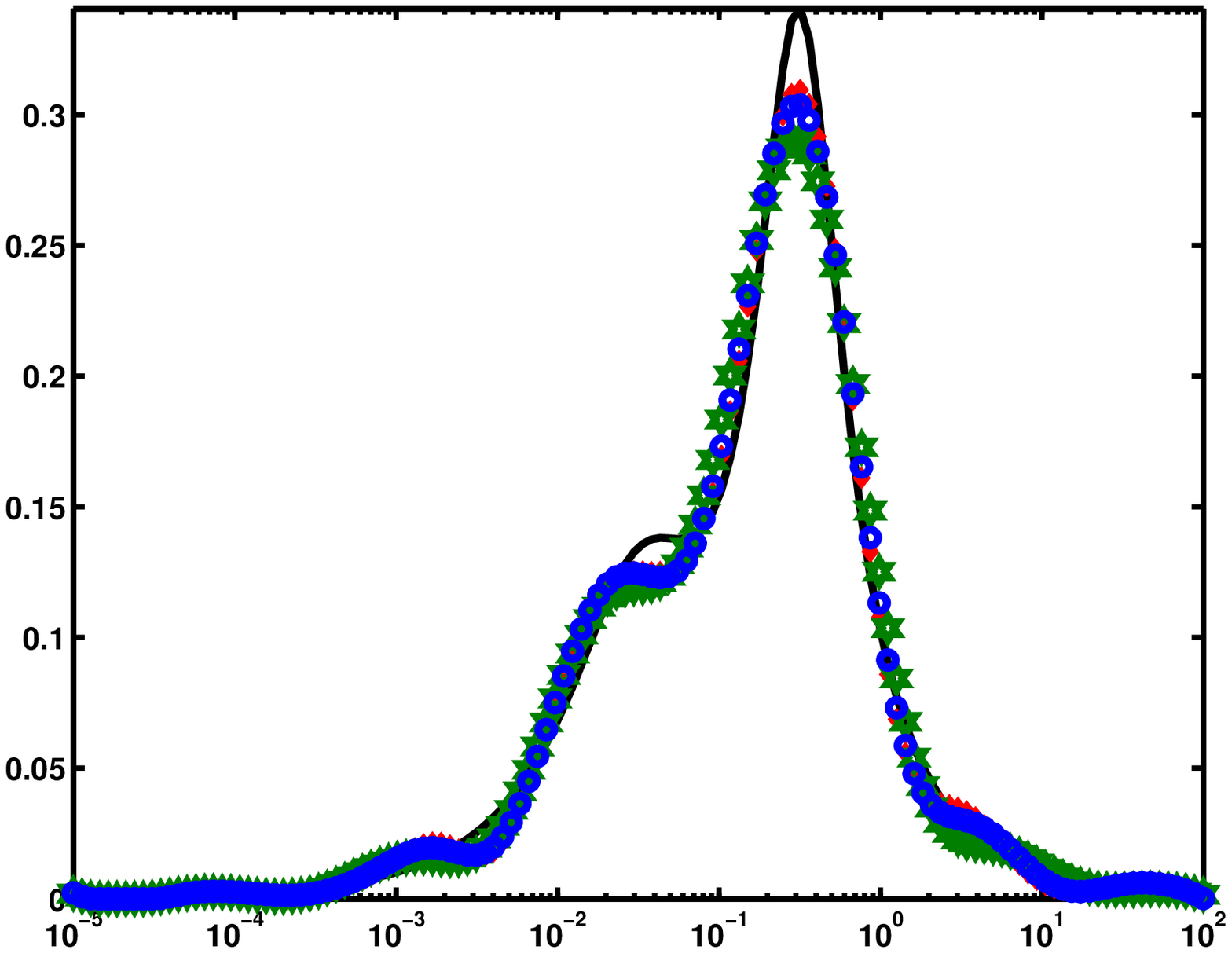}}
\subfigure[$L=L_2$]{\includegraphics[width=1.7in]{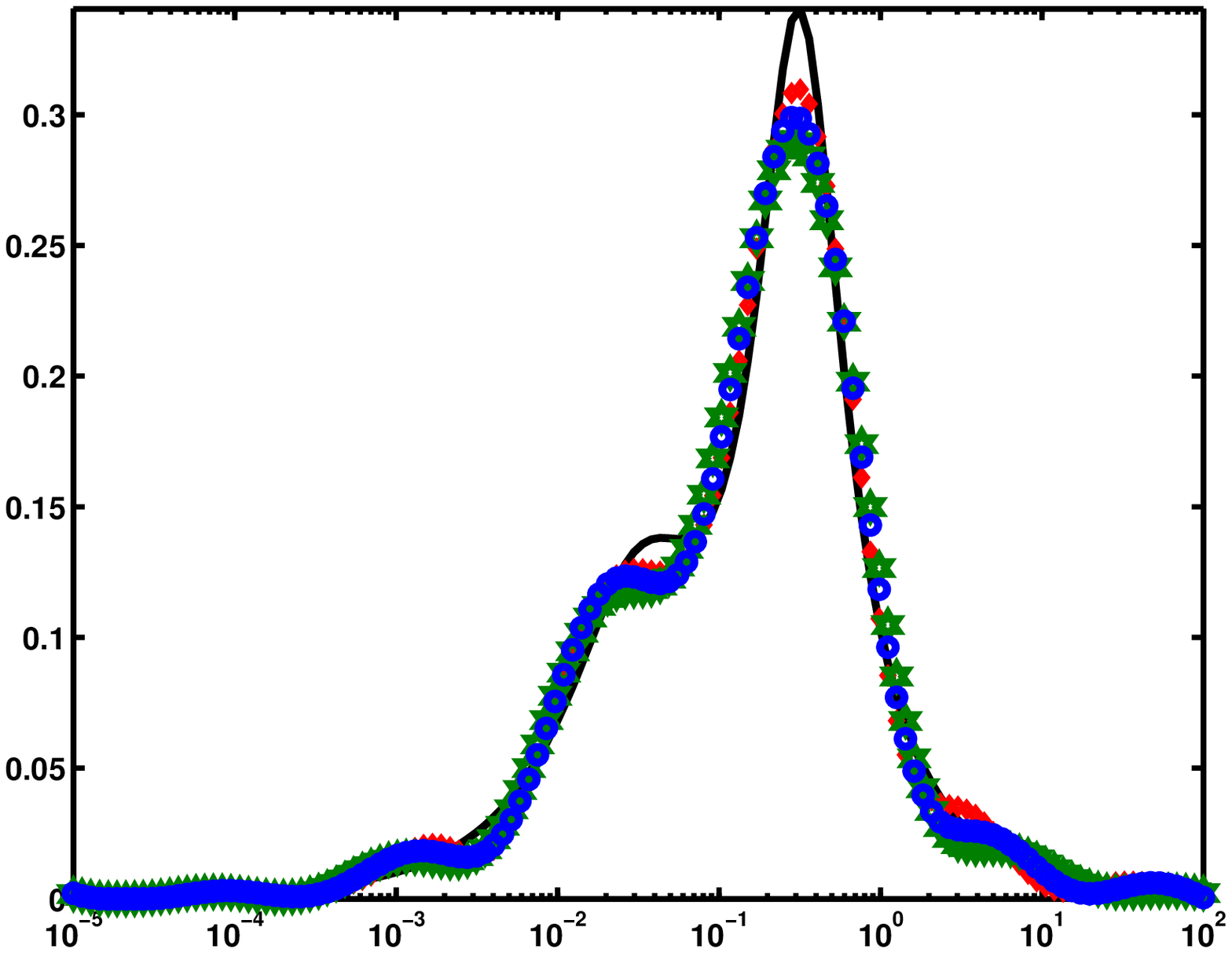}}
\caption{Mean error and example CVX NNLS solutions. $.1\%$ noise, RQ-B data set, matrix $A_4$.}
\label{fig-lambdachoiceRQ5A4LNCVX}
\end{figure}

\begin{figure}[!ht]
  \centering
\subfigure[$L=I$]{\includegraphics[width=1.7in]{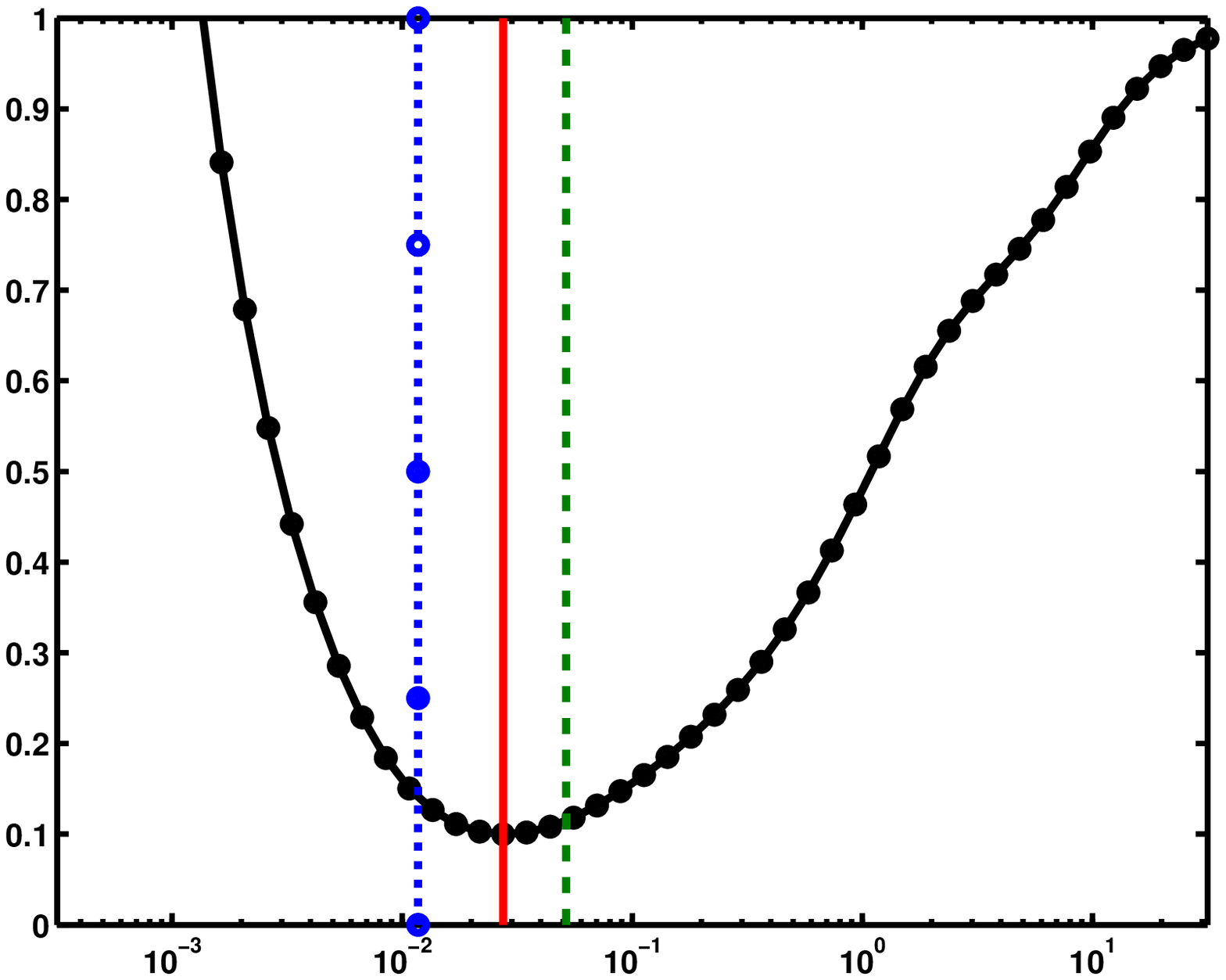}}
\subfigure[$L=L_1$]{\includegraphics[width=1.7in]{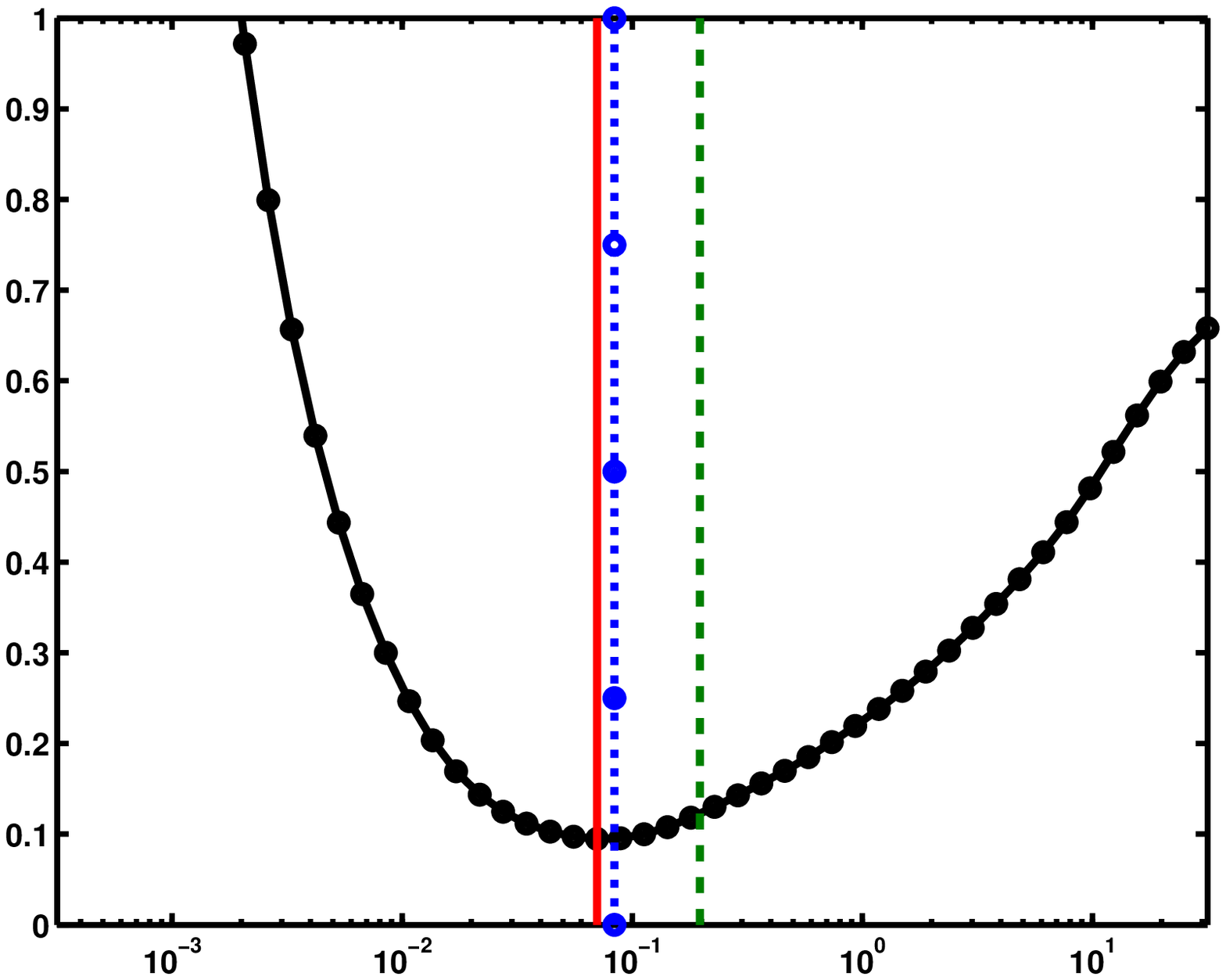}}
\subfigure[$L=L_2$]{\includegraphics[width=1.7in]{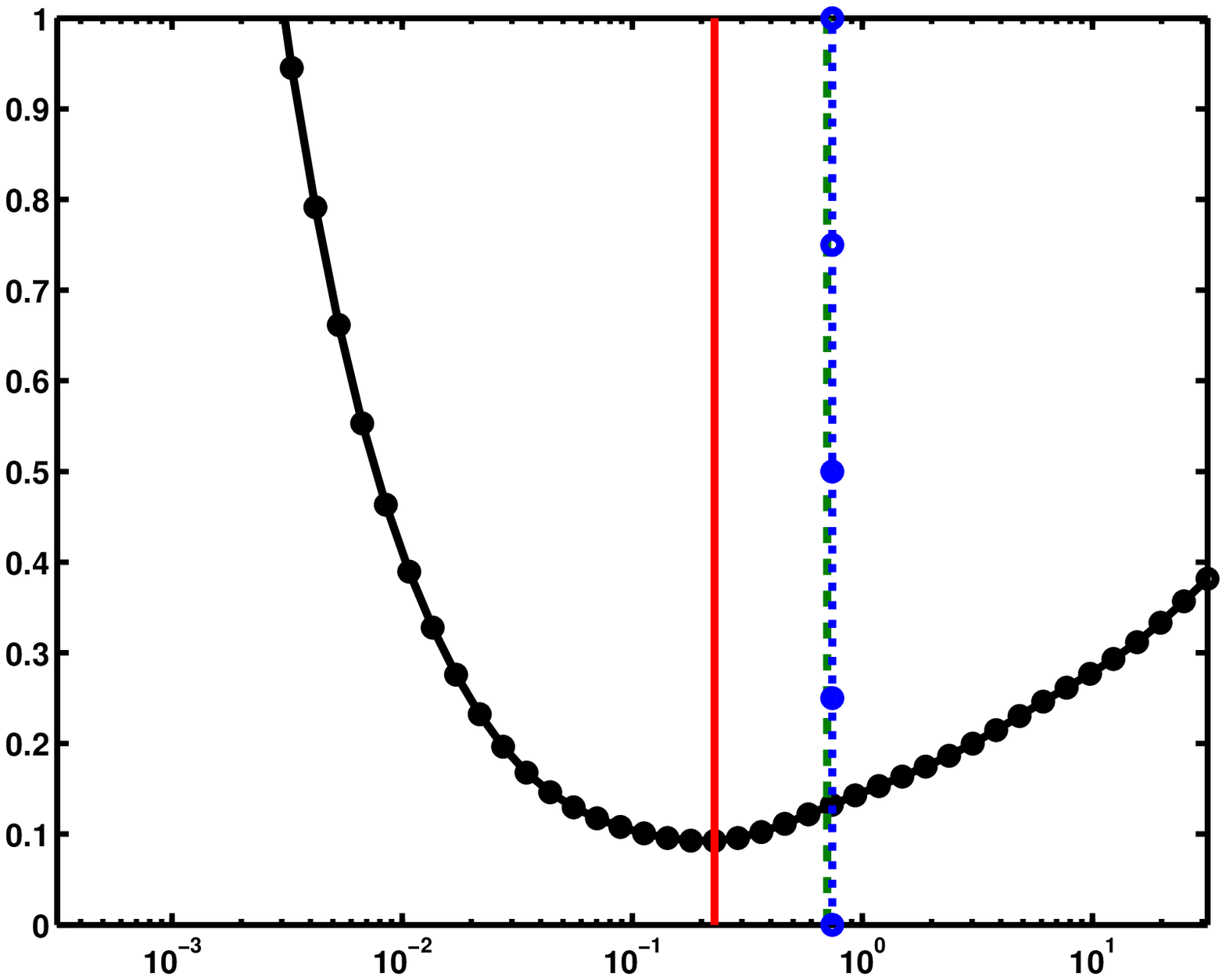}}
\subfigure[$L=I$]{\includegraphics[width=1.7in]{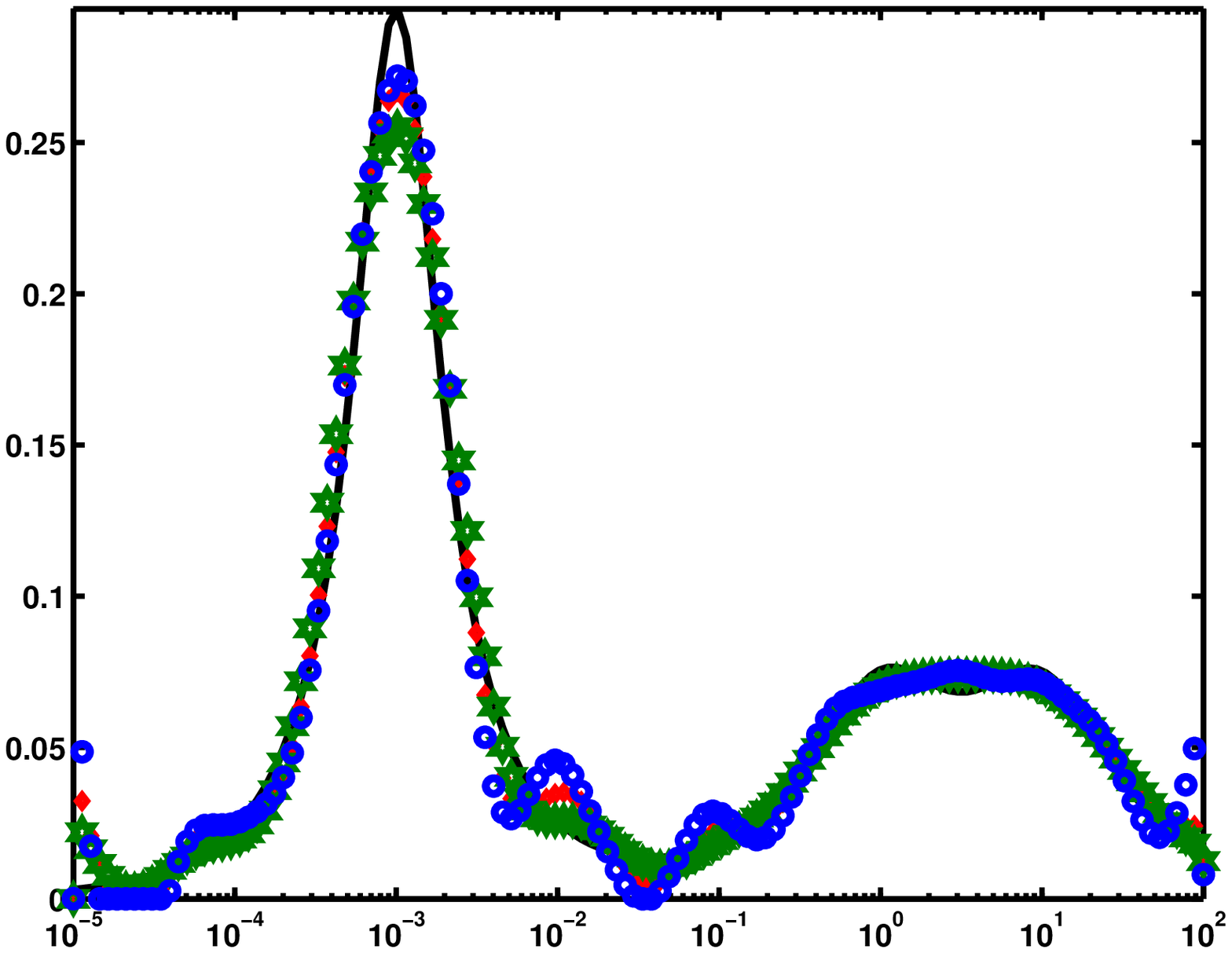}}
\subfigure[$L=L_1$]{\includegraphics[width=1.7in]{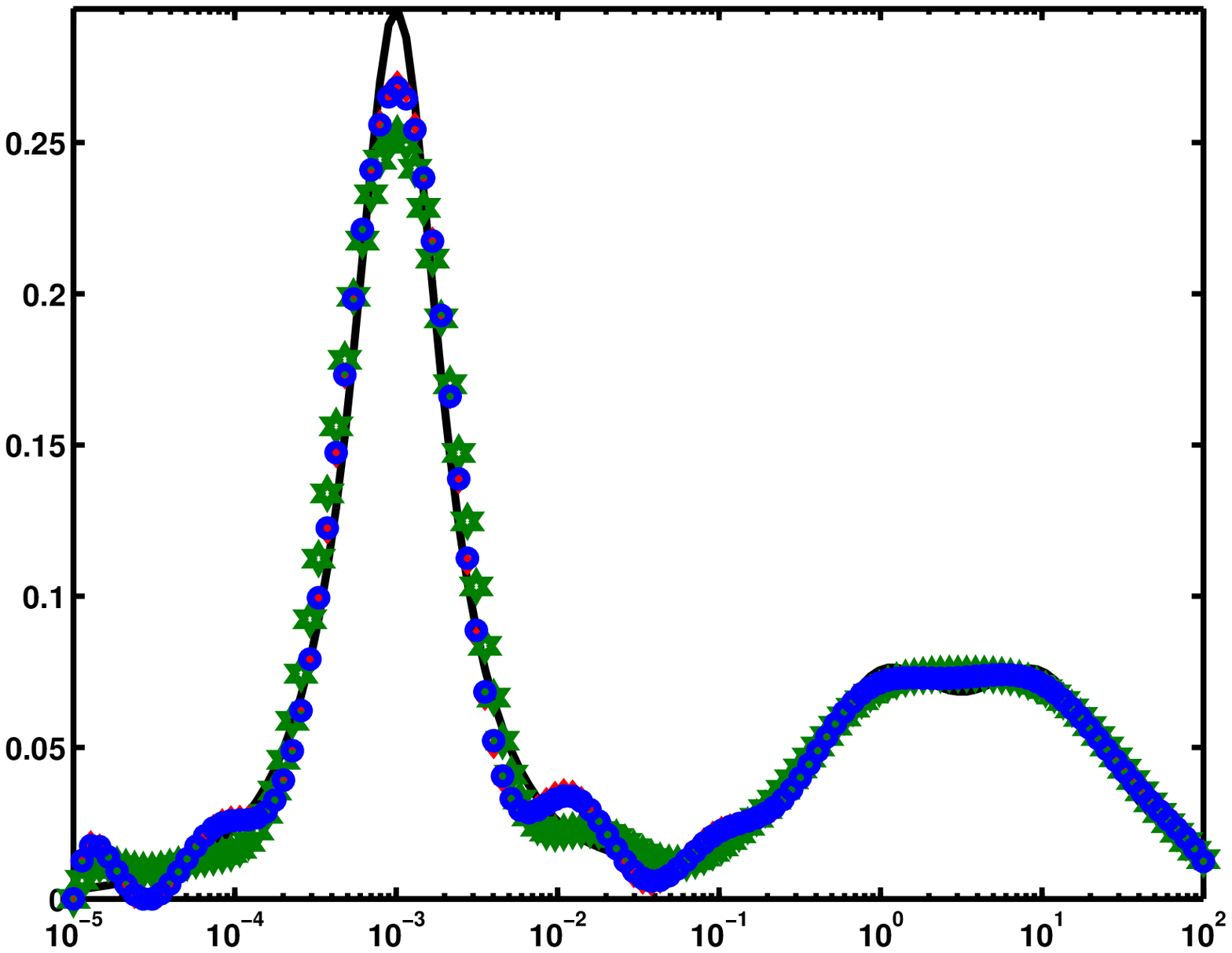}}
\subfigure[$L=L_2$]{\includegraphics[width=1.7in]{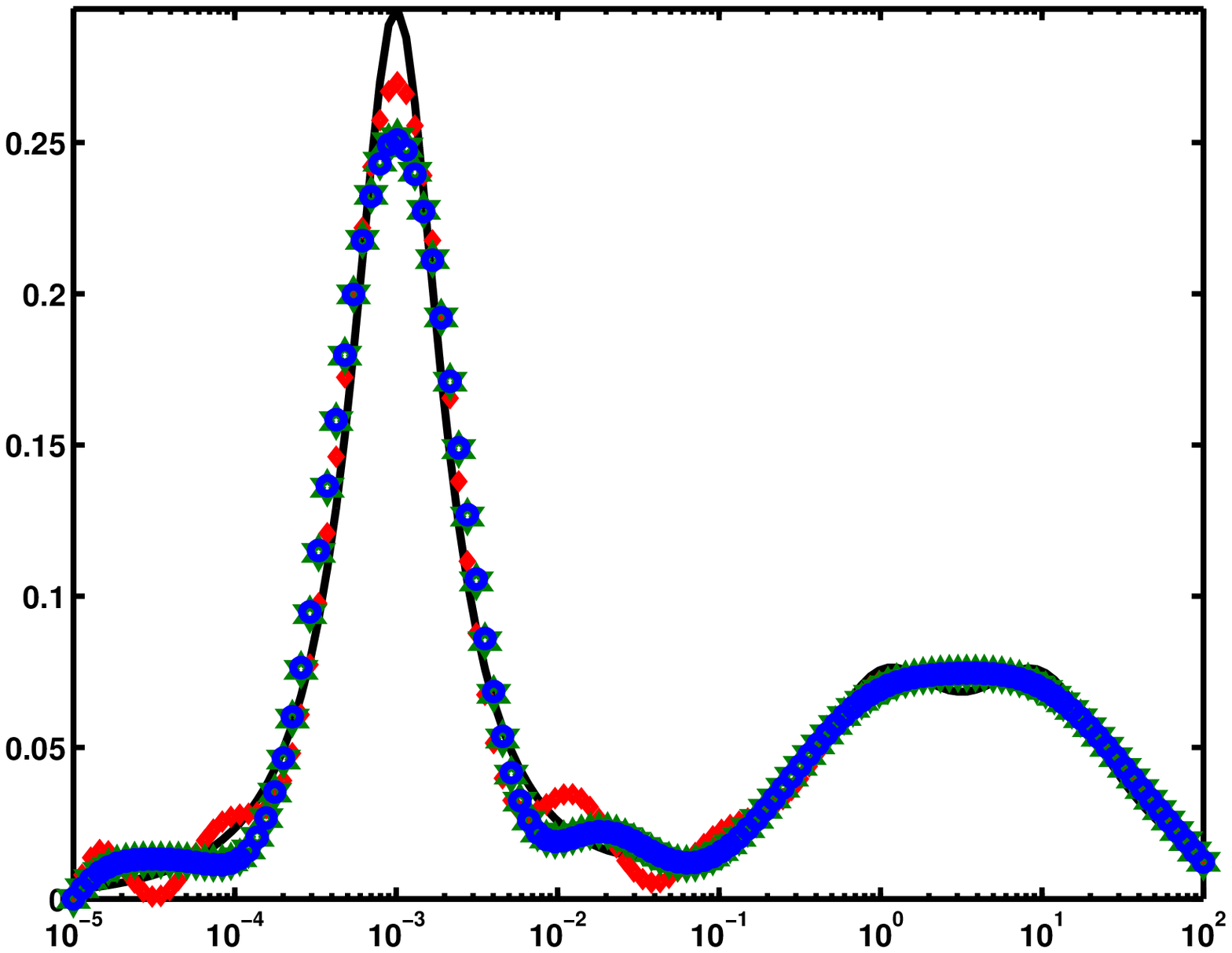}}
\caption{Mean error and example CVX NNLS solutions.  $.1\%$ noise, RQ-C data set, matrix $A_4$.}
\label{fig-lambdachoiceRQ6A4LNCVX}
\end{figure}

\begin{figure}[!ht]
  \centering
\subfigure[$L=I$]{\includegraphics[width=1.7in]{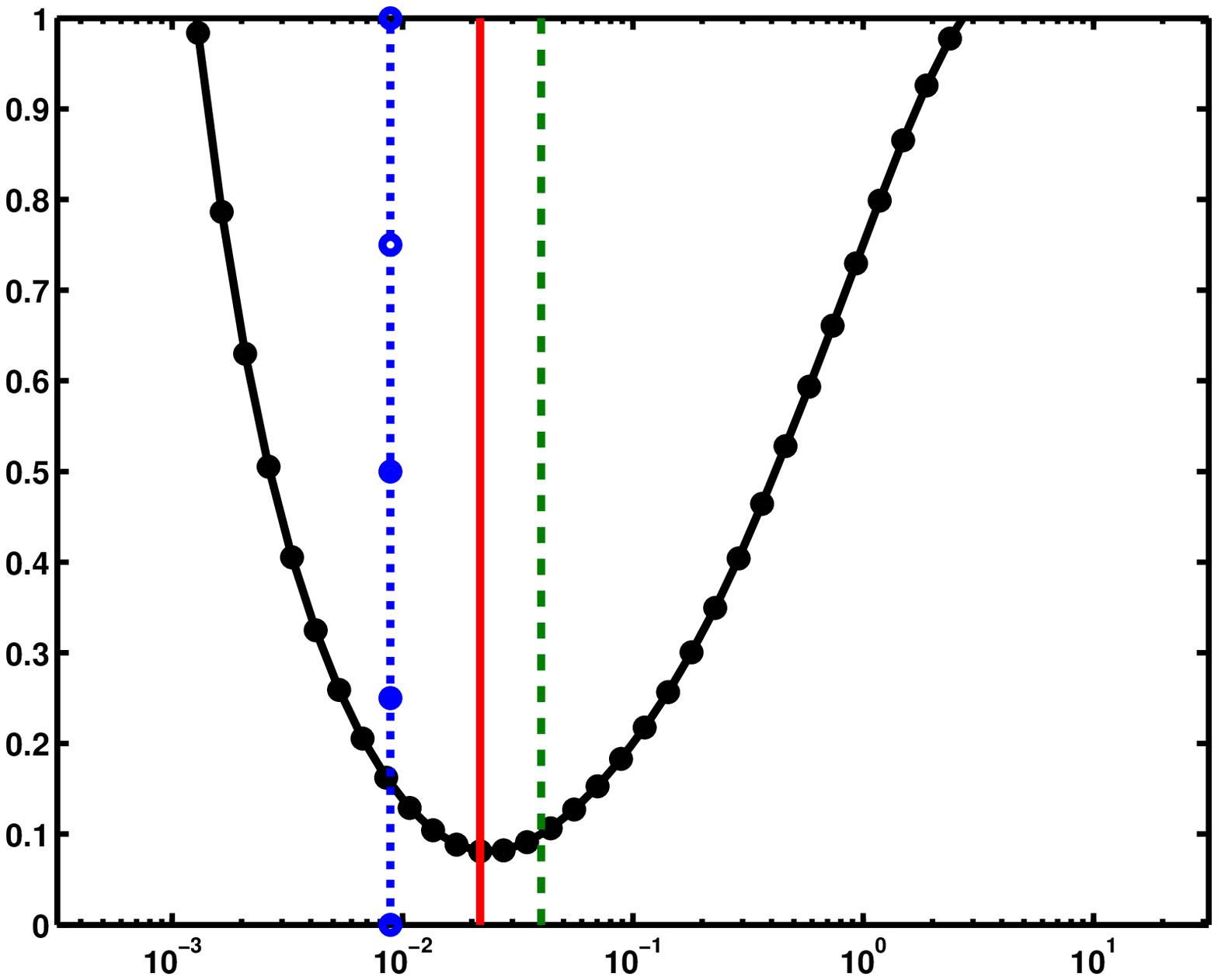}}
\subfigure[$L=L_1$]{\includegraphics[width=1.7in]{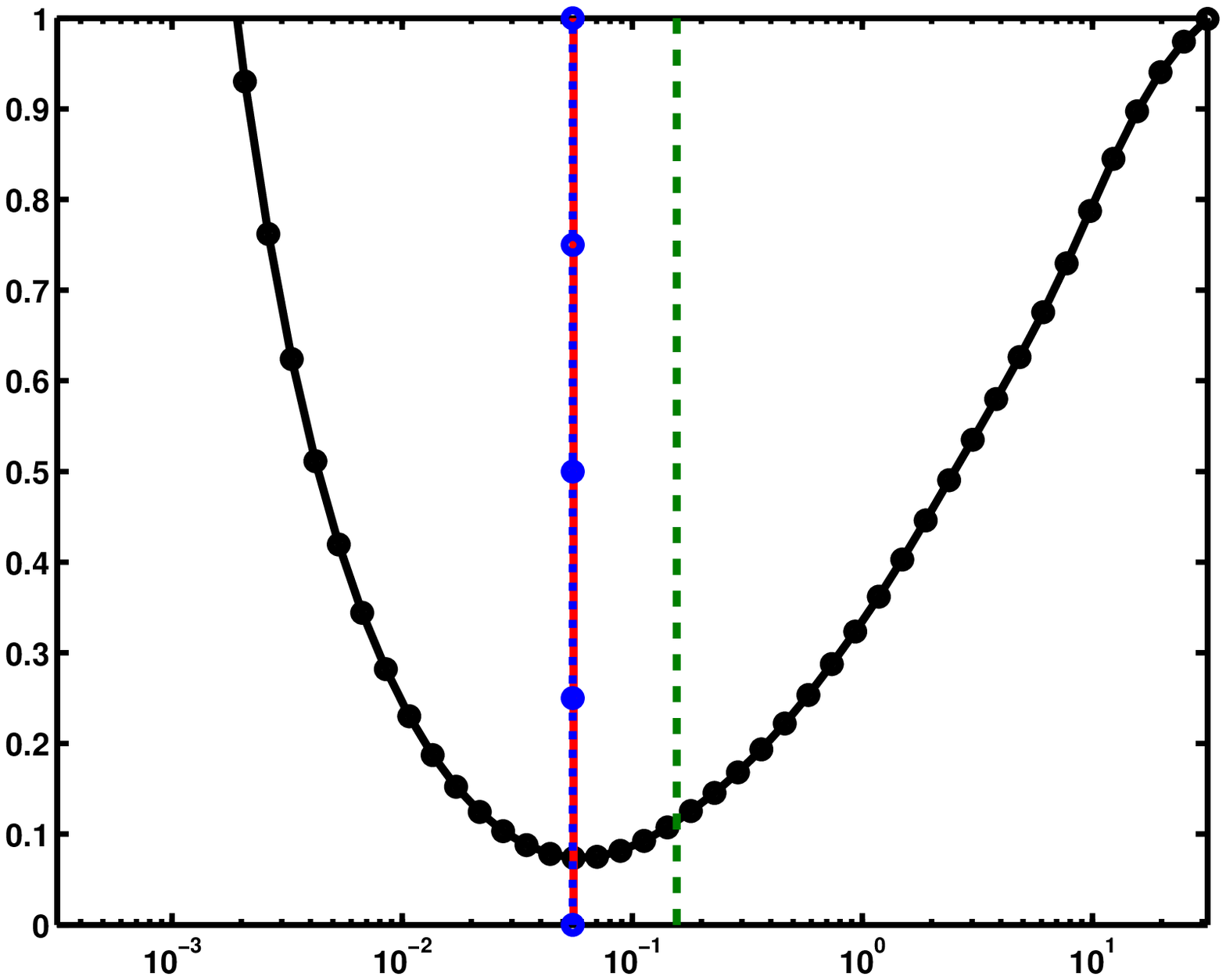}}
\subfigure[$L=L_2$]{\includegraphics[width=1.7in]{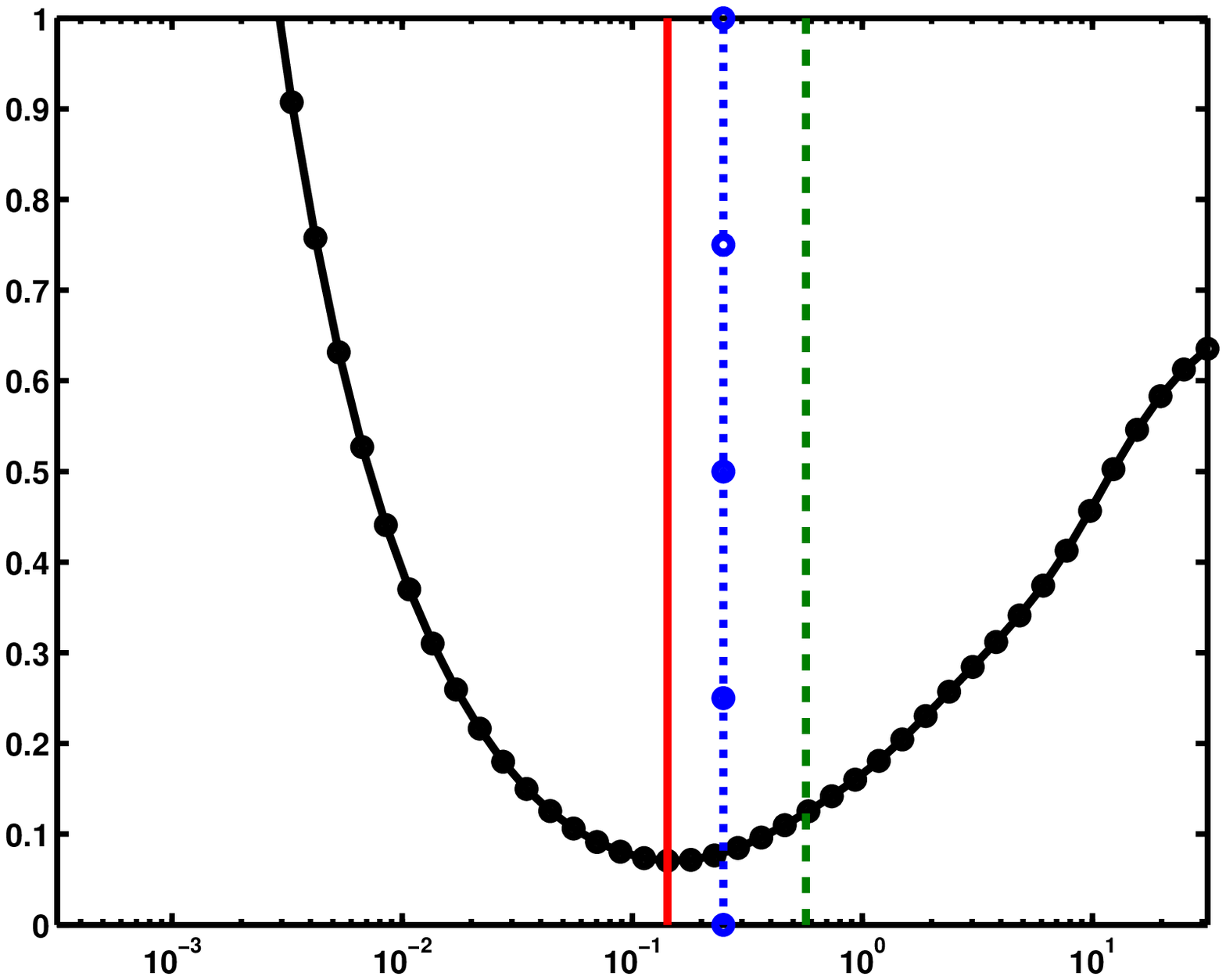}}
\subfigure[$L=I$]{\includegraphics[width=1.7in]{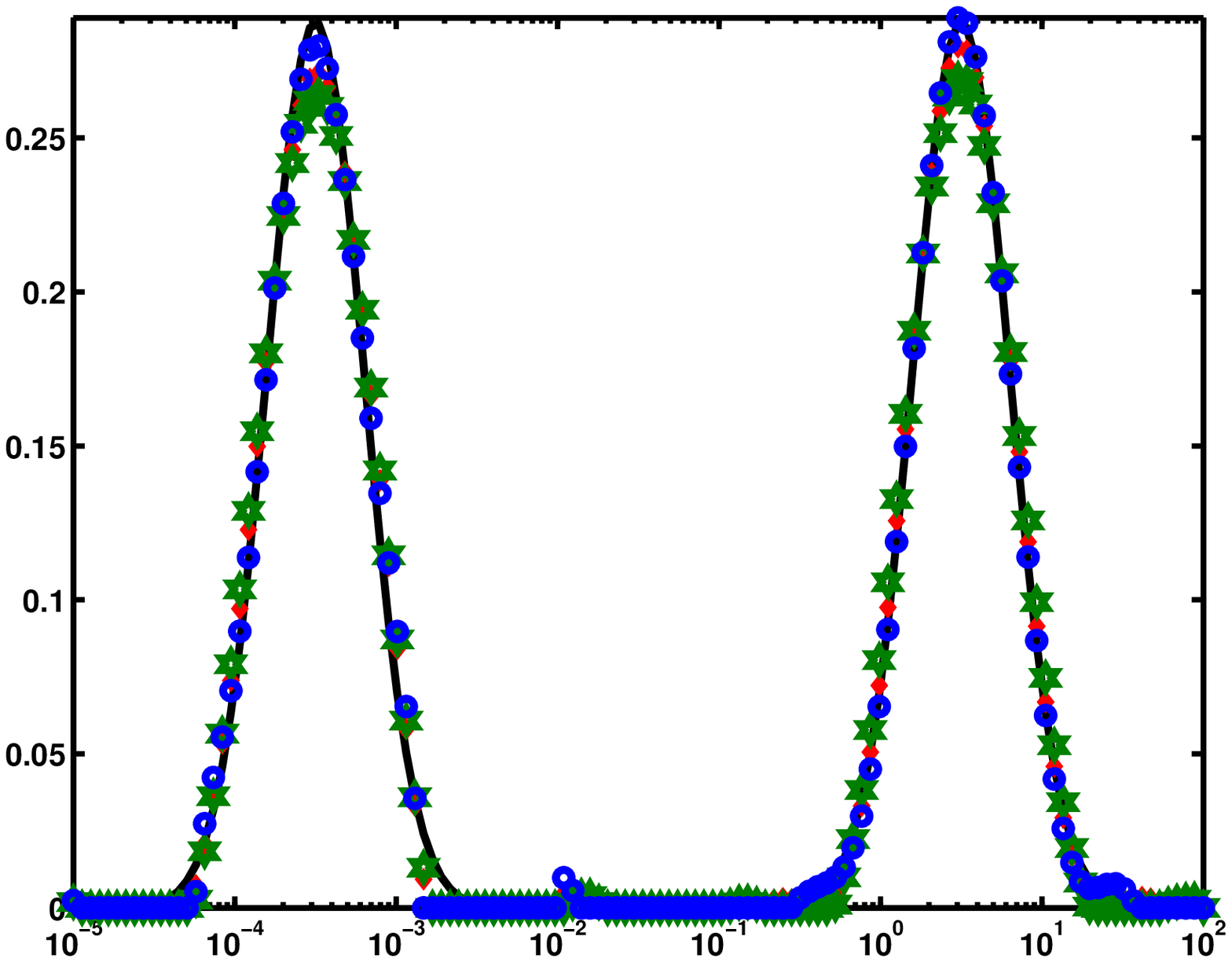}}
\subfigure[$L=L_1$]{\includegraphics[width=1.7in]{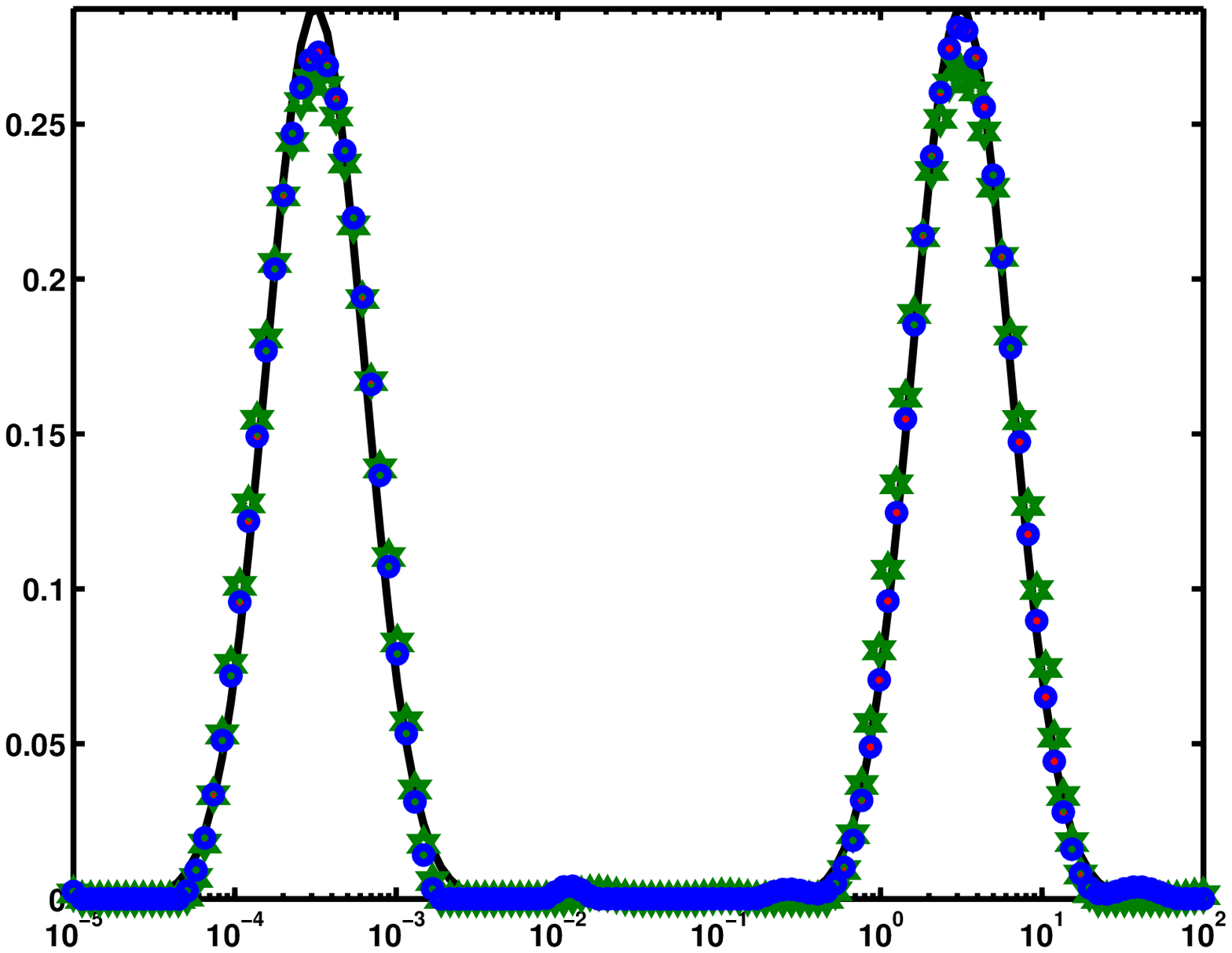}}
\subfigure[$L=L_2$]{\includegraphics[width=1.7in]{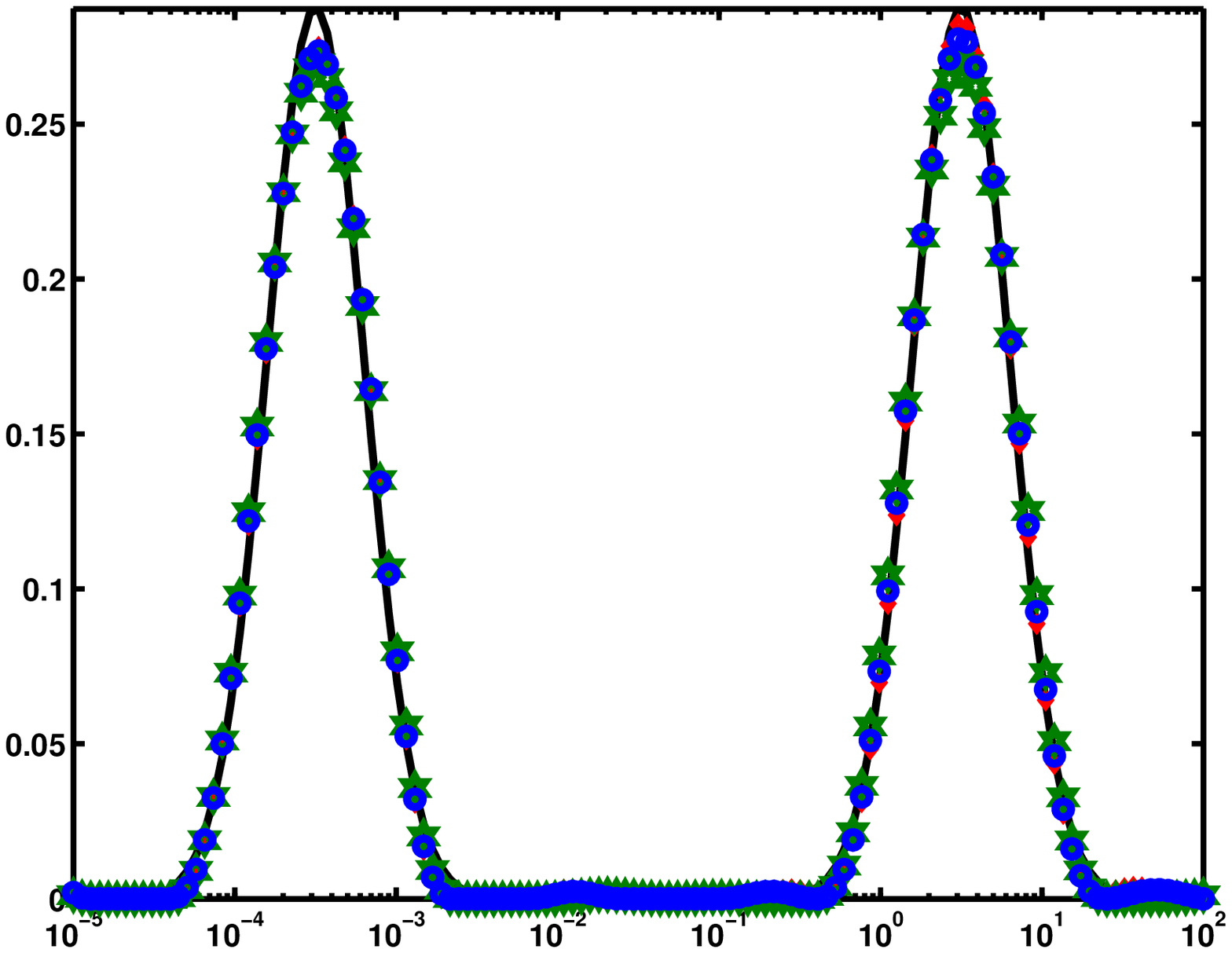}}
\caption{Mean error and example  CVX NNLS solutions.  $.1\%$ noise, LN-A data set, matrix $A_4$.}
\label{fig-lambdachoiceLN2A4LNCVX}
\end{figure}

\begin{figure}[!ht]
  \centering
\subfigure[$L=I$]{\includegraphics[width=1.7in]{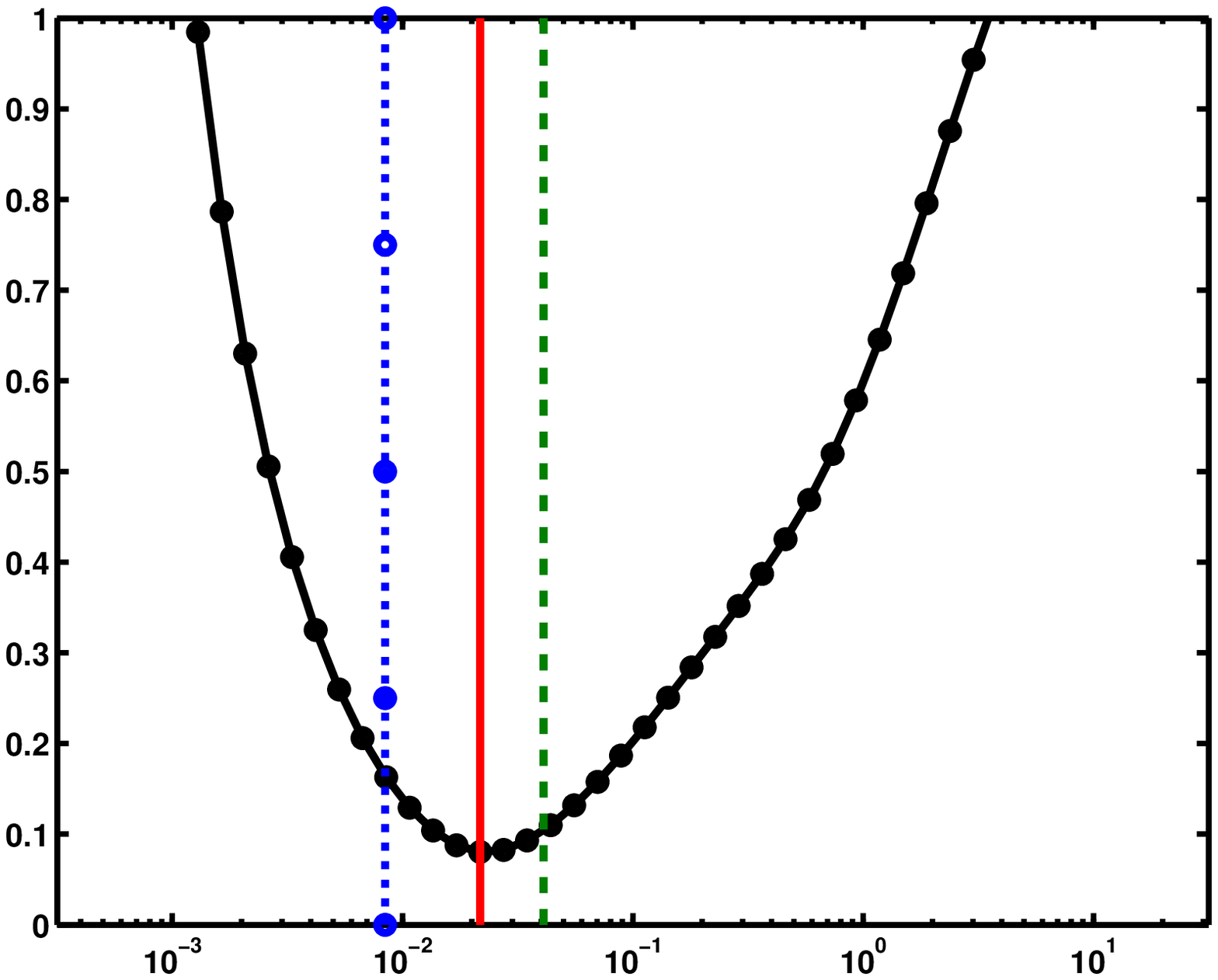}}
\subfigure[$L=L_1$]{\includegraphics[width=1.7in]{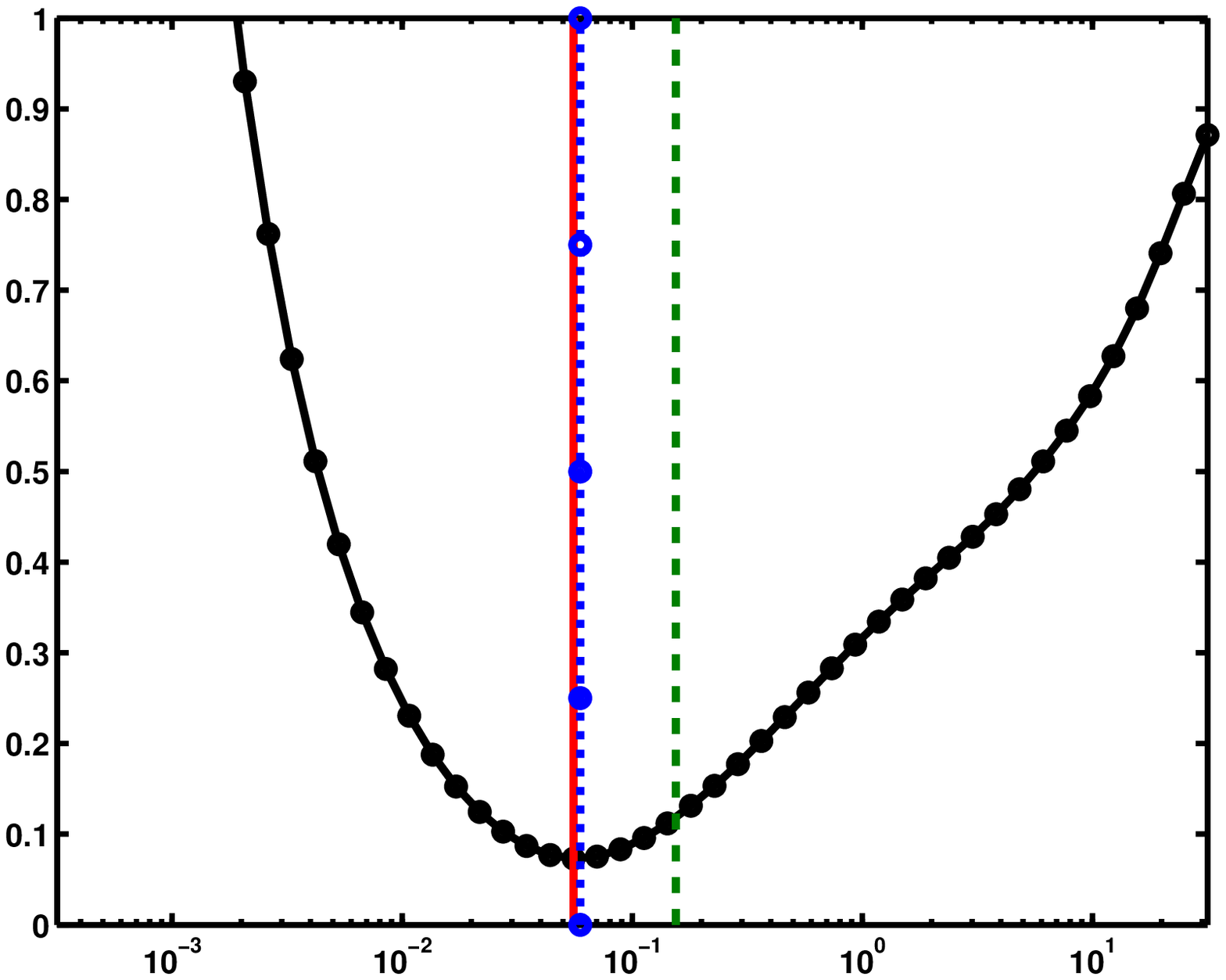}}
\subfigure[$L=L_2$]{\includegraphics[width=1.7in]{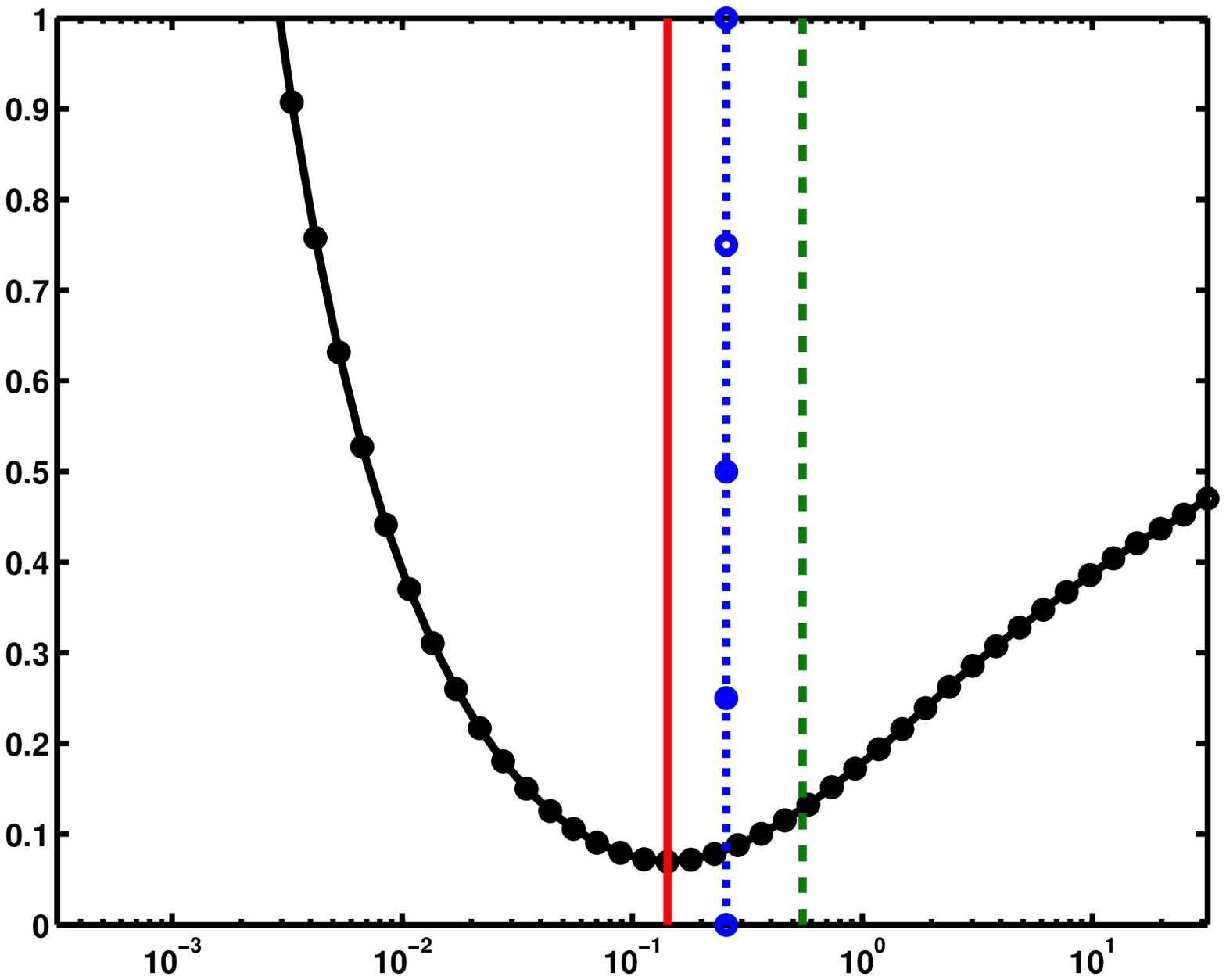}}
\subfigure[$L=I$]{\includegraphics[width=1.7in]{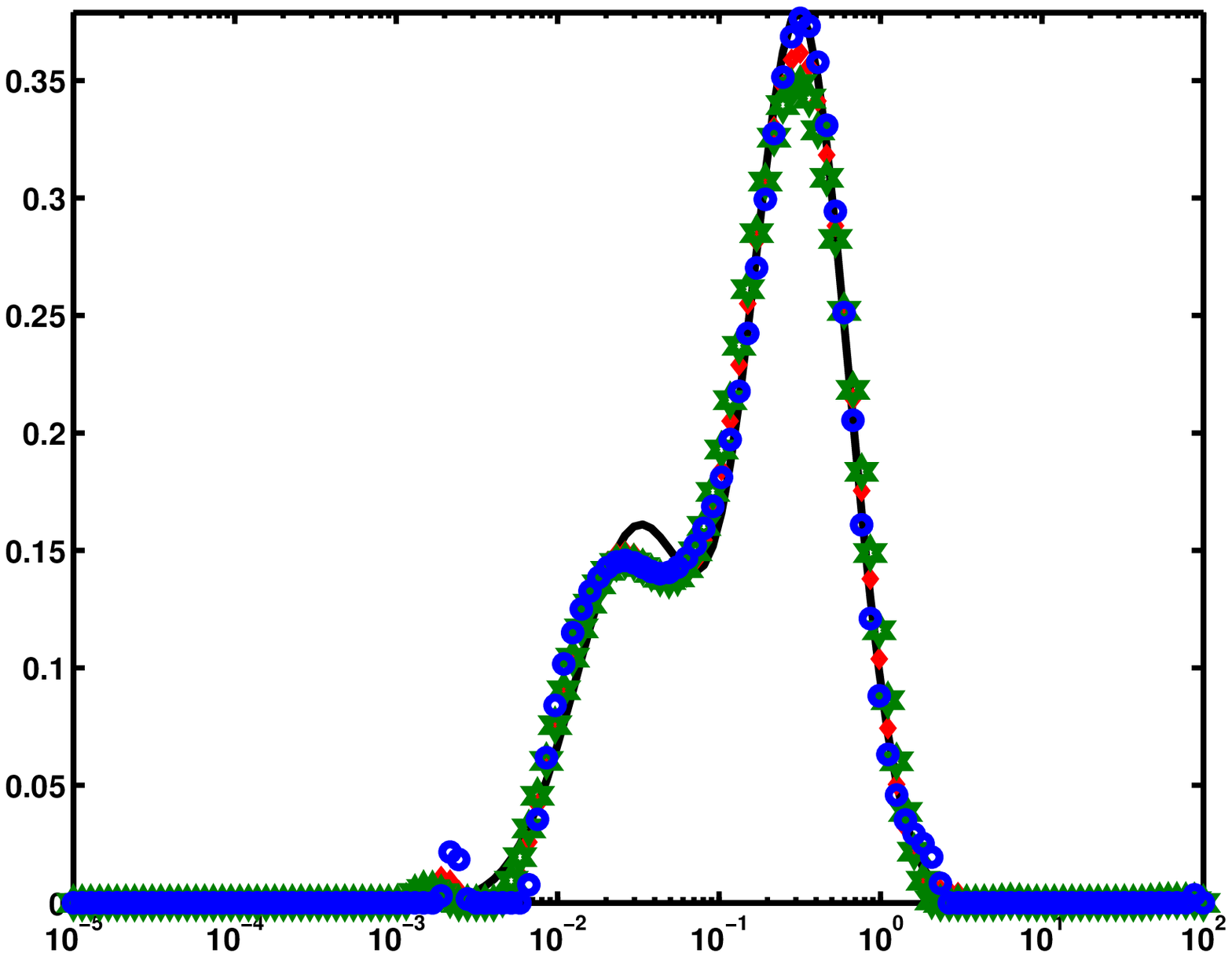}}
\subfigure[$L=L_1$]{\includegraphics[width=1.7in]{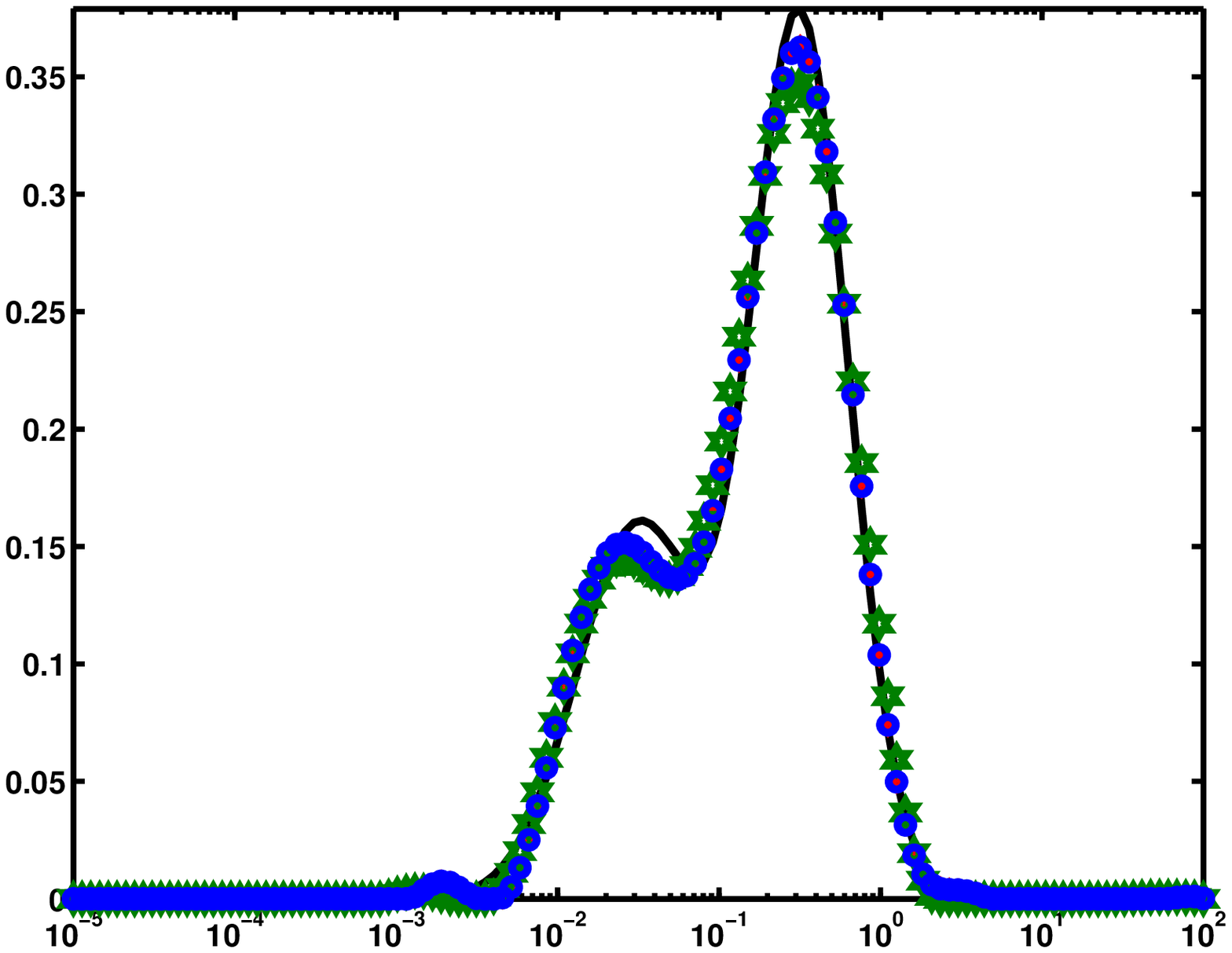}}
\subfigure[$L=L_2$]{\includegraphics[width=1.7in]{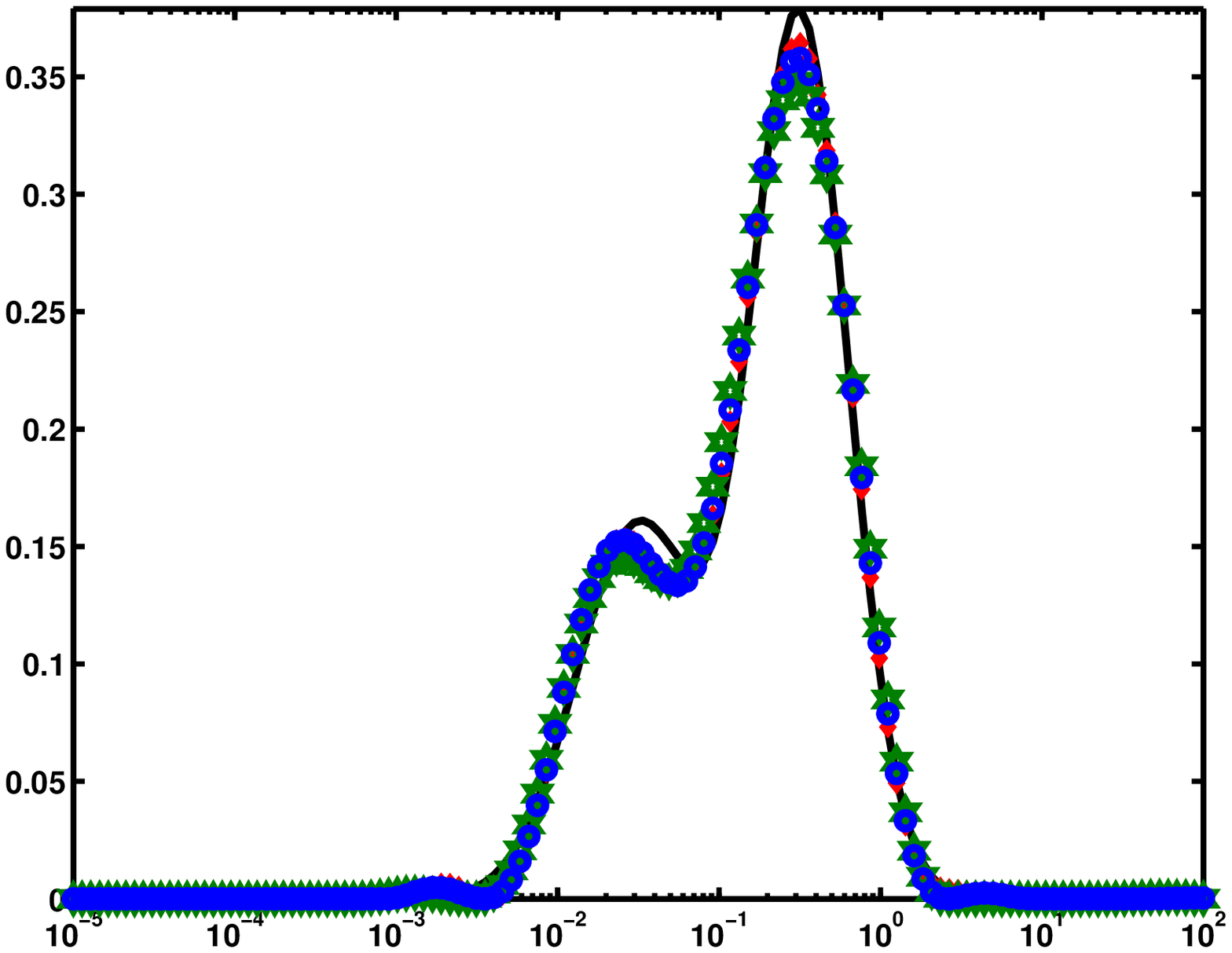}}
\caption{Mean error and example  CVX NNLS solutions.  $.1\%$ noise, LN-B data set, matrix $A_4$.}
\label{fig-lambdachoiceLN5A4LNCVX}
\end{figure}

\begin{figure}[!ht]
  \centering
\subfigure[$L=I$]{\includegraphics[width=1.7in]{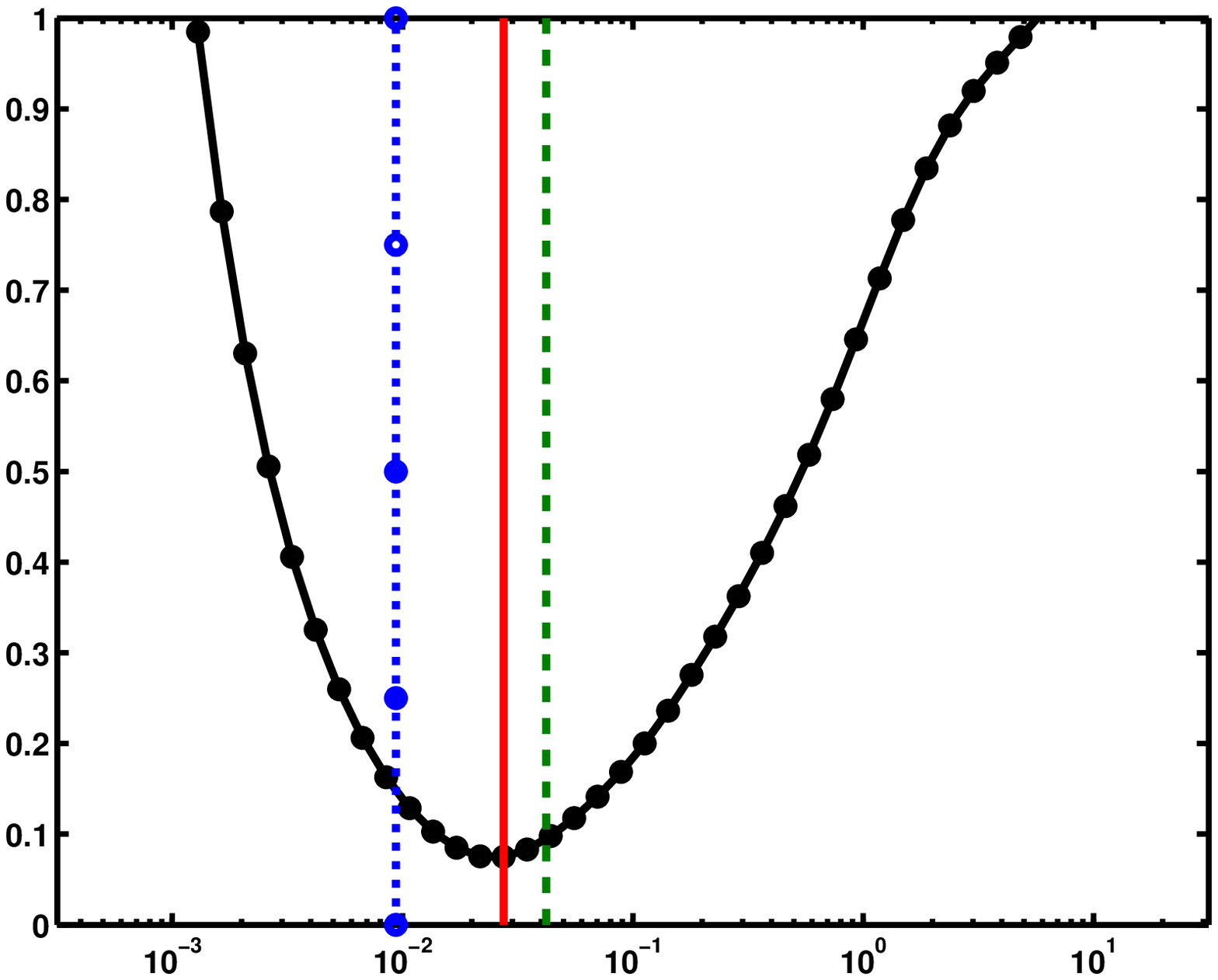}}
\subfigure[$L=L_1$]{\includegraphics[width=1.7in]{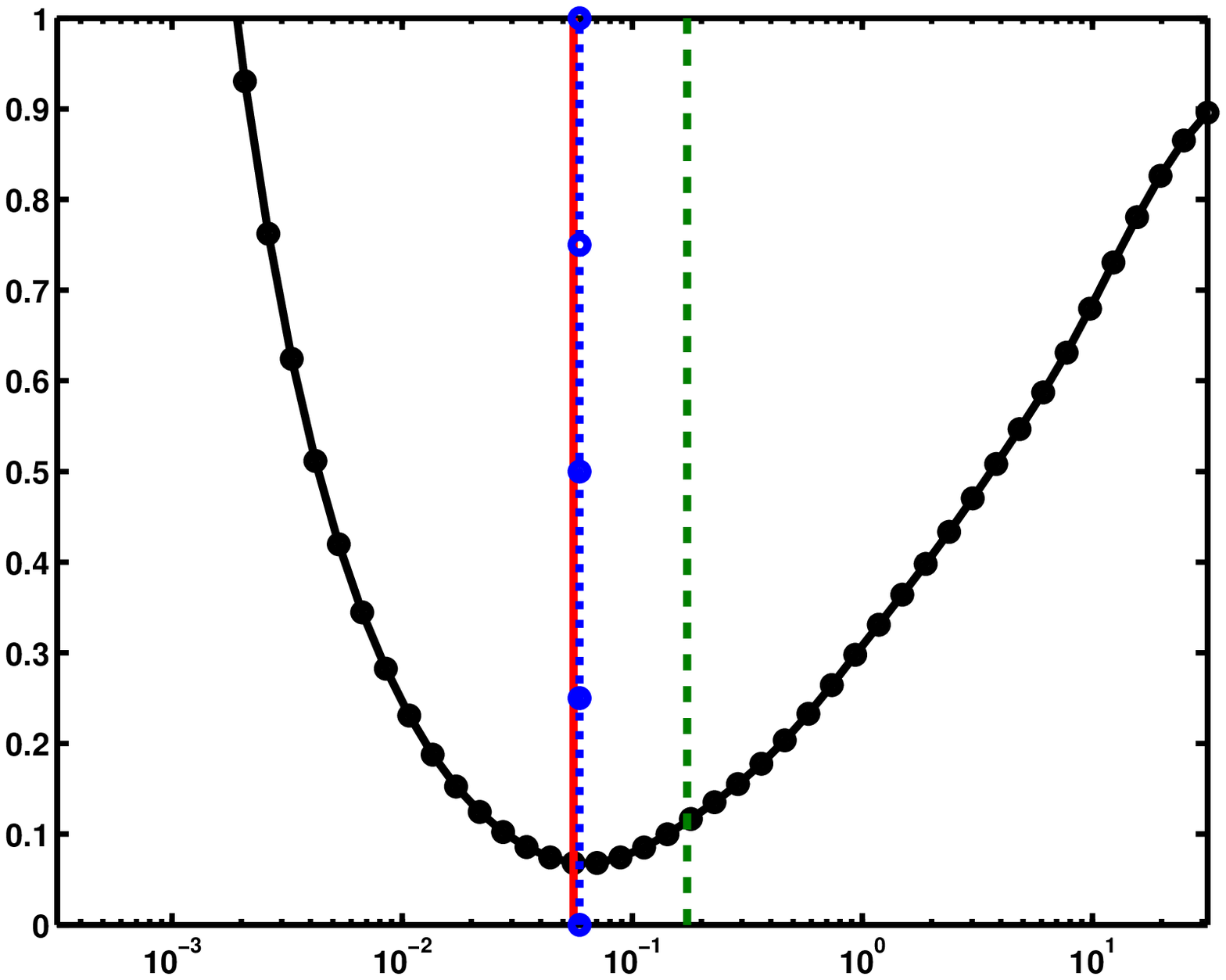}}
\subfigure[$L=L_2$]{\includegraphics[width=1.7in]{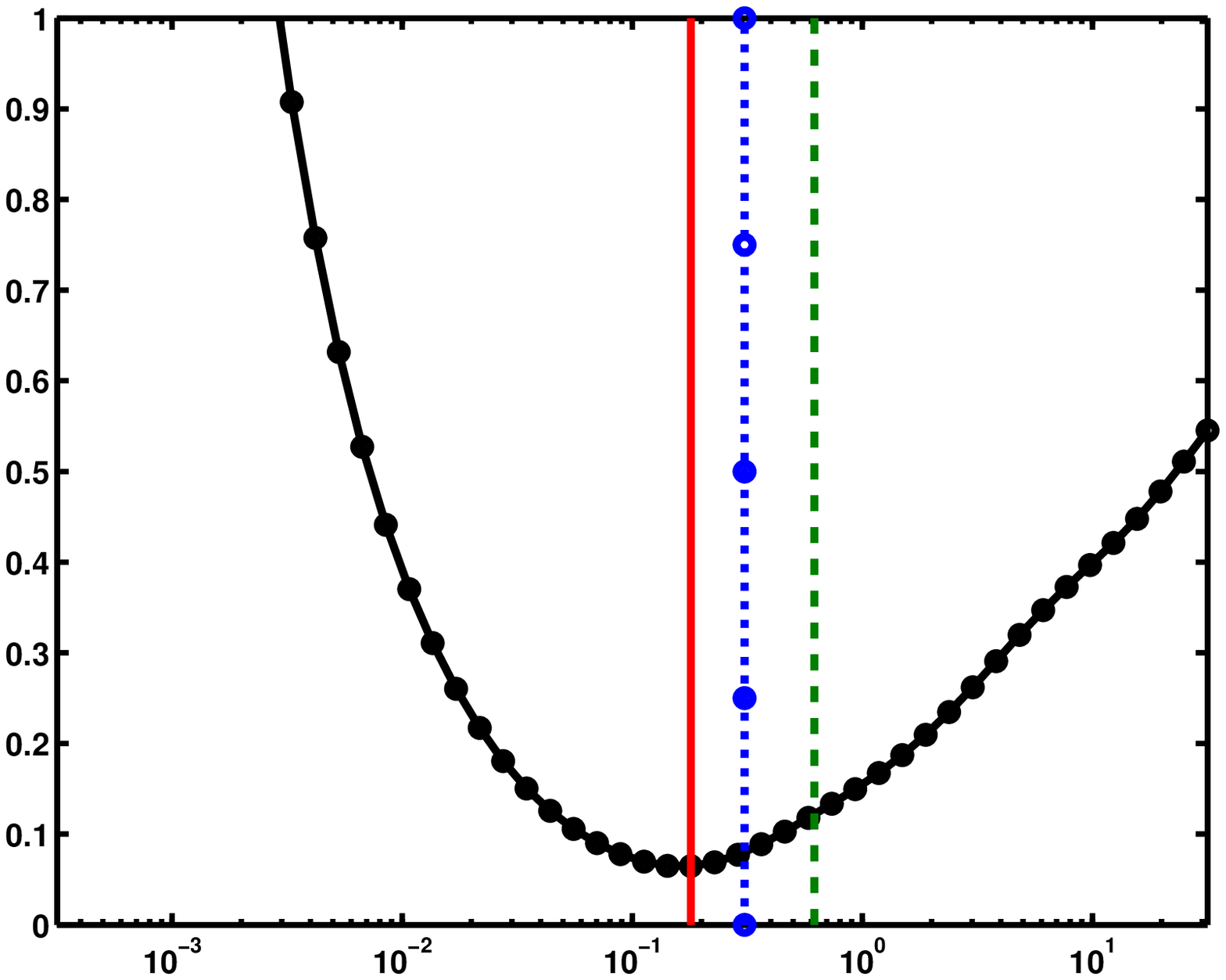}}
\subfigure[$L=I$]{\includegraphics[width=1.7in]{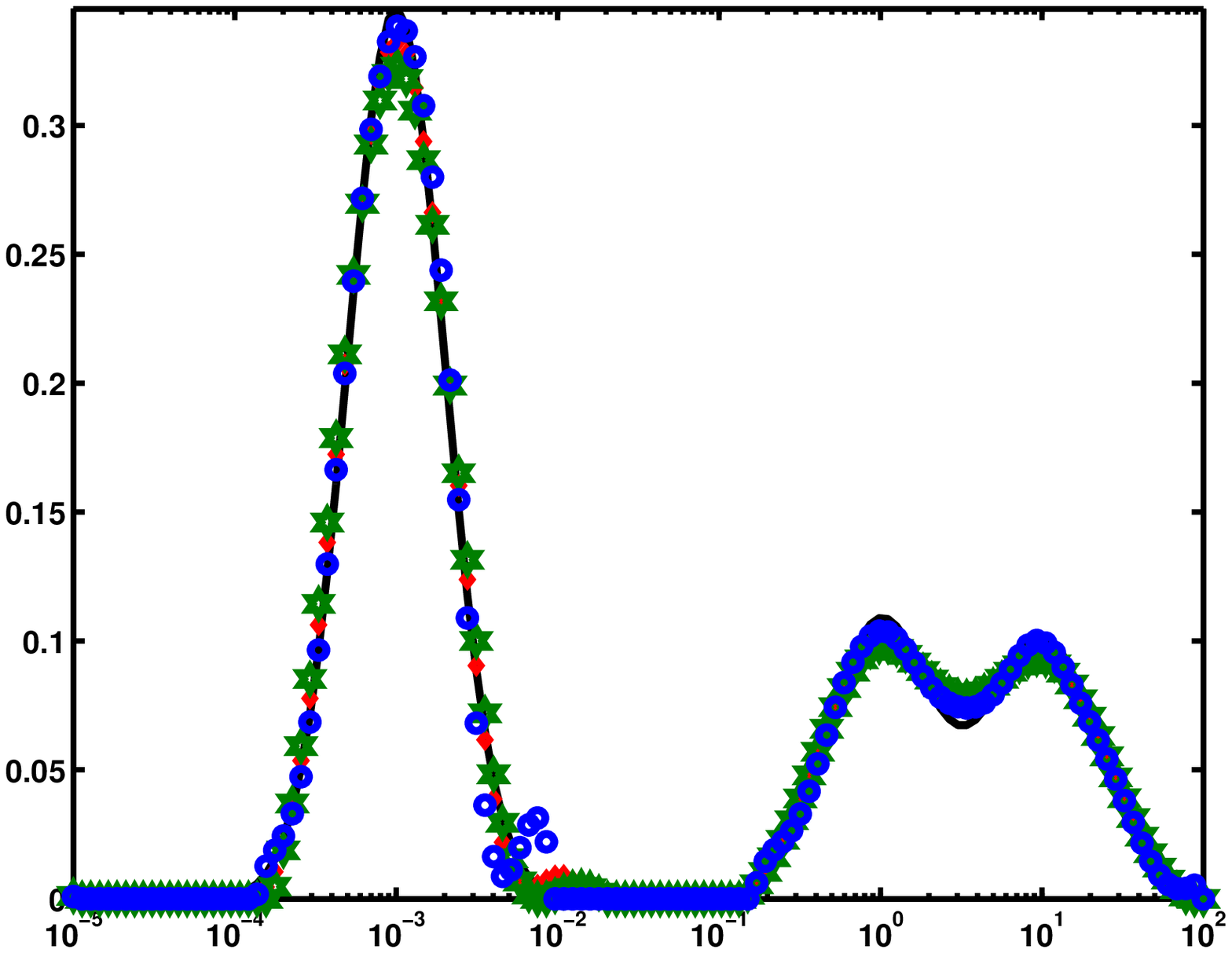}}
\subfigure[$L=L_1$]{\includegraphics[width=1.7in]{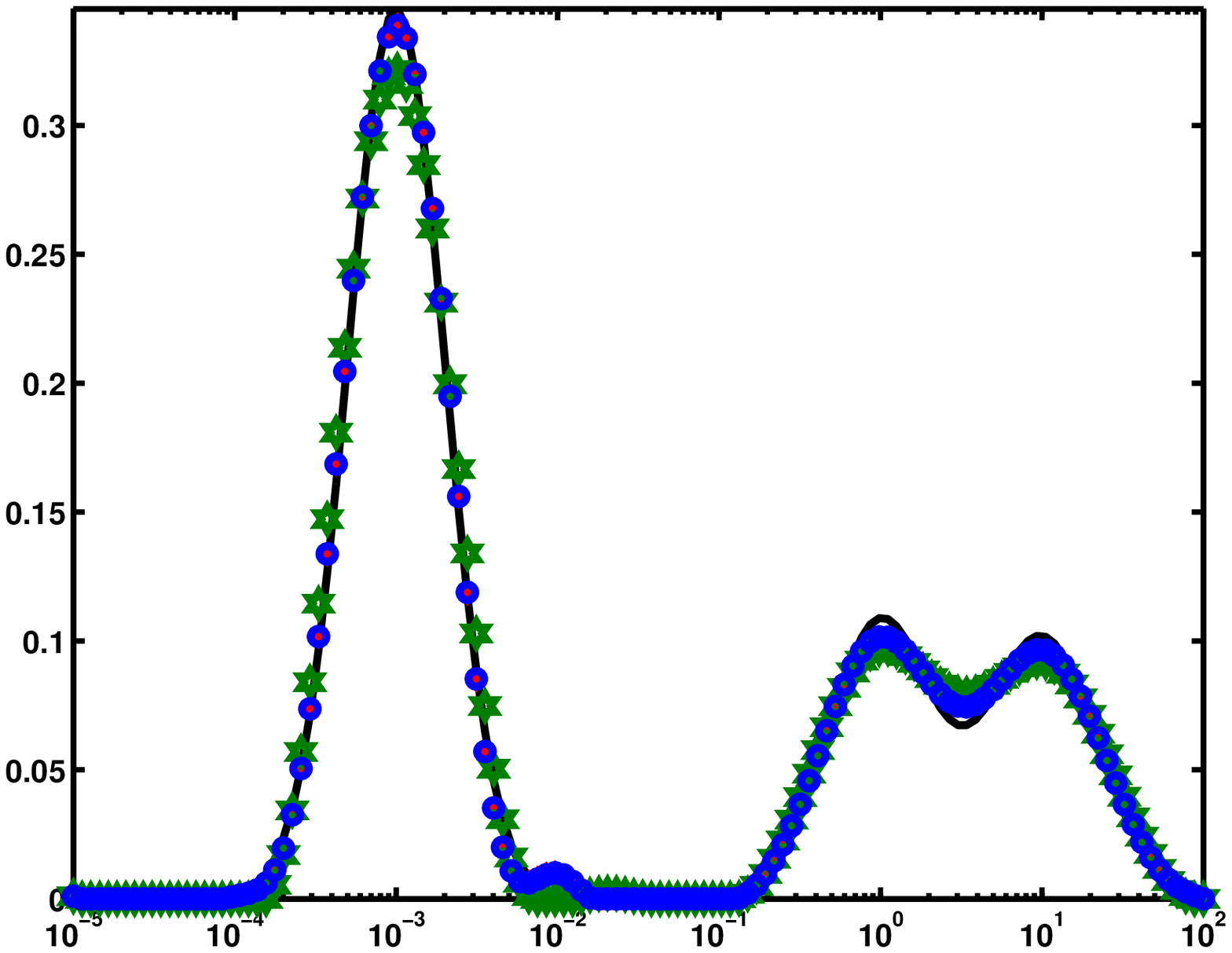}}
\subfigure[$L=L_2$]{\includegraphics[width=1.7in]{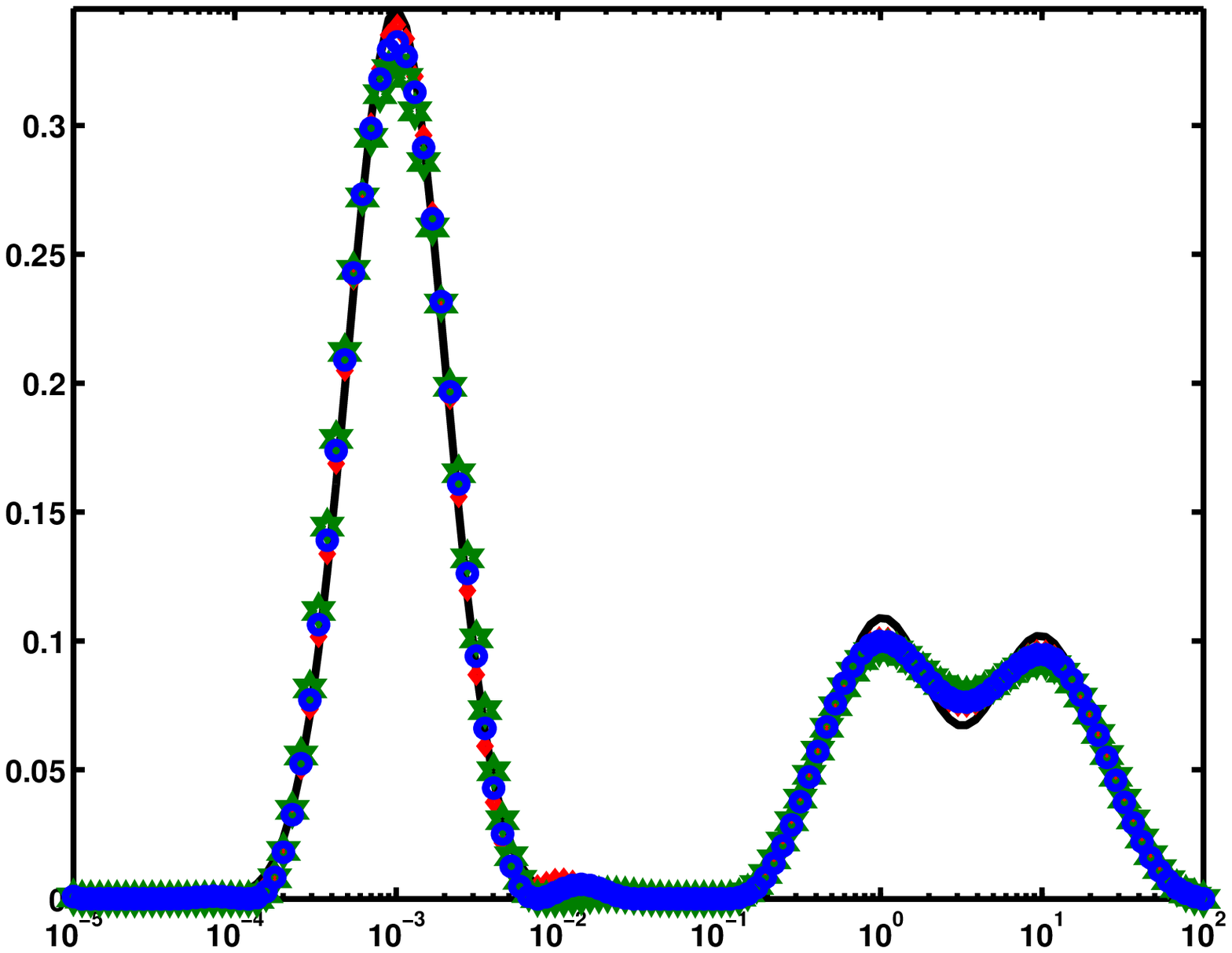}}
\caption{Mean error and example  CVX NNLS solutions.  $.1\%$ noise, LN-C data set, matrix $A_4$.}
\label{fig-lambdachoiceLN6A4LNCVX}
\end{figure}

To verify that these conclusions are not just a feature of the particular choice of the NNLS solver, we also repeated the experiments replacing the \texttt{lsqnonneg} solver for \eqref{nnls}. Several algorithms for NNLS are provided in the literature, see e.g. \cite{Bellavia, cvx, kimsra, Slawski}, but not all are relevant for dense but small ill-posed problems,. We chose CVX, which is easily implemented in Matlab, and is well--known as a robust package for specifying and solving convex programs, \cite{cvx,gb08}. These results are presented in Figures~
\ref{fig-lambdachoiceRQ1A4LNCVX}-\ref{fig-lambdachoiceLN6A4LNCVX} for exactly the same set of data. The results are remarkably consistent and demonstrate that the suggested technique to find the optimal regularization parameter for Tikhonov regularized non-negative least squares is robust to the underlying algorithm. Additional experiments, not reported here, also used the subspace Barzilai and Borwein (SBB) algorithm, presented in \cite{kimsra, sbbcode}. We found that this algorithm was overall less robust for solving the small-scale non-negative least squares problems described here, which does not contradict the results in \cite{kimsra}. Still, where the algorithm succeeded, for larger choices of $\lambda$, again NCP and L-curve estimates of the optimal choice for $\lambda$ are feasible. A selection of these results is presented in the supplementary materials, as are further experiments for higher noise levels using both CVX and \texttt{lsqnonneg}, \cite{Supp}.

\subsection{Comparison of LLS and NNLS}\label{lsnnls}
 
\begin{figure}[!ht]
 \centering
\subfigure[$L=I$]{\includegraphics[width=1.7in]{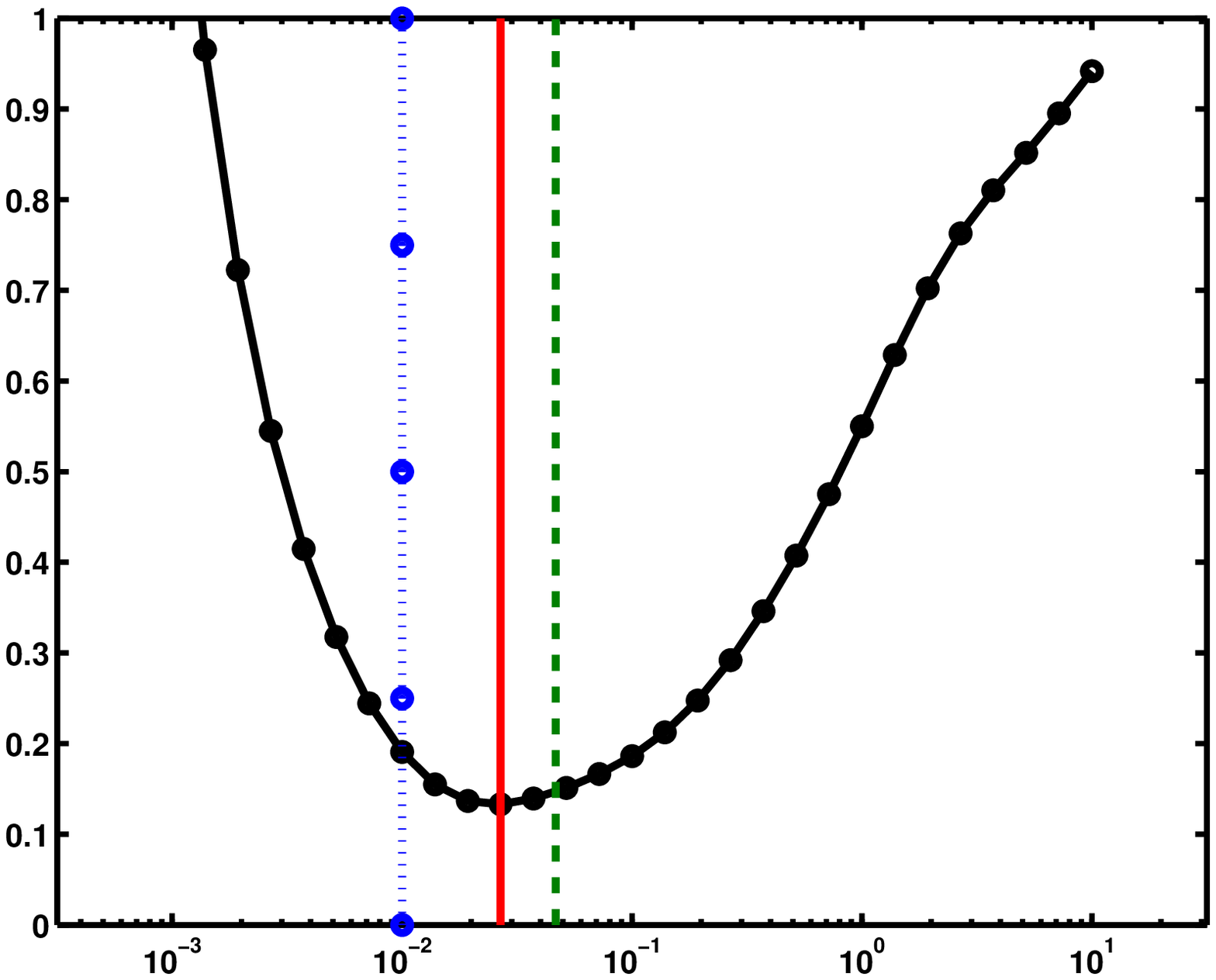}}
\subfigure[$L=L_1$]{\includegraphics[width=1.7in]{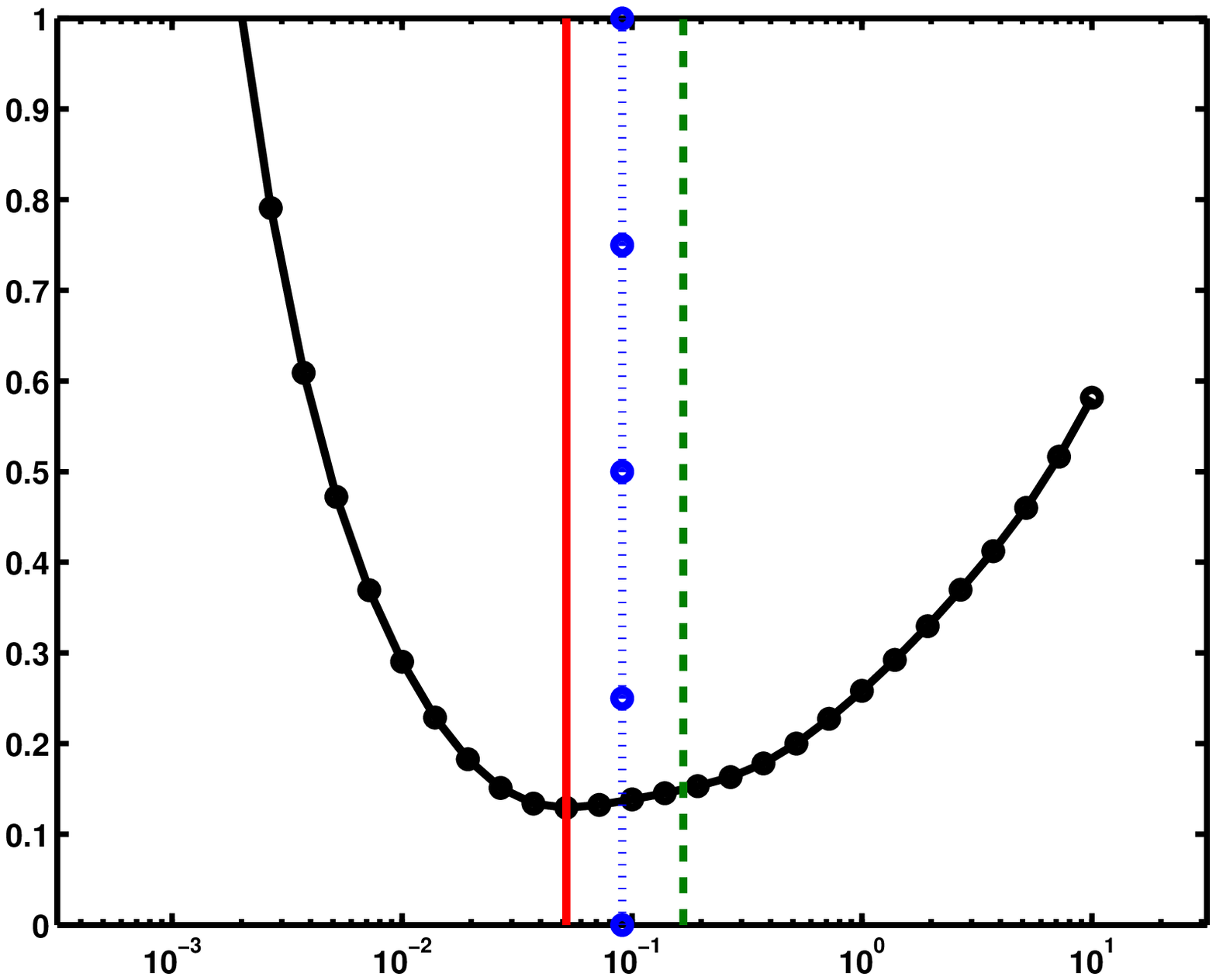}}
\subfigure[$L=L_2$]{\includegraphics[width=1.7in]{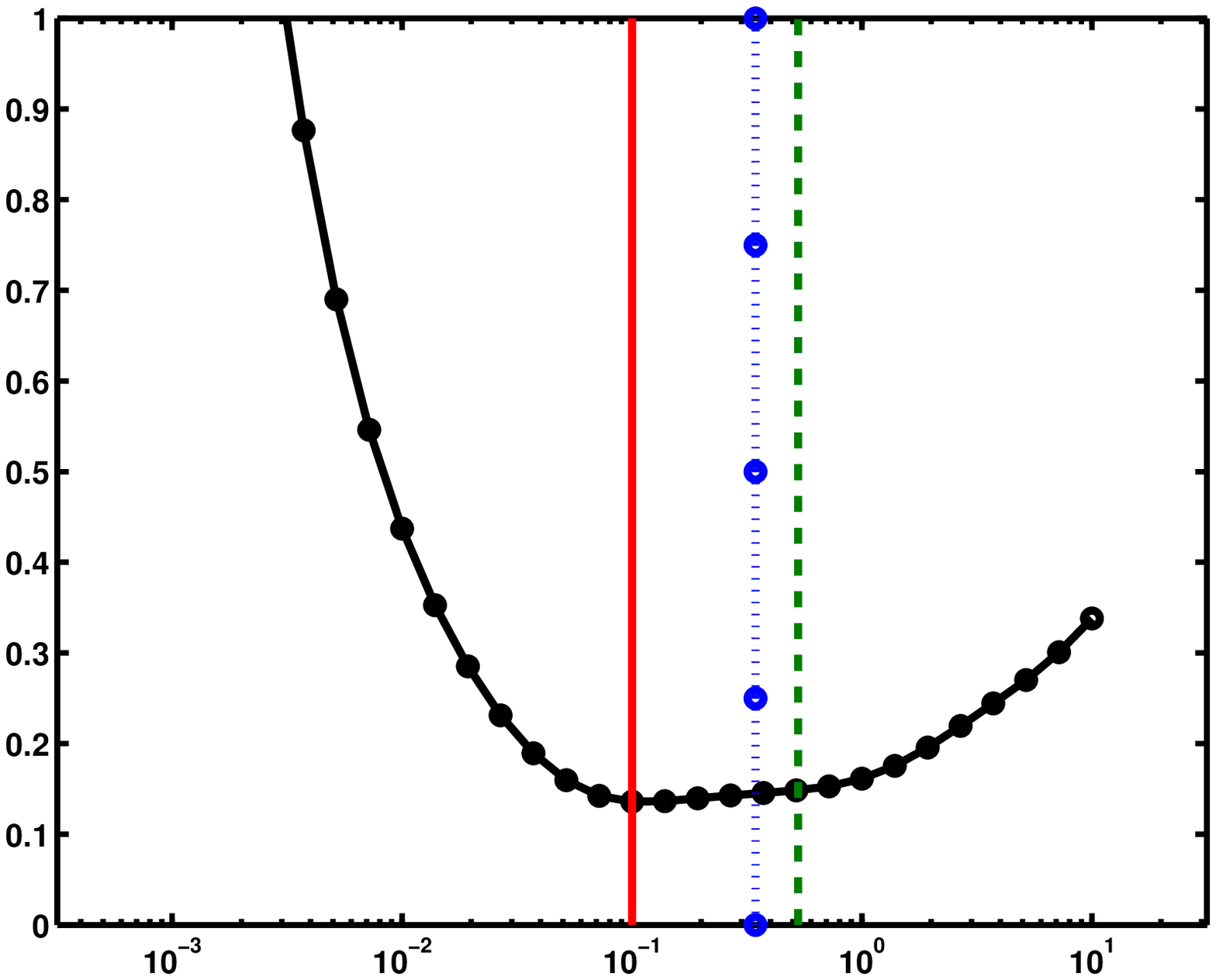}}
\subfigure[$L=I$]{\includegraphics[width=1.7in]{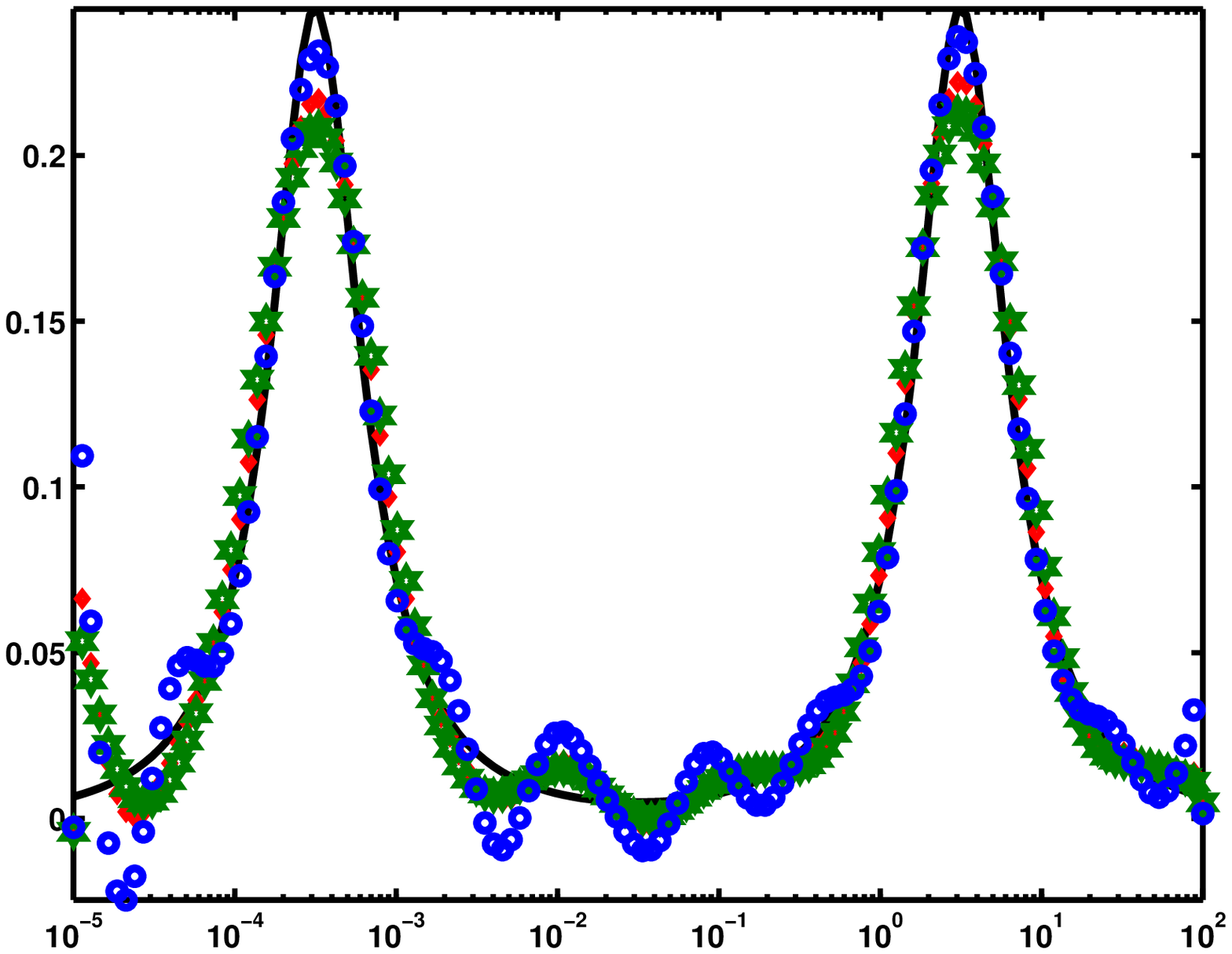}}
\subfigure[$L=L_1$]{\includegraphics[width=1.7in]{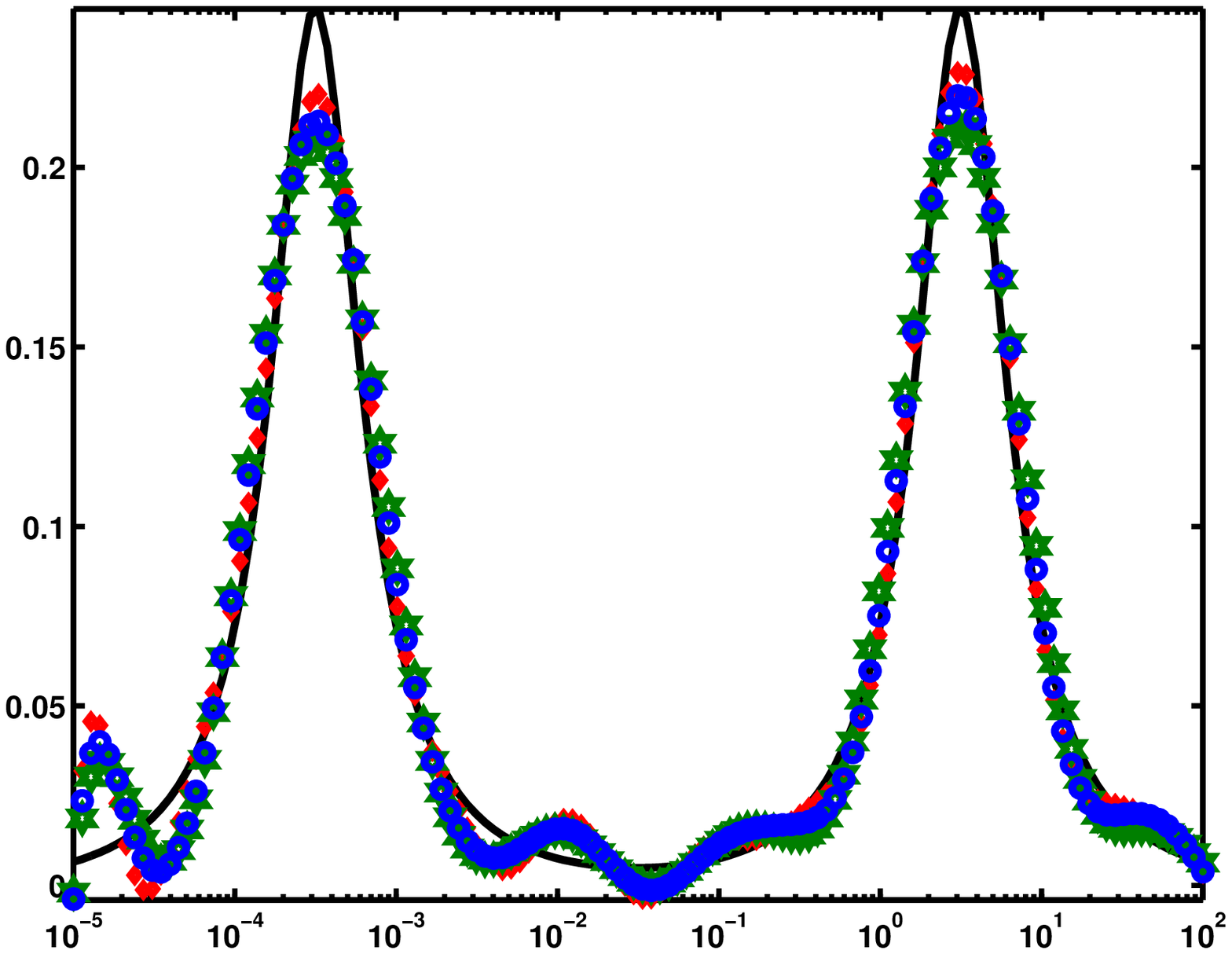}}
\subfigure[$L=L_2$]{\includegraphics[width=1.7in]{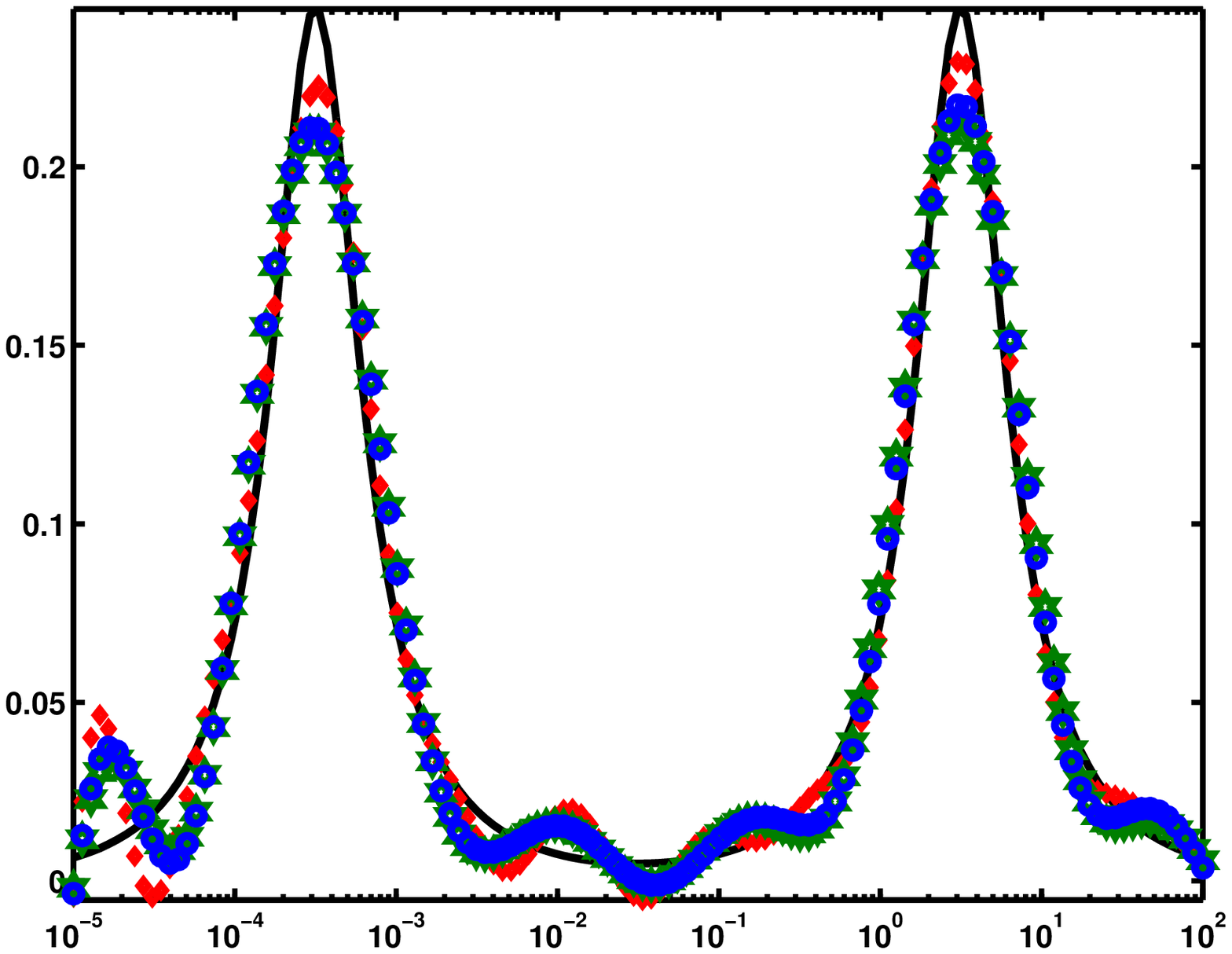}}
\caption{Mean error and example LLS solutions.  $.1\%$ noise, RQ-A data set, matrix $A_4$.}
\label{fig-lambdachoiceRQ1A4LNLS}
\end{figure}

\begin{figure}[!ht]
  \centering
\subfigure[$L=I$]{\includegraphics[width=1.7in]{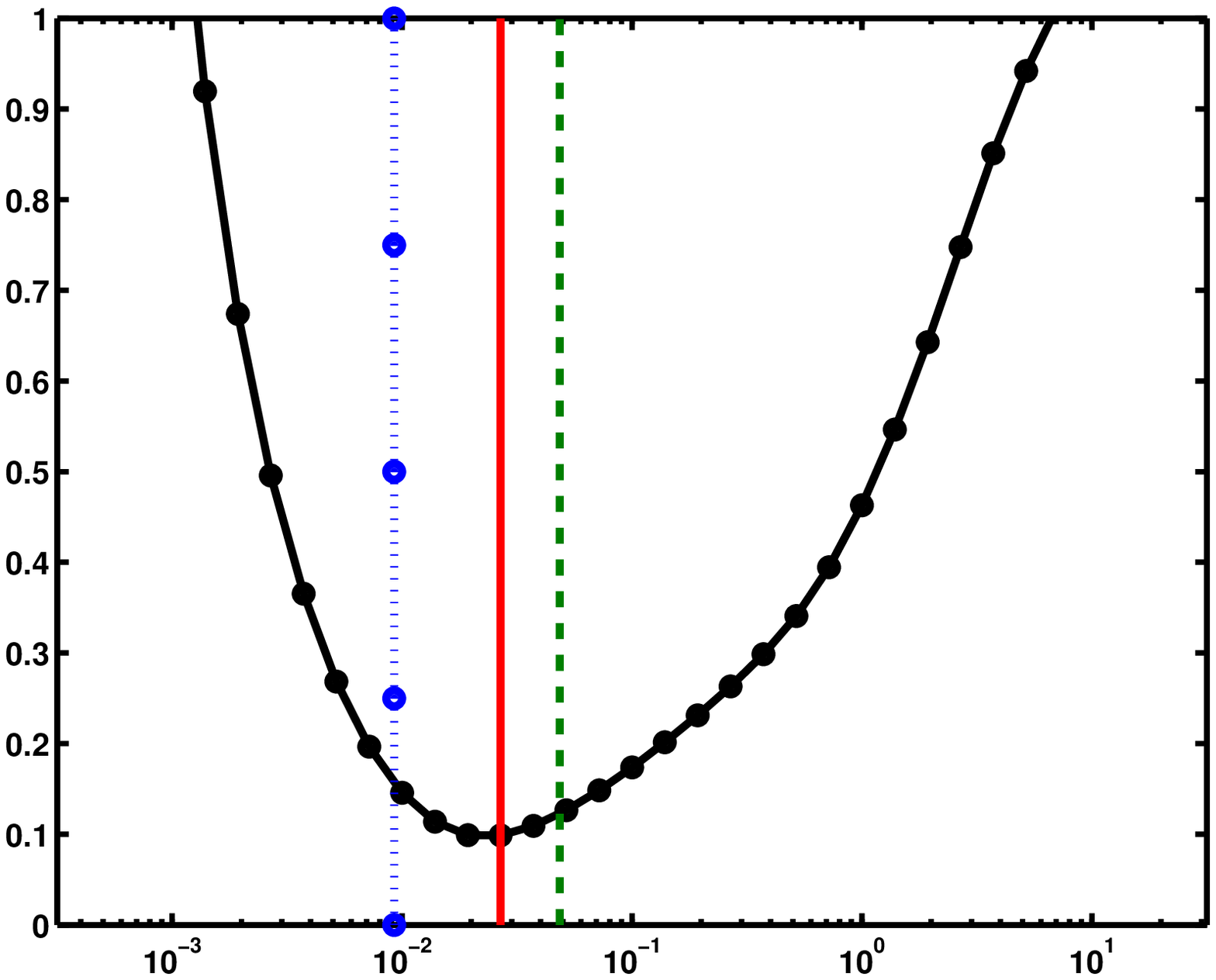}}
\subfigure[$L=L_1$]{\includegraphics[width=1.7in]{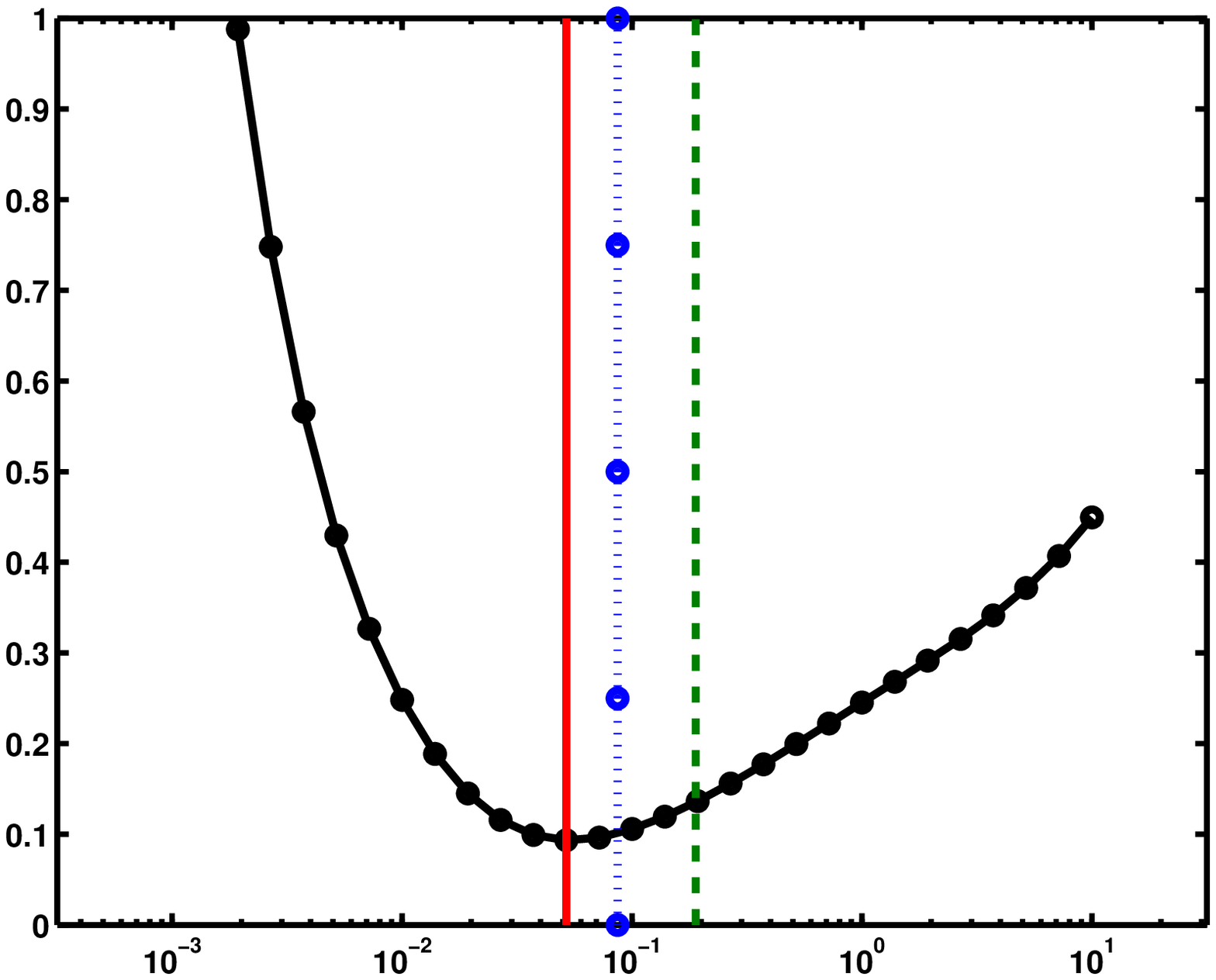}}
\subfigure[$L=L_2$]{\includegraphics[width=1.7in]{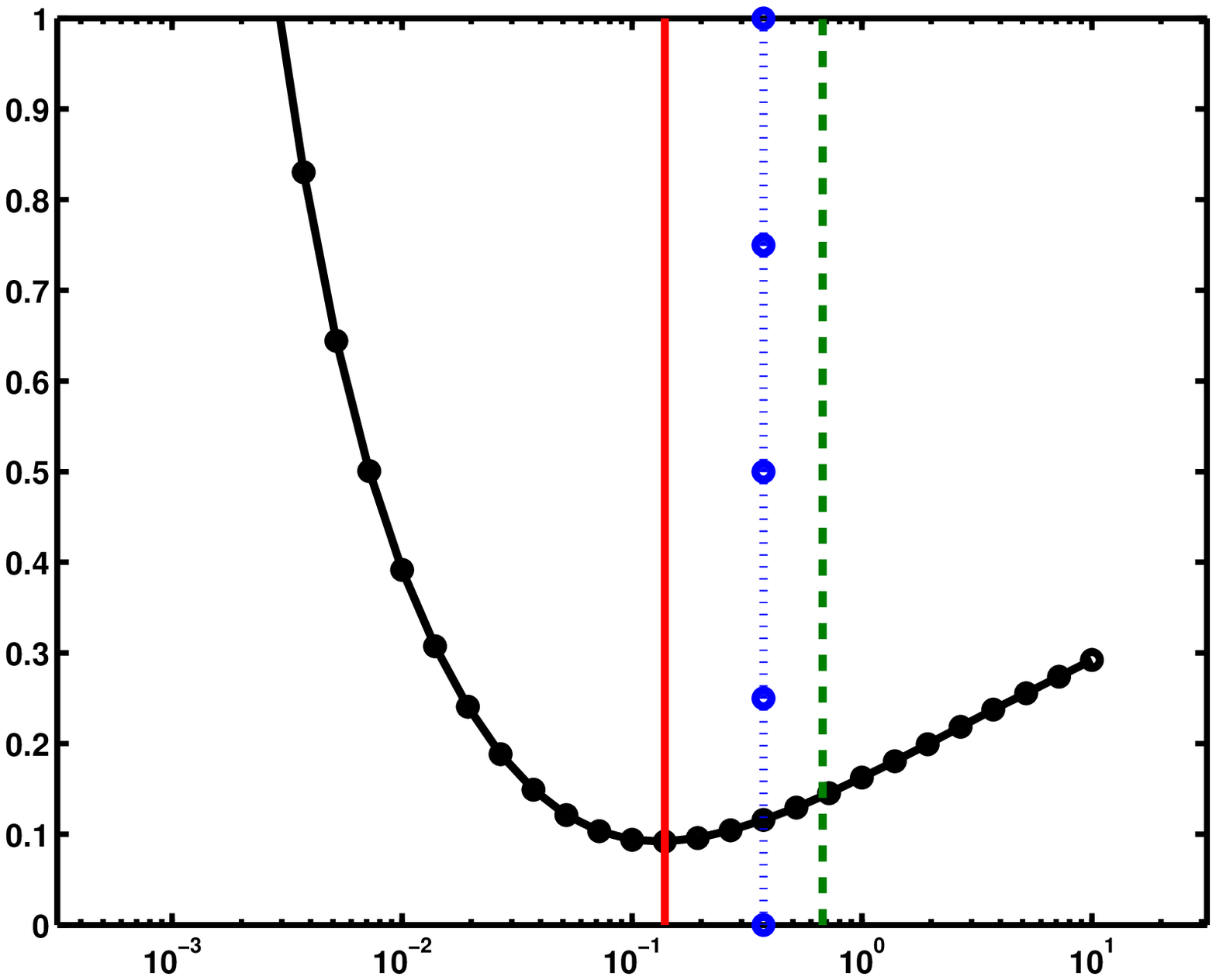}}
\subfigure[$L=I$]{\includegraphics[width=1.7in]{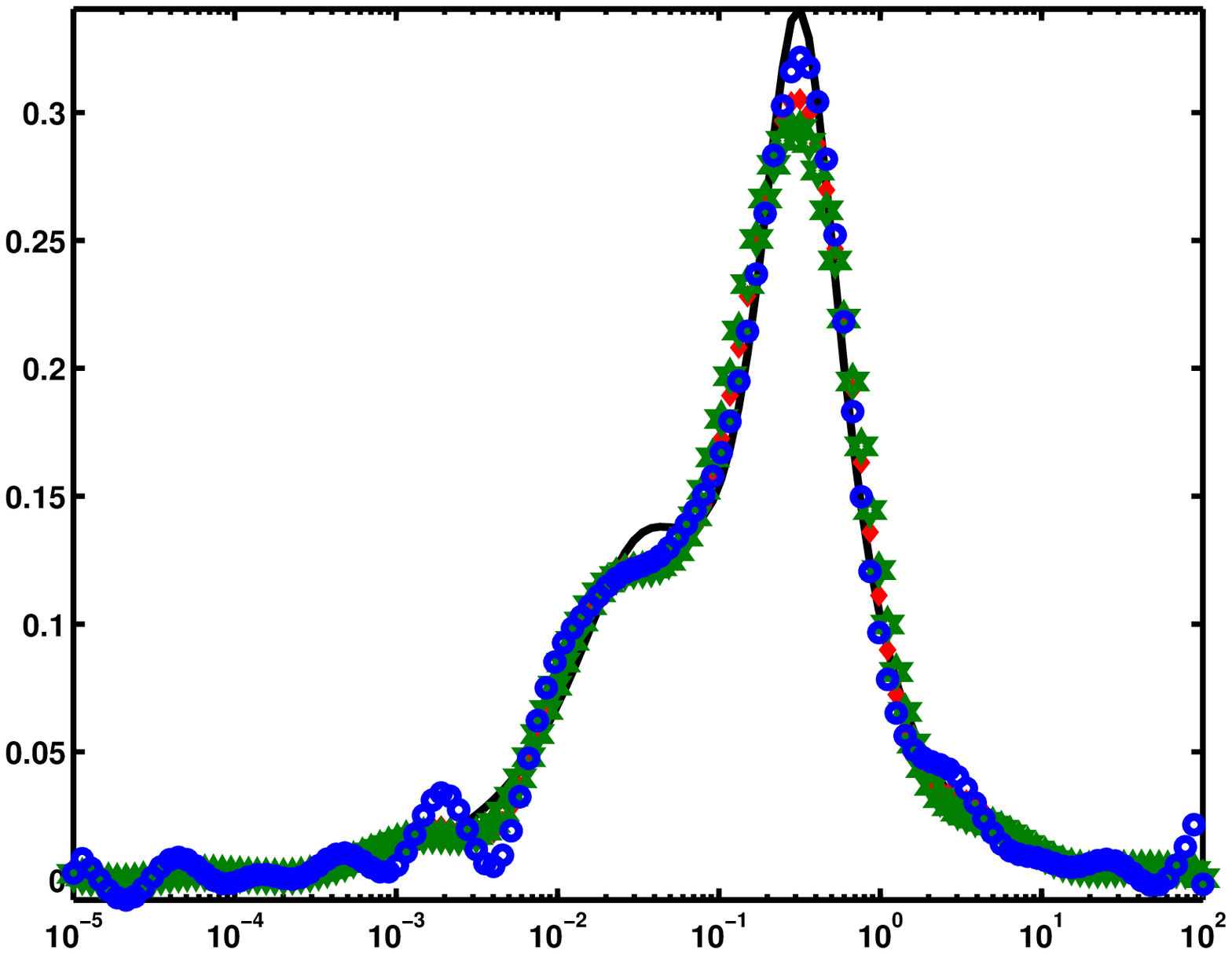}}
\subfigure[$L=L_1$]{\includegraphics[width=1.7in]{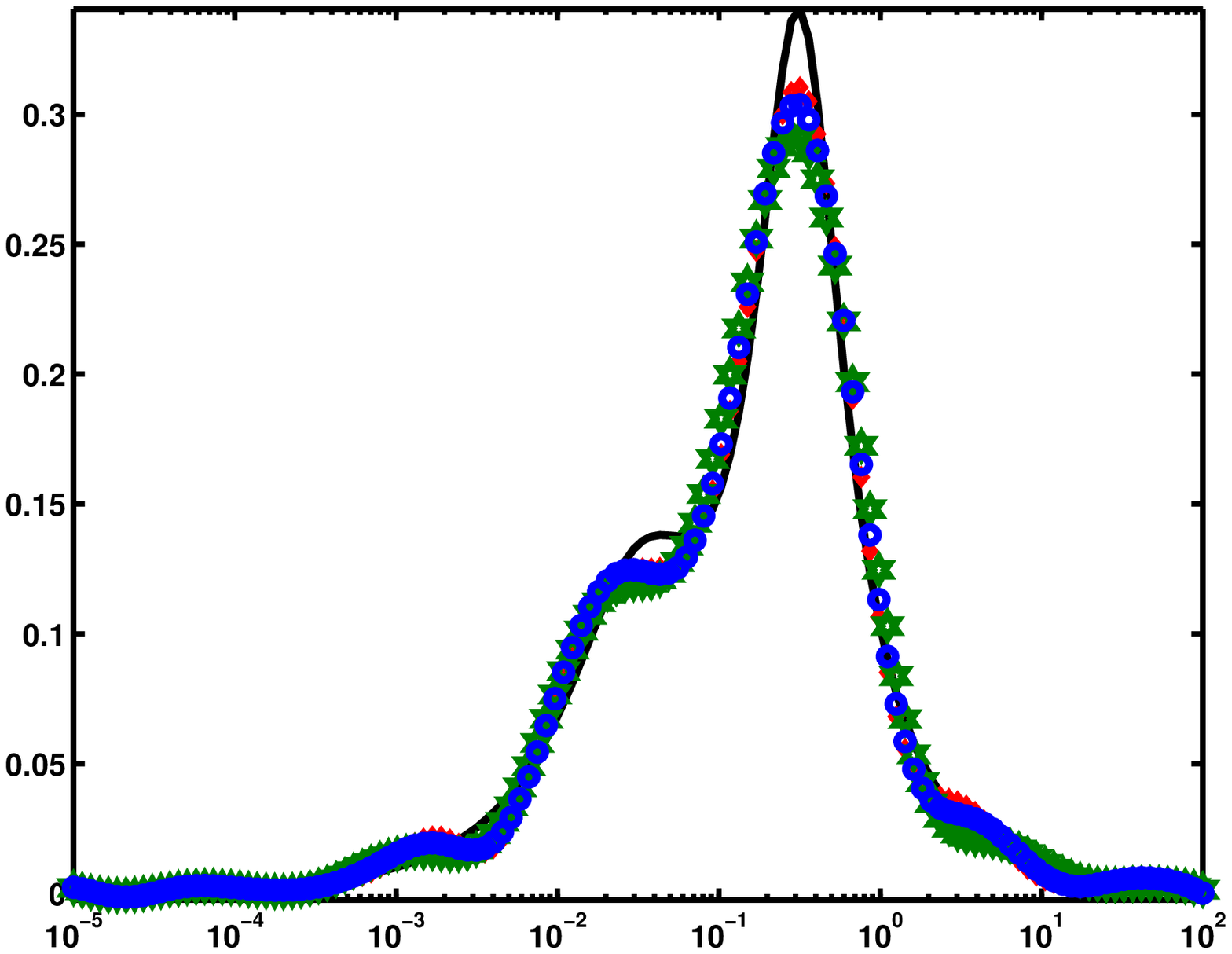}}
\subfigure[$L=L_2$]{\includegraphics[width=1.7in]{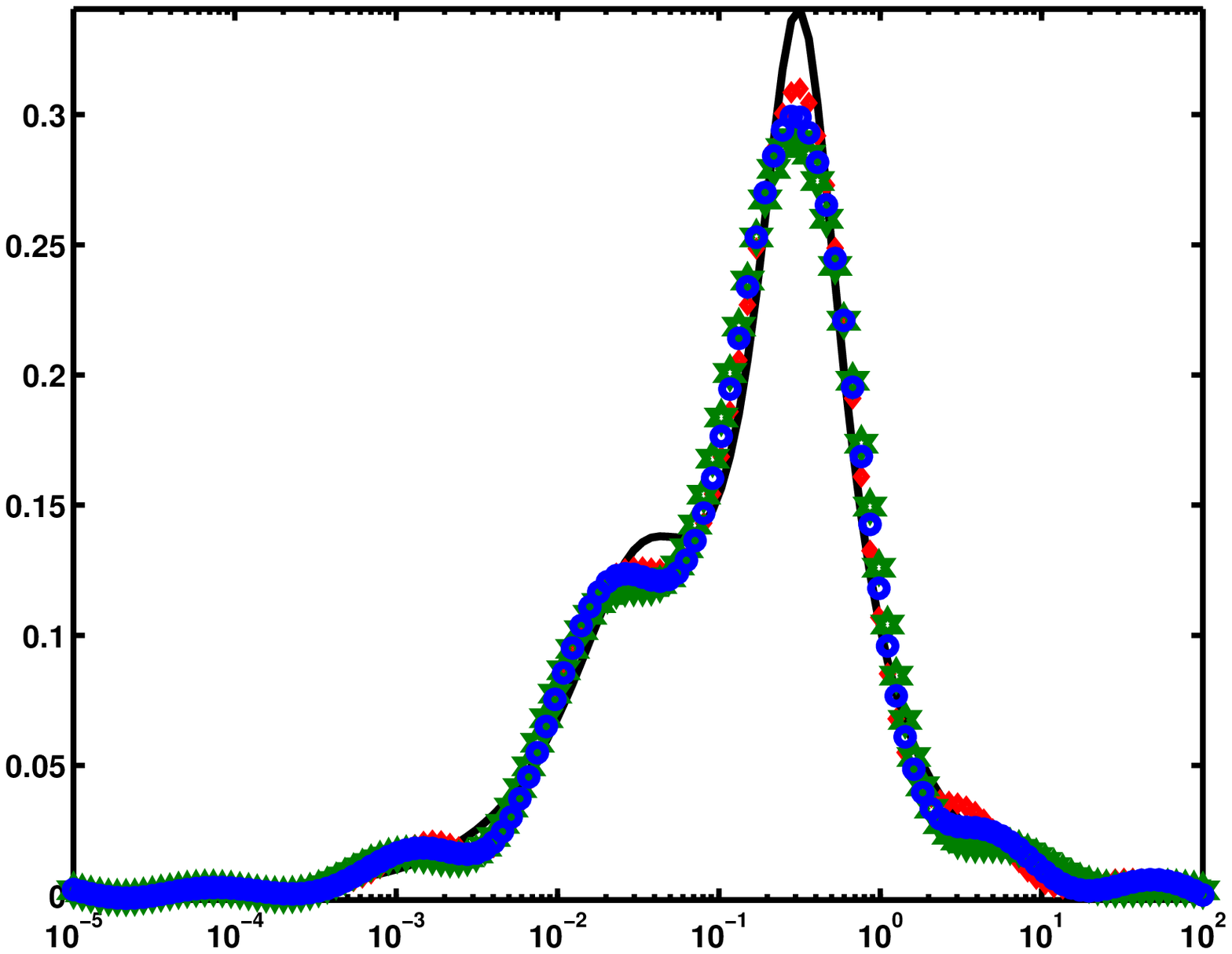}}
\caption{Mean error and example LLS solutions. $.1\%$ noise, RQ-B data set, matrix $A_4$.}
\label{fig-lambdachoiceRQ5A4LNLS}
\end{figure}

\begin{figure}[!ht]
  \centering
\subfigure[$L=I$]{\includegraphics[width=1.7in]{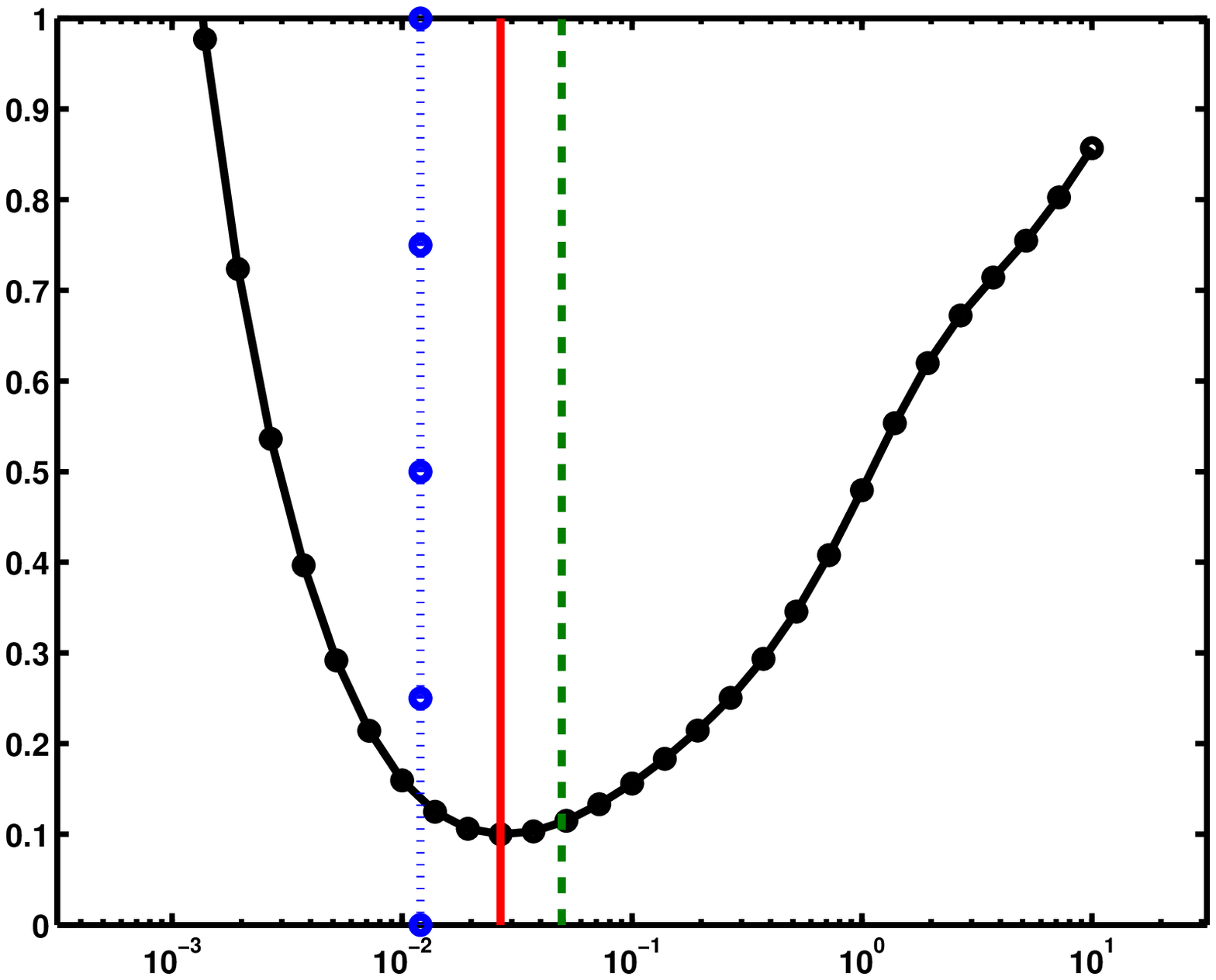}}
\subfigure[$L=L_1$]{\includegraphics[width=1.7in]{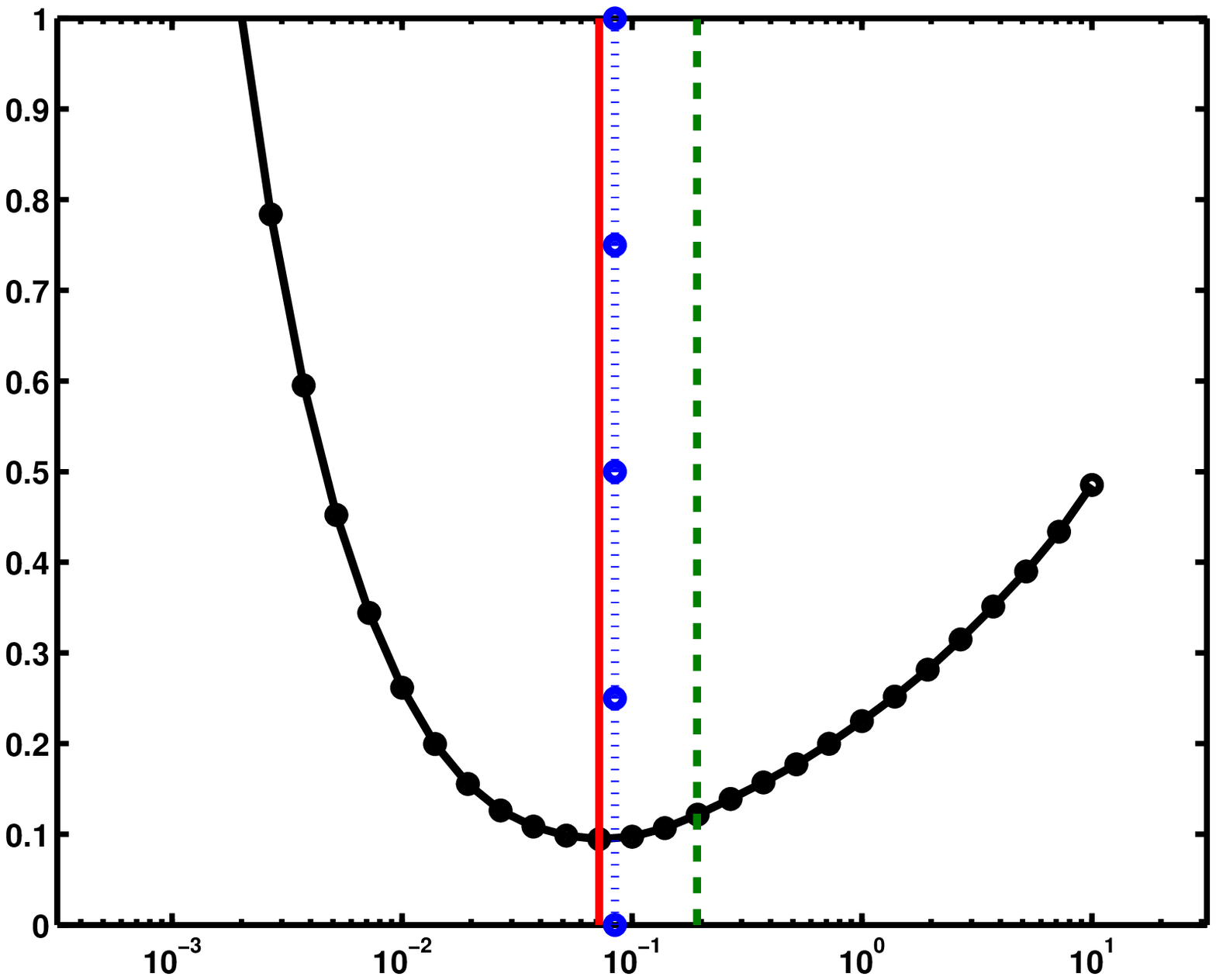}}
\subfigure[$L=L_2$]{\includegraphics[width=1.7in]{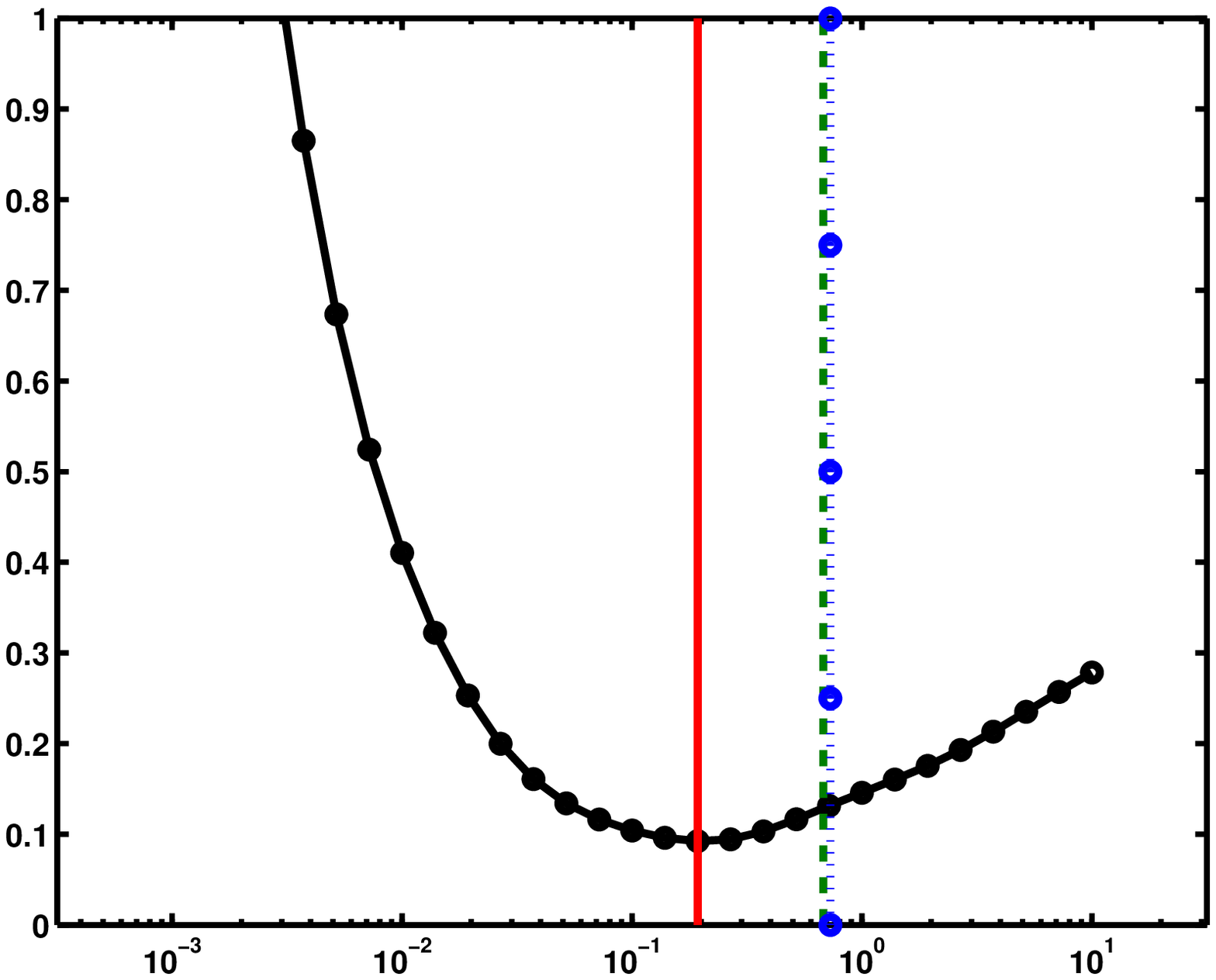}}
\subfigure[$L=I$]{\includegraphics[width=1.7in]{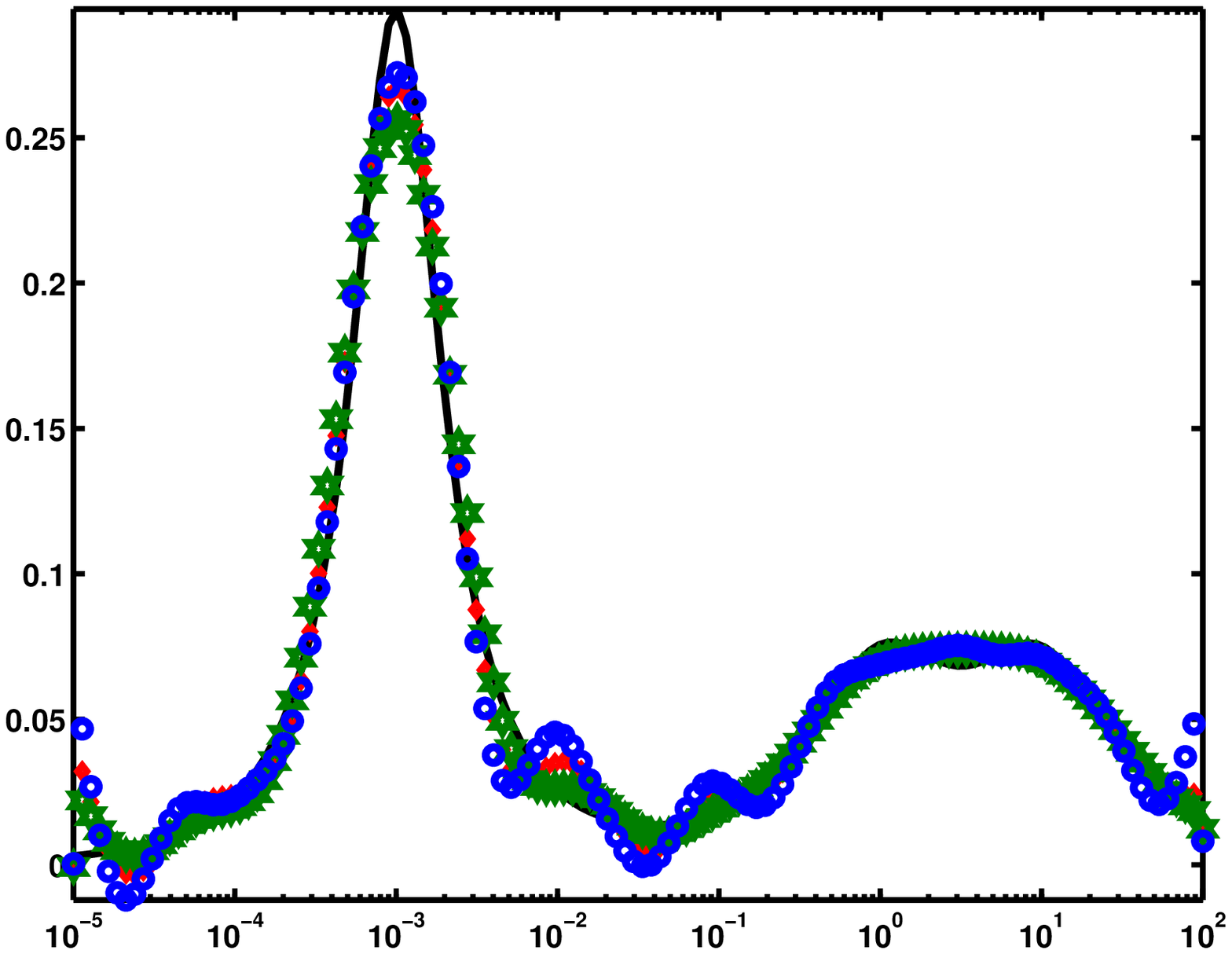}}
\subfigure[$L=L_1$]{\includegraphics[width=1.7in]{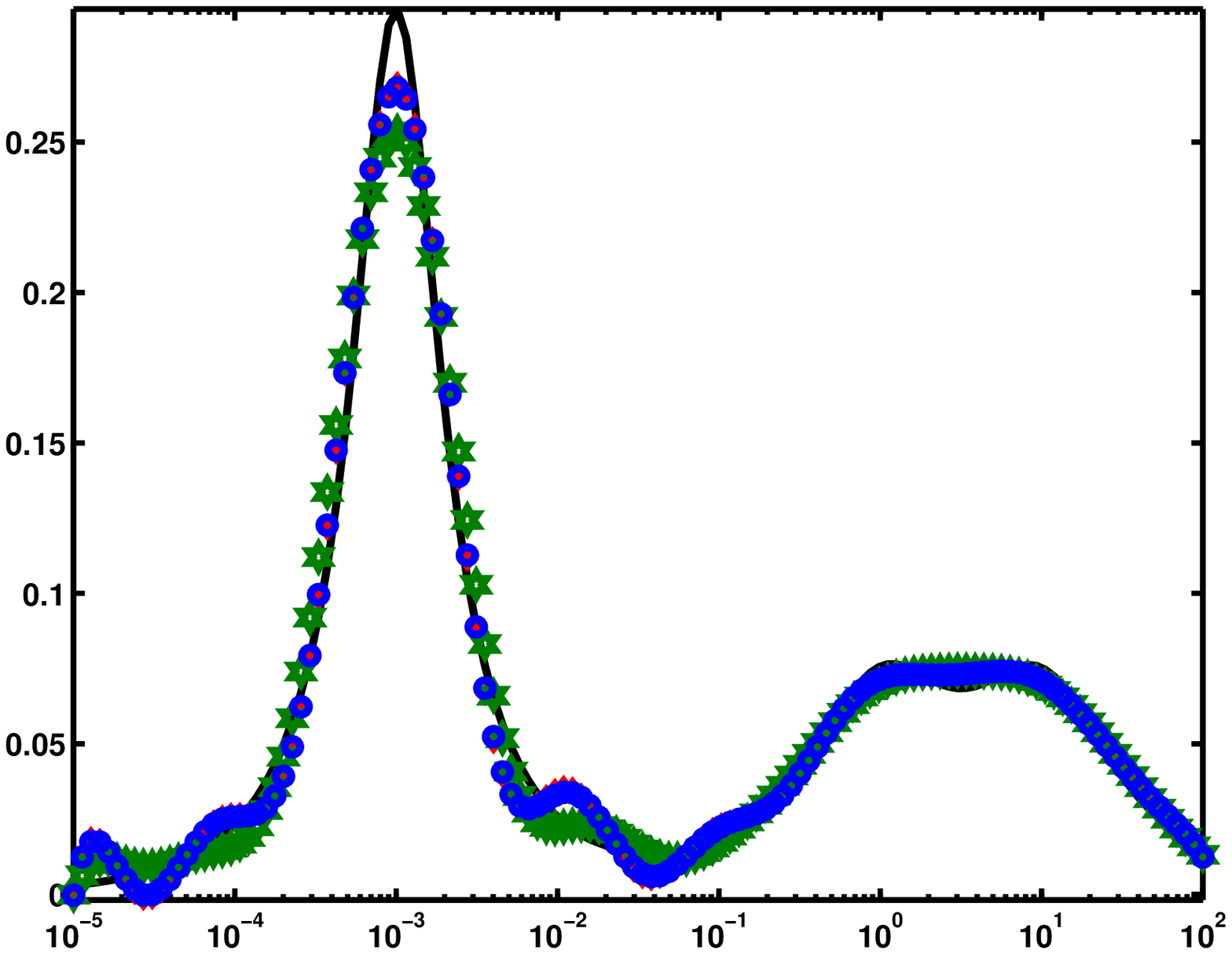}}
\subfigure[$L=L_2$]{\includegraphics[width=1.7in]{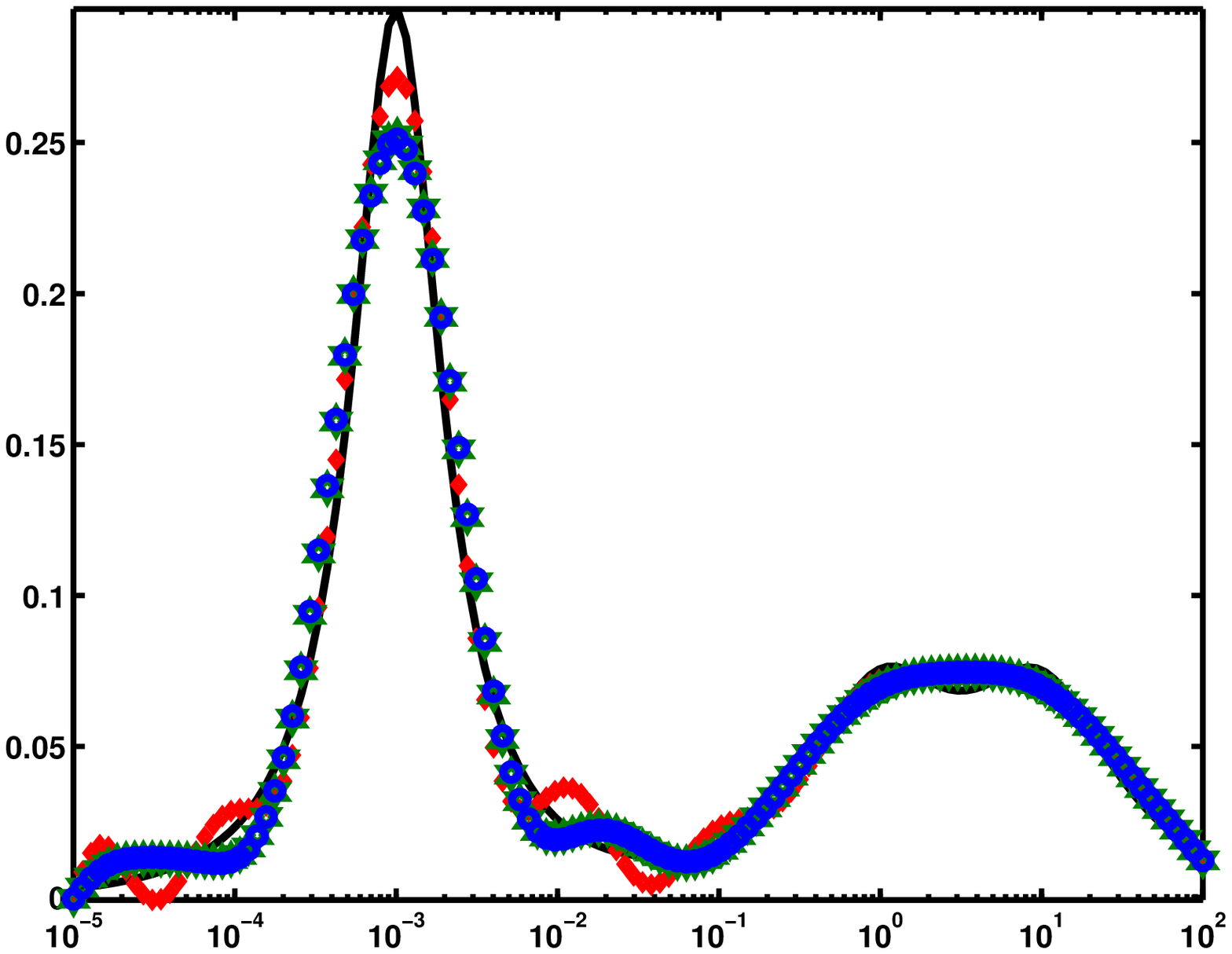}}
\caption{Mean error and example LLS solutions.  $.1\%$ noise, RQ-C data set, matrix $A_4$.}
\label{fig-lambdachoiceRQ6A4LNLS}
\end{figure}

\begin{figure}[!ht]
  \centering
\subfigure[$L=I$]{\includegraphics[width=1.7in]{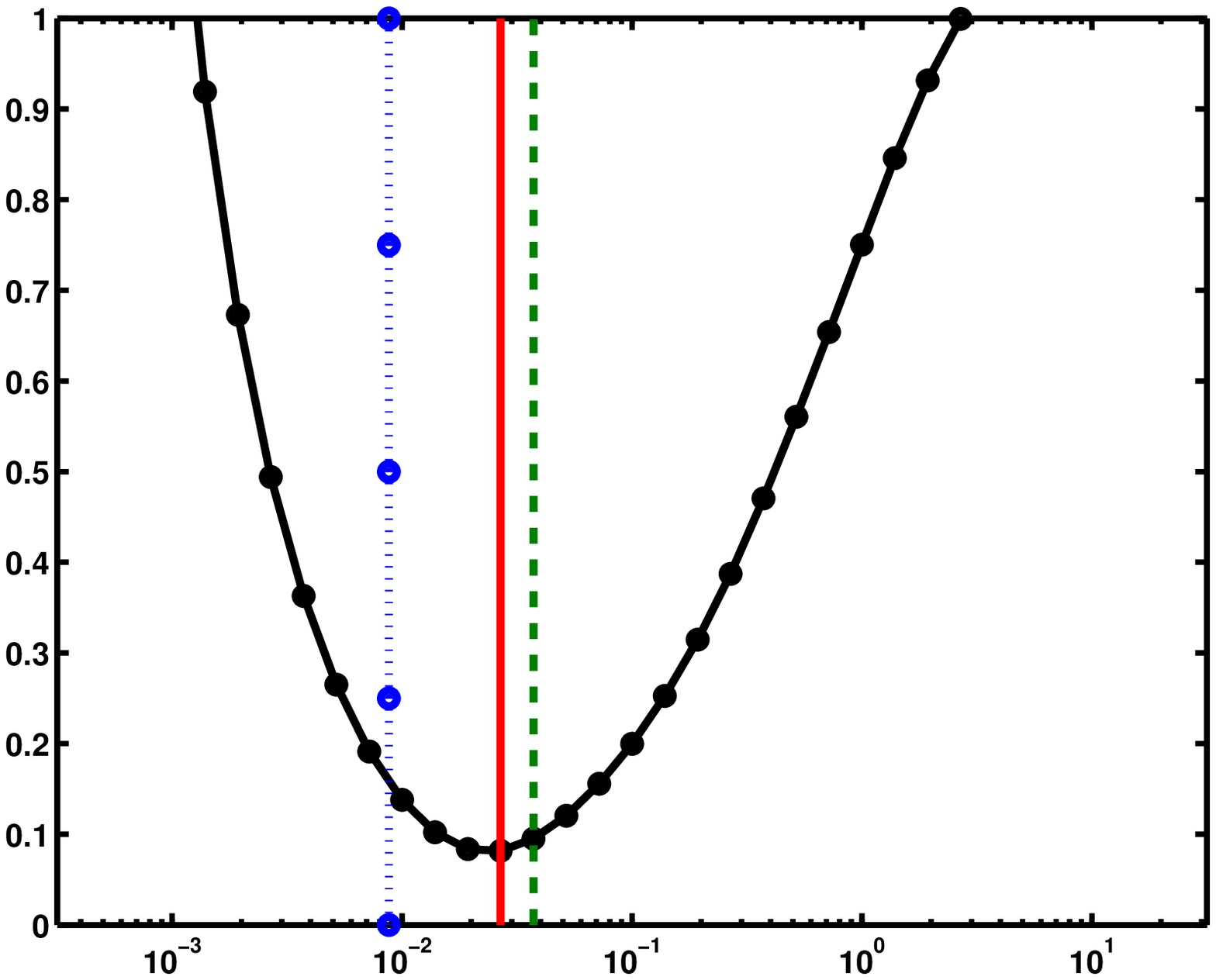}}
\subfigure[$L=L_1$]{\includegraphics[width=1.7in]{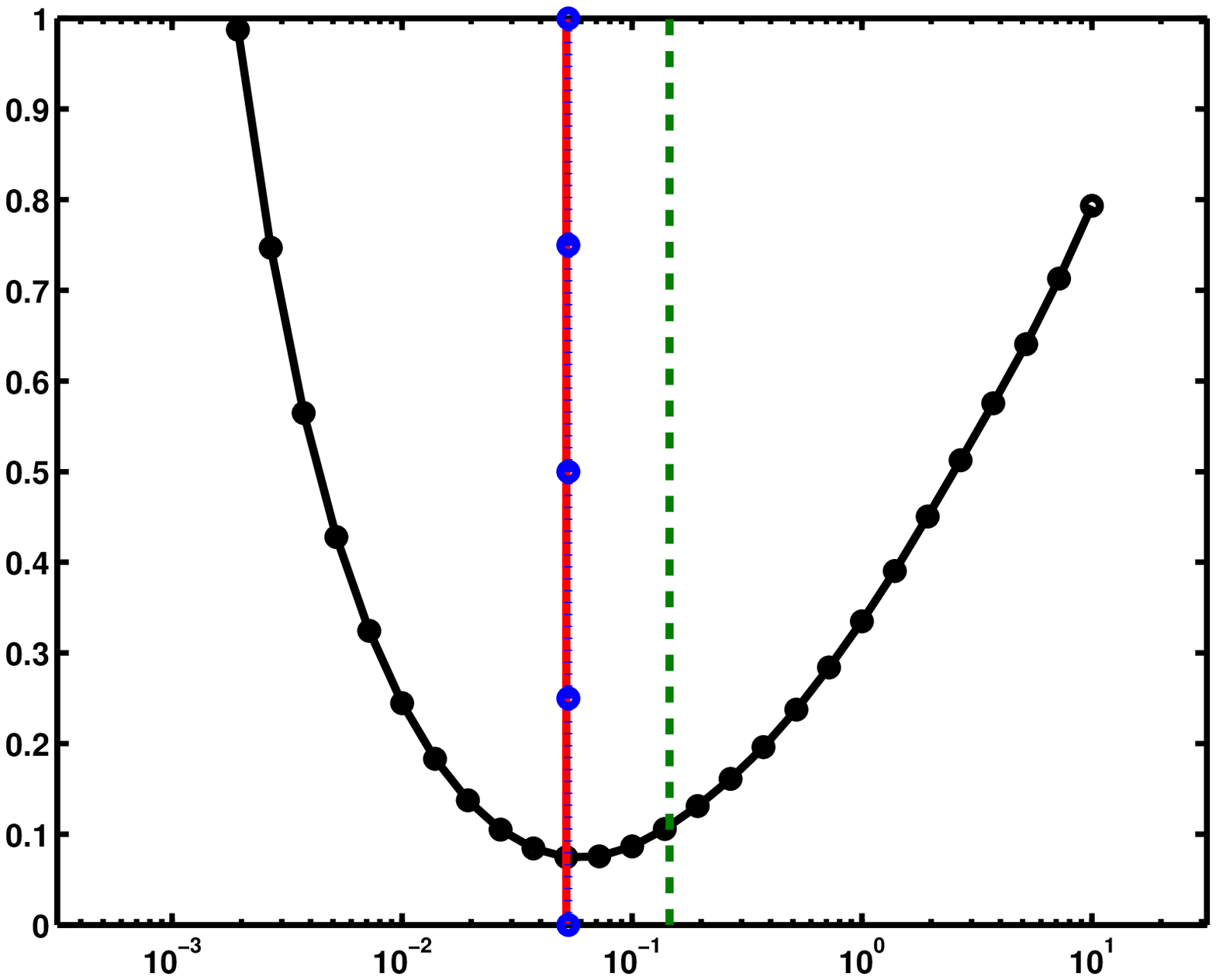}}
\subfigure[$L=L_2$]{\includegraphics[width=1.7in]{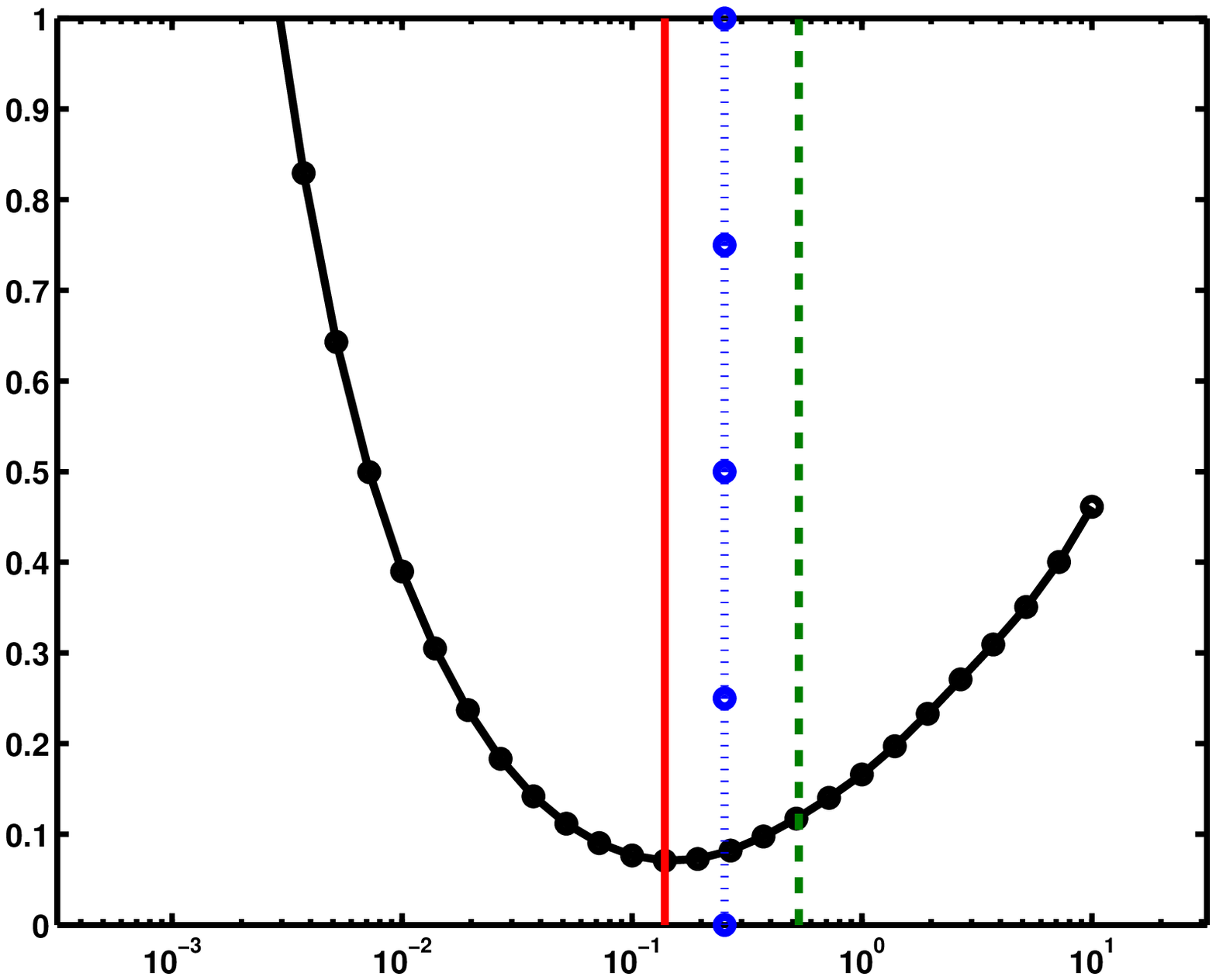}}
\subfigure[$L=I$]{\includegraphics[width=1.7in]{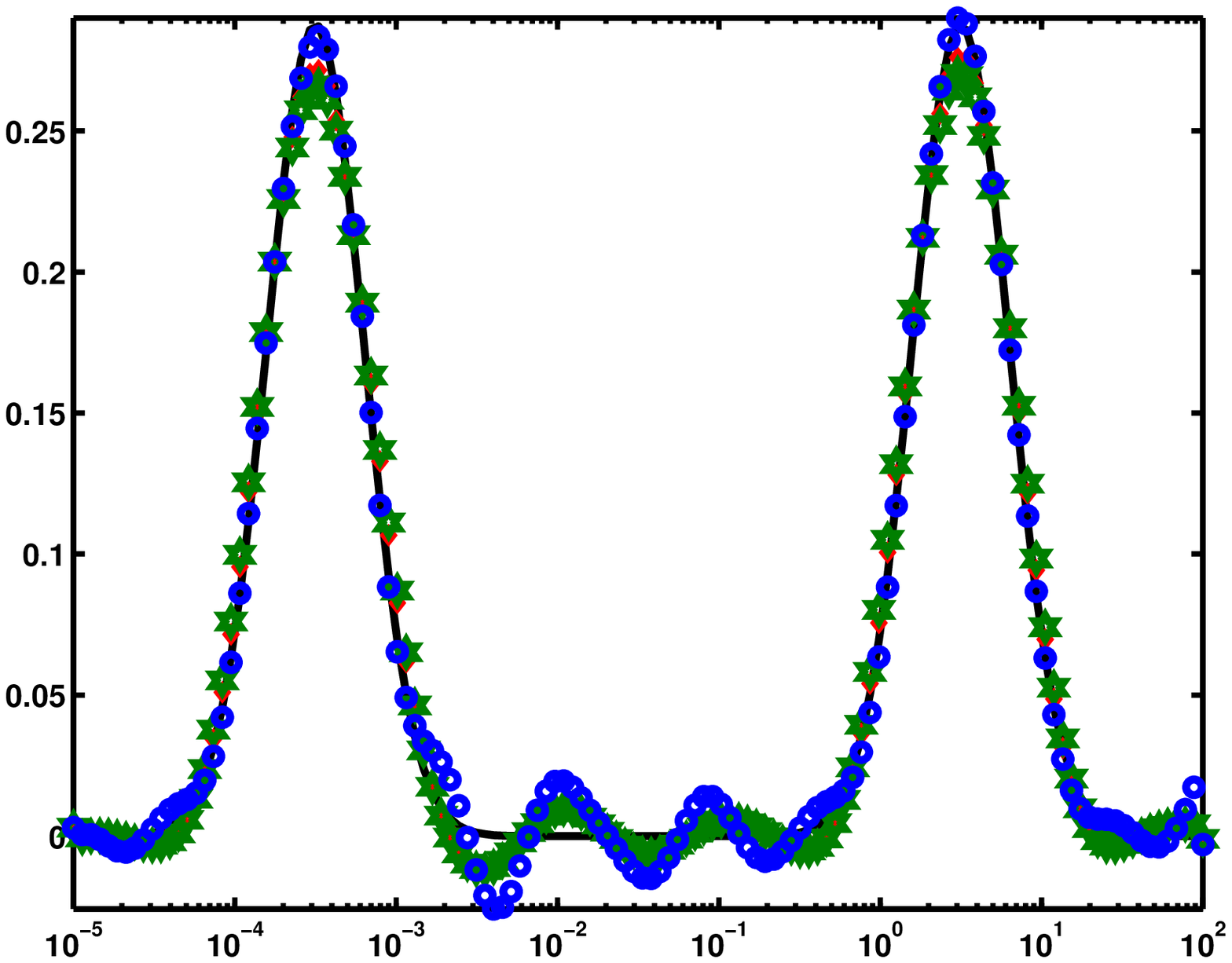}}
\subfigure[$L=L_1$]{\includegraphics[width=1.7in]{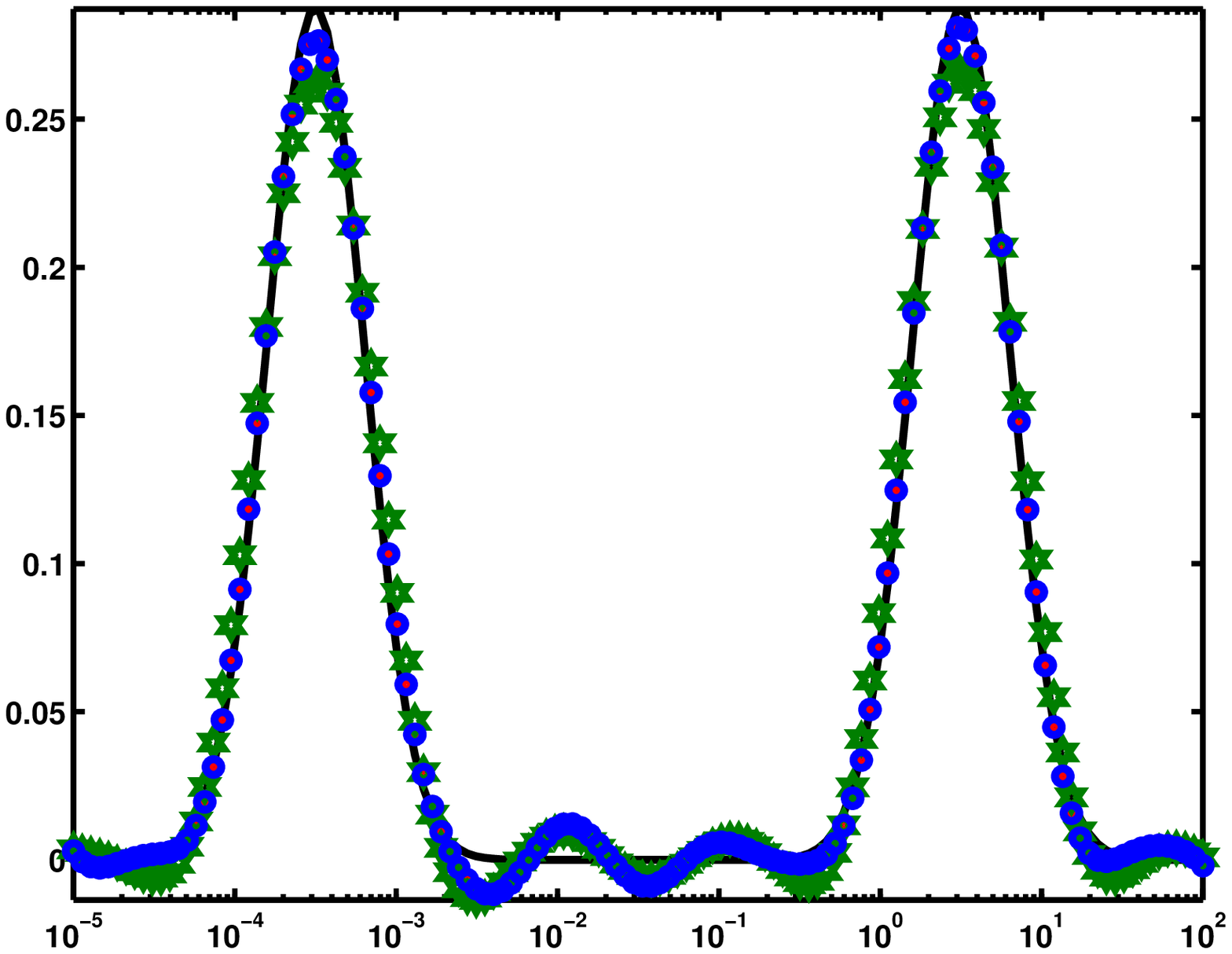}}
\subfigure[$L=L_2$]{\includegraphics[width=1.7in]{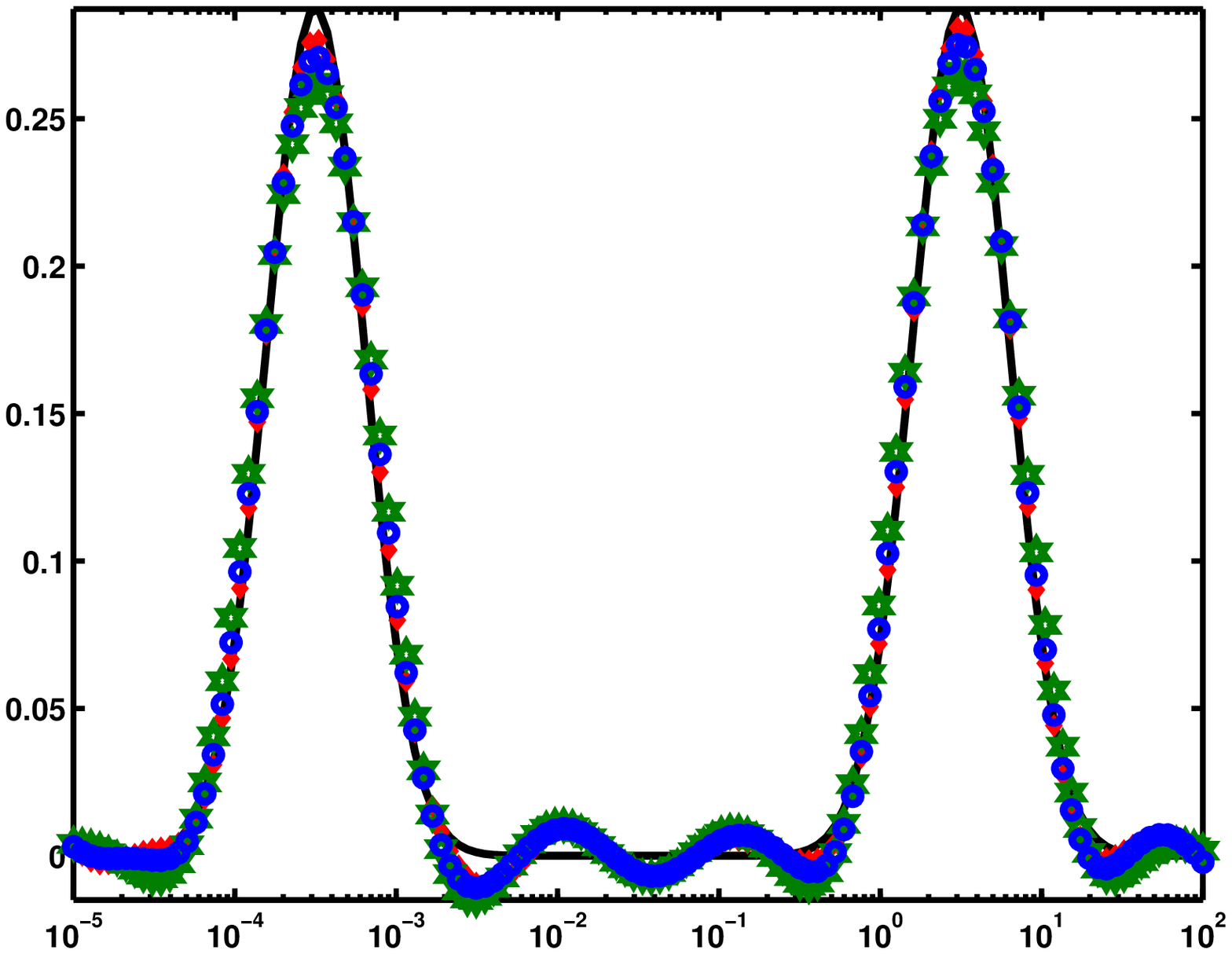}}
\caption{Mean error and example LLS solutions.  $.1\%$ noise, LN-A data set, matrix $A_4$.}
\label{fig-lambdachoiceLN2A4LNLS}
\end{figure}

\begin{figure}[!ht]
  \centering
\subfigure[$L=I$]{\includegraphics[width=1.7in]{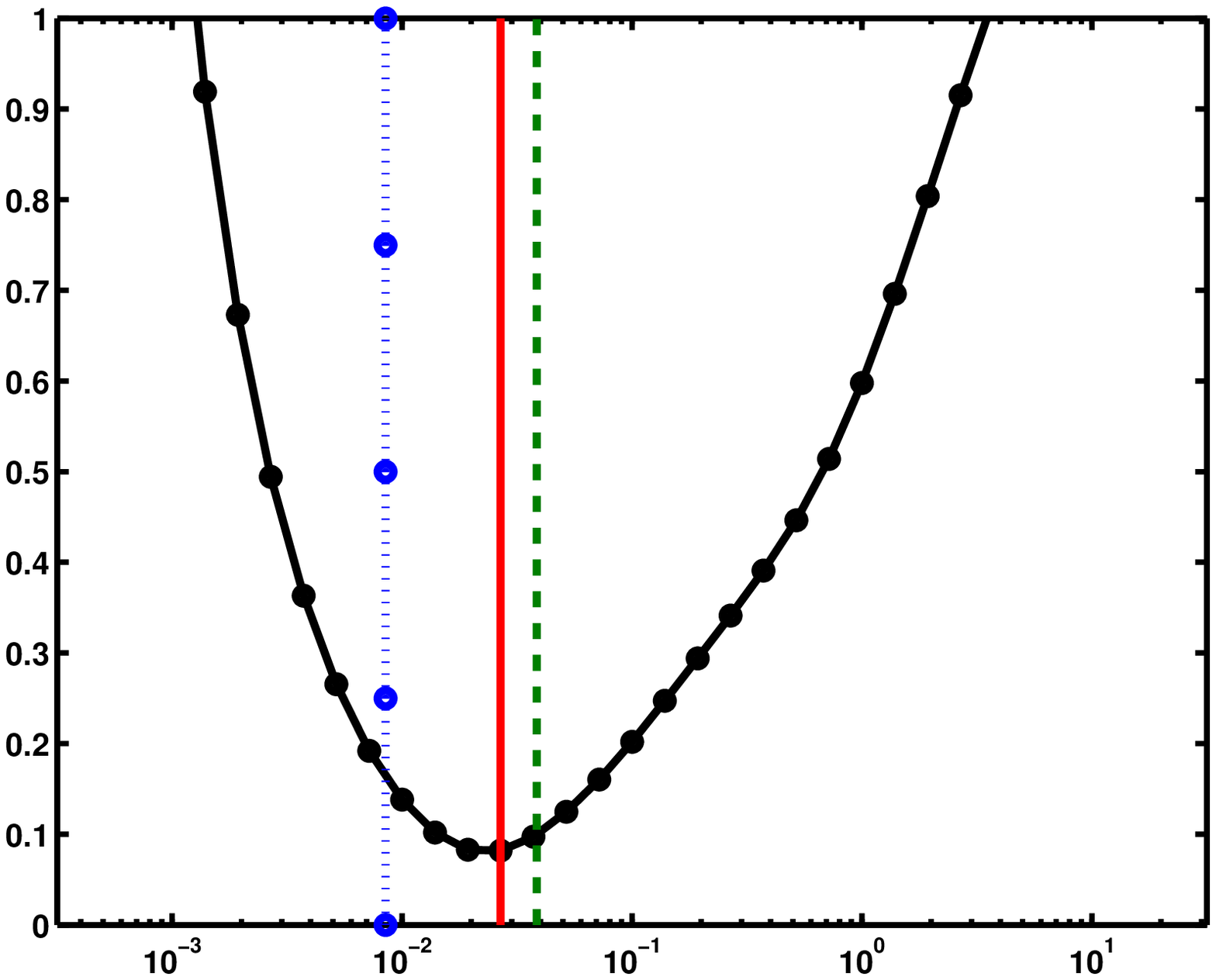}}
\subfigure[$L=L_1$]{\includegraphics[width=1.7in]{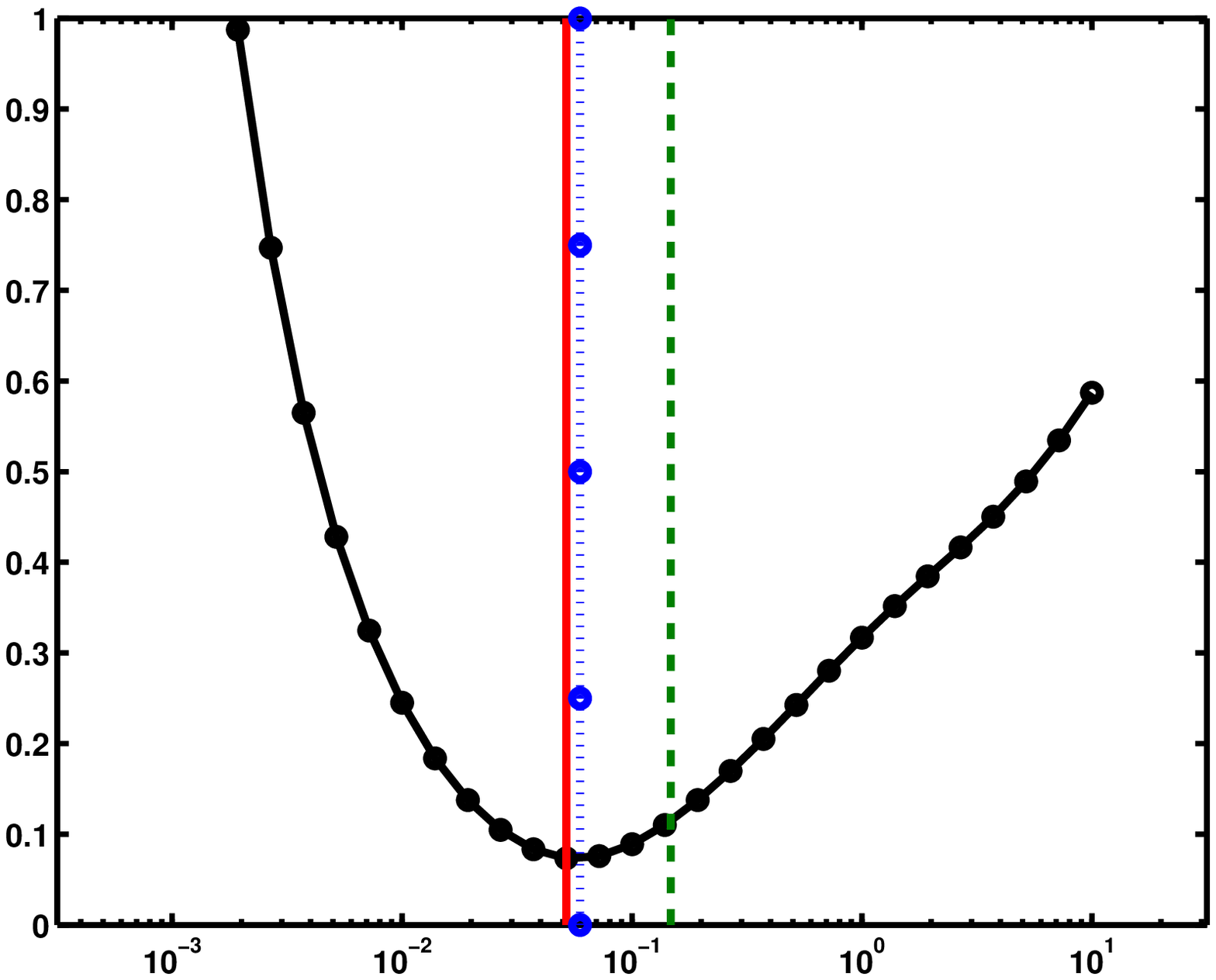}}
\subfigure[$L=L_2$]{\includegraphics[width=1.7in]{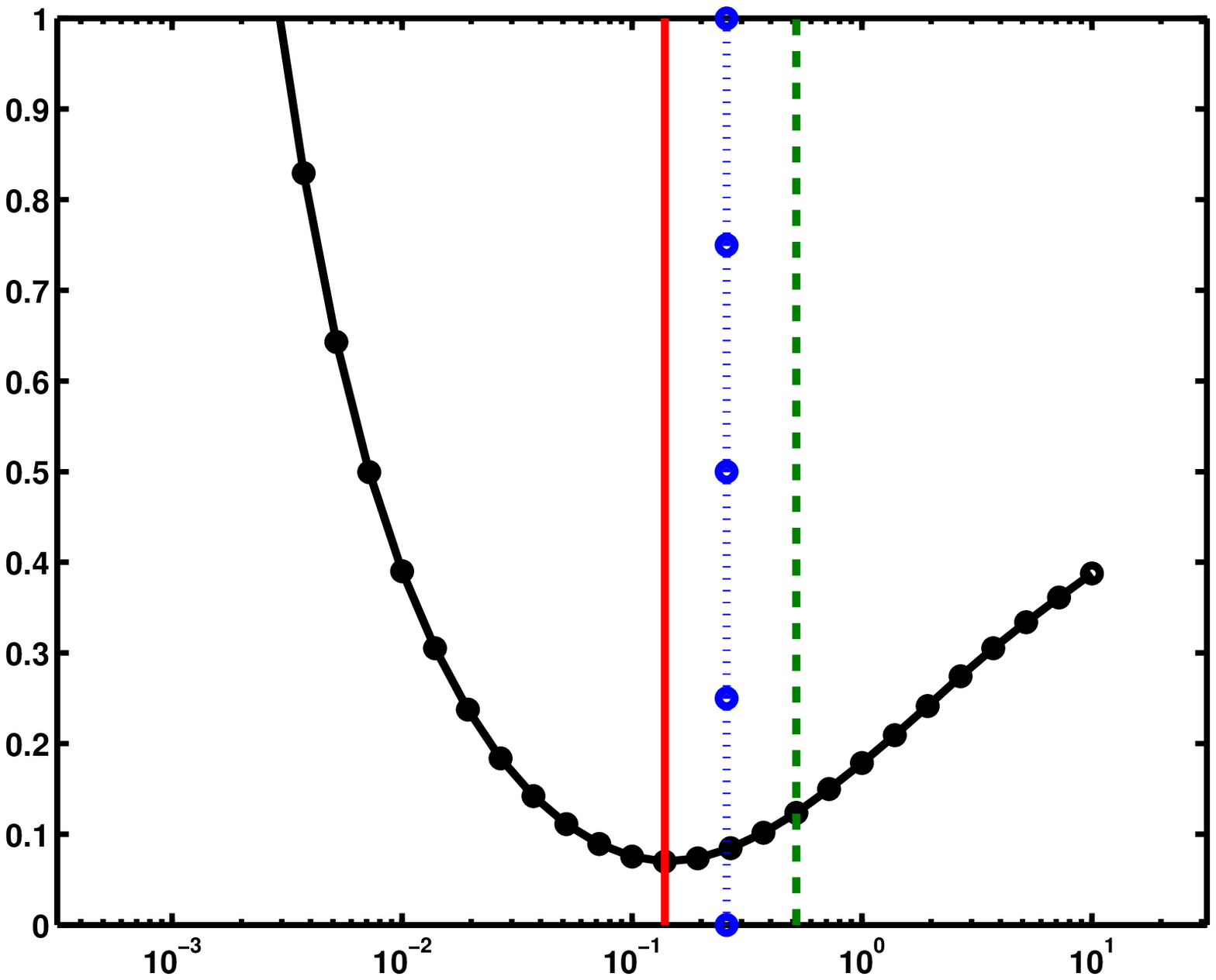}}
\subfigure[$L=I$]{\includegraphics[width=1.7in]{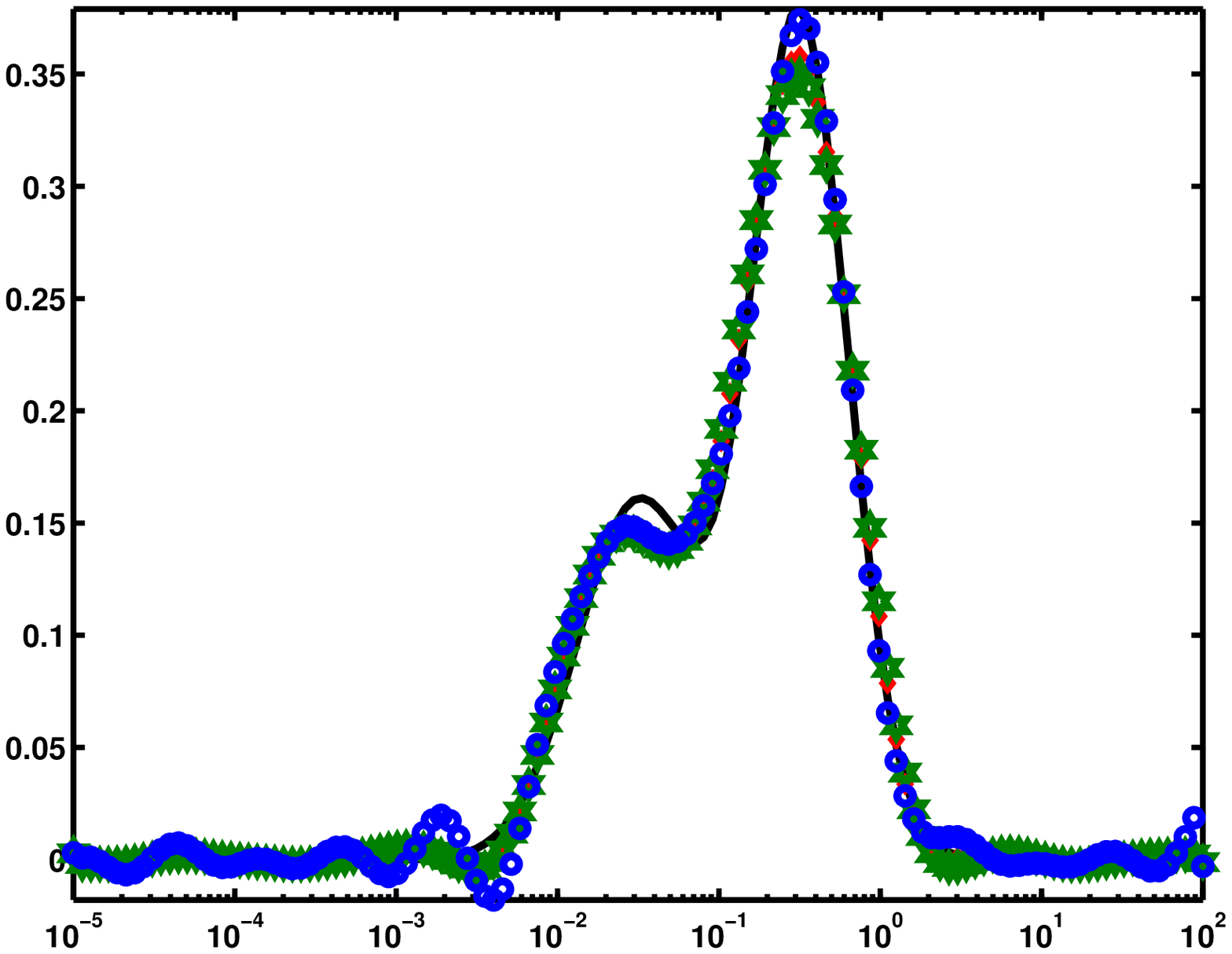}}
\subfigure[$L=L_1$]{\includegraphics[width=1.7in]{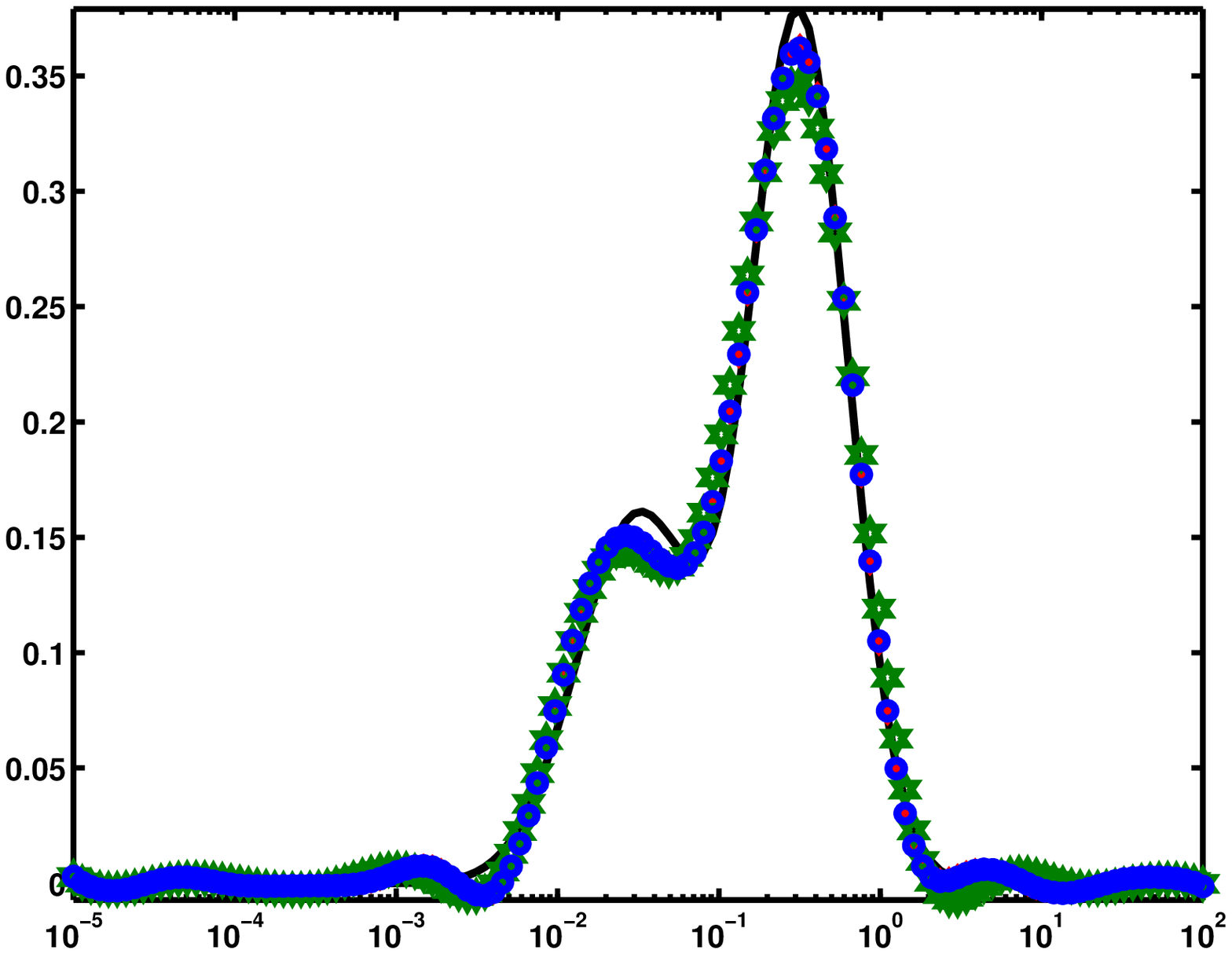}}
\subfigure[$L=L_2$]{\includegraphics[width=1.7in]{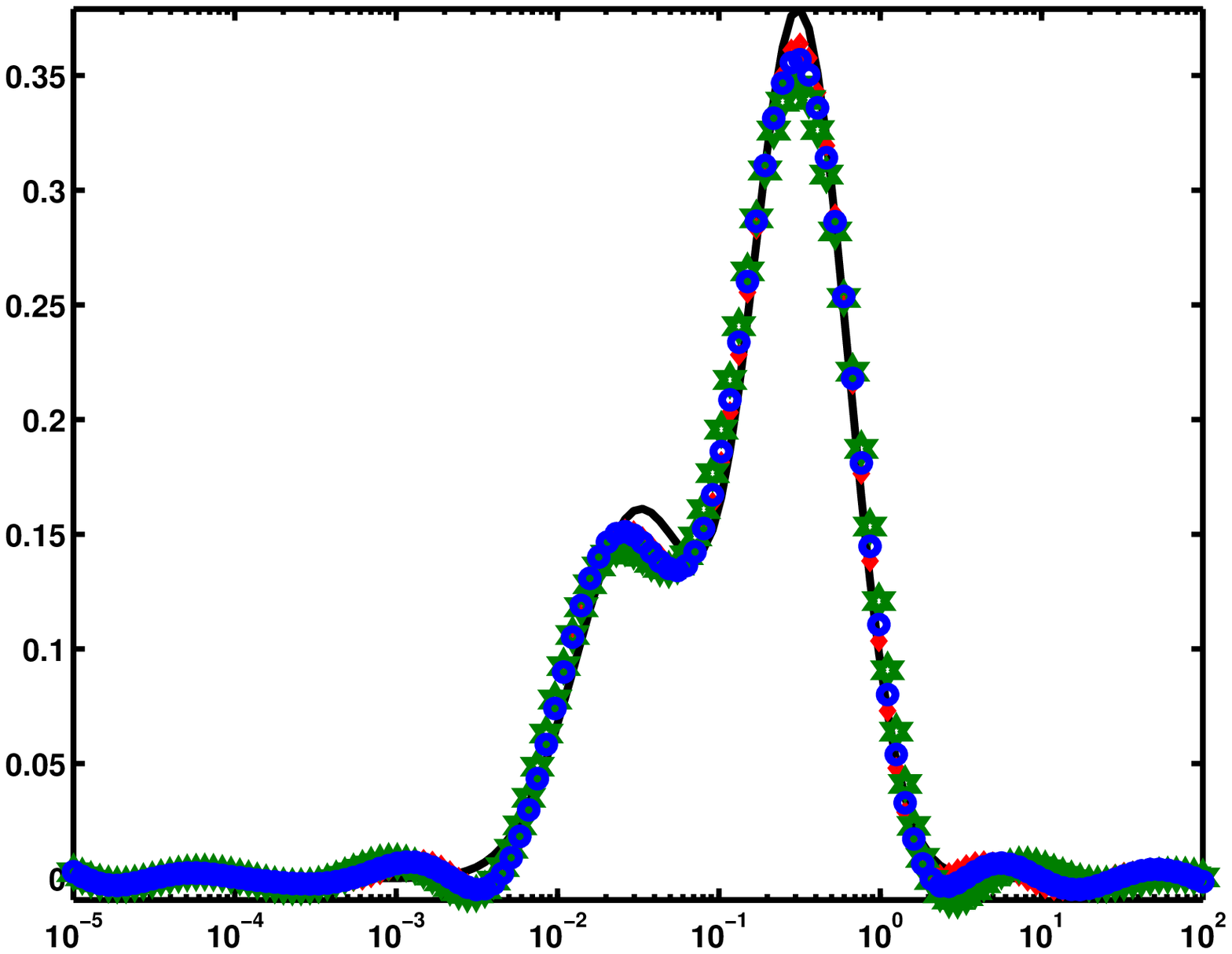}}
\caption{Mean error and example LLS solutions.  $.1\%$ noise, LN-B data set, matrix $A_4$.}
\label{fig-lambdachoiceLN5A4LNLS}
\end{figure}

\begin{figure}[!ht]
  \centering
\subfigure[$L=I$]{\includegraphics[width=1.7in]{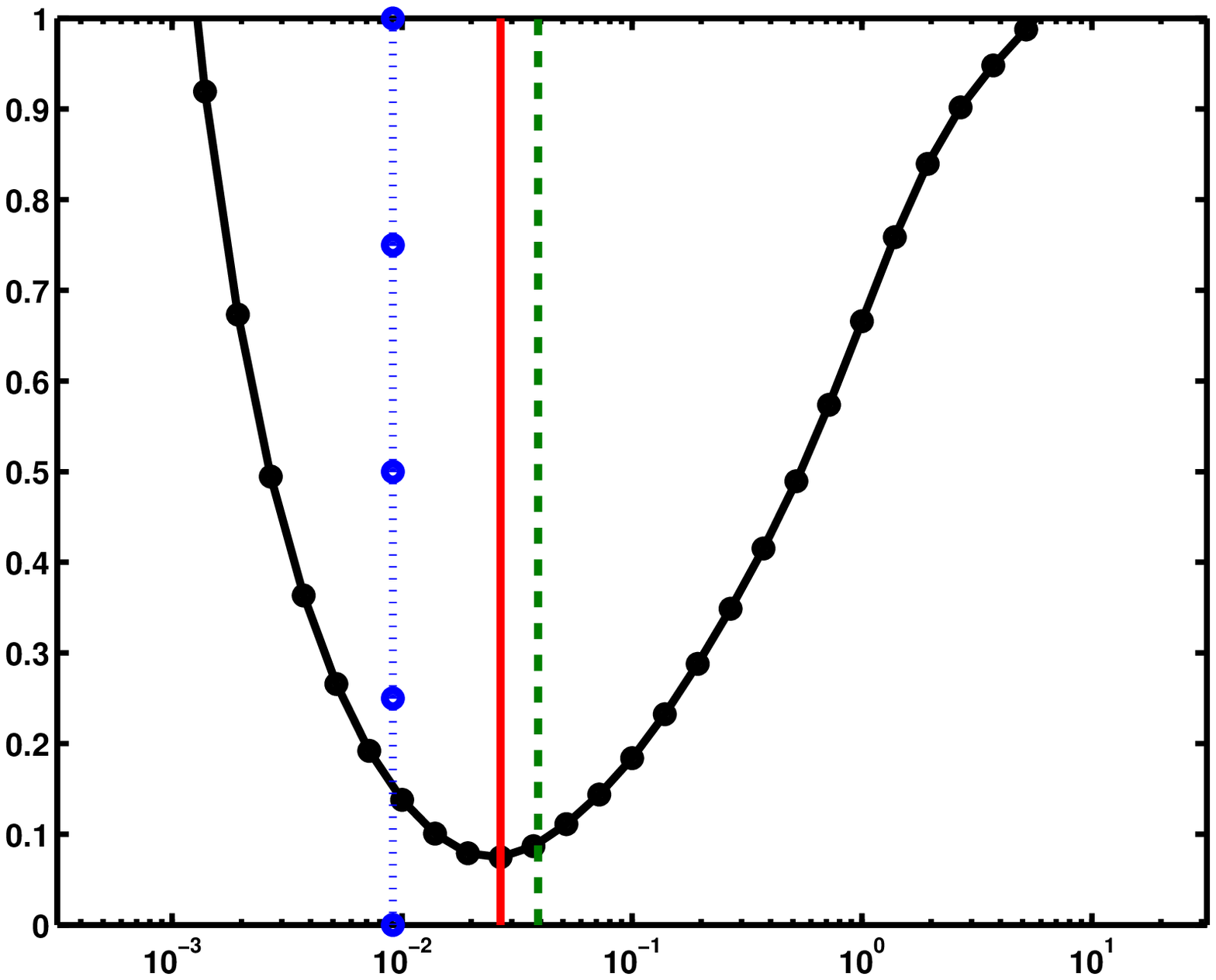}}
\subfigure[$L=L_1$]{\includegraphics[width=1.7in]{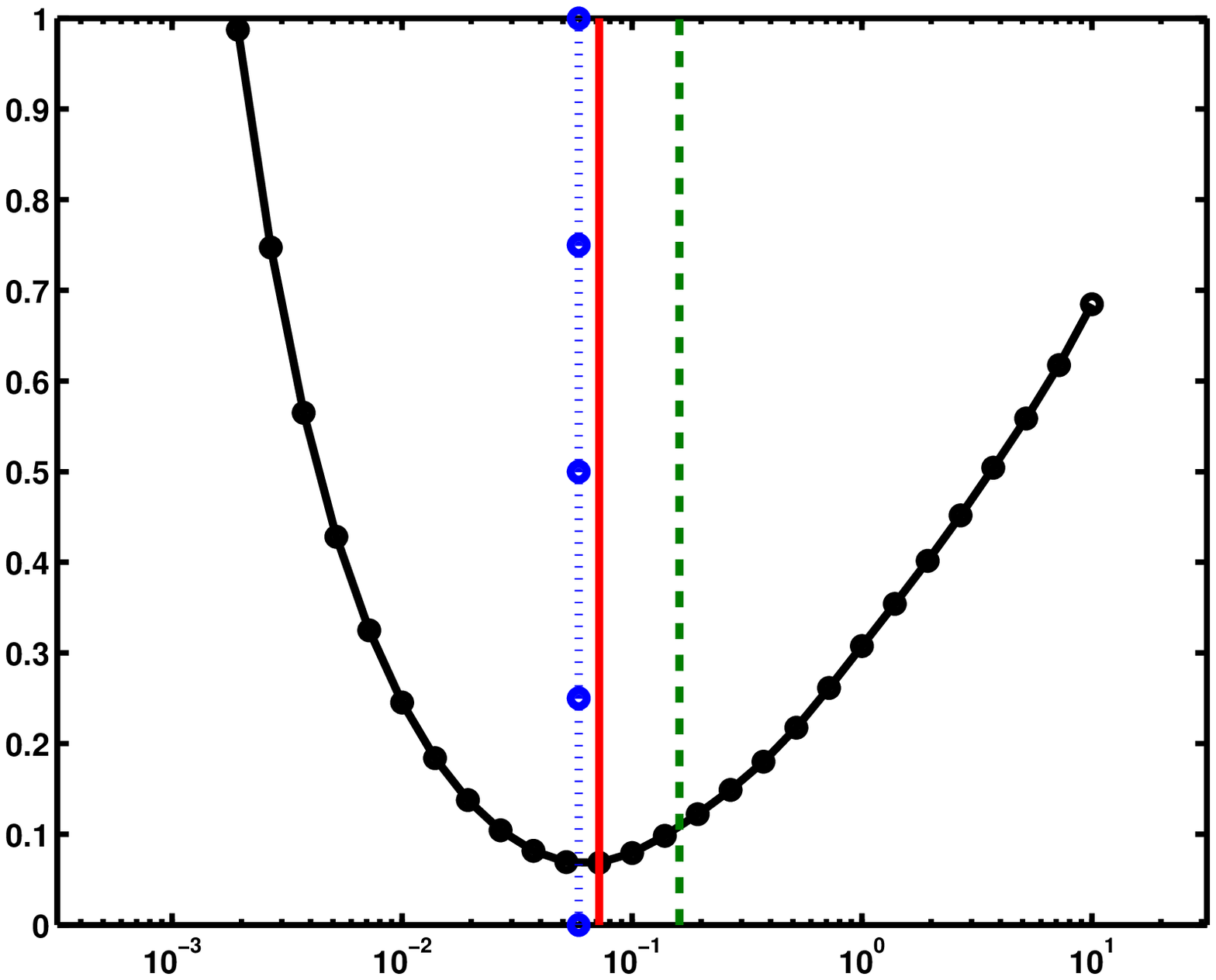}}
\subfigure[$L=L_2$]{\includegraphics[width=1.7in]{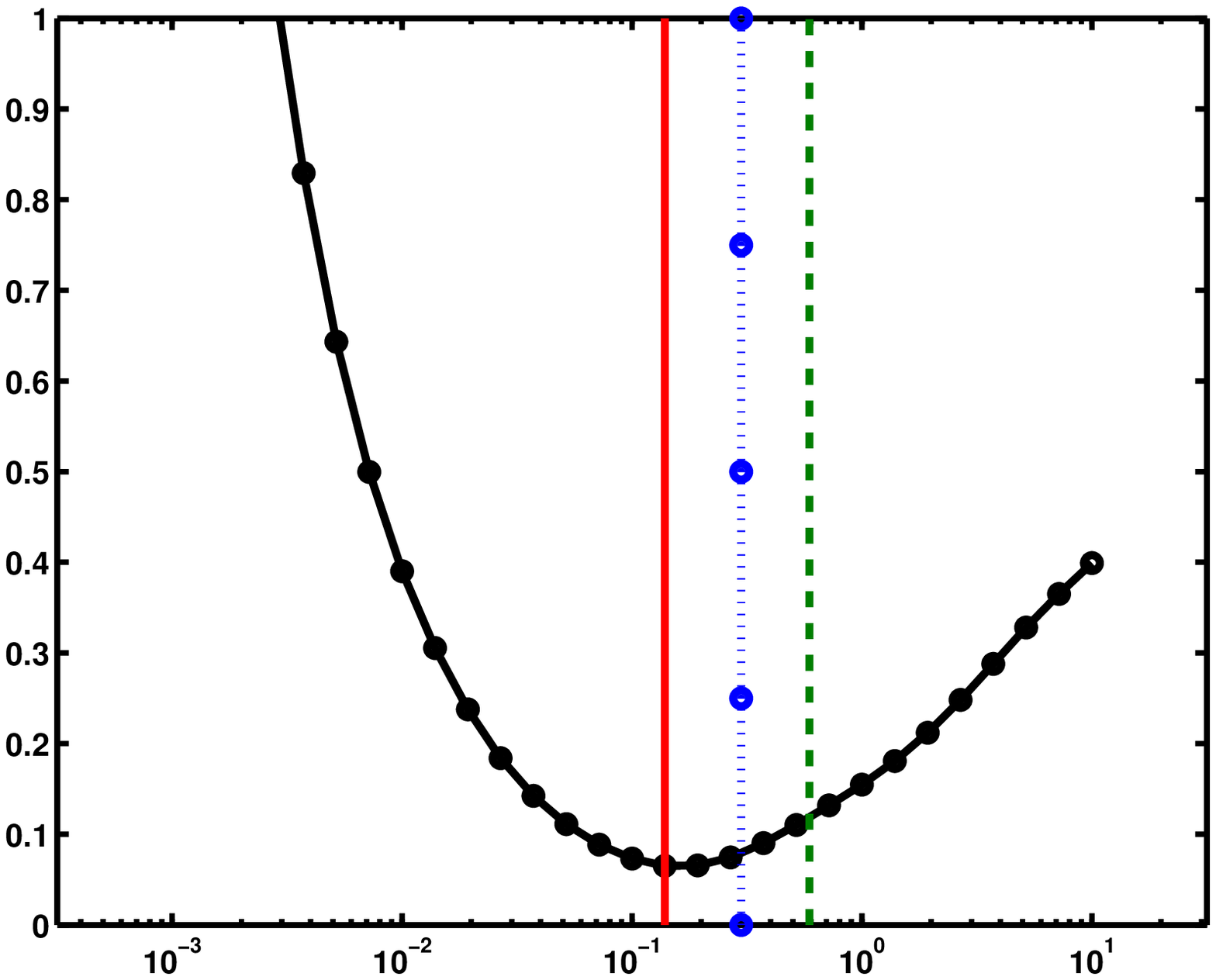}}
\subfigure[$L=I$]{\includegraphics[width=1.7in]{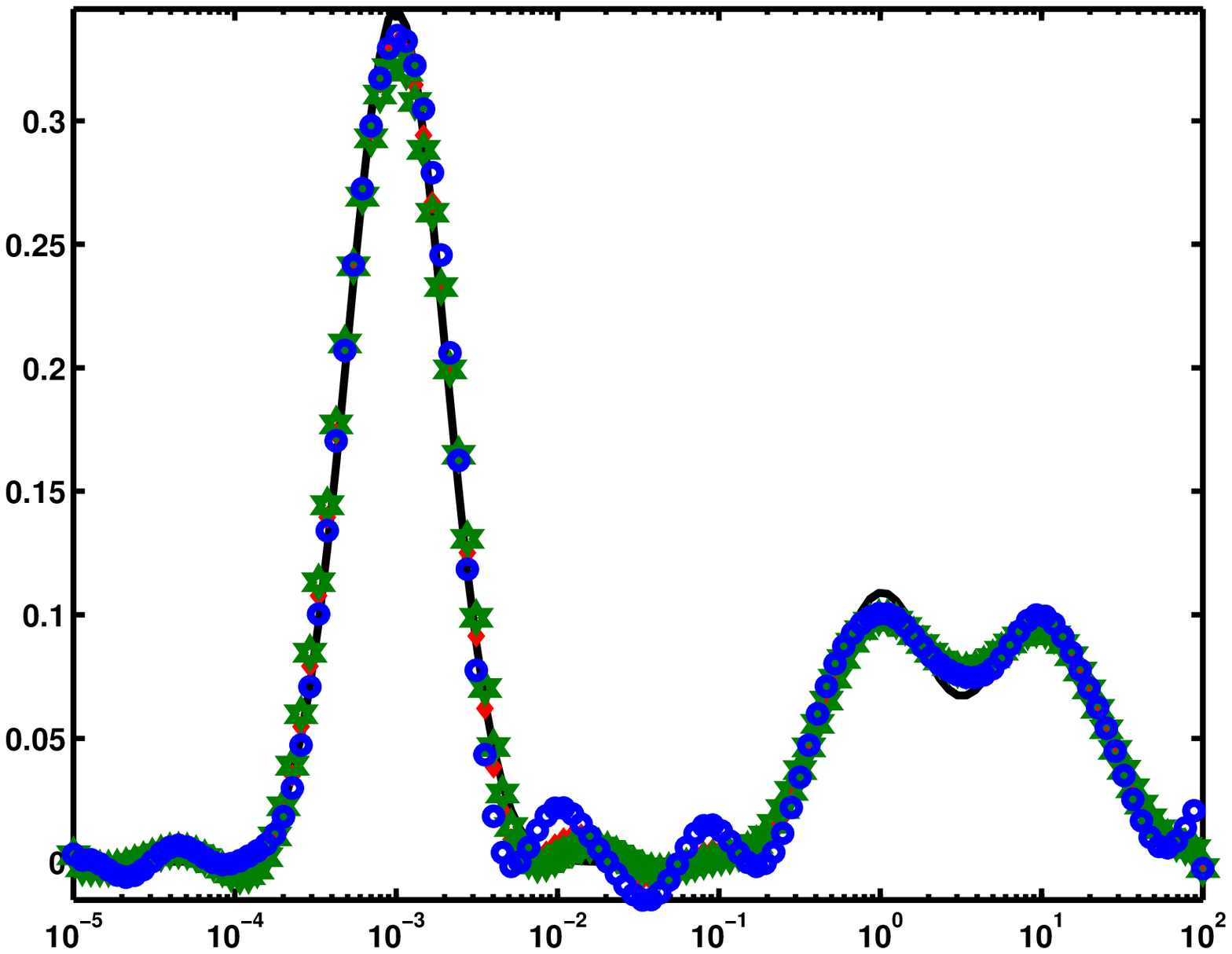}}
\subfigure[$L=L_1$]{\includegraphics[width=1.7in]{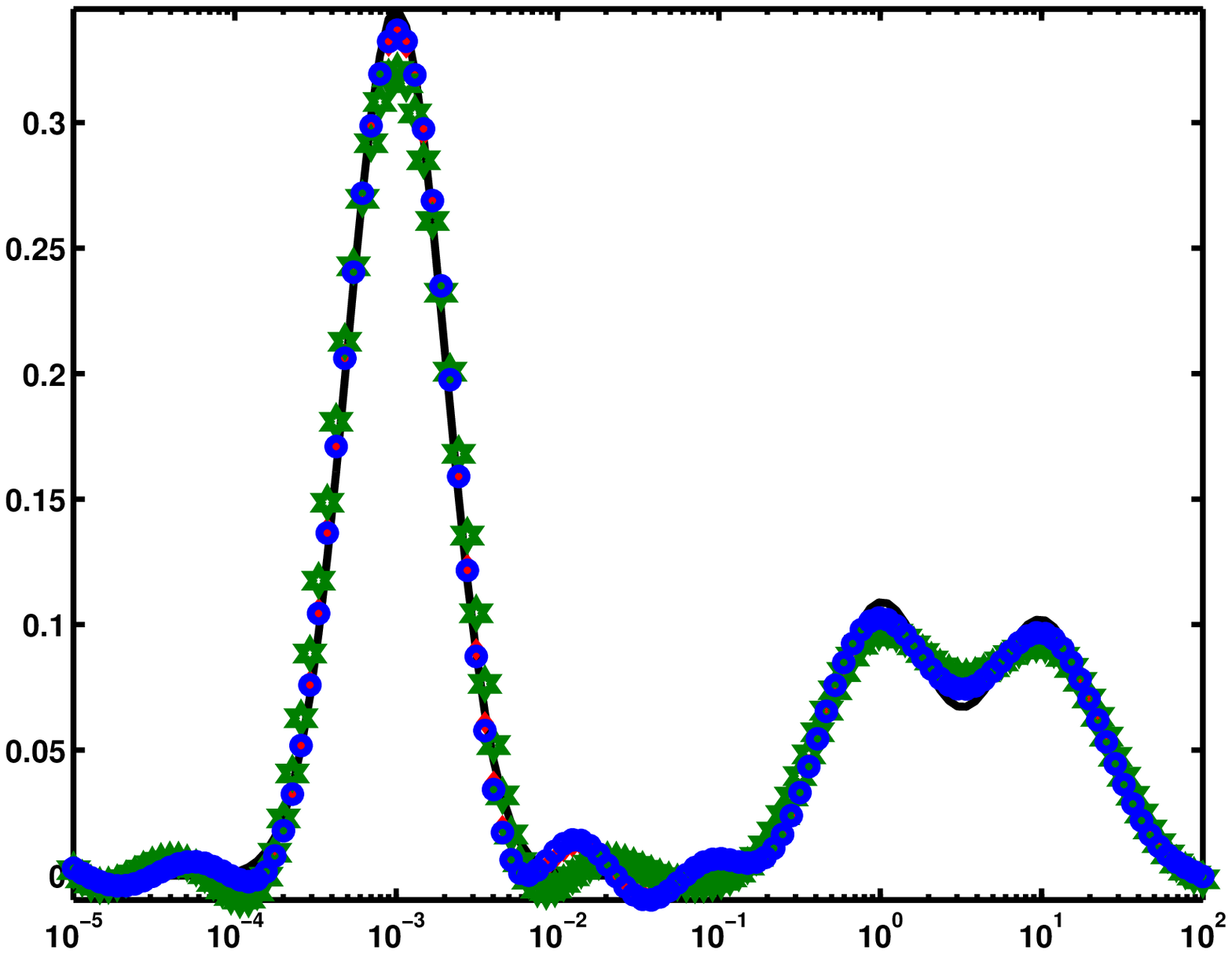}}
\subfigure[$L=L_2$]{\includegraphics[width=1.7in]{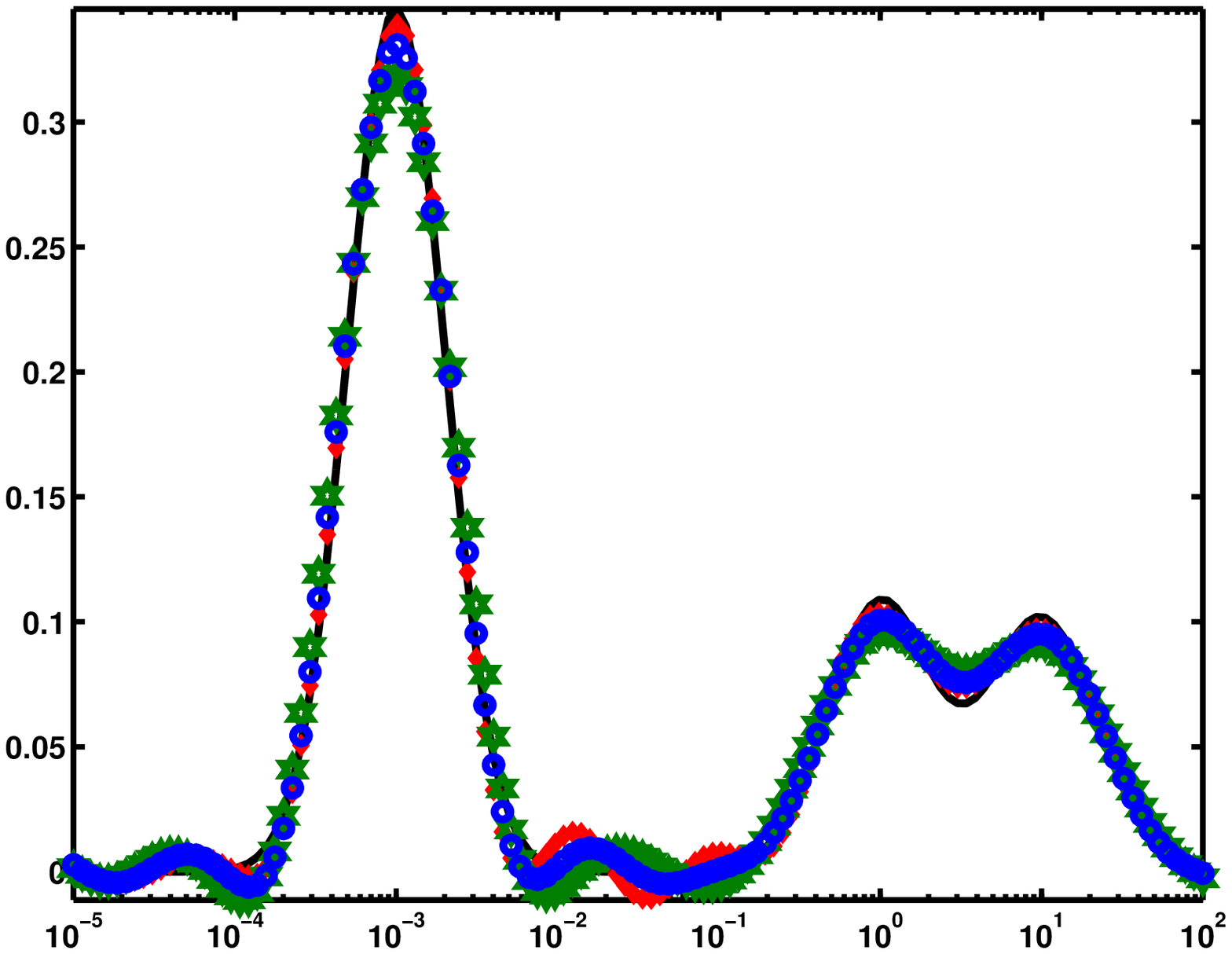}}
\caption{Mean error and example LLS solutions.  $.1\%$ noise, LN-C data set, matrix $A_4$.}
\label{fig-lambdachoiceLN6A4LNLS}
\end{figure}

Finally we present a comparison of the NNLS results with those that are obtained using the Tikhonov regularization \eqref{regsoln} for the same two parameter choice methods in order to assess whether the extra cost of the NN constraint is necessary. 
While extensive results were given in \cite{CSUMS12} for \eqref{regsoln}, these were only for  the $t-$quadrature matrices. The results for the $A_3$ matrix and other noise levels are given in \cite{Supp}, where comparative tables of relative errors are also given. Here we present the results visually in Figures~\ref{fig-lambdachoiceRQ1A4LNLS}-\ref{fig-lambdachoiceLN6A4LNLS} equivalent to the NNLS results in   Section~\ref{parameterchoice}.  We see again that we cannot always anticipate for the optimum choice of $\lambda$ chosen by a specific algorithm will provide a solution with the minimum error, as measured by sampling over multiple choices for the regularization parameter. On the other hand, the performance of these two parameter choice methods for the LS problem is good evidence that the performance for the NNLS problem is consistent. For the solutions it is evident that NNLS provides better control of oscillations around zero due to the positivity constraint. Moreover, because the NNLS does not need to control these oscillations, solutions are less smooth and provide better resolution of the peaks. 
Given the sample sizes for the particular microbial fuel cell application, and possibly other electrochemical applications with limited sampling,  the extra minimal cost associated with finding NNLS solutions offers significantly improved results. Indeed for higher noise levels, the LLS solutions offer little reliable information about the actual physical processes of the model.

\section{Conclusion}\label{conc}
The inverse problem associated with impedance spectroscopy of fuel cells has been discussed. Two models for the underlying distribution function  of relaxation times have been considered.  If the model for the DRT is known to be log-normal or RQ, then nonlinear fitting of the data using the right model can be done very consistently, while trying to fit to the wrong model consistently returns inaccurate results. Moreover, when the noise level is high, distinguishing which model to use in NLS fitting becomes problematic and may yield results that do not accurately describe the data. If the physical model for the process is not known, the inverse problem for estimating the DRT needs to be solved by finding a solution to the discrete linear system. For this system,  it has been shown that the model error can be made negligible by a  change of variables and by extending the effective range of quadrature. Moreover, the conditioning of the problem improves considerably when the  right-preconditioned matrices $A^{\mathrm{s}}$ are used.

To obtain feasible solutions to the discrete linear systems additional constraints are required. Simulations with artificial, but realistic, data demonstrate that the use of NNLS   with a smoothing norm provides higher quality solutions than those obtained without the NN constraint.   Using higher-order smoothing norms also reduces the error in the solutions. Moreover, the LC and NCP criteria are effective regularization parameter choice techniques in the context of the NNLS formulation.  Indeed,   the use of the NCP criterion for parameter choice with the NN  constraint is a novel development of more general use for NNLS in other applications, and has been validated for use both with \texttt{lsqnoneg} in Matlab and the CVX software, \cite{cvx,gb08}.

Although these results  have been verified within the context of the   analysis of fuel cells, there is no reason to suppose that they would not be relevant within the broader framework of solving Fredholm integral equations  for other applications. In particular,  an approach for optimal regularization parameter estimation in the context of non-negatively constrained Tikhonov least squares has been provided.

\section{Acknowledgements}
We would like to thank an anonymous referee who suggested applying the non-negative constrained Tikhonov regularization with algorithms other than \texttt{lsqnonneg}, hence leading to stronger conclusions on the relevance of our techniques. 
Authors Hansen, Hogue and Sander were supported by  NSF CSUMS grant DMS 0703587: ``CSUMS:
Undergraduate Research Experiences for Computational Math Sciences Majors at ASU". Renaut was supported by NSF MCTP grant DMS 1148771: ``MCTP: Mathematics Mentoring Partnership Between Arizona State University and the Maricopa County Community College District",  NSF grant DMS 121655: 
``Novel Numerical Approximation Techniques for Non-Standard Sampling Regimes",  and AFOSR grant 025717 ``Development and Analysis of Non-Classical Numerical Approximation Methods".
Popat was supported by ONR grant N000141210344: ``Characterizing electron
transport resistances from anode-respiring bacteria using electrochemical techniques". All authors wish to thank Professor C\'esar Torres from Arizona State
University for   discussions concerning the EIS modeling.
\appendix
\section{Truncation error}\label{appA}
In \eqref{quadt}-\eqref{quads} the semi-infinite (infinite) integrals are necessarily truncated. An analysis of the impact of the truncation for the lognormal DRT was presented in \cite{CSUMS12}. Here we investigate the extension of this result for the RQ model.  
Consider the single RQ process given by \eqref{DRTsCole}. The error from the upper truncation of the integral of this distribution at $s_{N}$ is given by
\begin{equation}\label{uppertrunc}
E_\mathrm{trunc}^\mathrm{u}(s_{N}) = \int_{s_{N}}^\infty f(s) \,ds = \frac{1}{2} -\frac{\tan ^{-1}\left(\tan \left(\frac{\pi  \beta}{2}\right) \tanh \left(\frac{\beta (s_{N}-\ln(t_0))}{2}\right)\right)}{\pi  \beta}.
\end{equation}
Thus, in order to have $E_\mathrm{trunc}^\mathrm{u}(s_{N}) < \delta$, we must have
\begin{equation}\label{uppertruncbd}s_{N} > \frac{2}{\beta}\tanh^{-1}\left(\frac{\tan\left(\frac{\pi\beta}{2}(1-2\delta)\right)}{\tan \frac{\pi\beta}{2}}\right) + \ln(t_0).\end{equation}
Similarly, because the error from the lower truncation is given by
\begin{equation}\label{lowertrunc}E_\mathrm{trunc}^\mathrm{l}(s_{1}) = \int_{-\infty}^{s_1} f(s)\,ds = \frac{1}{2} + \frac{\tan ^{-1}\left(\tan \left(\frac{\pi  \beta}{2}\right) \tanh \left(\frac{\beta (s_{1}-\ln(t_0))}{2}\right)\right)}{\pi  \beta},\end{equation}
keeping $E_\mathrm{trunc}^\mathrm{l} < \delta$ requires
\begin{equation}\label{lowertruncbd}s_{1} < -\frac{2}{\beta}\tanh^{-1}\left(\frac{\tan\left(\frac{\pi\beta}{2}(1-2\delta)\right)}{\tan \frac{\pi\beta}{2}}\right) +\ln(t_0).\end{equation}
Note that these error bounds are symmetric for $t_0=1$; then $s_{1} = -s_{N}$, $E_\mathrm{trunc}^\mathrm{l}(s_{1}) = E_\mathrm{trunc}^\mathrm{u}(s_{N})$. On the other hand, as $t_0$ moves to the right or left of $1$, the bounds for $s_{N}$, $s_{1}$ shift in tandem, and thus we obtain a requirement on the total range
\begin{equation}
\label{totaltruncE}
E_\mathrm{trunc} =E_\mathrm{trunc}^\mathrm{u}+E_\mathrm{trunc}^\mathrm{l}< 2\delta \,\, \mathrm{for} \,\,  s_{N} - s_{1} > \frac{4}{\beta}\tanh^{-1}\left(\frac{\tan\left(\frac{\pi\beta}{2}(1-2\delta)\right)}{\tan \frac{\pi\beta}{2}}\right).
\end{equation}  
Once $E_\mathrm{trunc}$ is bounded, Lemma 1 from \cite{CSUMS12} can be applied directly, without proof.
\begin{lem}
Suppose that $s_{1}$ and $s_{N}$ are such that the upper and lower truncation errors are each less than $\delta$, i.e.,
\[E_\mathrm{trunc}^\mathrm{u}(s_{N}) = \int_{-\infty}^\infty f(s) \, ds - \int_{-\infty}^{s_{N}} f(s) \, ds < \delta \,\mathrm{ and }\,
E_\mathrm{trunc}^\mathrm{l}(s_{1}) = \int_{-\infty}^\infty f(s) \, ds - \int_{s_1}^{\infty} f(s) \, ds < \delta.\]
Then 
\begin{align*}
E_\mathrm{trunc}^1(\omega) &= \int_{-\infty}^\infty \frac{f(s)}{1+\omega^2e^{2s}} \, ds - \int_{s_1}^{s_{N}} \frac{f(s)}{1+\omega^2e^{2s}} \, ds \leq \delta\left(1 + \frac{1}{1 + \omega^2 e^{2s_{N}}} \right) \leq 2\delta \\
E_\mathrm{trunc}^2(\omega) &= \int_{-\infty}^\infty \frac{\omega e^{2s} f(s)}{1+\omega^2e^{2s}} \, ds  - \int_{s_1}^{s_{N}} \frac{\omega e^{2s} f(s)}{1+\omega^2e^{2s}} \, ds \leq (E_2)_\mathrm{min} + (E_2)_\mathrm{max} <\delta,
\end{align*} 
where 
\begin{align*}
(E_2)_\mathrm{min} &= \int_{-\infty}^{s_1} \frac{\omega e^{s} f(s)}{1+\omega^2 e^{2s}} \, ds \leq 
\begin{cases}
 \delta \frac{\omega e^{s_1}}{1+\omega^2 e^{2s_1}} \le \frac{\delta}{2}&\quad  \omega e^{s_1} < 1\\
\frac{\delta}{2} &\quad   \omega e^{s_1} \geq 1
\end{cases}
 \\
(E_2)_\mathrm{max} &= \int_{s_{N}}^{\infty} \frac{\omega e^{s} f(s)}{1+\omega^2 e^{2s}} \, ds \leq  
\begin{cases}
\delta \frac{\omega e^{s_{N}}}{1+\omega^2 e^{2s_{N}}} \le \frac{\delta}{2}&\quad   \omega e^{s_{N}} \geq 1\\
\frac{\delta}{2} &\quad   \omega e^{s_{N}} < 1.
\end{cases}
\end{align*}
\end{lem}

Observe here that for the standard choice $s_{1}=\ln (T_{\mathrm{min}})$ and $s_{N}=\ln (T_{\mathrm{max}})$,  with $T_{\mathrm{min}}=1/\omega_{\mathrm{max}}$ and $T_{\mathrm{max}}=1/\omega_{\mathrm{min}}$,  then $ \omega e^{s_1} <1$, $\omega e^{s_{N}} >1$, and the bounds for $E_2$ simplify as given. 
Provided that the range for the integration moves with the location of $t_0$, either to the right or left, \eqref{totaltruncE}, the error is controlled appropriately. If on the other hand we always assume $t_0=1$ and pick the range symmetrically with respect to $0$, the error will depend on the actual location of $t_0$.

\end{document}